\documentclass{article}
\usepackage{amssymb}
\usepackage{amsmath}
\usepackage{amsfonts}
\usepackage{graphicx}

\newtheorem{lemma}{Lemma}[section]
\newtheorem{theorem}[lemma]{Theorem}
\newtheorem{corollary}[lemma]{Corollary}
\newtheorem{proposition}[lemma]{Proposition}
\newtheorem{definition}[lemma]{Definition}

\newtheorem{condition}[lemma]{Condition}

\newtheorem{remark}[lemma]{Remark}
\newenvironment{proof}[1][Proof]{\noindent\textbf{#1.} }{\ \rule{0.5em}{0.5em}}
\newcommand{\ZZ}{\mathbb{Z}}

\newcommand{\RR}{\mathbb{R}}

\renewcommand{\ll}{\lambda}

\begin{document}

\title{Small divisor problem in the theory of three-dimensional water gravity waves}
\author{G\'{e}rard Iooss$^{\dag }$, Pavel Plotnikov$^{\ddag }$ \\
\dag {\small \ IUF, INLN UMR 6618 CNRS - UNSA, 1361 rte des
Lucioles, 06560
Valbonne, France}\\
\ddag {\small Russian academy of Sciences, Lavryentyev pr. 15,
Novosibirsk
630090, Russia} \\
{\tiny gerard.iooss@inln.cnrs.fr, } {\tiny
plotnikov@hydro.nsc.ru}} \maketitle

\begin{abstract}
We consider doubly-periodic travelling waves at the surface of an
infinitely deep perfect fluid, only subjected to gravity $g$ and
resulting from the nonlinear interaction of two simply periodic
travelling waves making an angle $2\theta $ between them. \newline
Denoting by $\mu =gL/c^{2}$ the dimensionless bifurcation
parameter  ( $L$ is the wave length along the direction of the
travelling wave and $c$ is the velocity of the wave), bifurcation
occurs for $\mu =\cos \theta$. For non-resonant cases, we first
give a large family of formal three-dimensional gravity travelling
waves, in the form of an expansion in powers of the amplitudes of
two basic travelling waves. "Diamond waves" are a particular case
of such waves, when they are symmetric with respect to the
direction of propagation.\newline \emph{The main object of the
paper is the proof of existence} of such symmetric waves having
the above mentioned asymptotic expansion. Due to the
\emph{occurence of small divisors}, the main difficulty is the
inversion of the linearized operator at a non trivial point, for
applying the Nash Moser theorem. This operator is the sum of a
second order differentiation along a certain direction, and an
integro-differential operator of first order, both depending
periodically of coordinates. It is shown that for almost all
angles $\theta $, the 3-dimensional travelling waves bifurcate for
a set of "good" values of the bifurcation parameter having
asymptotically a full measure near the bifurcation curve in the
parameter plane $(\theta ,\mu ).$
\end{abstract}

\tableofcontents

\numberwithin{equation}{section}

\section{Introduction}

\subsection{Presentation and history of the problem}

We consider small-amplitude three-dimensional doubly periodic
travelling gravity waves on the free surface of a perfect fluid.
These \emph{unforced} waves appear in literature as steady
3-dimensional water waves, since they are steady in a suitable
moving frame. The fluid layer is supposed to be infinitely deep,
and the flow is irrotational only subjected to gravity. The
bifurcation parameter is the horizontal phase velocity, the
infinite depth case being not essentially different from the
finite depth case, except for very degenerate situations that we
do not consider here. The essential difficulty here, with respect
to the existing literature is that \emph{we assume the absence of
surface tension}. Indeed the surface tension plays a major role in
all existing proofs for three-dimensional travelling
gravity-capillary waves, and when the surface tension is very
small, which is the case in many usual situations, this implies a
reduced domain of validity of these results.

\begin{figure}[h]
\includegraphics[height=1.90in]{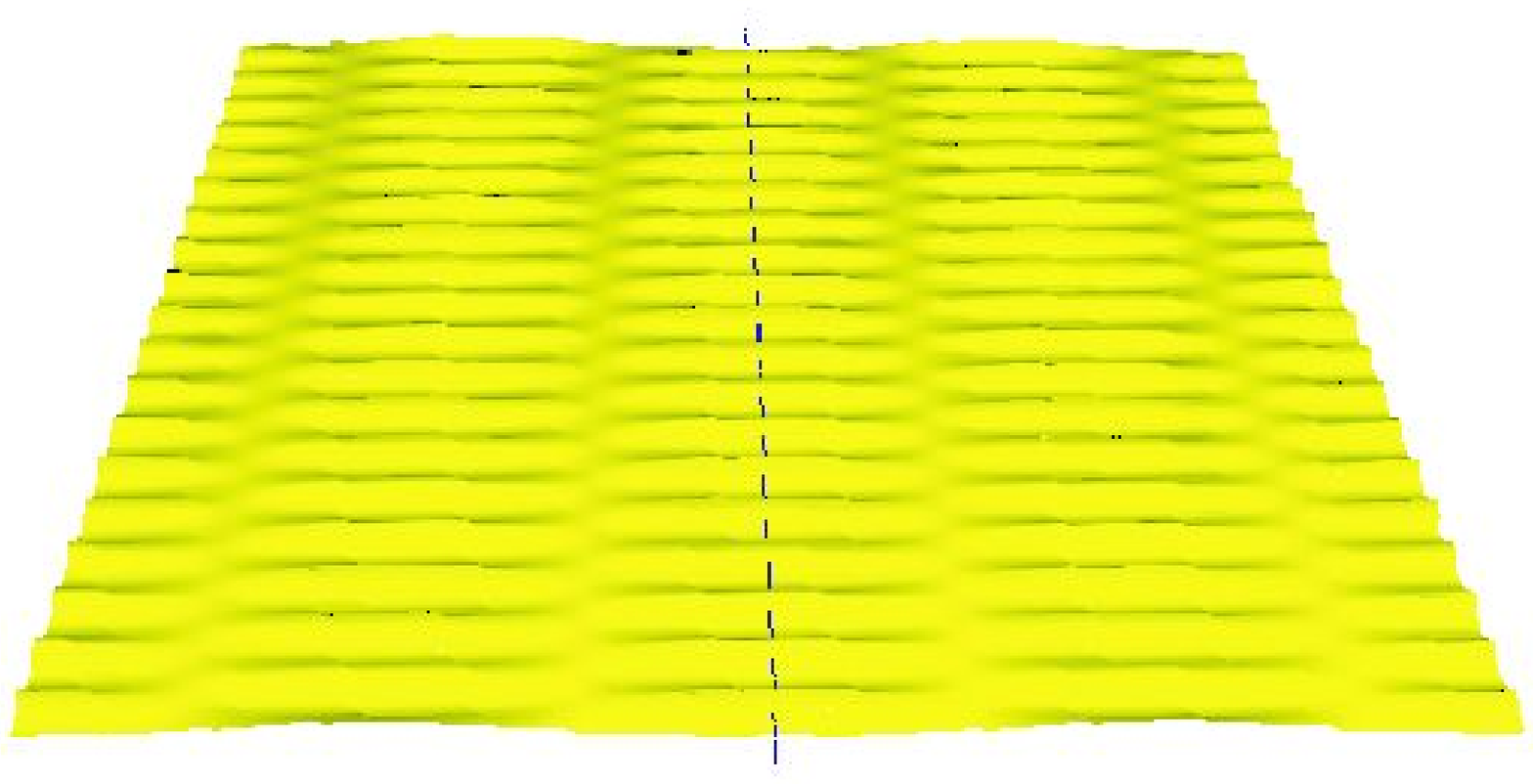} %
\includegraphics[height=1.90in]{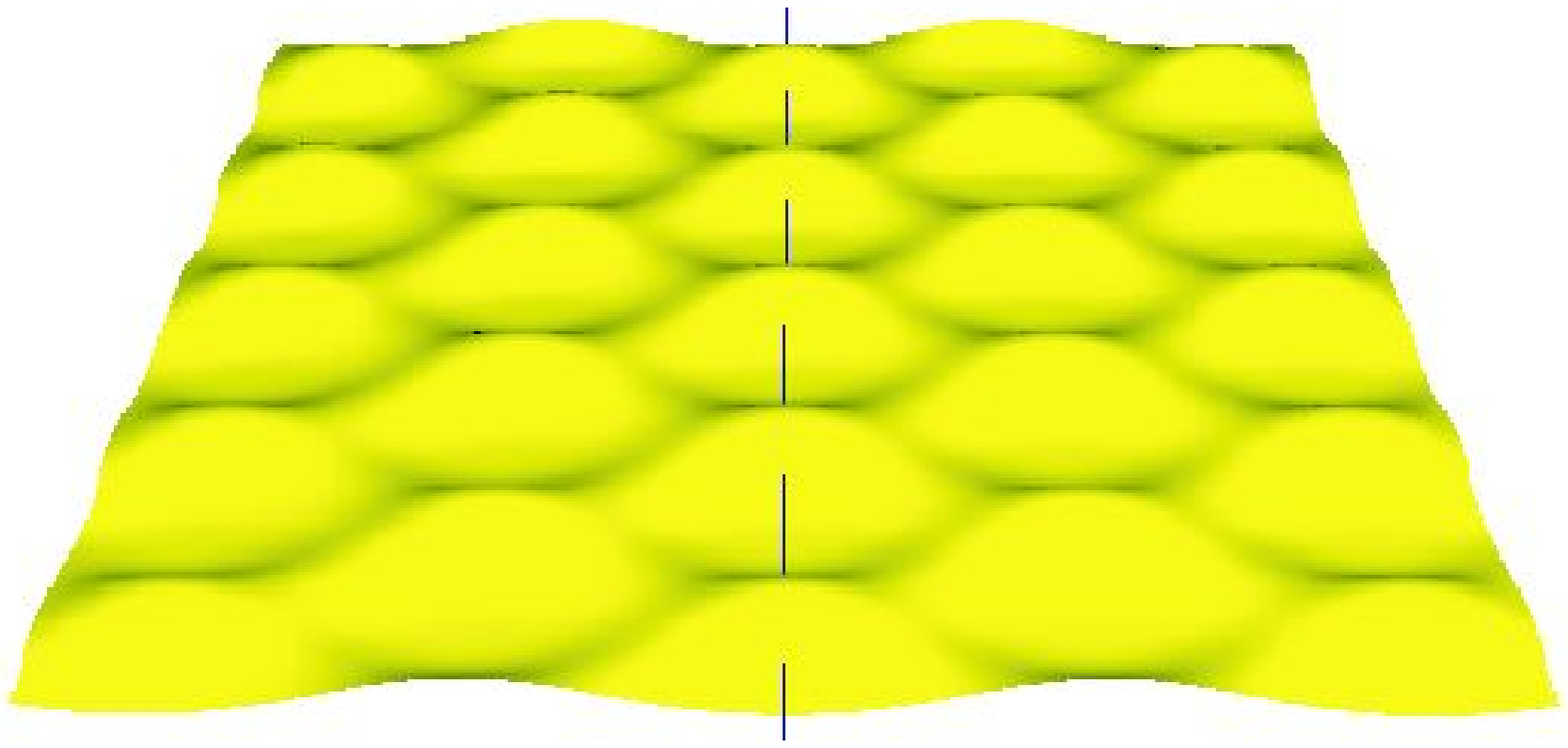} 
\caption{3-dim travelling wave, the elevation $\protect\eta_{\protect%
\varepsilon }^{(2)} $ is computed with formula (\protect\ref{UN}). Top: $%
\protect\theta=11.3^{o}, \protect\tau=1/5, \protect\varepsilon=0.8 \protect%
\mu_c$; bottom: $\protect\theta=26.5^{o}, \protect\tau=1/2,\protect%
\varepsilon=0.6\protect\mu_c$. The dashed line is the direction of
propagation of the waves. Crests are dark and troughs are grey.}
\label{fig1}
\end{figure}

In 1847 Stokes \cite{Stokes} gave a nonlinear theory of
\emph{two-dimensional} travelling gravity waves, computing the
flow up to the cubic order of the amplitude of the waves, and the
first mathematical proofs for such periodic
two-dimensional waves are due to Nekrasov \cite{Nekrasov}, Levi-Civita \cite%
{LeviCivita} and Struik \cite{Struik} about 80 years ago.
Mathematical progresses on the study of \emph{three-dimensional}
doubly periodic water waves came much later. In particular, to our
knowledge, first formal expansions in powers of the amplitude of
three-dimensional travelling waves can be found in papers
\cite{Fuchs} and \cite{Sret}. One can find many references and
results of researches on this subject in the review paper of Dias
and Kharif \cite{dia-khar} (see section 6).
The work of Reeder and Shinbrot (1981)\cite%
{Reed-Shin} represents a big step forward. These authors consider
symmetric diamond patterns, resulting from (horizontal) wave
vectors belonging to a lattice $\Gamma ^{\prime }$ (dual to the
spatial lattice $\Gamma $ of the doubly periodic pattern) spanned
by two wave vectors $K_{1}$ and $K_{2}$ with the \emph{same
length}, the velocity of the wave being in the direction of the
bissectrix of these two wave vectors, taken as the $x_{1}$
horizontal axis. We give at Figure \ref{fig1} two examples of
patterns for these waves
(see the detailed comment about these pictures at the end of subsection \ref%
{diamond}). These waves also appear in litterature as "short
crested waves"
(see Roberts and Schwartz \cite{Rob-Schw}, Bridges, Dias, Menasce \cite%
{B-D-M} for an extensive discussion on various situations and
numerical
computations). If we denote by $\theta $ the angle between $K_{1}$ and the $%
x_{1}-$ axis, Reeder and Shinbrot proved that bifurcation to
diamond waves occurs provided the angle $\theta $ is not too close
to 0 or to $\pi /2,$ and provided that the \emph{surface tension
is not too small}. In addition their result is only valid outside
a "bad" set in the parameter space, corresponding to resonances, a
quite small set indeed. This means that if one considers the
dispersion relation $\Delta (K,\mathbf{c})=0,$ where $K$ and
$\mathbf{c}\in
\mathbb{R}
^{2}$ are respectively a wave vector and the velocity of the
travelling
wave, then there is no resonance if for the critical value of the velocity $%
\mathbf{c}_{0}$ there are only the four solutions $\pm K_{1},\pm
K_{2}$ of the dispersion equation, for $K\in \Gamma ^{\prime }$
(i.e. for $K$ being any integer linear combination of $K_{1}$ and
$K_{2})$. The fact that the surface tension is supposed not to be
too small is essential for being able to use Lyapunov-Schmidt
technique, and the authors mention a \emph{small
divisor problem if there is no surface tension}, as computed for example in \cite%
{Rob-Schw}. Notice that the existence of  spatially bi-periodic
\emph{gravity} water waves was proved by Plotnikov   in
\cite{PIP1}, \cite{PIP2}
 in  the case of finite depth and for fixed rational values of $gL/c^{2}\tan\theta$, where $g,L,c$ are
  respectively the acceleration of gravity, the wave length in the direction of propagation,
  and the velocity of the wave.  Indeed, such a special choice of  parameters
avoids  resonances and the small divisor problem, because the
pseudo-inverse of the linearized operator is bounded.

Craig and Nicholls (2000) \cite{Craig-Nicholls2000} used the
hamiltonian formulation introduced by Zakharov
\cite{V.E.Zakharov}, in coupling the Lyapunov-Schmidt technique
with a variational method on the bifurcation equation. Still in
the presence of surface tension, they could suppress the
restriction of Reeder and Shinbrot on the "bad" resonance set in
parameter space, but they pay this complementary result in losing
the smoothness of the solutions. Among other results, the other
paper by Craig and Nicholls (2002) \cite{craig} gives the
principal parts of "simple" doubly periodic waves (i.e. in the non
resonant cases), expanded in Taylor series, taking into account
the two-dimensions of the parameter $\mathbf{c}$. They emphasize
the fact that this expansion is only formal in the absence of
surface tension.

Mathematical results of another type are obtained in using
"spatial dynamics", in which one of the horizontal coordinates
(the distinguished direction) plays the role of a time variable,
as was initiated by Kirchg\"assner \cite{Kirch} and extensively
applied to two-dimensional water wave problems (see a review in
\cite{Dias-Io}). The advantage of this method is that one does not
choose the behavior of the solutions in the direction of the
distinguished coordinate, and solutions periodic in this
coordinate are a particular case, as well as quasi-periodic or
localized solutions (solitary waves). In this framework one may a
priori assume periodicity in a direction transverse to the
distinguished direction, and a periodic solution in the
distinguished direction is automatically doubly periodic. The
first mathematical results obtained by this method, containing
3-dimensional doubly
periodic travelling waves, start with Haragus, Kirchg%
\"{a}ssner, Groves and Mielke (2001) \cite{Gr-Mi}, \cite{Gro}, \cite%
{Harag-Kirch}, generalized by Groves and Haragus (2003)
\cite{Grov-Harag}. They use a hamiltonian formulation and center
manifold reduction. This is essentially based on the fact that the
spectrum of the linearized operator is discrete and has only a
finite number of eigenvalues on the imaginary axis. These
eigenvalues are related with the dispersion relation mentioned
above. Here, one component (or multiples of such a component) of
the wave vector $K$ is imposed in a direction transverse to the
distinguished one, and there is no restriction for the component
of $K$ in the distinguished direction, which, in solving the
dispersion relation, gives the eigenvalues of the linearized
operator on the imaginary axis. The resonant situations, in the
terminology of Craig and Nicholls correspond here to more than one
pair of eigenvalues on the imaginary axis, (in addition to the
origin). In all cases it is known that the largest eigenvalue on
the imaginary axis leads to a family of periodic solutions, via
the Lyapunov center theorem (hamiltonian case), so, here again,
there is no restriction on the resonant set in the parameter space
at a fixed finite depth. The only restriction with this
formulation is that it is necessary to assume that the depth of
the fluid layer is finite. This ensures that the spectrum of the
linearized operator has a spectral gap near the imaginary axis,
which allows to use the center manifold reduction method. In fact
if we restrict the study to periodic solutions as here, the center
manifold reduction is not necessary, and the infinite depth case
might be considered in using an extension of the proof of
Lyapunov-Devaney  center theorem in the spirit of \cite{Io}, in
this case where 0 belongs to the continuous spectrum.
However, it appears that the \emph{%
number of imaginary eigenvalues becomes infinite when the surface
tension cancels}, which prevents the use of center manifold
reduction in the limiting case we are considering in the present
paper, not only because of the infinite depth.

\subsection{Formulation of the problem\label{basic}}

Since we are looking for waves travelling with velocity
$\mathbf{c}$, let us consider\emph{\ the system in the moving
frame} where the waves look steady.
Let us denote by $\varphi $ the potential defined by%
\begin{equation*}
\varphi =\phi -\mathbf{c}\cdot X,
\end{equation*}%
where $\phi $ is the usual velocity potential, $X=(x_{1},x_{2})$
is the 2-dim horizontal coordinate, $x_{3}$ is the vertical
coordinate, and the fluid region is
\begin{equation*}
\Omega =\{(X,x_{3}):-\infty <x_{3}<\eta (X)\},
\end{equation*}%
which is bounded by the free surface $\Sigma $ defined by%
\begin{equation*}
\Sigma =\{(X,x_{3}):x_{3}=\eta (X)\}.
\end{equation*}%
We also make a scaling in choosing $|\mathbf{c}|$ for the velocity
scale,
and $L$ for a length scale (to be chosen later), and we still denote by $%
(X,x_{3})$ the new coordinates, and by $\varphi ,\eta $ the
unknown functions. Now defining the parameter $\mu
=\frac{gL}{c^{2}}$ (the Froude number is $\frac{c}{\sqrt{gL}})$
where $g$ denotes the acceleration of gravity, and $\mathbf{u}$
the unit vector in the direction of $\mathbf{c},$
the system reads%
\begin{eqnarray}\label{equation1}
\Delta \varphi &=&0\text{ \ in }\Omega , \\\label{equation2}
\nabla _{X}\eta \cdot (\mathbf{u}+\nabla _{X}\varphi
)-\frac{\partial \varphi }{\partial x_{3}} &=&0\text{ \ on }\Sigma
, \\\label{equation3} \mathbf{u}\cdot \nabla _{X}\varphi
+\frac{(\nabla \varphi )^{2}}{2}+\mu \eta &=&0\text{ \ on }\Sigma
, \\\nonumber \nabla \varphi &\rightarrow &0\text{ as
}x_{3}\rightarrow -\infty .
\end{eqnarray}

\textbf{Hilbert spaces of periodic functions.} We specialize our
study to \emph{spatially periodic 3-dimensional travelling waves},
i.e. the solutions $\eta $ and $\varphi $ are \emph{bi-periodic in
}$X.$ This means that there are two independent wave vectors
$K_{1},K_{2}\in
\mathbb{R}
^{2}$ generating a lattice
\begin{equation*}
\Gamma ^{\prime }=\{K=n_{1}K_{1}+n_{2}K_{2}:n_{j}\in
\mathbb{Z}
\},
\end{equation*}%
and a dual lattice $\Gamma $ of periods in $%
\mathbb{R}
^{2}$ such that%
\begin{equation*}
\Gamma =\{\lambda =m_{1}\lambda _{1}+m_{2}\lambda _{2}:m_{j}\in
\mathbb{Z}
,\lambda _{j}\cdot K_{l}=2\pi \delta _{jl}\}.
\end{equation*}%
The Fourier expansions of $\eta $ and $\varphi $ are in terms of
$e^{iK\cdot X},$ where $K\in \Gamma ^{\prime }$ and $K\cdot
\lambda =2n\pi ,$ $n\in
\mathbb{Z}
$, for $\lambda \in \Gamma .$ \emph{The situation we consider in
the further analysis, is with a lattice }$\Gamma ^{\prime
}$\emph{\ generated by the symmetric wave vectors} $K_{1}=(1,\tau
),$ $K_{2}=(1,-\tau ).$ In such a
case the functions on $%
\mathbb{R}
^{2}/\Gamma $ are 2$\pi -$ periodic in $x_{1},$ 2$\pi /\tau -$ periodic in $%
x_{2},$ and invariant under the shift $(x_{1},x_{2})\mapsto
(x_{1}+\pi ,x_{2}+\pi /\tau )$ (and conversely). We define the
Fourier coefficients of
a bi- periodic function $u$ on such lattice $\Gamma _{\tau }$ by%
\begin{equation*}
\widehat{u}(k)=\frac{\sqrt{\tau }}{2\pi }\int_{[0,2\pi ]\times
\lbrack 0,2\pi /\tau ]}u(X)\exp (-ik\cdot X)dX.
\end{equation*}%
For $m\geq 0$ we denote by $H^{m}(%
\mathbb{R}
^{2}/\Gamma )$ the Sobolev space of bi-periodic functions of $X\in
\mathbb{R}
^{2}/\Gamma $ which are square integrable on a period, with their
partial
derivatives up to order $m$, and we can choose the norm as%
\begin{equation*}
||u||_{m}=\left( \sum_{k\in \Gamma ^{\prime }}(1+|k|)^{2m}|\widehat{u}%
(k)|^{2}\right) ^{1/2}.
\end{equation*}

\textbf{Operator equations.} Now, we reduce the above system for
$(\varphi ,\eta )$ to a system of two scalar equations in choosing
the new unknown
function%
\begin{equation*}
\psi (X)=\varphi (X,\eta (X)),
\end{equation*}%
and we define the Dirichlet-Neumann operator $\mathcal{G}_{\eta }$ by%
\begin{eqnarray}
\mathcal{G}_{\eta }\psi &=&\sqrt{1+(\nabla _{X}\eta )^{2}}\frac{d\varphi }{dn%
}|_{x_{3}=\eta (X)}  \label{Dir-Neu} \\
&=&\frac{\partial \varphi }{\partial x_{3}}|_{x_{3}=\eta
(X)}-\nabla _{X}\eta \cdot \nabla _{X}\varphi  \notag
\end{eqnarray}%
where $n$ is normal to $\Sigma $, exterior to $\Omega ,$ and
$\varphi $ is
the solution of the $\eta -$ dependent Dirichlet problem%
\begin{eqnarray*}
\Delta \varphi &=&0,\text{ \ }x_{3}<\eta (X) \\
\varphi &=&\psi ,\text{ \ }x_{3}=\eta (X), \\
\nabla \varphi &\rightarrow &0\text{ as }x_{3}\rightarrow -\infty
.
\end{eqnarray*}%
Notice that this definition of $\mathcal{G}_{\eta }$ follows
\cite{Lannes}
and insures the selfadjointness and positivity of this linear operator in $%
L^{2}(%
\mathbb{R}
^{2}/\Gamma )$ (see Appendix \ref{a0}). Our definition differs from another usual way
of defining the Dirichlet - Neumann operator without the square root in factor in (\ref{Dir-Neu}). 
Now we have the identity(\ref%
{Dir-Neu}) and the system to solve reads%
\begin{equation}
\mathcal{F}(U,\mu ,\mathbf{u})=0,\text{ \ }\mathcal{F}=(\mathcal{F}_{1},%
\mathcal{F}_{2}),  \label{F(U,mu)1}
\end{equation}%
where $U=(\psi ,\eta ),$ and%
\begin{eqnarray}
\mathcal{F}_{1}(U,\mu ,\mathbf{u}) &=&:\mathcal{G}_{\eta }(\psi )-\mathbf{u}%
\cdot \nabla _{X}\eta ,  \label{basic1} \\
\mathcal{F}_{2}(U,\mu ,\mathbf{u}) &=&:\mathbf{u}\cdot \nabla
_{X}\psi +\mu
\eta +\frac{(\nabla \psi )^{2}}{2}+  \label{basic2} \\
&&-\frac{1}{2(1+(\nabla _{X}\eta )^{2})}\{\nabla _{X}\eta \cdot
(\nabla _{X}\psi +\mathbf{u})\}^{2}.  \notag
\end{eqnarray}%
Let us define the 2-components function space
\begin{equation*}
\mathbb{H}^{m}(%
\mathbb{R}
^{2}/\Gamma )=H_{0}^{m}(%
\mathbb{R}
^{2}/\Gamma )\times H^{m}(%
\mathbb{R}
^{2}/\Gamma )
\end{equation*}%
We denote the norm of $U$ in $\mathbb{H}^{m}(%
\mathbb{R}
^{2}/\Gamma )$ by
\begin{equation*}
||U||_{m}=||\psi ||_{H^{m}}+||\eta ||_{H^{m}},
\end{equation*}%
where $H_{0}^{m}$ means functions with 0 average, and $U=(\psi
,\eta ).$ The 0 average condition comes from the fact that the
value $\psi $ of the potential is defined up to an additive
constant (easily checked in equations (\ref{basic1}),
(\ref{basic2})). Moreover, the average of the right hand side of
(\ref{basic1}) is 0 as it can be easily checked (this is proved
for instance in \cite{craig}). We have the following

\begin{lemma}
\label{F(U,mu)} For any fixed $m\geq 3$, the mapping%
\begin{equation*}
(U,%
{\mu}%
,\mathbf{u})\mapsto \mathcal{F}(U,%
{\mu}%
,\mathbf{u})\text{\ \ is \ }C^{\infty }:\mathbb{H}^{m}(%
\mathbb{R}
^{2}/\Gamma )\times
\mathbb{R}
\times \mathbb{S}_{1}\rightarrow \mathbb{H}^{m-1}(%
\mathbb{R}
^{2}/\Gamma )
\end{equation*}%
in the neighborhood of $\{0\}\times
\mathbb{R}
\times \mathbb{S}_{1}$. Moreover $\mathcal{F}(\cdot ,\mu
,\mathbf{u})$ is
equivariant under translations of the plane:%
\begin{equation*}
\mathcal{T}_{\mathbf{v}}\mathcal{F}(U,\mu ,\mathbf{u})=\mathcal{F}(\mathcal{T%
}_{\mathbf{v}}U,\mu ,\mathbf{u})
\end{equation*}%
where
\begin{equation*}
\mathcal{T}_{\mathbf{v}}U(X)=U(X+\mathbf{v}).
\end{equation*}%
In addition, there is $M_{3}>0,$ such that for $||U||_{3}\leq M_{3}$ and $%
|\mu |\leq M_{3},$ $\mathcal{F}$ satisfies for any $m\geq 3$ the
"tame"
estimate%
\begin{equation}
||\mathcal{F}(U,\mu ,\mathbf{u})||_{m-1}\leq
c_{m}(M_{3})||U||_{m}, \label{tame F}
\end{equation}%
where $c_{m}$ only depends on $m$ and $M_{3}.$
\end{lemma}

\begin{proof}
The $C^{\infty }$ smoothness of $(\psi ,\eta )\mapsto
\mathcal{G}_{\eta
}(\psi ):\mathbb{H}^{m}(%
\mathbb{R}
^{2}/\Gamma )\rightarrow H^{m-1}(%
\mathbb{R}
^{2}/\Gamma )$ comes from the study of the Dirichlet-Neumann operator, see (%
\ref{equtheta},\ref{Gdemi-plan}), and the properties of elliptic
operators.
This result is proved in particular by Craig and Nicholls in \cite%
{Craig-Nicholls2000}, and by D.Lannes in \cite{Lannes}. Notice that $H^{s}(%
\mathbb{R}
^{2}/\Gamma )$ is an algebra for $s>1$. Notice that it is proved
by Craig et
al \cite{Craig-Sch-Sul} that the mapping $(\psi ,\eta )\mapsto \mathcal{G}%
_{\eta }(\psi ):H^{m}(%
\mathbb{R}
^{2}/\Gamma )\times C^{m}(%
\mathbb{R}
^{2}/\Gamma )\rightarrow H^{m-1}(%
\mathbb{R}
^{2}/\Gamma )$ is analytic and the authors give the explicit
Taylor expansion near 0, with the same type of "tame" estimates
that we shall use in the following sections. We choose here to
stay with $(\psi ,\eta )\in
\mathbb{H}^{m}(%
\mathbb{R}
^{2}/\Gamma )$ and we just use the $C^{\infty }$ smoothness of the
mapping, in addition to the tame estimates (see \cite{Lannes}).

The equivariance of $\mathcal{F}$ under translations of the plane
is obvious.

We refer to \cite{Lannes} for the proof of the following "tame"
estimate, valid for any $k\geq 1$ (here simpler than in
\cite{Lannes} since we have periodic functions and since there is
no bottom wall)
\begin{equation}
||\mathcal{G}_{\eta }(\psi )||_{k}\leq c_{k}(||\eta
||_{3})\{||\eta ||_{k+1}||\psi ||_{3}+||\psi ||_{k+1}\},
\label{tame G}
\end{equation}%
necessary to get estimate (\ref{tame F}).
\end{proof}

\subsection{Results}

We are now in a position to formulate the main result of this
paper on the
existence of non-linear diamond waves satisfying operator equation %
\eqref{F(U,mu)1}. We find an explicit solution to \eqref{F(U,mu)1}
in the vicinity of an approximate solution $U_{\varepsilon
}^{(N)}$ which existence is stated in the following lemma
restricted to "diamond waves", i.e to solutions belonging to the
important subspace (still with $\Gamma
^{\prime }$ generated by $(1,\pm \tau )$)%
\begin{equation*}
\mathbb{H}_{(S)}^{k}=\{U=(\psi ,\eta )\in \mathbb{H}^{k}(%
\mathbb{R}
^{2}/\Gamma ):\psi \text{ odd in }x_{1},\text{ even in
}x_{2},\text{ }\eta \text{ even in }x_{1}\text{ and in }x_{2}\}.
\end{equation*}%
For these solutions the unit vector $\mathbf{u}_{0}=(1,0)$ is
fixed (see a more general statement at Theorem \ref{Lembifurc},
with non necessarily symmetric formal solutions).

\begin{lemma}
\label{Lembifurc0} Let $N\geq 3$ be an arbitrary positive number
and the critical value of parameter $\mu _{c}(\tau )=(1+\tau
^{2})^{-1/2}$ is such that the dispersion equation $n^{2}+\tau
^{2}m^{2}=\mu _{c}^{-2}n^{4}$ has only the solution $(n,m)=(1,1)$
in the circle $m^{2}+n^{2}\leq N^{2}$. Then approximate
3-dimensional diamond waves are given by
\begin{gather}
U_{\varepsilon }^{(N)}=(\psi ,\eta )_{\varepsilon
}^{(N)}=\sum_{1\leq p\leq N}\varepsilon ^{p}U^{(p)}\in
\mathbb{H}_{(S)}^{k},\text{ for any }k,
\label{UN} \\
U^{(1)}=(\sin x_{1}\cos \tau x_{2},\frac{-1}{\mu _{c}}\cos
x_{1}\cos \tau
x_{2}),\quad \mu _{\varepsilon }^{(N)}=\mu _{c}+\tilde{\mu},\quad \tilde{\mu}%
=\mu _{1}\varepsilon ^{2}+O(\varepsilon ^{4}),  \notag
\end{gather}%
where
\begin{equation*}
\mu _{1}=\big(\frac{1}{4\mu _{c}^{3}}-\frac{1}{2\mu
_{c}^{2}}-\frac{3}{4\mu _{c}}+2+\frac{\mu
_{c}}{2}-\frac{9}{4(2-\mu _{c})}\big),
\end{equation*}%
and where for any $k$,
\begin{equation*}
\mathcal{F}(U_{\varepsilon }^{(N)},\mu _{\varepsilon }^{(N)},\mathbf{u}%
_{0})=\varepsilon ^{N+1}Q_{\varepsilon },
\end{equation*}%
\ $Q_{\varepsilon }$ uniformly bounded in $\mathbb{H}_{(S)}^{k},$
with respect to $\varepsilon $. There is one critical value $\tau
_{c}$ of $\tau $ such that $\mu _{1}(\tau _{c})=0,$ and $\mu
_{1}(\tau )<0$ for $\tau <\tau _{c},$ $\mu _{1}(\tau )>0$ for
$\tau >\tau _{c}.$
\end{lemma}

\begin{proof}
The lemma is a particular case of the general Theorem
\ref{Lembifurc} in the symmetric case.
\end{proof}

The following theorem on existence of $3D$-diamond waves is the
main result of the paper (notice that $\tau =\tan \theta )$

\begin{theorem}
\label{Thmexistence0} Let us choose arbitrary integers $l\geq 23$,
$N\geq 3$ and a real number $\delta <1$. Assume that
\begin{equation*}
\tau \in (\delta ,1/\delta ),\quad \mu _{c}=(1+\tau ^{2})^{-1/2}.
\end{equation*}%
Then there is a set $\mathfrak{N}$ of full measure in (0,1) with
the following property. If $\mu _{c}\in \mathfrak{N}$ and $\tau
\neq \tau _{c}$, then there exists a positive $\varepsilon
_{0}=\varepsilon _{0}(\mu _{c},N,l,\delta )$ and a set
$\mathcal{E}=\mathcal{E}(\mu _{c},N,l,\delta )$ so that
\begin{equation*}
\lim\limits_{\varepsilon \rightarrow 0}\frac{2}{\varepsilon ^{2}}%
\int\limits_{\mathcal{E}\cap (0,\varepsilon )}s\,ds=1,
\end{equation*}%
and for every $\mu =\mu _{\varepsilon }^{(N)}$ with $\varepsilon \in \mathcal{E%
}$, equation \eqref{F(U,mu)1} has a "diamond wave" type solution $%
U=U_{\varepsilon }^{(N)}+\varepsilon ^{N}W_{\varepsilon }$ with $%
W_{\varepsilon }\in \mathbb{H}_{(S)}^{l}$. Moreover, $W:\mathcal{E}%
\rightarrow \mathbb{H}_{(S)}^{l}$ is a Lipschitz function cancelling at $%
\varepsilon =0$, and for $\tau <\tau _{c}$ (resp. $\tau >\tau
_{c}$), and when $\varepsilon $ varies in $\mathcal{E}$, the
parameter $\mu =\mu _{\varepsilon }^{(N)}$ runs over a measurable
set of the interval $(\mu _{\varepsilon _{0}}^{(N)},\mu _{c})$
(resp. $(\mu _{c},\mu _{\varepsilon _{0}}^{(N)})$ of
asymptotically full measure near $\mu _{c}$.
\end{theorem}

\begin{figure}[h]
\begin{center}
\includegraphics[height=1.90in]{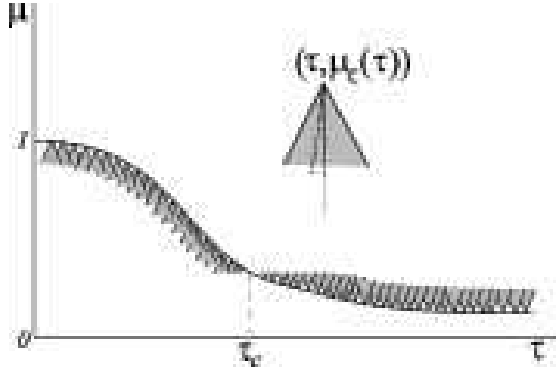}
\end{center}
\caption{Small sectors where 3-dimensional waves bifurcate. Their vertices
lie on the critical curve $\protect\mu=\protect\mu_c(\protect\tau)$. The
good set of points is asymptotically of full measure at the vertex on each
half line (see the detail above). In the paper we only give the proof for
each half line $\protect\tau=const$ (dashed line on the figure) }
\label{fig2}
\end{figure}

We can roughly express our result in considering the
two-dimensional parameter plane $(\tau ,\mu )$ where $\tau =\tan
\theta ,$ $2\theta $ being the angle between the two basic wave
vectors of same length generating the two-dimensional lattice
$\Gamma ^{\prime }$ dual of the lattice $\Gamma $ of
periods for the waves. The critical value $\mu _{c}(\tau )$ of $\mu $ $%
(=gL/c^{2})$, where $\mu _{c}(\tau )=(1+\tau ^{2})^{-1/2}=\cos
\theta $, corresponds to the solutions of the dispersion relation
we consider here (in particular 3-dimensional diamond waves
propagating in the direction of the
bisectrix of the wave vectors). We show that for $\tau <\tau _{c}$ ($%
\approx 2.48)$ the bifurcating (diamond) waves of size $O(|\mu
-\mu _{c}(\tau )|^{1/2})$ occur for $\mu <\mu _{c}(\tau ),$ while
for $\tau >\tau _{c}$ it occurs for $\mu >\mu _{c}(\tau ).$ We
prove that bifurcation of these 3-dimensional waves occurs on half
lines $\tau =const$ of the plane, with their origin on the
critical curve, for "good" values of $\tau $ (which appear to be
nearly all values of $\tau ).$ Moreover, we prove that on each
half line, these waves exist for "good" values of $\mu ,$ this set
of "good" values being asymptotically of full measure at the
bifurcation point $\mu =\mu _{c}(\tau )$ (see Figure \ref{fig2}).

Another way to describe our result is in terms of a bifurcation
from a non isolated eigenvalue in the spectrum of the linearized
operator at the origin. Indeed, for our critical values $(\tau
,\mu _{c}(\tau ))$ of the parameter, the differential at the
origin is a selfadjoint operator with in general a non isolated 0
eigenvalue (see Theorem \ref{141l}). Our result means that from
each point $(\tau ,\mu _{c}(\tau ))$ where $\tau $ is chosen in a
full measure set of $(0,\infty )$, a branch of solutions
bifurcates in the following sense. Every half line \ $\tau =const$
with origin at the point $(\tau ,\mu _{c}(\tau ))$ and on the good
side of this curve, contains a measurable set of points where the
bi-periodic gravity waves exist, with an amplitude $O(|\mu -\mu
_{c}|^{1/2})$, this set being asymptotically of full measure near
$\mu _{c}(\tau )$.

In fact, we can improve our result in replacing the half lines
mentioned above, by \emph{small sectors centered on these half
lines}. Each half line in each sector, with origin at the vertex
of the sector, contains a
measurable set where the bi-periodic gravity waves exist, with an amplitude $%
O(|\mu -\mu _{c}|^{1/2})$, this set being asymptotically of full
measure near $\mu _{c}(\tau ).$ The proof of such a result
introduces many technicalities, which are not essential for the
understanding of the paper. This complication is mainly due to the
fact that we then need to work with a lattice $\Gamma $ now
depending on $\varepsilon .$ We just mention in various places
what is really needed for such an extension of the result proved
here.

\subsection{Mathematical background}
There are some aspects of our method which deserve brief mention.
First we use  the Nash-Moser method, which is now an integral part
of nonlinear analysis \cite{Deimling}, for proving Theorem
\ref{Thmexistence0}. The crucial point for the Nash-Moser method
is to
 obtain a priori bounds on an approximate
right-inverse of the partial derivative
$\partial_U\mathcal{F}(U,\mu,\mathbf u)$. As it is shown in
Section 2, this problem is equivalent to the problem of
invertibility of a second-order selfadjoint pseudodifferential
operator with multiple characteristics. Second, we use the Moser
theory of foliation on a torus \cite{Moser} and the
\emph{invariant parametric representation} of the
Dirichlet-Neumann operator to reduce the linearized  equation  to
a canonical form with constant coefficients in the principal part.
We employ a modification of the Weil Theorem \cite{Weil} on
uniform distribution of numbers $\{\omega n^2\}$ modulo $1$ to
deduce the effective estimates of \emph{small divisors}, and hence
to prove the invertibility of the principal part of the linearized
operator. The most essential ingredient of our approach is the
algebraic \emph{descent method} \cite{IPT}, \cite{Plot-Tol} which
allows to reduce the canonical pseudodifferential equation on
$2$-dimensional torus to a Fredholm-type equation.

\subsection{Structure of the paper}

Now we can explain the organization of the paper. In section 2, we
prove Theorem \ref{Lembifurc} which establishes the existence of
approximate solutions under the form of power series of the
amplitudes of the two incident mono-periodic travelling waves,
corresponding to symmetric basic wave vectors. \emph{The parameter
is two-dimensional here}, due to the freedom in the direction of
propagation of the three-dimensional wave. To show this result, we
use a formal Lyapunov-Schmidt technique, assuming that the angle
$2\theta $ between the two basic wave vectors satisfies that $\tau
=\tan \theta $ is such that the equation for positive integers
$(n,m)$
\begin{equation*}
n^{2}+\tau ^{2}m^{2}=n^{4}(1+\tau ^{2})
\end{equation*}%
has the unique solution $(n,m)=(1,1)$ (non resonance property). In
playing with scales and parameters, this condition is not
restrictive among non resonant situations, which indeed represent
the general case. In such a case, the kernel of the linearized
operator at rest state (taken as the origin) is four-dimensional,
and in using extensively the symmetries of the system
(\ref{F(U,mu)1}), we obtain, for a fixed value of the bifurcation
parameter $(\mu ,\mathbf{u}),$ doubly-periodic formal travelling
gravity waves propagating in the direction $\mathbf{u}$. Limiting
cases are the mono-periodic travelling waves corresponding to one
of the basic wave vectors. The Lemma \ref{Lembifurc0} is a
particular case of the above theorem.

From now on, we restrict the study to solutions called "diamond
waves", which are \emph{symmetric with respect to the direction of
propagation},
here the $x_{1}-$ axis. In section 3 we consider the linear operator $%
\mathcal{L}(U,\mu )$ corresponding to the differential of
(\ref{F(U,mu)1}) at a non zero point in $\mathbb{H}_{(S)}^{k}$,
which we need to invert for using the Nash-Moser theorem. The
principal part of this operator is the symmetric sum
$$
-\mathcal{J}^{\ast }(\frac{1}{\mathfrak{a}}\mathcal{J}\cdot)+
\mathcal{G}_{\eta }, \quad \mathcal{J}=V\cdot \nabla,
$$
of a second order derivative in the direction of a periodic vector
field $V(X)$, and of the Dirichlet-Neumann operator which is
integro-differential of first order, both parts depending
periodically on coordinates. More precisely, $V=\mathbb
G^{-1}(X)(\mathbf u_0+\nabla\psi(X))$, where $\mathbb G(X)\,
dX\cdot dX$ is the first fundamental form of the free surface.
Recall that $\mathbb G$ is a covariant tensor field on $\Sigma$,
and
 for the standard
parametrization  $x_3=\eta(X)$, it is given by $\mathbb
G(X)=1+\nabla \eta\otimes \nabla\eta$. It follows from the
kinematic condition \eqref{equation2} that integral curves of the
vector field $V(X)$ coincide with trajectories of liquid particles
moving along $\Sigma$ and submitted to the vertical gravity $\mu$.

  Section 3 is concerned with the first step
of the long way towards the inversion of $\mathcal{L},$ which
consists in finding a diffeomorphism of the torus for which the
highest order terms of the operator $\mathcal{L}$ become constants
(depending on the linearization point). We begin (Lemma
\ref{changevar1}) with the construction of a diffeomorphism which
takes integral curves of the vector field $V$ onto straight lines
parallel to the abscissa axis. Being endowed with the Jacobi
metric $ ds^2=\big(1/2-\mu \eta(X)\big) \mathbb G(X)\, dX\cdot dX$
the free surface becomes a Riemannian manifold on which the
integral curves of $V$ coincide with geodesics (see Appendix
\ref{particles}). Hence, by Lemma \ref{changevar1}, they form a
geodesic foliation on $\Sigma$. Moreover, since the distance
between each of these curves and the abscissa axis is finite, the
foliation has a zero rotation number. It is at this point where
the restriction to symmetric solutions (diamond waves) is
necessary, since we don't know yet how to manage such a
diffeomorphism in the non symmetric case, see \cite{Moser} for
discussion. Recall that the Moser Theorem \cite{Moser} guarantees
the existence of at least one geodesic for any given rotation
number.

The second  result of Section 3 is Theorem \ref{thmDir-Neum} which
gives the parametric representation of the Dirichlet-Neumann
operator in
 arbitrary  coordinates $Y$ on $\Sigma$ so that a mapping
$X=X(Y)$ is a diffeomorphism of a torus. It follows from this
theorem that for any smooth periodic function $u(Y)$ and
$\check{u}(X)=u(Y(X)$, the Dirichlet-Neumann operator has the
decomposition
$$
\mathcal{G}_\eta \check{u}=\check{\mathcal G_1u}+ \check{\mathcal
G_0u}+\check{\mathcal G_{-1}u},
$$
in which $\mathcal G_i$ are pseudodifferential operators of order
$i$. We give an explicit expression for their symbols  in terms of
the first fundamental form  and  the principal curvatures of the
free surface. In particular, we show that, up to a positive
invariant multiplier, the symbol of the operator $\mathcal G_1$ is
equal to $\sqrt{\mathbb G(Y)^{-1} k\cdot k}$, and the real part of
the symbol of $\mathcal G_0$ coincides with the difference between
the sum of the principal curvatures and the normal curvature of
$\Sigma$ in the direction of $\mathbb G^{-1}k$. This leads to the
interesting conclusion: \emph{ the manifold $\Sigma$ is defined by
its Dirichlet-Neumann operator up to translation and rotation  of
the embedding space}.

Combining  Lemma  \ref{changevar1} and Theorem \ref{thmDir-Neum}
gives  the main result of  Section 3 -- Theorem
\ref{thmChangeCoord}. This theorem   ensures the existence of a
diffeomorphism  $X=X(Y)$ of the $2$-dimensional torus, which
brings the linearized operator to the  canonical form
$$
\mathfrak L+\mathfrak H=\mathfrak L+\mathfrak
A\partial_{y_1}+\mathfrak B+\mathfrak L_{-1},
$$
where the remainder $\mathfrak L_{-1}$ is of order $-1$,
$\mathfrak A$ and $\mathfrak B$ are zero-order pseudodifferential
operators, and the principal part
$$\mathfrak L=\nu
\partial_{y_1}^2+(-\Delta)^{1/2}, \quad
\Delta=\partial_{y_1}^2+\tau^2 \partial_{y_2}^2$$ is a selfadjoint
pseudodifferential operator. Here  the parameter  $\nu$ depends on
the point of linearization,  with
$\nu(0)=\nu_0=\mu_{c}(\tau)^{-1}$ (Lemma \ref{Lemcoef-nu}).

In Section 4
 we study the operator $\mathfrak L$ in many details, and give estimates on its
 resolvent in Sobolev spaces of bi-periodic functions which are odd in $y_1$,
  and even in $y_2$. We begin with the observation that
for $\nu_0=\mu_c(\tau)^{-1}$ and almost every positive $\tau$,
zero is a simple eigenvalue of the operator $\mathfrak L_0= \nu_0
\partial_{y_1}^2+(-\Delta)^{1/2}$ and
\begin{equation*}
\|\mathfrak L_0^{-1} u\|_s\leq
c(\tau,\alpha)\|u\|_{s+(1+\alpha)/2}
\end{equation*}
for all $u$ orthogonal to the kernel of $\mathfrak L_0$ and
$\alpha>0$. Next we study the perturbation of its resolvent
assuming that $\nu=\nu_0-\varepsilon^2\nu_1+O(\varepsilon^3)$ and
with a spectral parameter $\varkappa=O(\varepsilon^2)$, both being
Lipschitz functions of a small parameter $\varepsilon$. Here we
have a \emph{ small divisor} problem, and we meet  the necessity
to restrict the parameter values to "good ones", for being able to
find suitable estimates. Calculations (Lemma \ref{145l}) show that
the resolvent of $\mathfrak L$ satisfies the estimate
\begin{equation}\label{smallresolv}
\|(\mathfrak L-\varkappa)^{-1} u\|_s\leq c\|u\|_{s+1}, \quad
u\in(\text{~ker~} \mathfrak L_0)^{\perp},
\end{equation}
if  parameters $\nu$ and $\varkappa$ satisfy the quadratic
Diophantine inequalities
\begin{equation}\label{polinomialappr}
|\omega n^2-m-C|\geq c n^{-2}\text{~~for all positive
integers~~}n,m,
\end{equation}
where  $ \omega =\nu\tau^{-1}$ and $ C=(2\nu \tau )^{-1}-\varkappa
\tau^{-1}$. Note that there is a  difference between linear and
polynomial Diophantine approximations: in classic theory of linear
Diophantine forms, see \cite{Cassels} for general theory, the
integers for which "small divisors" are really small, form a
sparse set in the  integral lattice. This property was used in
pioneering works of  Siegel \cite{Siegel} and in the Arnold proof
of the Kolmogorov Theorem \cite{Arnold}. In contrast to the linear
case, the couples $(m,n)$, for which the left hand side of
inequality \eqref{polinomialappr} is small, can form clusters in
$\ZZ^2$, and the problem of obtaining   small divisors estimates
becomes more complicated. It turns out that the validity of
inequalities \eqref{polinomialappr} with a constant $c$
independent of the small parameter  is a consequence of estimate
\begin{equation}\label{numbers2}
N^{-1}\text{~card~}\{n:\, \omega_0 n^2\text{~modulo~} 1\leq
\varepsilon\text{~~and~~} 1\leq n\leq N\}\leq
c\varepsilon\text{~~for all~~}N\geq \varepsilon^{-\lambda}.
\end{equation}
 Recall that, by the Weil Theorem \cite{Weil}, \cite{Cassels},
for each fixed $\varepsilon$, the left hand side tends to
$\varepsilon$ as $N\to \infty$.  Hence the inequality holds true
for some $c$ depending on $\varepsilon$. In Appendix \ref{numbers}
we make this result more precise and prove the existence of
absolute constant $c$ such that inequality \eqref{numbers2} is
fulfilled for all $\lambda \geq 78$ and all intervals of length
$\varepsilon$. This leads to the main result of this section --
Theorem \ref{142t}, which shows that with a suitable choice of the
parameters, the resolvent operator provides a loss of one in the
degree of differentiability. Moreover, estimate
\eqref{smallresolv} holds true for all $\varepsilon^2$ in an
asymptotically full measure set on every half line $\theta =const$
of the parameter plane, the origin of which being chosen
arbitrarily in a full measure set, on the bifurcation curve $\mu
=\cos \theta $.

In Section 5 we take into account all remaining terms of the linear operator $%
\mathfrak L+\mathfrak{H}$ and  prove its invertibility with a loss
of differentiability.  The main difficulty  is that the operator
$\mathfrak L+\mathfrak H$ involves the principal part $\mathfrak
L$, which inverse   is unbounded, and arbitrary operators
$\mathfrak A$, $\mathfrak B$ with "variable coefficients".

Most, if not all, existing results related to such problems were
obtained by use of  \emph{the Fr\" ohlich-Spencer method} proposed
in \cite{FS83}, cf \cite{PIP1,PIP3}, and developed  by Craig and
Wayne  \cite{CW93,C00} and Bourgain  \cite{B95,B98}. The basic
idea of the method is  a representation of operators in the form
of infinite matrices with elements labelled by some lattice and
block decompositions of this lattice. Let us use the operator
$\mathfrak L+\mathfrak H$ to illustrate the main features of this
method. First we have to replace a periodic function $u$ by the
sequence of its Fourier coefficients $\{\widehat u(k)\}$, $k\in
\Gamma'$, and the operator $\mathfrak L$ by the diagonal matrix
with the elements $L(k)$. Then we have to split the lattice
$\Gamma'$ into a "regular" part which consists of all $k$ with
"large" $L(k)$, and an "irregular" part  which includes all $k$
corresponding to "small" values of $L(k)$. Using the contraction
mapping principle we can eliminate the "regular" component and
reduce the inversion of
 $\mathfrak L+\mathfrak H$  to the  inversion of an
infinite matrix on the "irregular" subspace. The existence of an
inverse to this matrix is established by using a special iteration
process which is the core of the  method. Note that the Fr\"
ohlich-Spencer method is working in our case only if $\mathfrak
A=0$.

Our approach is  based on \emph{the descent method} which was
proposed in \cite{Plot-Tol,IPT} and  dates back to the classic
Floquet-Lyapunov theory. The descent method is a pure algebraic
procedure which brings the canonical operator to an operator with
constant coefficients and does not depend on the structure and
spectral properties of the principal part $\mathfrak L$. The heart
of the method is the following identity (Theorem \ref{descent1})
\begin{equation*}
\big(\mathfrak{L}+\mathfrak{AD}_{1}+\mathfrak{B}\big)(1+\mathfrak{C})u=(1+\mathfrak{E})
(\mathfrak{L}-\varkappa )u+\mathfrak{F}u,
\end{equation*}%
which holds true for all functions $u\in H^2(\RR^2/\Gamma)$ odd in
$y_1$. Here $\mathfrak C$ and $\mathfrak E$ are bounded operators
in the Sobolev spaces of periodic functions $H^s(\RR^2/\Gamma)$;
the remainder $\mathfrak F$  is a bounded operator
$:H^{s-1}(\RR^2/\Gamma)\mapsto H^{s}(\RR^2/\Gamma)$; the  Floquet
exponent $\varkappa$ has an explicit expression in terms of
operators $\mathfrak A$ and $\mathfrak B$. Moreover, if
$\|\mathfrak A\|,\|\mathfrak B\|\sim \varepsilon$, then
$\|\mathfrak C\|,\|\mathfrak E\|,\|\mathfrak F\|\sim \varepsilon$
and $\varkappa=O(\varepsilon^2)$. The proof of these results
constitutes Section \ref{descentmethod} and Appendix
\ref{pseudodifferential}. The technique used is more general than
the one used in \cite{IPT}, since we use here general properties
of pseudodifferential operators, however taking into account of
the symmetry properties of $\mathcal{L}.$

The descent method of algebraic character presented here, might be
easily used for example on the one-dimensional KDV and
Schr\"odinger equations, avoiding the heavy technicalities of the
Fr\" ohlich-Spencer method.

Thus we reduce the problem of the inversion of the canonical
operator $\mathfrak L+\mathfrak H$ to the problem of the inversion
of operator $ \mathfrak L-\varkappa +\tilde{\mathfrak L}_{-1}, $
where  $\tilde{\mathfrak L}_{-1}$ is a smoothing remainder. It is
then possible to use the result of section 4 for inverting the
full operator and to prove  Theorem \ref{t121} -- the main result
on the existence and estimates of  $(\mathfrak L+\mathfrak
H)^{-1}$. In particular, this theorem implies that if $\mathfrak
A$ and $\mathfrak B$ are Lipschitz operator-valued functions of a
small parameter $\varepsilon$, which vanish for $\varepsilon=0$
and satisfy symmetry and metric conditions (Section 5), and  if
$\mathfrak L-\varkappa$ meets all requirements of Theorem
\ref{142t}, then for all $\varepsilon^2$ taken in an
asymptotically full measure set, the resolvent has the
representation
$$
(\mathfrak L+\mathfrak H)^{-1}(\varepsilon)=\frac{1}{\mathfrak c}
\mathfrak H_0(\varepsilon)+ \mathfrak H_1(\varepsilon),
$$
in which  operators $\mathfrak H_1(\varepsilon):
H^{s}(\RR^2/\Gamma)\mapsto H^{s-1}(\RR^2/\Gamma)$ are uniformly
bounded in $\varepsilon$, and $\mathfrak H_0(\varepsilon)$ are
bounded operators of  rank 1,  the coefficient $\mathfrak
c=\varepsilon^2 @+O(\varepsilon^3) $ being given by
\eqref{125aa}. We show at the end of the section (see Theorem
\ref{inv L}), that the results apply to the linear operator $
\mathcal{L}(U,\mu )$ corresponding to the differential of
(\ref{F(U,mu)1}) at a non zero point in $\mathbb{H}_{(S)}^{k}$. In
particular, we give the sufficient conditions which provide the
existence of the bounded inverse $\mathcal{L}(U,\mu )^{-1}:
H^{s+3}(\RR/\Gamma)\mapsto H^{s}(\RR/\Gamma)$.

Section 6 applies extensively the result proved in \cite{IPT}
concerning the Nash-Moser theorem with parameters in a Cantor set.
The main result, which is the main result of the paper is Theorem
\ref{Thmexistence0} establishing the existence of smooth
bi-periodic travelling gravity waves symmetric with respect to the
direction of propagation, in the region of the parameter space
mentioned above. Notice that, a part from the last section, which
heavily rests upon the self contained Appendix N of \cite{IPT},
the rest of the paper is self contained, with some details of
computations and basics on pseudodifferential operators put in
Appendix, for providing an easy reading.

\section{Formal solutions}

\subsection{Differential of $\mathcal{G}_{\protect\eta }$}

\label{Dir The difficulty here, is that we need to play on two
parameters, taking into account that the small divisor problem
already occurs at criticality (origin of
the half line).-Neum}In this subsection we study the structure of the operator $%
\mathcal{G}_{\eta }$, and we give useful formulas and estimates.

The following regularity property holds

\begin{lemma}
\label{Lemma diff G}The differential $h\mapsto \partial _{\eta }\mathcal{G}%
_{\eta }[h]$ of $G_{\eta }$ satisfies for $\eta ,$ $\psi ,$ $h$
smooth
enough bi-periodic functions%
\begin{eqnarray}
\partial _{\eta }\mathcal{G}_{\eta }[h]\psi &=&-\mathcal{G}_{\eta }(h\zeta
)+\nabla _{X}\cdot \{(\zeta \nabla _{X}\eta -\nabla _{X}\psi )h\},
\label{diffG_2} \\
\zeta &=&\frac{1}{1+(\nabla _{X}\eta )^{2}}\{\mathcal{G}_{\eta
}\psi +\nabla _{X}\eta \cdot \nabla _{X}\psi \}.  \label{diffG_3}
\end{eqnarray}%
Moreover, despite of the apparent loss of derivatives for $(\psi
,\eta )$ in (\ref{diffG_2},\ref{diffG_3}) we have for
$||U||_{3}\leq M_{3},$ the
following tame estimate%
\begin{equation*}
||\partial _{\eta }\mathcal{G}_{\eta }[h]\psi ||_{k}\leq
c_{k}(M_{3})\{||h||_{k+1}+||U||_{k+1}||h||_{2}\}.
\end{equation*}
\end{lemma}

\begin{proof}
We refer to Appendix \ref{a0} for the formula
(\ref{diffG_2},\ref{diffG_3})
already proved for instance in \cite{Lannes}, and we also refer to \cite%
{Lannes} for the tame estimate.
\end{proof}

From the formulas of the above Lemma \ref{Lemma diff G}, we are
now able to
compute successive derivatives of $\mathcal{F}$. Observe that in (\ref%
{diffG_2},\ref{diffG_3}) there is a loss of two derivatives for
$(\psi ,\eta ).$ In fact there is a compensation cancelling the
dependence into the second order derivatives and we have the
following Lemma which completes Lemma \ref{F(U,mu)}:

\begin{lemma}
\label{regul F}For $||U||_{3}\leq M_{3},$ and $|\mu |\leq M_{3}$
the following tame estimates hold (and analogous ones for higher
order
derivatives)%
\begin{eqnarray*}
||\partial _{U}\mathcal{F}(U,\mu ,\mathbf{u})[\delta U]||_{k}
&\leq
&c_{k}(M_{3})\{||\delta U||_{k+1}+||U||_{k+1}||\delta U||_{3}\}, \\
||\partial _{UU}^{2}\mathcal{F}(U,\mu ,\mathbf{u})[\delta
U_{1},\delta U_{2}]||_{k} &\leq &c_{k}(M_{3})\{||\delta
U_{1}||_{k+1}||\delta U_{2}||_{3}+
\\
&&+||\delta U_{2}||_{k+1}||\delta U_{1}||_{3}+||U||_{k+1}||\delta
U_{1}||_{3}||\delta U_{2}||_{3}\}.
\end{eqnarray*}
\end{lemma}

\subsection{Linearized equations at the origin and dispersion relation}

\label{lineqorigin}

The linearization at the origin of system (\ref{basic1}),
(\ref{basic2})
leads to%
\begin{eqnarray}
\mathcal{G}^{(0)}(\psi )-\mathbf{u}\cdot \nabla _{X}\eta &=&0,
\label{linsyst1} \\
\mathbf{u}\cdot \nabla _{X}\psi +\mu \eta &=&0,  \label{linsyst2}
\end{eqnarray}%
where the following operator
\begin{equation*}
\mathcal{G}^{(0)}=(-\Delta )^{1/2}
\end{equation*}%
is defined more precisely in Appendix \ref{a01}. Now expanding in Fourier series, we have%
\begin{equation*}
\psi (X)=\sum_{K\in \Gamma ^{\prime }}\psi _{K}e^{iK\cdot
X},\text{ \ }\eta (X)=\sum_{K\in \Gamma ^{\prime }}\eta
_{K}e^{iK\cdot X},
\end{equation*}%
and (\ref{linsyst1}), (\ref{linsyst2}) give for any $K\in \Gamma
^{\prime }$
\begin{eqnarray*}
|K|\psi _{K}-i(K\cdot \mathbf{u})\eta _{K} &=&0, \\
i(K\cdot \mathbf{u})\psi _{K}+\mu \eta _{K} &=&0.
\end{eqnarray*}%
Hence, the \emph{dispersion relation} reads%
\begin{equation}
\Delta (K,\mu ,\mathbf{u})\overset{def}{=}\mu |K|-(K\cdot
\mathbf{u})^{2}=0. \label{dispersion}
\end{equation}

The point now to discuss is the number of solutions $K\in \Gamma
^{\prime }$ of (\ref{dispersion}), for a fixed vector
$\mathbf{u}\in \mathbb{S}_{1},$
and a fixed parameter $\mu .$ We restrict our analysis to a \emph{lattice }$%
\Gamma ^{\prime }$\emph{\ generated by two vectors }$K_{1}$\emph{\ and }$%
K_{2}$\emph{\ symmetric with respect to the }$x_{1}-$\emph{\ axis,
taken in the direction of }$\mathbf{u}$, which is the situation if
one is looking for
\emph{short crested waves}:%
\begin{equation*}
K_{1} =(1,\tau ) , \text{ \ } K_{2} =(1,-\tau )
\end{equation*}%
where $\tau $ is positive. When $\tau $ is small, the lattice
$\Gamma $ of periods is formed with diamonds elongated in the
$x_{2}$ direction (see Figure \ref{fig1}). Taking 1 for the first
component of $K_{1}$ implies that we choose the length scale $L$
as the wave length in the $x_{1}-$ direction divided by $2\pi .$

We consider in what follows, the cases when the direction
$\mathbf{u}_{0}$ of the travelling waves at criticality is the
$x_{1}-$ axis, and the critical parameter $\mu _{c}=(1+\tau
^{2})^{-1/2}$ is such that the equation for $(m_{1},m_{2})\in
\mathbb{N}
^{2}$
\begin{equation}
\mu _{c}\sqrt{m_{1}^{2}+\tau ^{2}m_{2}^{2}}-m_{1}^{2}=0
\label{disp0}
\end{equation}%
has only the solution
\begin{equation*}
(m_{1},m_{2})=(1,1).
\end{equation*}%
In case we have a solution $(m_{1},m_{2})\neq (1,1),$ one can make
the
change $(\mu _{c},\tau )\mapsto (\frac{\mu _{0}}{m_{1}},\tau \frac{m_{2}}{%
m_{1}})$ to recover the case we study here. Moreover, changing $\mu $ into $%
\mu /m_{1}$ corresponds to changing the length scale $L$ into
$L/m_{1}$ which indeed corresponds to the new wave length in the
$x_{1}$ direction. So, it is clear that \emph{we do not restrict
the generality} in choosing the case of a solution
$(m_{1},m_{2})=(1,1).$

Notice that for any integer $l,$ when $\tau =l$ or $1/l,$ there is
an infinite number of solutions $(m_{1},m_{2})$ of (\ref{disp0}),
hence we need to avoid such choices for $\tau .$

\textbf{Remark.} We notice here the fundamental difference between
the present type of study and the works using spatial dynamics for
finding travelling waves, as for instance Groves and Haragus in
\cite{Grov-Harag}. Their study only consider cases with surface
tension, and cannot work without surface tension, since this would
lead to an infinite set of imaginary eigenvalues $\pm im_{1}$
(hence preventing the use of center manifold reduction), with no
restriction for $m_{1}$ to be an integer, while $m_{2}\in
\mathbb{N}
$ (this corresponds to fixing the length scale with the period in
$x_{2},$ transverse to the direction of the travelling waves.

\subsection{Formal computation of 3-dimensional waves in the simple case}

\label{comp-Diamond}

In this subsection we make a formal bifurcation analysis for the
simple case. We denote by $\mu _{c}$ the critical value of $\mu ,$
and we denote by $\mathbf{u}_{0}=(1,0)$ the critical direction for
the waves (this direction of propagation may be changed for
bifurcating travelling waves). The
lattice $\Gamma^{^{\prime}}$is generated by the two symmetric wave vectors $%
K_{1}=(1,\tau ),$ $K_{2}=(1,-\tau ),$ where
\begin{equation*}
\mu _{c}^{-2}=1+\tau ^{2}.
\end{equation*}%
We notice that we have the following Fourier series for $U=(\psi ,\eta ):$%
\begin{equation*}
U=\sum_{n=(n_{1},n_{2})\in
\mathbb{Z}
^{2}}U_{n}e^{i(n_{1}K_{1}\cdot X+n_{2}K_{2}\cdot X)},\text{ \
}U_{n}=(\psi _{n},\eta _{n}),\text{ \ }\psi _{0}=0,
\end{equation*}%
and we notice that%
\begin{eqnarray*}
n_{1}K_{1}\cdot X+n_{2}K_{2}\cdot X &=&m_{1}x_{1}+\tau m_{2}x_{2}, \\
m_{1} &=&n_{1}+n_{2},\text{ \ }m_{2}=n_{1}-n_{2},
\end{eqnarray*}%
which gives functions which are $2\pi -$ periodic in $x_{1}$ and
$2\pi /\tau -$ periodic in $x_{2}.$

We already noticed the equivariance of system (\ref{basic1}),
(\ref{basic2}) with respect to translations of the plane,
represented \ by the linear operator $\mathcal{T}_{\mathbf{v}},$
$\mathbf{v}$ being any vector of the plane. Let us complete the
symmetry properties of our system by the symmetries
$\mathcal{S}_{0}$ and $\mathcal{S}_{1}$ defined by the
representations of respectively the symmetry with respect to 0,
and the
symmetry with respect to $x_{1}$ axis%
\begin{eqnarray}
\mathcal{S}_{0}U &=&\sum_{n=(n_{1},n_{2})\in
\mathbb{Z}
^{2}}(SU_{n})e^{-i(n_{1}K_{1}\cdot X+n_{2}K_{2}\cdot X)},\text{ \
\ }SU_{n}=(-\psi _{n},\eta _{n}),  \label{symS} \\
\mathcal{S}_{1}U &=&\sum_{n=(n_{1},n_{2})\in
\mathbb{Z}
^{2}}U_{n}e^{i(n_{1}K_{2}\cdot X+n_{2}K_{1}\cdot X)}. \label{symR}
\end{eqnarray}%
The system (\ref{basic1}), (\ref{basic2}) is equivariant, under
the symmetry $\mathcal{S}_{0}$ in all cases, while it is
equivariant under $\mathcal{S}_{1}$ only if
\begin{equation*}
\mathbf{u}\cdot K_{1}=\mathbf{u}\cdot K_{2}.
\end{equation*}%
In particular, in such a case we have%
\begin{equation*}
\mathcal{L}_{0}\mathcal{S}_{1}=\mathcal{S}_{1}\mathcal{L}_{0},\text{ \ }\mathcal{L}_{0}%
\mathcal{S}_{0}=\mathcal{S}_{0}\mathcal{L}_{0},
\end{equation*}%
where we denote by $\mathcal{L}_{0}$ the symmetric linearized operator for $%
\mu =\mu _{c}$ $\ $and $\mathbf{u}=\mathbf{u}_{0}$
\begin{equation}
\mathcal{L}_{0}=\left(
\begin{array}{cc}
\mathcal{G}^{(0)} & -\mathbf{u}_{0}\cdot \nabla \\
\mathbf{u}_{0}\cdot \nabla & \mu _{c}%
\end{array}%
\right) .  \label{L_0}
\end{equation}%
Notice that the commutation property with the linear operator $\mathcal{S}%
_{0}$ is not trivial since the choice of writing our system in the
moving frame selects the direction $\mathbf{u}$ which breaks a
reflection symmetry. Indeed the symmetry property results from the
galilean invariance of the Euler equations.

In the following Lemma, we use the two parameters $\tilde{\mu}=\mu
-\mu _{c}$
and $\mathbf{\omega }=\mathbf{u}-\mathbf{u}_{0}$ and we notice that, since $%
\mathbf{u}$ is unitary, we have%
\begin{eqnarray*}
\mathbf{\omega } &=&(\omega _{1},\omega _{2}), \\
\omega _{2} &=&\frac{1}{2\tau }\mathbf{\omega }\cdot (K_{1}-K_{2}), \\
\omega _{1} &=&-\frac{\omega _{2}^{2}}{2}+O(\omega _{2}^{4}).
\end{eqnarray*}%
In the following subsections we prove the following

\begin{theorem}
\label{Lembifurc}Assume we are in the simple case, i.e. for $\tau
$ and the critical value of the parameter $\mu _{c}(\tau )=(1+\tau
^{2})^{-1/2}<1$ such that the equation $n^{2}+\tau ^{2}m^{2}=\mu
_{c}^{-2}n^{4}$ has only the solution $(n,m)=(1,1)$ in $%
\mathbb{N}
^{2}.$ Then, for any $N\geq 1,$ and any $\mathbf{v}\in
\mathbb{R}
^{2}$, approximate 3-dimensional waves are given by $\mathcal{T}_{\mathbf{v}%
}U_{\varepsilon }^{(N)}$($T^{2}-$ torus family of solutions) where $%
K_{1}=(1,\tau ),$ $K_{2}=(1,-\tau )$ are the wave vectors, and
\begin{eqnarray*}
U_{\varepsilon }^{(N)} &=&(\psi ,\eta )_{\varepsilon
}^{(N)}=\sum_{\left( p_{1},p_{2}\right) \in
\mathbb{N}
^{2},\text{ }p_{1}+p_{2}\leq N}\varepsilon _{1}^{p_{1}}\varepsilon
_{2}^{p_{2}}U^{(p_{1},p_{2})}\in \mathbb{H}^{k},\text{ for any }k, \\
U^{(1,0)} &=&\xi _{1}=(\sin (K_{1}\cdot X),\frac{-1}{\mu _{c}}\cos
(K_{1}\cdot X)), \\
U^{(0,1)} &=&\xi _{2}=(\sin (K_{2}\cdot X),\frac{-1}{\mu _{c}}\cos
(K_{2}\cdot X)),
\end{eqnarray*}%
\begin{eqnarray*}
U^{(2,0)} &=&\left( \frac{-1}{2\mu _{c}^{2}}\sin (2K_{1}\cdot X),\frac{1}{%
2\mu _{c}^{3}}\cos (2K_{1}\cdot X)\right) , \\
U^{(0,2)} &=&\left( \frac{-1}{2\mu _{c}^{2}}\sin (2K_{2}\cdot X),\frac{1}{%
2\mu _{c}^{3}}\cos (2K_{2}\cdot X)\right) , \\
U^{(1,1)} &=&\left( \frac{1-2\mu _{c}}{\mu _{c}(2-\mu _{c})}\sin
((K_{1}+K_{2})\cdot X),\frac{\mu _{c}^{2}+2\mu _{c}-2}{\mu
_{c}^{2}(2-\mu
_{c})}\cos ((K_{1}+K_{2})\cdot X)\right) + \\
&&+\left( 0,\frac{\tau ^{2}}{\mu _{c}}\cos ((K_{1}-K_{2})\cdot
X)\right) ,
\end{eqnarray*}%
\qquad \qquad
\begin{eqnarray*}
\tilde{\mu} &=&-\frac{\mu _{c}^{2}}{8}(\alpha _{0}+\beta
_{0})(\varepsilon _{1}^{2}+\varepsilon _{2}^{2})+O\{(\varepsilon
_{1}^{2}+\varepsilon
_{2}^{2})^{2}\}, \\
\omega \cdot (K_{1}-K_{2}) &=&(\varepsilon _{1}^{2}-\varepsilon
_{2}^{2})\left( \frac{\mu _{c}}{8}(\alpha _{0}-\beta
_{0})+O\{(\varepsilon _{1}^{2}+\varepsilon _{2}^{2})\}\right) ,
\end{eqnarray*}%
with%
\begin{eqnarray*}
\alpha _{0}+\beta _{0} &=&\frac{4}{\mu _{c}^{2}}\left( -\frac{1}{\mu _{c}^{3}%
}+\frac{2}{\mu _{c}^{2}}+\frac{3}{\mu _{c}}-8-2\mu _{c}+\frac{9}{2-\mu _{c}}%
\right) \\
\beta _{0}-\alpha _{0} &=&\frac{4}{\mu _{c}^{2}}\left( -\frac{3}{\mu _{c}^{3}%
}+\frac{2}{\mu _{c}^{2}}+\frac{3}{\mu _{c}}-8-2\mu _{c}+\frac{9}{2-\mu _{c}}%
\right) ,
\end{eqnarray*}%
and where for any $N$ and $k$%
\begin{equation*}
\mathcal{F}(U_{\varepsilon }^{(N)},\mu _{c}+\tilde{\mu},\mathbf{u}_{0}+%
\mathbf{\omega })=(\varepsilon _{1}^{2}+\varepsilon _{2}^{2})^{\frac{N+1}{2}%
}Q_{\varepsilon },
\end{equation*}%
\ $Q_{\varepsilon }$ uniformly bounded in $\mathbb{H}^{k},$ with respect to $%
\varepsilon _{1},\varepsilon _{2}$. There are critical values
$\tau _{c}$ and $\tau _{c}^{\prime }$ of $\tau $ such that
$(\alpha _{0}+\beta _{0})(\tau _{c})=0,$ and $(\alpha _{0}-\beta
_{0})(\tau _{c}^{\prime })=0,$ and $\alpha _{0}+\beta _{0}$ is
positive for $\tau \in (0,\tau _{c}),$ negative for $\tau >\tau
_{c}$, $\tau _{c}\approx 2.48,$ while $\alpha _{0}-\beta _{0}$ is
negative for $\tau \in (0,\tau _{c}^{\prime }),$ positive for
$\tau >\tau _{c}^{\prime }$, $\tau _{c}^{\prime }\approx 0.504.$
Moreover, for $\varepsilon _{1}=\varepsilon _{2}$ we have "diamond
waves" where the direction of propagation $\mathbf{u}$ is along
the $x_{1}-$ axis; and for $\varepsilon _{2}=0$ (resp.
$\varepsilon _{1}=0)$ we obtain 2-dimensional travelling waves of
wave vector $K_{1}$ (resp. $K_{2}).$ All
solutions are invariant under the shift $\mathcal{T}_{\mathbf{v}_{0}}:$ $%
X\mapsto X+(\pi ,\pi /\tau )$ and $\psi _{\varepsilon }^{(N)}$ is
odd in $X,$ while $\eta _{\varepsilon }^{(N)}$ is even in $X.$
\end{theorem}

\textbf{Remark 1}: In this lemma we assume that equation (\ref{disp0}) for $%
(m_{1},m_{2})\in
\mathbb{Z}
^{2}$ has only the four solutions $(m_{1},m_{2})=(\pm 1,\pm 1),$
corresponding to the four wave vectors $K=\pm K_{1}$ and $\pm
K_{2}.$ The corresponding pattern of the waves for $\varepsilon
_{1}=\varepsilon _{2}$ is in diamond form, and for $\tau $ close
to 0, the diamonds are flattened in the $x_{1}$ direction$,$ and
elongated in $x_{2}$ direction, looking like flattened hexagons or
flattened rectangles, because of the elongated shape of crests and
troughs. This last case is indeed observed experimentally for deep
fluid layers. Nearly all (in the measure sense) values of $\tau $
are indeed such that we are in the simple case.

\textbf{Remark 2}: The above Lemma is stated differently in
Theorem 4.1 of \cite{craig}; indeed we prove here that the
manifold of solutions has a simple formulation in terms of the two
parameters.

\textbf{Remark 3:} Since $\mu =gL/c^{2},$ the result of Theorem \ref%
{Lembifurc} about the sign of $\alpha _{0}+\beta _{0}$ shows that
for $\tau <\tau _{c}$ the bifurcation of "diamond waves" (i.e.
$\mathbf{\omega }=0)$
occurs for a velocity $c$ of the waves larger than the critical velocity $%
c_{0}$ corresponding to $\mu _{c},$ while for $\tau >\tau _{c}$
the bifurcation occurs for $c<c_{0}.$ This is in accordance with
the numerical results of Bridges et al \cite{B-D-M} (see p.
166-167 with $A_{1}=A_{2}$ real, $T_{10}=1,\tau =0$ (no surface
tension)). Notice that $\tau _{c}^{2}\approx 6.15,$ i.e. this
corresponds to a critical angle $\theta $ between the wave number
$K_{1}$ and the direction of the travelling wave, such that
$\theta \approx 68^{0},$ which is very large, and not easy to
reach experimentally.

\textbf{Remark 4:} Notice that for $\tau $ near $\tau _{c}^{\prime
}$ we still have "diamond waves' (even for $\tau =\tau
_{c}^{\prime })$ and it may exist other bifurcating 3-dimensional
waves, as noticed in \cite{B-D-M}. However to confirm this, we
need to compute at least coefficients of order 4 in
(\ref{bifurcsyst}).

\textbf{Proof}: the proof is made in Appendix \ref{a02}.

\subsection{Geometric pattern of diamond waves}

\label{diamond}

Diamond waves are obtained for $\varepsilon _{1}=\varepsilon
_{2},$ they propagate along the $x_{1}-$ axis, and possess the
symmetry $\mathcal{S}_{1}.$
In making $\varepsilon _{1}=\varepsilon _{2}=\varepsilon /2$ in Lemma \ref%
{Lembifurc} we obtain%
\begin{eqnarray}
U_{\varepsilon }^{(N)} &=&(\psi ,\eta )_{\varepsilon
}^{(N)}=\sum_{n\in
\mathbb{N}
,\text{ }n\leq N}\varepsilon ^{n}U^{(n)}\in \mathbb{H}^{k},\text{
for any }k,
\label{Diamond} \\
U^{(1)} &=&\xi _{0}=(\sin x_{1}\cos \tau x_{2},-\frac{1}{\mu
_{c}}\cos
x_{1}\cos \tau x_{2}),  \notag \\
U^{(2)} &=&\frac{1}{4}(U^{(2,0)}+U^{(0,2)}+U^{(1,1)}),  \notag
\end{eqnarray}%
where $\psi _{\varepsilon }^{(N)}\in H_{o,e}^{k},\eta
_{\varepsilon
}^{(N)}\in H_{e,e}^{k}$ meaning that $\psi _{\varepsilon }^{(N)}$ is odd in $%
x_{1},$ even in $x_{2},$ and $\eta _{\varepsilon }^{(N)}$ is even
in both
coordinates. Moreover, we have%
\begin{equation*}
\mu -\mu _{c}=\tilde{\mu}=\varepsilon ^{2}\mu _{1}+O(\varepsilon
^{4})
\end{equation*}%
with%
\begin{eqnarray}
\mu _{c} &=&(1+\tau ^{2})^{-1/2},  \notag \\
\mu _{1} &=&-\frac{\mu _{c}^{2}}{16}(\alpha _{0}+\beta _{0})  \label{mu_1} \\
&=&\frac{1}{4\mu _{c}^{3}}-\frac{1}{2\mu _{c}^{2}}-\frac{3}{4\mu _{c}}+2+%
\frac{\mu _{c}}{2}-\frac{9}{4(2-\mu _{c})}.  \notag
\end{eqnarray}%
For $\tau $ close to 0 (corresponds to some of the experiments shown in \cite%
{Ham-Hend-Seg}), one has%
\begin{equation*}
\eta =\tilde{\varepsilon}\cos x_{1}\cos \tau x_{2}+\frac{\tilde{\varepsilon}%
^{2}}{4}\cos 2x_{1}(1+\cos 2\tau x_{2})+O(\tilde{\varepsilon}^{2}\tau ^{2}+|%
\tilde{\varepsilon}|^{3})
\end{equation*}%
with
\begin{equation*}
\tilde{\varepsilon}=-\varepsilon /\mu _{c}\sim -\varepsilon .
\end{equation*}%
One can assume that $\tilde{\varepsilon}>0$ since
$\tilde{\varepsilon}<0$ would correspond to a shift by $\pi $ of
$x_{1}.$ Then the above formula shows that crests (maxima) and
troughs (minima) are elongated in the $x_{2}$ direction, with
crests sharper than the troughs, and there are "nodal" lines
($\eta \approx 0)$ (as noticed in experiments \cite{Ham-Hend-Seg})
at
\begin{equation*}
x_{2}=\pi /2\tau +n\pi /\tau ,\text{ \ }n\in
\mathbb{Z}
,
\end{equation*}%
where $\eta $ is of order $O(\tilde{\varepsilon}^{2}\tau ^{2}+|\tilde{%
\varepsilon}|^{3})$. So, the pattern roughly looks asymptotically like \emph{%
rectangles elongated in }$x_{2}$ direction, \emph{narrow around
the crests, wide around the troughs}, \emph{organized in staggered
rows}. Notice that when $\tau \rightarrow 0,$ we have $\mu
_{1}\sim -3/4,$ hence $\varepsilon
^{2}\sim (4/3)(\mu _{c}-\mu ),$ i.e.%
\begin{equation*}
\frac{c-c_{0}}{c_{0}}\sim
\frac{3}{8}F_{0}^{2}\tilde{\varepsilon}^{2},\text{ \
}F_{0}=\frac{c_{0}}{\sqrt{gL}},
\end{equation*}%
where $F_{0}$ is the Froude number built with the short wave
length $L.$

For larger, still small, values of $\tau ,$ the nodal lines
disappear and \emph{the pattern looks like hexagons} where two
sides of crests, parallel to the $x_{2}$ axis, are connected to
the nearest tip of two analogue crests, shifted by half of the
wave length in $x_{1}$ and $x_{2}$ directions.

For values of $\tau $ near 1, the pattern of the surface looks
like juxtaposition of squares.

For large $\tau ,$ i.e. in particular $\tau >\tau _{c},$ we have
for the
free surface truncated at order $\varepsilon ^{2}$%
\begin{equation*}
\eta =\tilde{\varepsilon}\cos x_{1}\cos \tau x_{2}+\frac{\tau \tilde{%
\varepsilon}^{2}}{4}\cos 2\tau x_{2}(1+\cos 2x_{1})+O(\tilde{\varepsilon}%
^{2})
\end{equation*}%
where%
\begin{equation*}
\tilde{\varepsilon}=-\varepsilon /\mu _{c}\sim -\varepsilon \tau .
\end{equation*}%
The above formula shows that crests (maxima) and troughs (minima)
are elongated in the $x_{1}$ direction, with crests sharper than
the troughs, and there are "nodal" lines ($\eta \approx 0)$ at
\begin{equation*}
x_{1}=\pi /2+n\pi ,\text{ \ }n\in
\mathbb{Z}
.
\end{equation*}%
Moreover, in the above formula, we see that $\tau $ is in factor of $\tilde{%
\varepsilon}^{2}$, which means that the second order term in the
expansion influences much sooner the shape of the surface as
$\tilde{\varepsilon}$ increases. In particular for small values of
$\tilde{\varepsilon}$ there are local maxima between two minima in
the troughs. This phenomenon is seen in experiments (see
\cite{Ham-Hend-Seg}) in the case when $\tau $ is small,
however the values of $\tilde{\varepsilon}$ allowing such a phenomenon for $%
\tau $ small are $O(1)$ and cannot be justified mathematically. So, when $%
\tau $ is large, the pattern roughly looks asymptotically like \emph{%
rectangles elongated in }$x_{1}$ direction, \emph{narrow around
the crests, wide around the troughs}, \emph{organized in staggered
rows, and where local maxima in the middle of the troughs may
occur for a large enough amplitude}. Notice in addition that when
$\tau \rightarrow \infty ,$ then $\mu _{1}\sim \frac{1}{4}\tau
^{7/2},$ hence $\varepsilon ^{2}\tau ^{2}\sim \frac{4}{\tau
^{3/2}}(\mu -\mu _{c}),$ i.e.%
\begin{equation*}
\frac{c_{0}-c}{c_{0}}\sim \frac{1}{8}\tau ^{3/2}F_{0}^{2}\tilde{\varepsilon}%
^{2},\text{ \ }F_{0}=\frac{c_{0}}{\sqrt{gL}},
\end{equation*}%
where we observe that in this last formula $L$ is the physical
wave length along $x_{1},$ which is in this case the long wave
length $(=\tau L_{1},$ if
$L_{1}$ denotes the short one). We plot at Figure \ref{fig1} the elevation $%
\eta _{\varepsilon }^{(2)}$ a) for $\tau =1/5$ ($\theta \sim 11.3^{o}),$ $%
\varepsilon =0.8\mu _{c},$ b) for $\tau =1/2$ ($\theta \sim 26.5^{o}),$ $%
\varepsilon =0.6\mu _{c}.$ These cases correspond to $\tau $ very
small or
moderately small, currently observed in experiments (see \cite{Ham-Hend-Seg}%
). Observe however that in both cases we need to consider $\tau $
not exactly 1/5 or 1/2 since both cases are not "simple cases" as
required at Theorem \ref{Lembifurc}. Indeed, if we consider
solutions $(n,m)\in
\mathbb{N}
^{2}$ different from $(1,1)$ for the critical dispersion relation%
\begin{equation}
n^{2}+\tau ^{2}m^{2}=(1+\tau ^{2})n^{4},  \label{criticdisp}
\end{equation}%
then the smallest values are $(n,m)=(61,18971)$ for $\tau =1/5,$ and $%
(13,377)$ for $\tau =1/2.$ This means that the computation fails
for coefficients of $\varepsilon ^{18971}$ in the first case, and
$\varepsilon ^{377}$ in the second case. Notice that the formal
computations made by Roberts and Schwartz \cite{Rob-Schw} (1983)
correspond to diamond waves with $\tau =1$ ($\theta =45^{o})$ and
$\tau =\sqrt{3}$ ($\theta =60^{o}).$ They observed numerically the
flattening of troughs and sharpening of crests.
However, both cases are not "simple cases" in the sense of Theorem \ref%
{Lembifurc}, and these formal computation should break at order
$\varepsilon ^{35}$ for $\tau =1,$ and at order $\varepsilon
^{195}$ for $\tau =\sqrt{3}$
(due to solutions of (\ref{criticdisp}) $(n,m)=(5,35)$ for $\tau =1,$ and $%
(n,m)=(13,195)$ for $\tau =\sqrt{3}).$

\section{Linearized operator}

In this section we study the linearized problem at a non zero
$U=(\psi ,\eta )$ for the system (\ref{basic1}), (\ref{basic2}).
We restrict our study to
diamond waves, i.e. the direction of the waves we are looking for is $%
\mathbf{u}=\mathbf{u}_{0}=(1,0)$ and the system possesses the symmetries $%
\mathcal{S}_{0}$ and $\mathcal{S}_{1}$ (see (\ref{symS}),
(\ref{symR})). The purpose is to invert the linearized operator,
for being able to use the Newton method, as it is required in
Nash-Moser theorem.

\subsection{Linearized system in $(\protect\psi ,\protect\eta )\neq 0$}

\label{lin syst}

Let us write the nonlinear system (\ref{basic1}), (\ref{basic2})
under the
form (\ref{F(U,mu)1})%
\begin{equation*}
\mathcal{F}(U,\mu )=0,
\end{equation*}%
where%
\begin{equation*}
U=(\psi ,\eta ),
\end{equation*}%
and we omit the argument $\mathbf{u}_{0}$ since it is now fixed.
Then, for
any given $(f,g)$ the linear system%
\begin{equation*}
\partial _{U}\mathcal{F}(U,\mu )[\delta U]=F,\quad F:=(f,g)
\end{equation*}%
can be written as follows%
\begin{gather*}
\partial _{\eta }\mathcal{G}_{\eta }[\delta \eta ]\psi +\mathcal{G}_{\eta
}(\delta \psi )-\mathbf{u}_{0}\cdot \nabla (\delta \eta )=f, \\
V\cdot \nabla (\delta \psi )+\mu \delta \eta
+\big(\mathfrak{b}^{2}\nabla \eta -\mathfrak{b}(\nabla \psi
+\mathbf{u}_{0})\big)\cdot \nabla (\delta \eta )=g
\end{gather*}%
where%
\begin{equation}
V=\nabla \psi +\mathbf{u}_{0}-\mathfrak{b}\nabla \eta ,\quad \mathfrak{b}=%
\frac{1}{1+|\nabla \eta |^{2}}\{\nabla \eta \cdot
(\mathbf{u}_{0}+\nabla \psi )\}.  \label{Z}
\end{equation}%
Now defining%
\begin{equation}
\delta \phi =\delta \psi -\mathfrak{b}\delta \eta ,  \notag
\end{equation}%
and after using (\ref{diffG_2},\ref{diffG_3}), we obtain the new system%
\begin{equation}
\mathcal{L}(U,\mu )[\delta \phi ,\delta \eta ]=F+\mathcal{R}(\mathcal{F}%
,U)[\delta U],  \label{defR}
\end{equation}%
where the linear \emph{symmetric} operator $\mathcal{L}(U,\mu )$
is defined
by%
\begin{equation}
\mathcal{L}(U,\mu )=\left(
\begin{array}{cc}
\mathcal{G}_{\eta } & \mathcal{J}^{\ast } \\
\mathcal{J} & \mathfrak{a}%
\end{array}%
\right) ,  \label{basiclinop}
\end{equation}%
\begin{equation}
\mathcal{J}=V\cdot \nabla (\cdot ),\text{ \ \ \
}\mathfrak{a}=V\cdot \nabla \mathfrak{b}+\mu .  \label{a}
\end{equation}%
The rest $\mathcal{R}$ has the form%
\begin{gather*}
\mathcal{R}(\mathcal{F},U)[\delta U]=(R_{1}(\mathcal{F},U)[\delta U],0), \\
R_{1}(\mathcal{F},U)[\delta U]=\mathcal{G}_{\eta }\left( \frac{\mathcal{F}%
_{1}\delta \eta }{1+(\nabla \eta )^{2}}\right) -\nabla \cdot \left( \frac{%
\mathcal{F}_{1}\delta \eta }{1+(\nabla \eta )^{2}}\nabla \eta
\right) ,
\end{gather*}%
and cancels when $U$ is a solution of $\mathcal{F}(U,\mu )=0.$ We
also
notice that for $U=(\psi ,\eta )\in \mathbb{H}_{(S)}^{m}$ where%
\begin{equation*}
\mathbb{H}_{(S)}^{m}=\big\{U\in
\mathbb{H}^{m}(\mathbb{R}^{2}/\Gamma ):\psi
\text{ odd in }x_{1},\text{ even in }x_{2},\text{ }\eta \text{ even in }x_{1}%
\text{ and in }x_{2}\big\},
\end{equation*}%
then%
\begin{gather*}
V=(V_{1},V_{2})\in H_{e,e}^{m-1}(\mathbb{R}^{2}/\Gamma )\times H_{o,o}^{m-1}(%
\mathbb{R}^{2}/\Gamma ), \\
\mathfrak{b}\in H_{o,e}^{m-1}(\mathbb{R}^{2}/\Gamma ),\quad
\mathfrak{a}\in H_{e,e}^{m-2}(\mathbb{R}^{2}/\Gamma ),
\end{gather*}%
all these functions being invariant under the shift%
\begin{equation*}
\mathcal{T}_{\mathbf{v}_{0}}:(x_{1},x_{2})\mapsto (x_{1}+\pi
,x_{2}+\pi /\tau ).
\end{equation*}%
Moreover we have the following "tame" estimates

\begin{lemma}
\label{tame coef} Let $U\in \mathbb{H}_{(S)}^{m},$ $m\geq 3.$
Then, there exists $M_{3}>0$ such that for $||U||_{3}\leq M_{3},$
one has
\begin{equation*}
||V-\mathbf{u}||_{m-1}+||\mathfrak{a}-\mu ||_{m-2}\leq
c_{m}(M_{3})||U||_{m},
\end{equation*}%
\begin{equation*}
||R_{1}(\mathcal{F},U)[\delta U]||_{m-2}\leq c_{s}(M_{3})\{||\mathcal{F}%
_{1}||_{2}(||\eta ||_{m}||\delta \eta ||_{2}+||\delta \eta ||_{m-1})+||%
\mathcal{F}_{1}||_{m-1}||\delta \eta ||_{2}\}.
\end{equation*}
\end{lemma}

\begin{proof}
The tame estimates on $V-\mathbf{u}$ and $\mathfrak{a}-\mu $
result directly from their definitions, from the following
inequality, valid for any $f,g\in
H^{m}(%
\mathbb{R}
^{2}/\Gamma ),$ $m\geq 2$%
\begin{equation}
||fg||_{m}\leq c_{m}\{||f||_{2}||g||_{m}+||f||_{m}||g||_{2}\},
\label{tame prod}
\end{equation}%
and from interpolation estimates like%
\begin{equation}
||f||_{\lambda \alpha +(1-\lambda )\beta }\leq c||f||_{\alpha
}^{\lambda }||f||_{\beta }^{1-\lambda },  \label{interp-estim}
\end{equation}%
which leads to%
\begin{equation*}
||f||_{\alpha _{2}}||f||_{\beta _{2}}\leq c||f||_{\alpha
_{1}}||f||_{\beta _{1}}
\end{equation*}%
when
\begin{equation*}
0\leq \alpha _{1}\leq \alpha _{2}\leq \beta _{2}\leq \beta _{1},\text{ \ }%
\alpha _{1}+\beta _{1}=\alpha _{2}+\beta _{2}.
\end{equation*}%
The tame estimate on $R_{1}(\mathcal{F},U)[\delta U]$ follows from
the tame estimates (\ref{tame prod}), (\ref{tame G}) and from the
following
interpolation estimate deduced from (\ref{interp-estim}):%
\begin{equation*}
||f||_{3}||g||_{k+1}\leq
c_{k}\{||f||_{2}||g||_{k+2}+||f||_{k+1}||g||_{3}\}.
\end{equation*}
\end{proof}

The main problem in using the Nash Moser theorem, is to invert the
approximate linearized system, i.e. invert the linear system%
\begin{equation*}
\mathcal{L[}\delta \phi ,\delta \eta ]=(f,g),
\end{equation*}%
which leads to the scalar equation%
\begin{equation}
\mathcal{G}_{\eta }(\delta \phi )-\mathcal{J}^{\ast }(\frac{1}{\mathfrak{a}}%
\mathcal{J}(\delta \phi ))=h  \label{linequa}
\end{equation}%
with%
\begin{equation*}
h=f-\mathcal{J}^{\ast }(\frac{1}{\mathfrak{a}}g)\in H_{o,e}^{s}(%
\mathbb{R}
^{2}/\Gamma )
\end{equation*}%
and where we look for $\delta \phi $ in some $H_{o,e}^{s-r}(%
\mathbb{R}
^{2}/\Gamma ).$

\subsection{Pseudodifferential operators and diffeomorphism of the torus}

In this section we use a diffeomorphism: $\mathbb{R} ^{2}\rightarrow \mathbb{%
R}^{2},$ such that the principal part of the symbol of the linear
operator occurring in (\ref{linequa}) has a simplified structure.
Its new structure will allow us to use further a suitable descent
method for obtaining, at the end of the process, a
pseudodifferential operator equation with constant coefficients,
plus a perturbation operator of "small" order.

Let us denote the change of coordinates by $X=X(Y)$, where
$X(\cdot )$ is a
not yet determined diffeomorphism of $\mathbb{R}^{2}$ such that%
\begin{equation}
X(Y)=\mathbb{T}Y+\widetilde{\mathcal{V}}(Y),\text{ \ }\mathbb{T}%
Y=(y_{1},y_{2}/\tau ),  \label{varchange2}
\end{equation}%
We assume that $\widetilde{\mathcal{V}}(Y)$ and
$\tilde{\eta}(Y)=\eta (X(Y))$ satisfy the following

\begin{condition}
\label{diffeomorphism} Functions $\tilde\eta$ and
$\widetilde{\mathcal{V}}$
are doubly $2\pi$-periodic, $\widetilde \eta$ is even in $y_1$ and $y_2$, $%
\widetilde{\mathcal{V}}_1$ is odd in $y_1$ and even in $y_2$, $\widetilde{%
\mathcal{V}}_2$ is odd in $y_2$ and even in $y_1$, and
\begin{equation*}
\tilde\eta (y_1+\pi, y_2+\pi)=\tilde \eta(y_1,y_2),\quad \widetilde{\mathcal{%
V}} (y_1+\pi, y_2+\pi)=\widetilde{ \mathcal{V}}(y_1,y_2).
\end{equation*}
\end{condition}

In particular, $X(Y)$ takes diffeomorphically
$\mathbb{R}^{2}/\Gamma _{1}$ onto $\mathbb{R}^{2}/\Gamma _{\tau
}$. The lattices of periods $\Gamma _{1}$ and $\Gamma _{\tau }$
are respectively the dual of lattices $\Gamma _{1}^{\prime }$ and
$\Gamma _{\tau }^{\prime }$ generated by the wave vectors $(1,\pm
1)$ and $(1,\pm \tau )$. In new coordinates the free surface has
the parametric representation
\begin{equation*}
x=\mathbf{r}(Y):=(X_{1}(Y),X_{2}(Y),\widetilde{\eta }(Y))^{t},
\end{equation*}%
with 2x2 matrix $\mathbb{G}(Y)$ of the\emph{\ first fundamental
form of the free surface} defined by $g_{ij}=\partial
_{y_{i}}\mathbf{r}\cdot \partial _{y_{j}}\mathbf{r}$. We denote by
$J$ the determinant of the Jacobian matrix $\mathbb{B}(Y)=\nabla
_{Y}X(Y)$.

Our aim is to simplify the structure of the operators involved in
the basic equation \eqref{linequa} by choosing an appropriate
change of coordinates. The most suitable tool for organizing such
a process is the theory of pseudodifferential operators, and we
begin with recalling the definition of a pseudodifferential
operator . We consider the class of integro-differential operators
on a two-dimensional torus having the representation
\begin{equation*}
\mathfrak{A}u(Y)=\frac{1}{2\pi }\sum\limits_{k\in \mathbb{Z}%
^{2}}e^{ikY}A(Y,k)\widehat{u}(k),\quad \widehat{u}(k)=\frac{1}{2\pi }%
\int\limits_{{T}^{2}}e^{-ikY}u(Y)dY,
\end{equation*}%
which properties are completely characterised by the function $A:\mathbf{T}%
^{2}\times \mathbb{R}^{2}\mapsto \mathbb{C}$ named the symbol of $\mathfrak{A%
}$. We say that $\mathfrak{A}$ is a pseudodifferential operator,
if its symbol satisfy the condition

\begin{condition}
\label{condA}There are integers $l>0$, $m\geq 0$ and a real $r$
named the order of the operator $\mathfrak{A}$ so that
\begin{equation*}
|\mathfrak{A}|_{m,l}^{r}=\Vert A(\cdot ,0)\Vert
_{C^{l}}+\sup\limits_{k\in \mathbb{Z}^{2}\setminus
\{0\}}\sup\limits_{|\alpha |\leq m}|k|^{|\alpha |-r}\Vert \partial
_{k}^{\alpha }A(\cdot ,k)\Vert _{C^{l}}<\infty .
\end{equation*}
\end{condition}

Pseudodifferential operators enjoy many remarkable properties
including
explicit formulae for compositions and commutators (see Appendix \ref%
{pseudodifferential} for references and more details). Important
examples of
such operators are the first-order pseudodifferential operator $\mathcal{G}%
^{(0)}=(-\Delta )^{1/2}$ with the symbol $\big|\mathbb{T}^{-1}
k\big|$ and the second-order pseudodifferential operator
\begin{equation}  \label{mathfrakL}
\mathfrak{L}=\nu \mathfrak{D}_{1}^{2}+(-\Delta )^{1/2}, \text{~~where~~}%
\mathfrak{D}_{1}=\partial _{y_{1}},
\end{equation}%
with the symbol
\begin{equation}
L(k)=-\nu k_{1}^{2}+\big|\mathbb{T}^{-1}k\big|.  \label{LL}
\end{equation}%
On the other hand, integro-differential operators
$\mathfrak{D}_{1}^{j}$, defined by
\begin{equation}
\mathfrak{D}_{1}^{j}u(Y)=\frac{1}{2\pi }\sum\limits_{k_{1}\neq 0}e^{ikY}\Big(%
{ik_{1}}\Big)^{j}\widehat{u}_{k},\quad j\in \mathbb{Z}.
\label{1018}
\end{equation}%
are not pseudodifferential for $j\leq 0,$ which easy follows from
the formulae
\begin{equation*}
\mathfrak{D}_{1}^{j}u=\partial _{y_{1}}^{j}u\text{~~for~~}j>0,\text{~~and~~}%
\mathfrak{D}_{1}^{0}u=\Pi _{1}u=u-\frac{1}{2\pi }\int\limits_{-\pi
}^{\pi }u(s,y_{2})ds.
\end{equation*}

Further we will consider also the special class of zero-order
pseudodifferential operators $\mathfrak{A}$ with symbols having
the form of composition $A(Y,\xi (k))$, where the vector field
$\xi (k)=(\xi _{1}(k),\xi _{2}(k))$ is defined by
\begin{equation*}
\xi (k)=\mathbb{T}^{-1}k\big/|\mathbb{T}^{-1}
k|\text{~~for~~}k\neq 0,\quad \xi (0)=0.
\end{equation*}%
The metric properties of such operators are characterized by the
norm
\begin{equation}
|\mathfrak{A}|_{m,l}=\sup\limits_{|\alpha |\leq
m}\sup\limits_{|\xi |\leq 1}\Vert \partial _{\xi }^{\alpha
}A(\cdot ,\xi )\Vert _{C^{l}}<\infty , \label{103}
\end{equation}%
which is equivalent to the norm $|\mathfrak{A}|_{m,l}^{0}$. By
abuse of notation, further we will write simply $\xi $ instead of
$\xi (k)$ and use both the notations $A(Y,\xi _{1},\xi _{2})$ and
$A(Y,\xi )$ for $A(Y,\xi )$.

We are now in a position to formulate the main result of this
section.

\begin{theorem}
\label{thmChangeCoord}For any integers $\rho ,m\geq 14$, real
$\tau \in
(\delta ,\delta ^{-1})$, and $U=(\psi ,\eta )\in $ $\mathbb{H}_{(S)}^{m}(%
\mathbb{R}^{2}/\Gamma _{\tau })$, there exists $\varepsilon
_{0}>0$ so that for
\begin{equation*}
\Vert U\Vert _{\rho }\leq \varepsilon \text{~~with~~}\varepsilon
\in \lbrack 0,\varepsilon _{0}]
\end{equation*}%
there are a diffeomorphism of the torus of the form
\eqref{varchange2}, satisfying Condition \ref{diffeomorphism},
zero-order pseudodifferential operators $\mathfrak{A}$,
$\mathfrak{B}$ and an integro-differential operator
$\mathfrak{L}_{-1}$ of order $-1$ such that:

\begin{itemize}
\item[(\textit{i})] The identity
\begin{equation}  \label{equchange1}
\Big[\mathcal{G}_\eta \check u-\mathcal{J}^*\big(\frac{1}{\mathfrak{a}}%
\mathcal{J }\check u\big)\Big]\circ (\mathbb{T}+\tilde{\mathcal{V}})=\kappa%
\Big[ \mathfrak{L}u+\mathfrak{AD}_{1}u+\mathfrak{B}u
+\mathfrak{L}_{-1}u\Big]
\end{equation}
holds true for any $u \in H^2_{o,e}(\mathbb{R}^2/\Gamma_1)$ and
$\check u(X)=u(Y(X))$.

\item[(\textit{ii})] The operators $\mathfrak{A}$, $\mathfrak{B}$ and $%
\mathfrak{L}_{-1}$ have the bounds
\begin{gather}
|\mathfrak{A}|_{4,m-6}+|\mathfrak{B}|_{4,m-6}\leq c_{m}(\Vert
U\Vert
_{6})\Vert U\Vert _{m},  \notag \\
||\mathfrak{L}_{-1}u||_{r}\leq c\varepsilon
||u||_{r-1},\text{~~for~~}1\leq
r\leq \rho -13,  \notag \\
||\mathcal{L}_{-1}u||_{s}\leq c(\varepsilon
||u||_{s-1}+||U||_{s+13}||u||_{0}).  \notag
\end{gather}

\item[(\textit{iii})] Operators $\mathfrak{A}$, $\mathfrak{B}$ and $%
\mathfrak{L}_{-1}$ are invariant with respect to the symmetries $%
Y\rightarrow \pm Y^{\ast }$, $Y^{\ast }=(-y_{1},y_{2})$ which is
equivalent to the equivariant property
\begin{equation}
\begin{split}
\mathfrak{AD}_{1}u(\pm Y^{\ast })& =\mathfrak{AD}_{1}u^{\ast }(\pm
Y),\quad
\mathfrak{B}u(\pm Y^{\ast })=\mathfrak{B}u^{\ast }(\pm Y), \\
\mathfrak{L}_{-1}u(\pm Y^{\ast })& =\mathfrak{L}_{-1}u^{\ast }(\pm
Y),\quad u^{\ast }(Y)=u(Y^{\ast });
\end{split}
\label{irr1}
\end{equation}%
they are also invariant with respect to transform $Y\rightarrow
Y+(\pi ,\pi ) $.

\item[(\textit{iv})] Diffeomorphism \eqref{varchange2} of the
torus can be inverted as $Y=\mathbb{T}^{-1}(X-\mathcal{V}(X))$,
\begin{equation*}
y_{1}=x_{1}+d(x_{1},x_{2}),\quad y_{2}=\tau x_{2}+\tau
e(x_{1},x_{2}).
\end{equation*}%
Functions $d\in C_{o,e}^{m-4}(\mathbb{R}^{2}/\Gamma _{\tau })$,
$e\in
C_{e,o}^{m-4}(\mathbb{R}^{2}/\Gamma _{\tau })$, $\kappa ,J\in C_{e,e}^{m-4}(%
\mathbb{R}^{2}/\Gamma _{1})$ and parameter $\nu $ satisfy the inequalities%
\begin{gather}
\Vert d\Vert _{C^{m-4}}+\Vert e\Vert _{C^{m-4}}\leq c_{m}(\Vert
U\Vert _{4})\Vert U\Vert _{m},\quad |\nu -1/\mu |\leq c(\Vert
U\Vert _{4})\Vert
U\Vert _{4},  \label{denu} \\
\Vert \kappa -1\Vert _{C^{m-5}}+\Vert J-1/\tau \Vert
_{C^{m-5}}\leq c_{m}(\Vert U\Vert _{5})\Vert U\Vert _{m}.
\label{kappaJ}
\end{gather}
\end{itemize}
\end{theorem}

The proof is based on two propositions, the first of which gives
the representation of the Dirichlet-Neumann operator in the
invariant parametric form, and the second shows that trajectories
of liquid particles on the free surface forms a foliation of the
two-dimensional torus.

In order to formulate them it is convenient to introduce the
notations
$$
\mathbf G_1(Y,k)=\sqrt{\mathbb G^{-1}k\cdot k},\quad
\mathbf{div}\,\mathbf {q}(Y)=\frac{1}{\sqrt{\text{det~} \mathbb
G}}\text{~div}_Y\,\big(\sqrt{\text{det~} \mathbb G}\mathbf
{q}(Y)\big).
$$
Recall that $\mathbb Gk\cdot k$ is the first fundamental form of
the surface $\Sigma$.
\begin{theorem}
\label{thmDir-Neum} Suppose that functions $\tilde{\eta}$ and $\widetilde{%
\mathcal{V}}$ satisfy Condition \ref{diffeomorphism} and there are integers $%
\rho ,$ $l$ such that
\begin{eqnarray*}
||\widetilde{\eta }||_{C^{\rho
}}+||\widetilde{\mathcal{V}}||_{C^{\rho }}
&\leq &\varepsilon ,\text{ \ \ }10\leq \rho \leq l, \\
||\widetilde{\eta }||_{C^{s}}+||\widetilde{\mathcal{V}}||_{C^{s}}
&\leq &E_{l},\text{ \ \ \ }s\leq l.
\end{eqnarray*}%
Then there exists $\varepsilon _{0}>0$ depending on $\rho $ and
$l$ only such that for $0\leq \varepsilon \leq \varepsilon _{0}$
and $2\pi $-periodic sufficiently smooth function $u$, the
operator $\mathcal{G}_{\eta }$ has the
representation%
\begin{equation}
\mathcal{G}_{\eta }\check{u}\circ (\mathbb{T}+\widetilde{\mathcal{V}})=%
\mathcal{G}_{1}u+\mathcal{G}_{0}u+\mathcal{G}_{-1}u,\quad \check{u}%
(X)=u(Y(X)).  \label{representationG}
\end{equation}%
Here $\mathcal{G}_{1}$ is a first order pseudodifferential
operator with
symbol%
\begin{equation}
G_{1}(Y,k)=\frac{\sqrt{\det \mathbb{G}}}{J}\mathbf{G}_1(Y,k),%
\text{ \ }Y\in
\mathbb{R}
^{2},\text{ \ }k\in
\mathbb{Z}
^{2},  \label{defG_1}
\end{equation}%
$\mathcal{G}_{0}$ is a zero order pseudodifferential operator with
 symbol
\begin{gather}\nonumber
G_0=\text{\rm Re~} G_0+i\text{\rm Im~} G_0,\\
\label{cor11} \text{\rm Re~}G_0(Y,k)=\frac{\text{\rm det~}\mathbb
G}{2J^2}\,\big[\,\frac{1}{\mathbf G_1^2}\,Q(Y,k)+
\mathbf{div}(\mathbb G^{-1}\nabla_Y\tilde\eta)\,\big],
\\\label{cor12} \text{\rm Im~} G_0(Y,k)=-\frac{\sqrt{\text{\rm
det~}\mathbb G}}{2J}\mathbf{div}(\nabla_k \mathbf G_1).
\end{gather}
Here the quadratic form $Q(Y,.)$ is given by
\begin{equation}\label{cor13}
Q(Y,k)=\frac{1}{2}\nabla_Y(\mathbb G^{-1}k\cdot k)\cdot (\mathbb
G^{-1}\nabla_Y\tilde\eta)-\mathbb G^{-1}k\cdot\nabla_Y(\mathbb
G^{-1}\nabla_Y\tilde\eta\cdot k),
\end{equation}
and the operator  $\mathcal G_0$  satisfies the
estimates%
\begin{equation}
|\mathcal{G}_{0}u|_{4,\rho -2}^{0}\leq c\varepsilon ,\text{ \ }|\mathcal{G}%
_{0}u|_{4,l-2}^{0}\leq cE_{l},  \label{estimG0}
\end{equation}%
while the linear operator $\mathcal{G}_{-1}$ satisfies%
\begin{eqnarray}
||\mathcal{G}_{-1}u||_{r} &\leq &c\varepsilon ||u||_{r-1},\text{ \ for \ }%
1\leq r\leq \rho -9,  \notag \\
||\mathcal{G}_{-1}u||_{s} &\leq &c(\varepsilon ||u||_{s-1}+E_{l}||u||_{0}),%
\text{ \ for \ }s\leq l-9.  \label{estimG-1}
\end{eqnarray}%
Moreover, operators $\mathcal{G}_{1},$ $\mathcal{G}_{0}$ and $\mathcal{G}%
_{-1}$ satisfy the symmetry properties%
\begin{equation*}
\mathcal{G}_{j}u(\pm Y^{\ast })=\mathcal{G}_{j}u^{\ast }(\pm Y),\text{ \ }%
j=1,0,-1,\text{ \ }u^{\ast }(Y)=u(Y^{\ast }).
\end{equation*}
\end{theorem}
\begin{proof}
The proof is given in Appendix \ref{a03}.
\end{proof}

The formula for the principal term $G_1$ is a classic result of
the theory of pseudodifferential operators \cite{Hormander}. The
expression for the second term in  local Riemann coordinates
 was given  in \cite{AntonineBarucq}. It seems that the general formulae
\eqref{cor11}, \eqref{cor12} are new. Note  that the ratio
$\text{det~}\mathbb G/J^2=1+|\nabla_X\eta(X)|^2$ is a scalar
invariant and the real part of  $G_0$ can be rewritten in  the
invariant form
\begin{equation}\label{geometry0}
\frac{J}{\sqrt{\text{det~}\mathbb G}}\text{Re~}
G_0=\frac{1}{2}\frac{LG-2MF+NF}{EG-F^2}-\frac{1}{2}\frac{L\xi_1^2+2M
\xi_1\xi_2+N\xi_2^2} {E\xi_1^2+2F \xi_1\xi_2+G\xi_2^2}.
\end{equation}
Here we use the standard notations for the second fundamental form
$L\xi_1^2+2M \xi_1\xi_2+N\xi_2^2$ and the first fundamental form
$E\xi_1^2+2F \xi_1\xi_2+G\xi_2^2$ of the surface $x_3=\eta(X)$,
the vector $\xi$ is connected with the covector $k$ by the
relation $k=\mathbb G \xi$. The right hand side of
\eqref{geometry0} is the difference between the mean curvature of
$\Sigma$ and half of the normal curvature of $\Sigma$ in the direction
$\xi$. Note also that the conclusion of Theorem \ref{thmDir-Neum}
holds true without assumption on the smallness of $\varepsilon$,
but the proof  becomes more complicated and goes far beyond the
scope of the paper.

\begin{lemma}
\label{changevar1} For $m\geq 4,$ and $U\in \mathbb{H}_{(S)}^{m}(\mathbb{R}%
^{2}/\Gamma _{\tau })$ with $||U||_{4}$ small enough, there exists
a unique
function $\mathcal{Z\in C}^{m-3}(%
\mathbb{R}
^{2})$ such that%
\begin{equation}
\frac{\partial \mathcal{Z}}{\partial z_{1}}=\frac{V_{2}}{V_{1}}(z_{1},%
\mathcal{Z}),\quad {\Pi }_{1}\mathcal{Z}=z_{2},  \label{charact}
\end{equation}%
where ${\Pi }_{1}$ denotes the average over a period in $z_{1}$,
and $V_{i}$
are the components of the vector field $V$ defined by \eqref{Z}. Moreover, $%
\mathcal{Z}$ is even in $z_{1},$ odd in $z_{2}$,
\begin{equation*}
\mathcal{Z}(Z)=\mathcal{Z}(Z+(2\pi ,0))=\mathcal{Z}(Z+(0,2\pi
/\tau ))-2\pi /\tau =\mathcal{Z}(Z+(\pi ,\pi /\tau ))-\pi /\tau ,
\end{equation*}%
and the shifted function $\mathcal{T}_{\delta
}\mathcal{Z}=\mathcal{Z}(\cdot +\delta ,\cdot )$ is solution of
the system (\ref{charact}) where $V=V(\cdot +\delta ,\cdot )$.
Moreover, the mapping
\begin{equation}
x_{1}=z_{1},\quad x_{2}=z_{2}-\widetilde{d_{1}}(z_{1},z_{2})=\mathcal{Z}%
(z_{1},z_{2}),  \label{U1}
\end{equation}%
with his inverse
\begin{equation*}
z_{1}=x_{1},\quad z_{2}=x_{2}+d_{1}(x_{1},x_{2})
\end{equation*}%
define automorphisms of the torus $\mathbb{R}/\Gamma _{\tau }$: $X\mapsto Z=%
\mathcal{U}_{1}(X)$ , $X=\mathcal{U}_{1}^{-1}(Z)$. The functions
$d_{1}$ and $\widetilde{d_{1}}$ have symmetry $(e,o),$ as above
and we have the
following tame estimates%
\begin{equation*}
||d_{1}||_{C^{m-3}}+||\widetilde{d_{1}}||_{C^{m-3}}\leq
c_{m}(||U||_{4})||U||_{m}.
\end{equation*}%
The automorphism $\mathcal{U}_{1}$ takes integral curves of the
vector field
$V$, which coincide with the bicharacteristics of the operator $\mathcal{G}%
_{\eta }-\mathcal{J}^{\ast }(\mathfrak{a}^{-1}\mathcal{J}\cdot )$,
onto straight lines $\{z_{2}=\text{const.}\}.$ In other words,
bicharacteristics form a foliation of the torus with a rotation
number equal to $0$.
\end{lemma}

\begin{proof}
The proof is made in Appendix \ref{a2}.
\end{proof}

Let us turn to the proof of Theorem \ref{thmChangeCoord}. We look
for the
desired diffeomorphism $Y\rightarrow X$ in the form of the composition $(%
\mathcal{U}_{2}\circ \mathcal{U}_{1})^{-1}$,
\begin{equation*}
X\overset{\mathcal{U}_{1}}{\longrightarrow }Z\overset{\mathcal{U}_{2}}{%
\longrightarrow }Y,
\end{equation*}%
where the diffeomorphism $\mathcal{U}_{1}:\mathbb{R}^{2}/\Gamma
_{\tau
}\mapsto \mathbb{R}^{2}/\Gamma _{\tau }$ is completely defined by Lemma \ref%
{changevar1}, and the diffeomorphism
$\mathcal{U}_{2}:\mathbb{R}^{2}/\Gamma _{\tau }\mapsto
\mathbb{R}^{2}/\Gamma _{1}$ is unknown. We look for it in
the form $Y=\mathcal{U}_{2}(Z)$ with%
\begin{equation}
y_{1}=z_{1}+d_{2}(Z),\quad y_{2}=\tau (z_{2}+e_{2}(z_{2})),
\label{U_2}
\end{equation}%
where functions $d_{2}$ and $e_{2}$ will be specified below. Our
first task
is to make the formal change of variables in the left hand side of %
\eqref{equchange1}. We begin with the consideration of the
second-order differential operator $\mathcal{J}^{\ast
}(\mathfrak{a}^{-1}\mathcal{J}\cdot )$. It follows from the
equality $\mathcal{J}=V\cdot \nabla _{X}$ that
\begin{equation*}
-\mathcal{J}^{\ast }\big(\frac{1}{\mathfrak{a}}\mathcal{J}\check{u}\big)%
)\circ \mathcal{U}_{1}^{-1}\circ \mathcal{U}_{2}^{-1}\equiv \frac{1}{J}%
\nabla _{Y}\cdot
\Big\{\frac{J}{\mathfrak{a}}\big(\mathbb{B}^{-1}V\cdot \nabla
_{Y}u\big)\mathbb{B}^{-1}V\Big\},
\end{equation*}%
where
\begin{equation*}
V=V\big(X(Y)\big),\quad
\mathfrak{a}=\mathfrak{a}\big(X(Y)\big),\quad \mathbb{B}(Y)=\nabla
_{Y}X(Y),J(Y)=\text{det~}\mathbb{B}(Y).
\end{equation*}%
In particular, we have
\begin{equation*}
\mathbb{B}^{-1}\big(Y(Z)\big)=\nabla _{Z}Y(Z)\;\big[\nabla _{Z}X(Z)\big]%
^{-1}.
\end{equation*}%
On the other hand, Lemma \ref{changevar1} and formulae \eqref{U_2}
imply
\begin{equation*}
\big[\nabla _{Z}X(Z)\big]^{-1}\,V\big(X(Z)\big)=V_{1}\big(X(Z)\big)\mathbf{e}%
_{1},\quad \nabla _{Z}Y(Z)\,\mathbf{e}_{1}=\big(1+\partial _{z_{1}}d_{2}(Z)%
\big)\mathbf{e}_{1},
\end{equation*}%
where $\mathbf{e}_{1}=(1,0)$. Thus we get
\begin{equation*}
-\mathcal{J}^{\ast }\big(\frac{1}{\mathfrak{a}}\mathcal{J}\check{u}\big)%
)\equiv \frac{1}{J}\partial
_{y_{1}}\big\{\frac{J}{\mathfrak{a}}(1+\partial
_{z_{1}}d_{2})^{2}V_{1}^{2}\partial _{y_{1}}u\big\}.
\end{equation*}%
From this and parametric representation \eqref{representationG} of
the Dirichlet-Neumann operators we conclude that the left hand
side of the desired identity \eqref{equchange1} is equal to
\begin{equation*}
\mathfrak{p}\big(Z(Y)\big)\partial
_{y_{1}}^{2}u(Y)+\mathfrak{s}(Y))\partial
_{y_{1}}u(Y)+\mathcal{G}_{1}u(Y)+\mathcal{G}_{0}u(Y)+\mathcal{G}_{-1}u(Y),
\end{equation*}%
where
\begin{gather}
\mathfrak{p}(Z)=\frac{1}{\mathfrak{a}\big(X(Z)\big)}\big(1+\partial
_{z_{1}}d_{2}(Z)\big)^{2}V_{1}\big(X(Z)\big)^{2},  \label{coefficientsJ} \\
\mathfrak{s}(Y)=\frac{1}{J(Y)}\partial _{y_{1}}\big\{J(Y)\mathfrak{p}(Z(Y))%
\big\}.  \notag
\end{gather}%
Note that (\ref{equchange1}) is a second order operator with
respect to the
variable $y_{1}$ and a first- order operator with respect to variable $y_{2}$%
. Its principal part is the pseudodifferential operator $\mathfrak{p}%
\partial y_{1}^{2}+\mathcal{G}_{1}$ which symbol reads
\begin{equation*}
-\mathfrak{p}k_{1}^{2}+J^{-1}\big(\text{adj}\,\mathbb{G}\;k\cdot k\big)%
^{1/2},
\end{equation*}%
which we write in the form
\begin{equation*}
-\mathfrak{p}k_{1}^{2}+(\tau J)^{-1}\sqrt{g_{11}}\Big(\tau
^{2}k_{2}^{2}-2\tau ^{2}g_{11}^{-1}g_{12}k_{1}k_{2}+\tau
^{2}g_{11}^{-1}g_{22}k_{1}^{2}\Big)^{1/2}.
\end{equation*}%
Hence operator (\ref{equchange1}) can be reduced to the operator
with constant coefficients at principal derivatives if for some
constant $\nu $ ,
\begin{equation}
\nu ^{-1}\mathfrak{p}\big(Z(Y)\big)=\big(\tau
J(Y)\big)^{-1}\sqrt{g_{11}(Y)}. \label{PG11}
\end{equation}%
This equality can be regarded as first order differential equation
for functions $d_{2}(Z)$, $e_{2}(z_{2})$ and a constant $\nu $. It
becomes clear if we write the right hand side as a function of the
variable $Z$. To this end note that since $\partial
_{y_{1}}z_{2}(Y)=\partial _{z_{2}}x_{1}(Z)=0$,
\begin{gather*}
g_{11}=\big|\nabla _{Z}X\,\partial _{y_{1}}Z\big|^{2}+\big(\nabla
_{X}\eta
\cdot \nabla _{Z}X\cdot \partial _{y_{1}}Z\big)^{2}= \\
\lbrack \partial _{y_{1}}z_{1}(Y)]^{2}\big[1+(\partial _{z_{1}}\mathcal{Z}%
)^{2}+(\partial _{x_{1}}\eta +\partial _{x_{2}}\eta \partial _{z_{1}}%
\mathcal{Z})^{2}\big],
\end{gather*}%
which along with \eqref{charact} gives
\begin{equation}
g_{11}=\frac{V^{2}+(V\cdot \nabla _{X}\eta
)^{2}}{V_{1}^{2}(1+\partial _{z_{1}}d_{2})^{2}}.  \label{g_11_2}
\end{equation}%
On the other hand, we have
\begin{equation*}
J\big(Y(Z)\big)^{-1}=\text{det~}\nabla
_{Z}Y(Z)\big(\text{det~}\nabla _{Z}X(Z)\big)^{-1}=\newline
\tau (1+\partial _{z_{1}}d_{2})(1+\partial _{z_{2}}e_{2})(\partial _{z_{2}}%
\mathcal{Z})^{-1}.
\end{equation*}%
Substituting these equalities into \eqref{PG11} we obtain the
differential equation for $d_{2}$ and $e_{2}$,
\begin{equation}
\partial _{z_{1}}d_{2}(Z)=\big[\nu \mathfrak{q}(Z)\,(1+e_{2}^{\prime
}(z_{2}))\,\big]^{1/2}-1,  \label{d_2}
\end{equation}%
where $\mathfrak{q}$ is equal to
\begin{equation}
V_{1}(X(Z))^{-3}\Big(|V(X(Z))|^{2}+\big(V(X(Z))\cdot \nabla _{X}\eta (X(Z))%
\big)^{2}\Big)^{1/2}\mathfrak{a}(X(Z))\partial _{z_{2}}\big(\mathcal{Z}(Z)%
\big)^{-1},  \label{q}
\end{equation}%
has symmetry $(e,e)$, and is close to $\mu $. Note that it is
completely defined by equalities \eqref{Z}, \eqref{a}, and Lemma
\ref{changevar1}. This
equation can be solved as follows. We first note that it gives a unique $%
d_{2}$ with symmetry $(o,e)$ provided that%
\begin{equation}
\nu ^{-1/2}(1+e_{2}^{\prime }(z_{2}))^{-1/2}=\frac{1}{2\pi
}\int_{-\pi }^{\pi }\mathfrak{q}^{1/2}(Z)dz_{1}.  \label{e_2}
\end{equation}%
This now determines a unique periodic odd function $e_{2}$ provided that%
\begin{equation}
\nu =2\pi \tau \int_{-\pi /\tau }^{\pi /\tau }\left( \int_{-\pi }^{\pi }%
\mathfrak{q}^{1/2}(Z)dz_{1}\right) ^{-2}dz_{2},  \label{nu}
\end{equation}%
which is just $1/\mu $ when $(\psi ,\eta )=0.$ We observe that the
invariance of $\mathfrak{q}(Z)$ under the shift
$T_{\mathbf{v}_{0}},$ leads
to the invariance of the right hand side of (\ref{e_2}) by the change $%
z_{2}\mapsto z_{2}+\pi /\tau .$ This means that $e_{2}$ is indeed
$\pi /\tau -$ periodic. It follows from Appendix G (see Lemma G.3)
of \cite{IPT} and \ Lemma \ref{changevar1} that
\begin{eqnarray*}
||\mathfrak{q}-\mu ||_{C^{m-4}} &\leq &c_{m}(||U||_{4})||U||_{m}, \\
|\nu -1/\mu | &\leq &c||U||_{4}, \\
||d_{2}||_{C^{m-4}}+||e_{2}||_{C^{m-4}} &\leq
&c_{m}(||U||_{4})||U||_{m}.
\end{eqnarray*}%
From this, the identities
\begin{equation*}
d(X)=d_{2}(x_{1},x_{2}+d_{1}(X)),\quad
e(X)=d_{1}(X)+e_{2}(x_{2}+d_{1}(X))),
\end{equation*}%
and from Appendix G (see Lemma G.1) of \cite{IPT} we deduce estimates %
\eqref{denu} for the functions $d$, $e$, and also estimate
\eqref{kappaJ} for the Jacobian $J$. Hence the diffeomorphism
$Y\rightarrow X$ is well
defined and meets all requirements of assertion (\textit{iv}) of Theorem \ref%
{thmChangeCoord}. Next set
\begin{equation*}
\kappa (Y)=(\tau J(Y))^{-1}\sqrt{g_{11}(Y)},\text{ \ }b(Y)=-\tau
g_{11}(Y)^{-1}g_{12}(Y),\text{ \ }1+2a=\tau
^{2}g_{11}(Y)^{-1}g_{12}(Y).
\end{equation*}%
It follows from \eqref{denu} that
\begin{equation}
||g_{11}-1||_{C^{m-5}}+||g_{22}-1/\tau
^{2}||_{C^{m-5}}+||g_{12}||_{C^{m-5}}\leq
c_{m}(||U||_{5})||U||_{m}, \label{estimg}
\end{equation}%
which leads to estimate \eqref{kappaJ} for $\kappa $ and the
following estimates for the functions $a$ and $b$
\begin{equation}
||a||_{C^{m-5}}+||b||_{C^{m-5}}\leq c_{m}(||U||_{5})||U||_{m}.
\label{estimab}
\end{equation}%
We are now in a position to define the operators $\mathfrak{A}$, $\mathfrak{B%
}$ and $\mathfrak{L}_{-1}$ defined by the left hand side of identity %
\eqref{equchange1}. It follows from \eqref{coefficientsJ} and
\eqref{PG11}, that its coefficients satisfies the inequalities.
\begin{equation*}
\mathfrak{p}=\nu \kappa ,\quad
\mathfrak{p}^{-1}\mathfrak{s}=\partial _{y_{1}}(\ln \kappa J)
\end{equation*}%
From this we conclude that
\begin{equation*}
\Big[\mathcal{G}_{\eta }\check{u}-\mathcal{J}^{\ast }\big(\frac{1}{\mathfrak{%
a}}\mathcal{J}\check{u}\big)\Big]\circ (\mathcal{U}_{2}\circ \mathcal{U}%
_{1})^{-1}=\kappa
\Big[\mathfrak{L}+\widetilde{\mathcal{G}}_{1}u+\nu
\partial _{y_{1}}(\ln \kappa J)\mathfrak{D}_{1}u+\mathfrak{B}u+\mathfrak{L}%
_{-1}u\Big]
\end{equation*}%
where
\begin{equation}
\mathfrak{B}=\frac{1}{\kappa }\mathcal{G}_{0},\text{ \ \ }\mathfrak{L}_{-1}=%
\frac{1}{\kappa }\mathcal{G}_{-1},\quad \widetilde{\mathcal{G}}_{1}=\frac{1}{%
\kappa }\mathcal{G}_{1}-(-\Delta )^{1/2}.  \label{identifoper2}
\end{equation}%
By construction, the operator $\widetilde{\mathcal{G}}_{1}$ has
the
following symbol%
\begin{equation*}
\left\{ (1+2a(Y))k_{1}^{2}+2b(Y)\tau k_{1}k_{2}+\tau
^{2}k_{2}^{2}\right\} ^{1/2}-\left\{ k_{1}^{2}+\tau
^{2}k_{2}^{2}\right\} ^{1/2},
\end{equation*}%
which implies that $\widetilde{\mathcal{G}}_{1}=\mathfrak{A}_{0}\mathfrak{D}%
_{1}$ where the zero order pseudodifferential operator
$\mathfrak{A}_{0}$
has the symbol%
\begin{equation*}
\frac{-2i(a(Y)k_{1}+\tau b(Y)k_{2})}{\left\{
(1+2a(Y))k_{1}^{2}+2b(Y)\tau k_{1}k_{2}+\tau ^{2}k_{2}^{2}\right\}
^{1/2}+\left\{ k_{1}^{2}+\tau ^{2}k_{2}^{2}\right\} ^{1/2}}.
\end{equation*}%
which can be written in the standard form
\begin{equation}
A_{0}(Y,\xi )=\frac{-2i\big(a(Y)\xi _{1}+b(Y)\xi _{2}\big)}{1+\big(%
1+2a(Y))\xi _{1}^{2}+2b(Y)\xi _{1}\xi _{2}+\xi
_{2}^{2}\big)^{1/2}}. \label{symbolazero}
\end{equation}%
Hence the basic identity \eqref{equchange1} holds with the operators $%
\mathfrak{B}$, $\mathfrak{L}_{-1}$ and the operator $\mathfrak{A}=\mathfrak{A%
}_{0}+\nu \partial _{y_{1}}(\ln g_{11})$ which yields
(\textit{i}).
Estimates (\textit{ii}) follow from formula \eqref{symbolazero}, estimates %
\eqref{estimab}, \eqref{estimg} \eqref{kappaJ} and Theorem
\ref{thmDir-Neum}.

Note that the functions $\kappa ,$ $J$, $a$ are even and the
function $b$ is odd both in $y_{1}$ and $y_{2}$. Moreover, they
are invariant under the shift $Y\rightarrow Y+(\pi ,\pi /\tau )$.
Hence the operator $\mathfrak{A}$
satisfies the symmetry conditions (\textit{iii}). Formulae %
\eqref{identifoper2} along with Theorem \ref{thmDir-Neum} imply
that the
operators $\mathfrak{B}$ and $\mathfrak{L}_{-1}$ also satisfy the conditions %
\eqref{irr1}, which completes the proof of Theorem
\ref{thmChangeCoord}.

\subsection{Main orders of the diffeomorphism and coefficient $\protect\nu $}

In using later the Nash-Moser theorem, we need to set%
\begin{equation*}
U=U_{\varepsilon }^{(N)}+\varepsilon ^{N}W
\end{equation*}%
where $U_{\varepsilon }^{(N)}$ is an approximate solution, up to order $%
\varepsilon ^{N}$of the system (\ref{basic1}), (\ref{basic2}), and
$W$ is the unknown perturbation. At Theorem \ref{Lembifurc}, for
$\varepsilon _{1}=\varepsilon _{2}=\varepsilon /2$ (Diamond
waves), we showed how to compute explicitly any order of
approximation $U_{\varepsilon }^{(N)},$ with $\mu -\mu _{c}(\tau
)=\varepsilon ^{2}\mu _{1}(\tau )+O(\varepsilon ^{4}),$ $\mu
_{1}(\tau )\neq 0$ for $\tau \neq \tau _{c}$. Then we need to
know in particular the principal part of the coefficient $\nu $ in (\ref{nu}%
) and (\ref{mathfrakL}), which implies the knowledge of the diffeomorphisms $%
\mathcal{U}_{1}$ and $\mathcal{U}_{2}.$ We show the following two
Lemmas:

\begin{lemma}
\label{Lemcoef-nu}For $N\geq 3,$ we have%
\begin{equation}
\nu (\varepsilon ,\tau ,W)=\mu _{c}^{-1}-\varepsilon ^{2}\nu
_{1}(\tau )+O(\varepsilon ^{3}),  \label{nu-dev}
\end{equation}%
where $\mu _{c}=(1+\tau ^{2})^{-1/2},$ and
\begin{eqnarray}
\nu _{1}(\tau ) &=&\frac{\mu _{1}(\tau )}{\mu _{c}^{2}}-\frac{3}{16\mu _{c}}+%
\frac{5}{4\mu _{c}^{3}}  \notag \\
&=&\frac{1}{4\mu _{c}^{5}}-\frac{1}{2\mu _{c}^{4}}+\frac{1}{2\mu _{c}^{3}}+%
\frac{7}{8\mu _{c}^{2}}-\frac{1}{4\mu _{c}}-\frac{9}{16(2-\mu
_{c})}, \label{nu_1}
\end{eqnarray}%
which is positive for any $\tau >0,$ $(\nu _{1}(\tau )>5/16\mu
_{c}).$
\end{lemma}

\begin{lemma}
\label{Lemmcoordchange} For $N\geq 3,$ the diffeomorphism
(\ref{varchange2}) of the torus, which allows to obtain the form
(\ref{equchange1}) for the
linear operator $\mathcal{G}_{\eta }-\mathcal{J}^{\ast }(\frac{1}{\mathfrak{a%
}}\mathcal{J}(\cdot ))$, and the principal part of
$\tilde{\eta}(Y),$ are
given by%
\begin{equation*}
X=(x_{1},x_{2}),\text{ \ \ \ }Y=(y_{1},y_{2}),\text{ \ \ }%
X(Y)=(X_{1}(Y),X_{2}(Y)),
\end{equation*}%
\begin{eqnarray*}
X_{1}(Y) &=&y_{1}+\frac{\varepsilon }{2}\sin y_{1}\cos
y_{2}+\varepsilon ^{2}(\xi _{11}\sin 2y_{1}+\xi _{12}\sin
2y_{1}\cos 2y_{2})+O(\varepsilon
^{3}), \\
X_{2}(Y) &=&\frac{1}{\tau }y_{2}+\varepsilon \tau \cos y_{1}\sin
y_{2}+\varepsilon ^{2}(\xi _{21}\sin 2y_{2}+\xi _{22}\cos
2y_{1}\sin
2y_{2})+O(\varepsilon ^{3}), \\
\tilde{\eta}(Y) &=&-\frac{\varepsilon }{\mu _{c}}\cos y_{1}\cos
y_{2}+\varepsilon ^{2}\{\eta _{00}+\eta _{01}\cos 2y_{1}+\eta
_{02}\cos 2y_{2}(1+\cos 2y_{1})\}+O(\varepsilon ^{3}),
\end{eqnarray*}%
where $\varepsilon >0$ is defined by the main order $\varepsilon
\xi _{0}$
of $U_{\varepsilon }^{(N)}$ and%
\begin{equation*}
\mu =\mu _{c}+\varepsilon ^{2}\mu _{1}(\tau )+O(\varepsilon ^{4}),\text{ \ }%
\mu _{c}=(1+\tau ^{2})^{-1/2},
\end{equation*}%
and where%
\begin{eqnarray*}
\xi _{11} &=&-\frac{1}{4}\left\{ \frac{3}{4(2-\mu _{c})}-\frac{3}{16}-\frac{5%
}{4\mu _{c}}+\frac{3}{4\mu _{c}^{2}}+\frac{1}{\mu
_{c}^{3}}-\frac{1}{\mu
_{c}^{4}}\right\} , \\
\xi _{12} &=&-\frac{1}{64},\text{ \ \ \ }\xi _{21}=-\frac{\tau }{8}+\frac{17%
}{32\tau },\text{ \ \ \ }\xi _{22}=-\frac{\tau }{8}, \\
\eta _{00} &=&\frac{1}{4\mu _{c}^{3}}-\frac{1}{8\mu _{c}},\text{ \
\ \ }\eta
_{02}=\frac{1}{8\mu _{c}}, \\
\eta _{01} &=&-\frac{1}{4\mu _{c}}-\frac{1}{4\mu
_{c}^{2}}+\frac{1}{4\mu _{c}^{3}}+\frac{3}{8(2-\mu _{c})}.
\end{eqnarray*}
\end{lemma}

\textbf{Proof:} see Appendix \ref{a3}.

\section{Small divisors. Estimate of $\mathfrak{L}-$ resolvent}

It follows from Theorem \ref{thmChangeCoord} that the
pseudodifferential operator $\mathfrak{L}=\mathcal{G}^{(0)}+\nu
\partial _{y_{1}}^{2}$ with the symbol
\begin{equation}
L(k)=-\nu k_{1}^{2}+(k_{1}^{2}+\tau ^{2}k_{2}^{2})^{1/2},\quad k\in \mathbb{Z%
}^{2}.  \label{140}
\end{equation}%
forms the principal part of the linear problem. In this section we
describe the structure of the spectrum of $\mathfrak{L}$ and
investigate in details the dependence of its resolvent on
parameters $\nu $ and $\tau $. The first result in this direction
is the following theorem which constitutes generic properties of
$\mathfrak{L}.$

\begin{theorem}
\label{141l} For every positive $\nu $ and $\tau $, the operator $\mathfrak{L%
}$ is selfadjoint in $H^{0}(%
\mathbb{R}
^{2}/\Gamma )$ and has the natural domain of definition
\begin{equation*}
D(\mathfrak{L})=\{u\in H^{0}(%
\mathbb{R}
^{2}/\Gamma ):\,\,\sum\limits_{k\in \mathbb{Z}^{2}}L(k)^{2}|\widehat{u}%
(k)|^{2}<\infty \}.
\end{equation*}%
The spectrum of $\mathfrak{L}$ coincides with the closure of the
discrete
spectrum which consists of all numbers $\lambda _{k}=L(k)$, $k\in \mathbb{N}%
^{2}$; the corresponding eigenfunctions are defined by
\begin{equation*}
(2\pi )^{-1}\exp (\pm ik\cdot Y),\quad (2\pi )^{-1}\exp (\pm
k^{\ast }\cdot Y)\text{~~where~~}k^{\ast }=(-k_{1},k_{2}).
\end{equation*}%
For each real $\varkappa \neq \lambda _{k}$, $k\in
\mathbb{Z}^{2}$, the resolvent $(\mathfrak{L}-\varkappa )^{-1}$ is
a selfadjoint unbounded operator defined in terms of the Fourier
transform by the formula
\begin{equation}
\widehat{(\mathfrak{L}-\varkappa )^{-1}u}(k)=(L(k)-\varkappa )^{-1}\widehat{u%
}(k),\quad k\in \mathbb{Z}^{2}.  \label{141}
\end{equation}
\end{theorem}

\begin{proof}
The operator $\mathfrak{L}$ is unitary equivalent to a
multiplication operator in the Hilbert space $l_{2}$, which can be
represented by the
infinite diagonal matrix with diagonal elements $L(k)$, $k\in \mathbb{Z}^{2}$%
. It remains to note that for such a matrix the spectrum coincides
with the closure of diagonal, and the discrete spectrum coincides
with the diagonal.
\end{proof}

Note that zero is nontrivial eigenvalue of $\mathfrak{L}$ if and
only if the dispersion equation $L(k)=0$ has a nontrivial solution
$k=k_{0}$. The number
of such solutions depends on the arithmetic nature of parameters $\nu $, $%
\tau $. We restrict our attention to the simplest case when
$k_{0}=(1,1)$. It is easy to see that in this case the point $(\nu
,\tau )=(\nu _{0},\tau )$ belongs to the positive branch of the
hyperbola $\nu _{0}^{2}-\tau ^{2}=1$.

Our aim is to investigate in details the dependence of $\mathfrak{L}$%
-resolvent on parameter $\nu .$ With further applications to the
Nash-Moser theory in mind we take the perturbed values of
parameter $\nu $, and a spectral parameter $\varkappa $ in the
form
\begin{equation}
\nu _{j}(\varepsilon )=\nu _{0}-\varepsilon ^{2}\nu _{1}+\varepsilon ^{3}%
\tilde{\nu}_{j}(\varepsilon ),\quad \varkappa _{j}(\varepsilon
)=\varepsilon ^{2}\tilde{\varkappa}_{j}(\varepsilon ),
\label{perturbation}
\end{equation}%
where $\nu _{0}=\nu _{0}(\tau )$ and $\nu _{1}=\nu _{1}(\tau ).$
Here functions $\tilde{\nu}_{j}$ and $\tilde{\varkappa}_{j}$,
$j\geq 1$, are defined on the segment $[0,r_{0}]$ and satisfy the
inequalities
\begin{gather}
|\tilde{\nu}_{j}(\varepsilon )|+|\tilde{\varkappa}_{j}(\varepsilon
)|\leq R,
\notag \\
|\tilde{\nu}_{j}(\varepsilon ^{\prime
})-\tilde{\nu}_{j}(\varepsilon
^{\prime \prime })|+|\tilde{\varkappa}_{j}(\varepsilon ^{\prime })-\tilde{%
\varkappa}_{j}(\varepsilon ^{\prime \prime })|\leq R|\varepsilon
^{\prime
}-\varepsilon ^{\prime \prime }|,  \label{perturbationb} \\
|\tilde{\nu}_{j+1}-\tilde{\nu}_{j}|+|\tilde{\varkappa}_{j+1}-\tilde{\varkappa%
}|\leq R(2^{-j}).  \notag
\end{gather}%
Denoting by $\mathfrak{L}_{0}$ the operator $\mathfrak{L}$ for
$\nu =\nu _{0}(\tau ),$ the next theorem establishes the basic
estimates for a
resolvent $(\mathfrak{L}-\varkappa )^{-1}$ on the orthogonal complement to $%
\text{ker~}\mathfrak{L}_{0}$, which are stable with respect to
perturbations of parameters from a suitable Cantor set.

\begin{theorem}
\label{142t}

\begin{itemize}
\item[(\textit{a})] For each $\alpha \in (0,1]$, there is a set of
full measure $\mathfrak{N}_{\alpha }\subset (1,\infty )$ so that
whenever $\nu _{0}\in \mathfrak{N}_{\alpha }$ and $\nu =\nu
_{0}(\tau )=(1+\tau ^{2})^{1/2} $, zero is a non-trivial
eigenvalue of $\mathfrak{L}_{0}$ and
\begin{gather}
\text{\textrm{ker~}}\mathfrak{L}_{0}=\text{\textrm{span~}}\big\{1,\exp
(\pm
ik_{0}\cdot Y),\quad \exp (\pm k_{0}^{\ast }\cdot Y)\big\},  \notag \\
\Vert \mathfrak{L}_{0}^{-1}u\Vert _{s-(1+\alpha )/2}\leq c\Vert u\Vert _{s},%
\text{~~when~~}u\in (\text{\textrm{ker~}}\mathfrak{L}_{0})^{\perp
}\cap H^{s}(\mathbb{R}^{2}/\Gamma ).  \label{142}
\end{gather}

\item[(\textit{b})] Suppose that $\nu _{0}\in \mathfrak{N}_{\alpha }$ with $%
\alpha \in (0,1/78)$, then there is a set $\mathcal{E}\subset
\lbrack 0,r_{0})$ so that
\begin{equation}
\frac{2}{r^{2}}\int\limits_{[0,r]\cap \mathcal{E}}\varepsilon
\,d\varepsilon \rightarrow 1\text{~~as~~}r\rightarrow 0,
\label{1425}
\end{equation}%
and for all $\varepsilon \in \mathcal{E}$, $j,s\geq 1$ and $\nu
=\nu _{j}(\varepsilon )$, $\varkappa =\varkappa _{j}(\varepsilon
)$,
\begin{equation}
\big\|\big(\mathfrak{L}-\varkappa \big)^{-1}u\big\|_{s-1}\leq
c\Vert u\Vert _{s}\text{~~when~~}u\in
(\text{\textrm{ker~}}\mathfrak{L}_{0})^{\perp }\cap
H^{s}(\mathbb{R}^{2}/\Gamma ),  \label{1426}
\end{equation}%
where the positive constant $c$ depends on $\alpha $, $\nu _{0}$
and $R$ only.
\end{itemize}
\end{theorem}

\begin{remark}
The fact proved at Lemma \ref{Lemcoef-nu} that $\nu _{1}(\tau
)>0$, allows the freedom to take $\tau $ as a function of
$\varepsilon $ as $\tau =\tau _{0}+\varepsilon ^{2}\tau _{1}.$
Then, for $|\tau _{1}|$ small enough, we have
\begin{equation*}
\omega _{1}:=\omega _{0}(\nu _{0}^{-1}\nu _{1}-\tau _{0}^{-1}\tau
_{1})\neq 0
\end{equation*}%
and Theorem \ref{142t} still holds true.
\end{remark}

It follows from formula \eqref{141} and the Parseval identity that Theorem %
\ref{142t} will be proved once we prove the following

\begin{theorem}
\label{142at}

\begin{itemize}
\item[(\textit{a})] For any $\alpha \in (0,1]$ there is a set of
full measure $\mathfrak{N}_{\alpha }\subset (1,\infty )$ so that
for all $\nu _{0}\in \mathfrak{N}_{\alpha }$ and $\nu =\nu
_{0}(\tau )$,
\begin{equation}
|L(k)|\geq c|k|^{-(1+\alpha )/2}\text{~~when~~}k\neq 0,\pm
k_{0},\pm k_{0}^{\ast }.  \label{142ab}
\end{equation}

\item[(\textit{b})] Suppose that $\nu _{0}\in \mathfrak{N}_{\alpha }$ with $%
\alpha \in (0,1/78),$ then, there is a set $\mathcal{E}\subset
\lbrack
0,r_{0})$ satisfying \eqref{1425} so that for all $\varepsilon \in \mathcal{E%
}$, $j\geq 1$ and $\nu =\nu _{j}(\varepsilon )$, $\varkappa
=\varkappa _{j}(\varepsilon )$,
\begin{equation}
\big|(L(k)-\varkappa )^{-1}\big|\geq c|k|^{-1}\text{~~for
all~~}k\neq 0,\pm k_{0},\pm k_{0}^{\ast },  \label{1426b}
\end{equation}%
where the positive constant $c$ depends on $\nu _{0}$, $\alpha $
and $R$ only.
\end{itemize}
\end{theorem}

The proof naturally falls into four steps.

\paragraph{First step.}

We begin with proving two auxiliary lemmas which establish a
connection between the symbol $L$ and the Diophantine function
defined by the formula
\begin{equation}
(m,n)\rightarrow \omega -\frac{m}{n^{2}}-\frac{C}{n^{2}},\quad
(m,n)\in \mathbf{N}^{2},  \label{1430}
\end{equation}%
where constants $\omega $, $C$ are connected with parameters $\nu
$, $\tau $ and $\varkappa $ by the relations
\begin{equation*}
\omega =\nu /\tau ,\quad C=1/(2\nu \tau )-\varkappa /\tau .
\end{equation*}

\begin{lemma}
\label{145l} Let for some $\varrho \geq 2$ and $\alpha \in \lbrack
0,1]$, parameters $\nu $, $\tau $ and a vector $k\in
\mathbb{Z}^{2}$ satisfy the inequalities
\begin{gather}
0<\varrho ^{-1}\leq \nu ,\tau <\varrho ,\quad |\varkappa |<\varrho
,\quad
|k_{1}|\geq 2\varrho ^{2},  \label{146a} \\
\big|\omega k_{1}^{2}-|k_{2}|-C\big|\geq 5(\varrho ^{2}+\varrho
^{10})|k_{1}|^{-1-\alpha },  \label{147}
\end{gather}%
Then $|L(k)-\varkappa |\geq |k|^{-(1+\alpha )/2}$.
\end{lemma}

\begin{proof}
Since $L(k)$ is even in $k$, it suffices to prove the lemma for
\begin{equation*}
k=(n,m)\text{~~with integers~~}n>2\varrho ^{2},\quad m>0.
\end{equation*}%
We begin with the observation that for $n^{2}\geq 2\varrho ^{2}m$,
\begin{gather*}
-L(k)+\varkappa \geq \varrho ^{-1}n^{2}-\sqrt{\varrho ^{2}m^{2}+n^{2}}%
-\varrho \geq \\
4\varrho ^{4}\Big(\frac{1}{\varrho }-\frac{1}{2}\big(\frac{1}{\varrho ^{2}}+%
\frac{1}{\varrho ^{4}}\big)^{1/2}\Big)-\varrho \geq
\frac{4}{3}\varrho ^{3}-\varrho >1.
\end{gather*}%
Hence it suffices to prove the lemma for
\begin{equation}
n^{2}\leq 2\varrho ^{2}m.  \label{148}
\end{equation}%
Since for all non-negative $\sigma $, $0\leq 1+\sigma /2-(1+\sigma
)^{1/2}\leq \sigma ^{2}/4$, we have
\begin{equation*}
-\tau ^{-1}L(k)=\omega n^{2}-m-(2\tau ^{2}m)^{-1}n^{2}+(4\tau
^{4}m^{3})^{-1}n^{4}\ o_{1},\quad o_{1}\in \lbrack 0,1],
\end{equation*}%
which along with the identity
\begin{equation*}
(2\tau ^{2}m)^{-1}n^{2}=(2\nu \tau )^{-1}+(2\nu \tau
m)^{-1}(\omega n^{2}-m)
\end{equation*}%
yields
\begin{equation}
-\tau ^{-1}e(k)\big(L(k)-\varkappa \big)=\omega n^{2}-m-C+o_{2}.
\label{149}
\end{equation}%
Here
\begin{equation*}
o_{2}=C(1-e(k))+e(k)(4\tau ^{4}m^{3})^{-1}n^{4}o_{1},\quad e(k)=\big(%
1-(2\tau \nu m)^{-1}\big)^{-1}.
\end{equation*}%
Since
\begin{equation*}
(2\nu \tau m)^{-1}\leq 2^{-1}\varrho ^{2}m^{-1}\leq \varrho
^{4}n^{-2}\leq 2^{-1},
\end{equation*}%
we have
\begin{equation*}
1\leq e(k)\leq 2,\quad e(k)-1\leq 2\varrho ^{4}n^{-2},
\end{equation*}%
which along with the inequality $|C|\leq 3\varrho ^{2}/2$ leads to
the estimate
\begin{equation*}
|o_{2}|\leq \Big(3\varrho ^{6}+\frac{\varrho ^{4}}{2}\frac{n^{6}}{m^{3}}\Big)%
\frac{1}{n^{2}}\leq 5\varrho ^{10}n^{-2}.
\end{equation*}%
From this, \eqref{147} and \eqref{149} we conclude that
\begin{gather*}
2\varrho |L(k)-\varkappa |\geq |\omega n^{2}-m-C|-2\varrho ^{10}n^{-2}\geq \\
5(\varrho ^{2}+\varrho ^{10})n^{-1-\alpha }-5\varrho
^{10}n^{-2}\geq 5\varrho ^{2}n^{-1-\alpha }.
\end{gather*}%
It remains to note that due \eqref{148} the right hand side is larger than $%
|k|^{-(1+\alpha )/2}$ and the lemma follows.
\end{proof}

\begin{lemma}
\label{149l} Suppose that parameters $\nu_0, \tau_0$ and
$\nu,\tau,\varkappa$ satisfy the following conditions.

\begin{itemize}
\item[(\textit{i})] There is $\varrho \geq 2$ so that
\begin{equation*}
0<\varrho ^{-1}\leq \nu ,\tau ,\nu _{0}<\varrho ,\quad |\varkappa
|<\varrho .
\end{equation*}

\item[(\textit{ii})] The dispersion equation
\begin{equation}
\nu _{0}k_{1}^{2}-\sqrt{k_{1}^{2}+\tau ^{2}k_{2}^{2}}=0
\label{dispersion0}
\end{equation}%
has the unique positive solution $k=k_{0}$.

\item[(\textit{iii})] There are positive $N$ and $q$ so that
\begin{gather*}
N\geq 2\varrho ^{2},\quad q\geq 5(\varrho ^{2}+\varrho ^{10}), \\
|\omega n^{2}-m-C|\geq qn^{-1-\alpha }\text{~~ for all
integers~~}n\geq N,\,m\geq 0
\end{gather*}
\end{itemize}

Then there are positive $\delta $ and $\gamma _{0}$, depending on
$N$, $q$ and $\nu _{0}$ only, so that
\begin{equation}
|L(k)-\varkappa |\geq \gamma _{0}|k|^{-(1+\alpha )/2}\text{~~for all ~~}%
k\neq 0,\pm k_{0},\pm k_{0}^{\ast }  \label{144a}
\end{equation}%
when
\begin{equation}
|\nu -\nu _{0}|+|\varkappa |\leq \delta .  \label{144ab}
\end{equation}
\end{lemma}

\begin{proof}
It follows from (\textit{iii}) that $\nu $, $\tau $ and $\varkappa
$ meet al requirements of Lemma \ref{145l} which implies that
$|L(k)-\varkappa |\geq |k|^{-(1+\alpha )/2}$ for all $|k_{1}|\geq
N$. On the other hand, since
\begin{equation*}
\lim\limits_{|k_{2}|\rightarrow \infty }|L(k)-\varkappa |\geq
\lim\limits_{|k_{2}|\rightarrow \infty }\sqrt{\varrho
^{-2}k_{2}^{2}+k_{1}^{2}}-\varrho k_{1}^{2}-\varrho =\infty ,
\end{equation*}%
there is $N^{\ast }$ depending on $\varrho $ and $N$ only so that $%
|L(k)-\varkappa |\geq |k|^{-(1+\alpha )/2}$ for all $|k|\geq
N^{\ast }$. Since the number of wave vectors $k$ in the circle
$|k|\leq N$ is finite, the existence $\delta $ and $\gamma _{0}$
follows from continuity of $L(k)$ as a function of parameters $\nu
,\tau $.
\end{proof}

\paragraph{Second step.}

Next we prove that condition (\textit{iii}) from the previous
lemma is the
generic property of function \eqref{1430} which leads to assertion (\textit{%
a)} of Theorem \ref{142at}. In order to formulate the results let
us introduce the important function $d:\mathbb{N}^{2}\mapsto
\mathbb{R}$ defined by the formula
\begin{equation}
d(m,n)=\omega _{0}-\frac{m}{n^{2}}-\frac{C_{0}}{n^{2}},\text{~~where~~}%
\omega _{0}=\nu _{0}/\tau ,C_{0}=(2\nu _{0}\tau )^{-1}. \label{DI}
\end{equation}

\begin{lemma}
\label{1410l} For each $\alpha \in (0,1]$ and $q>0$ there is a set
of a full measure $\mathfrak{M}_{\alpha }$ in $(1,\infty )$ so
that for all $\nu
_{0}\in \mathfrak{M}_{\alpha }$, there exists $N>0$, depending on $\nu _{0}$%
, $\alpha $ and $q$ only, such that
\begin{equation}
|d(m,n)|\geq qn^{-3-\alpha }\text{~~when~~}n\geq N,m\geq 0.
\label{1411}
\end{equation}
\end{lemma}

\begin{proof}
Noting that $C_{0}=(\omega _{0}-\omega _{0}^{-1})/2,$ we rewrite inequality %
\eqref{1411} in the equivalent form
\begin{equation*}
|d_{0}(\omega _{0},m,n)|\geq qn^{-3-\alpha
}\text{~~where~~}d_{0}(\omega _{0},m,n)=\omega
_{0}-\frac{1}{n^{2}}\big(m-2^{-1}\omega _{0}+2^{-1}\omega
_{0}^{-1}\big).
\end{equation*}%
Without loss of generality we can assume also that $\omega _{0}\in
\lbrack 1,a]$, where $a$ is an arbitrary positive number. The set
of points $\omega _{0}$ for which $|d_{0}(\omega _{0},m,n)|\leq
qn^{-3-\alpha }$ can be covered by the system of the intervals
$\iota (m,n)$ labelled by integers $m$ and $n$. It is easy to see
that their extremities $\omega _{0}^{\pm }(m,n)$
satisfy $\Big(n^{2}+\frac{1}{2}\Big)\omega _{0}^{\pm }(m,n)\geq m-\frac{q}{%
n^{1+\alpha }}.$ On the other hand, since
\begin{equation*}
\partial _{\omega }d_{0}(\omega _{0},m,n)=1+(2n^{2})^{-1}+(2\omega
_{0}^{2}n^{2})^{-1}\in \lbrack 1,2],
\end{equation*}%
the length of each intervals is less than $2qn^{-3-\alpha }.$
Hence for fixed $n$, the number of intervals $\iota (m,n)$ having
nonempty
intersections with a segment $[1,a]$ is less than $a(n^{2}+1/2)+q/{%
n^{1+\alpha }}$. From this we conclude that
\begin{equation*}
\sum\limits_{n\geq N}\sum\limits_{m:\iota (m,n)\cap \lbrack
1,a]\neq
\emptyset }\text{\textrm{~meas~}}\iota (m,n)\leq 2q(\frac{3}{2}%
a+q)\sum\limits_{n\geq N}n^{-1-\alpha }\leq cq(a+q)N^{-\alpha }.
\end{equation*}

Hence

\begin{equation*}
\begin{split}
\text{\textrm{meas~}}\big\{\omega _{0}& \in \lbrack
1,a]:\,|d_{0}(\omega
,m,n)|\geq qn^{-3-\alpha }\text{ ~~for all~~}m\geq 0,n\geq N\big\}\geq \\
a-cq(a+q)N^{-\alpha }& \rightarrow a\text{~~as~~}N\rightarrow
\infty .
\end{split}%
\end{equation*}%
It remains to note that the mapping $\nu _{0}\rightarrow \omega
_{0}$ takes diffeomorphically the interval $[1,\infty )$ onto
itself and the lemma follows.
\end{proof}

\noindent Next lemma shows that almost each point of
$\mathfrak{M}_\alpha$ satisfies the following \noindent

\textbf{Condition M. }\textit{There are positive constant $c_{0}$,
$c_{1}$ such that for all integers $n>0$, $m\geq 0$,
\begin{gather}
|d(m,n)|\geq c_{0}n^{-3-\alpha },  \label{DS} \\
\left\vert e^{2\pi i\omega _{0}n}-1\right\vert ^{-1}\leq
c_{1}n^{2} \label{DSb}
\end{gather}%
}

\noindent

\begin{lemma}
\label{conditm} For each $\alpha \in (0,1)$ there is a set $\mathfrak{N}%
_{\alpha }\subset \mathfrak{M}_{\alpha }$ such that $\text{\textrm{~meas~}}%
\big(\mathfrak{M}_{\alpha }\setminus \mathfrak{N}_{\alpha })=0$
and Condition M holds for each $\nu _{0}\in \mathfrak{N}_{\alpha
}$.
\end{lemma}

\begin{proof}
We begin with the observation that inequality \eqref{DSb} holds
for almost every $\omega _{0}$ with a constant $c_{1}$ depending
on $\omega _{0}$ only.
Hence it suffices to show that inequality \eqref{DS} holds true for each $%
\nu _{0}\in \mathfrak{M}_{\alpha }$. Fix $q=1$, $\nu \in \mathfrak{M}%
_{\alpha }$ and note that by lemma \eqref{1410l}, $|d(m,n)|\geq
n^{-3-\alpha }$ for all $n\geq N.$ Since
$\lim\limits_{m\rightarrow \infty }d(m,n)=-\infty $, we can choose
$N$ such that this inequality is fulfilled for all $(m,n)$ outside
of the simplex $n>0,m\geq 0$, $n+m\leq N$. It remains to note that
this simplex contains only finite number of integer points and
$d(m,n)\neq 0$ for irrational $\omega _{0}$ different from the sum
of a rational and the square root of a rational.
\end{proof}

\paragraph{ Third step}

We intend now to study the robustness of estimate (\ref{DS}) when
one adds a small perturbation to $d(m,n)$. In order to formulate
the corresponding result we introduce some notations. For each
$\varepsilon \in \lbrack 0,r_{0}]$ and $j\geq 1$ set
\begin{equation*}
\tilde{\omega}_{j}(\varepsilon )=\nu _{j}(\varepsilon )/\tau ,\quad \tilde{C}%
_{j}(\varepsilon )=\big(2\nu _{j}(\varepsilon )\tau
\big)^{-1}-\varkappa _{j}(\varepsilon )\tau ^{-1}
\end{equation*}%
It follows from \eqref{perturbation} that they have the
representation
\begin{equation*}
\tilde{\omega _{j}}(\varepsilon )=\omega _{0}-\varepsilon
^{2}\omega
_{1}+\varepsilon ^{3}\tilde{\Omega}_{j}(\varepsilon ),\tilde{C_{j}}%
(\varepsilon )=C_{0}-\varepsilon ^{2}\tilde{\varphi
_{j}}(\varepsilon )\quad \varepsilon \in \lbrack 0,r_{0}],
\end{equation*}%
in which $\omega _{1}=\tau ^{-1}\nu _{1}$. Our task is to obtain
the
estimates for the function $\tilde{D_{j}}:[0,r_{0}]\times \mathbb{N}%
^{2}\rightarrow \mathbb{R}$ defined by
\begin{equation}
\tilde{D}(\varepsilon ,m,n,)=\tilde{\omega _{j}}(\varepsilon )-\frac{m}{n^{2}%
}-\frac{C_{0}}{n^{2}}+\varepsilon ^{2}\frac{\tilde{\varphi
_{j}}(\varepsilon )}{n^{2}}.  \label{DC}
\end{equation}%
For technical reason it is convenient to formulate the problem in
terms of a new small parameter $\lambda $. Set
\begin{gather}
\lambda =\varepsilon ^{2},\quad \varphi _{j}(\lambda )=\tilde{\varphi _{j}}(%
\sqrt{\lambda }),\quad \Omega _{j}(\lambda )=\tilde{\Omega _{j}}(\sqrt{%
\lambda }),  \notag \\
{\omega _{j}}(\lambda )=\omega _{0}-\lambda \omega _{1}+\lambda
^{3/2}\Omega _{j}(\lambda ),\quad \lambda \in \lbrack 0,\rho
_{0}],  \label{001}
\end{gather}%
where $\rho _{0}=r_{0}^{2}$. It follows from \eqref{perturbationb}
that for suitable choice $R$ and $r_{0}$,
\begin{equation}
|\Omega _{j}|+|\varphi _{j}|\leq R,\text{ \ }|\Omega _{j+1}-\Omega
_{j}|+|\varphi _{j+1}-\varphi _{j}|\leq R(2^{-j}),\quad |\Omega
_{j}^{\prime
}(\lambda )|+|\varphi _{j}^{\prime }(\lambda )|\leq \frac{R}{2\sqrt{\lambda }%
}.  \label{HYPSA}
\end{equation}%
Set
\begin{equation}
D_{j}(\lambda ,m,n,):=\tilde{D_{j}}(\sqrt{\lambda },m,n)=\omega
_{j}(\lambda
)-\frac{m}{n^{2}}-\frac{C_{0}}{n^{2}}+\lambda \frac{\varphi _{j}(\lambda )}{%
n^{2}}.  \label{DCL}
\end{equation}

\begin{definition}
For positive $N$, $q$, $r$ denote by $\mathcal{H}_j(N,q,r)$ the
set defined by
\begin{equation}  \label{30}
\mathcal{H}_j(N,q,r)=\Bigl\{\ll\in (0,r): |D_j(\ll,m,n)|\geq qn^{-4}\text{%
~for all integers~}n\geq N, \quad m\geq 0\Bigr\}.
\end{equation}
\end{definition}

\begin{theorem}
\label{thm1} Assume that $\nu _{0}\in \mathfrak{N}_{\alpha }$ with
$\alpha \in (0,1/78).$ Then, for each $q>0$, there is $N>0$ such
that
\begin{equation}
\frac{1}{r}\text{\textrm{~meas~}}\bigcap\limits_{j\geq 1}\{\mathcal{H}%
_{j}(N,q,r)\}\rightarrow 1\text{~~as~~}r\rightarrow 0,  \label{31}
\end{equation}%
and there exists $c^{\ast }$ such that
\begin{equation}
\frac{1}{r}\text{\textrm{~meas~}}\bigcap\limits_{j\geq 1}\{\mathcal{H}%
_{j}(N,q,r)\}\geq 1-c^{\ast }r^{\varpi }\text{~~for~~}0<r<r^{\ast
}, \label{311}
\end{equation}%
where $0<\varpi <\min \{\alpha /(3-\alpha ),(78^{-1}-\alpha
)/(3+\alpha )\}$.
\end{theorem}

\begin{proof}
Subsection \ref{small2} is devoted to the proof.
\end{proof}

\paragraph{Fourth step.}

We are now in a position to complete the proof of Theorem
\ref{142at}. First note that by Lemmas \ref{1410l} and
\ref{conditm} for each $\nu _{0}\in \mathfrak{N}_{\alpha }$,
parameters $\nu =\nu _{0}$($\tau )$ and $\varkappa =0$ meet all
requirements of Lemma \ref{149l} with $\delta =0$. Applying
this lemma we obtain \eqref{142ab}. It remains to note that for almost all $%
\nu _{0}$ the dispersion equation has the only positive solution
$k_{0}$ and assertion (\textit{a}) follows. In order to prove
(\textit{b}) choose an arbitrary $\nu _{0}\in \mathfrak{N}_{\alpha
}$, $q>0$ and set
\begin{equation*}
\mathcal{E}=\big\{\varepsilon >0:\,\varepsilon ^{2}\in
\bigcap\limits_{j\geq 1}\mathcal{H}_{j}(N,q,r)\big\},
\end{equation*}%
where $N$ is given by Theorem \ref{thm1}. It follows from this
theorem that the set $\mathcal{E}$ satisfies density condition
\eqref{1425}. On the other
hand, \eqref{30} implies that for $\varepsilon \in \mathcal{E}$, parameters $%
\nu =\nu _{j}(\varepsilon )$, $\varkappa =\varkappa
_{j}(\varepsilon )$ meet all requirements of Lemma \ref{149l}.
Since
\begin{equation*}
|\nu _{j}(\varepsilon )-\nu _{0}|+|\varkappa _{j}(\varepsilon )|\rightarrow 0%
\text{~~as~~}\varepsilon \rightarrow 0
\end{equation*}%
uniformly with respect to $j$, there exists $\varepsilon _{2}>0$
depending
on $\nu _{0}$ and $R$ only such that for all $\varepsilon \in \mathcal{E}%
\cap (0,\varepsilon _{2})$, parameters $\nu =\nu _{j}(\varepsilon )$ and $%
\varkappa =\varkappa _{j}(\varepsilon )$ satisfy inequality
\eqref{144ab} which yields \eqref{144a} (where $\alpha =1),$ and
the theorem follows.

\subsection{Proof of Theorem \protect\ref{thm1}}

\label{small2}

Our approach is based on standard methods of metric theory of
Diophantine
approximations and Weil Theorem on the uniform distribution of a sequence $%
\{\omega _{0}n^{2}\}$. Without loss of generality we can assume
that $\omega _{1}>0$, and $\rho _{0}=r_{0}^{2}$ satisfies the
inequality
\begin{equation}
\rho _{0}<\frac{\omega _{0}}{4}(\omega _{1}+R)^{-1}\Rightarrow
\omega (\lambda )\geq \frac{3\omega _{0}}{4}\text{~~for ~~}\lambda
\leq \rho _{0}. \label{rr}
\end{equation}%
It follows from \eqref{HYPSA} that the sequences $\Omega _{j}$and
$\varphi _{j}$ converge uniformly on $[0,\rho _{0}]$ to functions
$\Omega _{\infty }$ and $\varphi _{\infty }$ such that
\begin{equation}
|\Omega _{\infty }|+|\varphi _{\infty }|\leq R,|\Omega _{\infty
}-\Omega _{j}|+|\varphi _{\infty }-\varphi _{j}|\leq
2^{1-j}R,|\Omega _{\infty }^{\prime }(\lambda )|+|\varphi _{\infty
}^{\prime }(\lambda )|\leq 2^{-1}\lambda ^{-1/2}R.  \label{infty}
\end{equation}%
Our task is to estimate the measure of intersection of the sets $\mathcal{H}%
_{j}$. We begin with the investigation of their structure.

\paragraph{Structure of a set $\mathcal{H}_j(N,q,r)$.}

The main result of this paragraph is the following covering lemma.
Fix an arbitrary positive $q$ and set
\begin{equation}
r_{2}=\min \Bigl\{\frac{\rho _{0}}{2},\,\frac{\omega _{1}^{2}}{128R^{2}}%
\Bigr\},\text{ \ }N_{3}(q)=\max \Bigl\{\left( \frac{16R}{\omega
_{1}}\right) ^{1/2},\Bigl(\frac{2q}{c_{0}}\Bigr)^{1/(1-\alpha
)},\Bigl(\frac{4q}{\omega _{1}r_{2}}\Bigr)^{1/4}\Bigr\},
\label{32}
\end{equation}%
where $c_{0}$ is the constant from Condition $M$.

\begin{lemma}
\label{LO} For $N>N_{3}(q)$, $r<r_{2}$ and $1\leq j\leq \infty $, the set $%
[0,r]\setminus \mathcal{H}_{j}(N,q,r)$ is covered by the system of
the intervals
\begin{equation*}
I_{j}(m,n)=(\lambda _{j}^{-}(m,n),\lambda _{j}^{+}(m,n)),\quad
n\geq N\quad I_{j}(m,n)\cap \lbrack 0,r]\neq \emptyset ,
\end{equation*}%
such that

\begin{itemize}
\item[\textit{(i)}] $\lambda_j^{\pm}(m,n)$ are solutions of the
equations
\begin{equation}
\lambda_j^{\pm}\Bigl(\omega _{1}-\frac{\varphi_j (\lambda_j^{\pm})}{n^{2}}%
\Bigr)-(\lambda_j^{\pm})^{3/2}\Omega (\lambda_j^{\pm})=d(m,n)\pm \frac{q}{%
n^{4}}  \label{42}
\end{equation}%
with $d(m,n)>0$. They satisfy the inequalities
\begin{gather}
0<\lambda_j^{-}(m,n)<\lambda_j^{+}(m,n)<2r_{2},  \label{50} \\
\frac{2}{5\omega _{1}}d(m,n)\leq \lambda_j^{\pm}\leq \frac{2}{\omega _{1}}%
d(m,n)  \label{41} \\
\frac{4q}{3\omega _{1}}\frac{1}{n^{4}}\leq
\lambda_j^{+}(m,n)-\lambda _j^{-}(m,n)\leq \frac{4q}{\omega
_{1}}\frac{1}{n^{4}}.  \label{050}
\end{gather}

\item[\textit{(ii)}] For a fixed $n>N$, the left extremities
$\lambda _{-}^{j}(m,n)$ strongly decreases in $m$,
\begin{gather}
\lambda _{j}^{-}(M_{j,n}(r),n)\leq \lambda
_{j}^{-}(M_{j,n}(r)-1,n)\leq
...\leq \lambda ^{-}(m_{j,n}(r),n),  \label{451} \\
\lambda _{j}^{-}(k-1,n)-\lambda _{j}^{-}(k,n)\geq \frac{2}{3\omega _{1}n^{2}}%
,  \label{452}
\end{gather}%
where
\begin{eqnarray*}
m_{j,n}(r) &=&\min \{m>0:I_{j}(m,n)\cap \lbrack 0,r]\neq \emptyset \}, \\
M_{j,n}(r) &=&\max \{m>0:I_{j}(m,n)\cap \lbrack 0,r]\neq \emptyset
\}.
\end{eqnarray*}

\item[\textit{(iii)}] For each such interval with $I_j(m,n)\cap
\lbrack 0,r]\neq \emptyset $ ,
\begin{equation}
d(m,n)\leq \frac{5}{2}\omega _{1}r,\quad n\geq \left(
\frac{2c_{0}}{5\omega _{1}}\right) ^{1/(3+\alpha )}r^{-1/(3+\alpha
)},\lambda_j^{+}(m,n)\leq 2r. \label{45}
\end{equation}

\item[\textit{(iv)}] If intervals $I_j(m,n)$ and $I_\infty$ have
nonempty intersections with $(0,r]$, then
\begin{equation}  \label{LINF}
|\lambda_j^\pm(m,n)-\lambda_\infty^\pm(m,n)|\leq 2^{-j-4}\omega_1.
\end{equation}
\end{itemize}
\end{lemma}

\begin{proof}
By abuse of notations, we suppress the index $j$. Proof of
(\textit{i}).
\begin{equation*}
-\partial _{\lambda }D(\lambda ,m,n)=\omega _{1}-\frac{\varphi
+\lambda \varphi ^{\prime }}{n^{2}}-\frac{3}{2}\lambda
^{1/2}\Omega (\lambda )-\lambda ^{3/2}\Omega ^{\prime }(\lambda ),
\end{equation*}%
it follows from \eqref{32}and \eqref{HYPSA} that for $\lambda \leq
2r_{2}$
\begin{equation*}
\omega _{1}-\frac{2R}{n^{2}}-2R\lambda ^{1/2}\leq -\partial
_{\lambda }D(\lambda ,m,n)\leq \omega
_{1}+\frac{2R}{n^{2}}+2R\lambda ^{1/2},
\end{equation*}%
which along with \eqref{32} implies the inequalities
\begin{equation}
\frac{\omega _{1}}{2}\leq -\partial _{\lambda }D(\lambda ,m,n)\leq \frac{%
3\omega _{1}}{2}\text{~~for~~}\lambda \in \lbrack 0,2r_{2}],\quad
n>N_{3}. \label{48}
\end{equation}%
Hence the function $-D(\lambda ,m,n)$ is strongly monotone on the interval $%
(0,2r_{2})$. Therefore for $r<r_{2}$,$N\geq N_{3}$ the set
$(0,r)\setminus \mathcal{H}(N,q,r)$ can be covered by the system
of the intervals
\begin{equation*}
J(m,n)=\bigl(\beta (m,n),\gamma (m,n)\bigr),\quad 0\leq \beta
(m,n)<\gamma (m,n)\leq r,
\end{equation*}%
such that
\begin{gather}
-D(\beta (m,n),m,n)=-\frac{q}{n^{4}}\text{~~if~~}\beta (m,n)>0,
\label{51}
\\
-D(\gamma (m,n),m,n)=+\frac{q}{n^{4}}\text{~~if~~}\gamma (m,n)<r,
\label{52}
\\
-\frac{q}{n^{4}}\leq -D(\beta
(m,n),m,n)<\frac{q}{n^{4}}\text{~~if~~}\beta
(m,n)=0  \label{53} \\
-\frac{q}{n^{4}}<-D(\gamma (m,n),m,n)\leq
\frac{q}{n^{4}}\text{~~if~~}\gamma (m,n)=r,  \label{54}
\end{gather}%
Note that, by condition $M$,
\begin{equation*}
\frac{q}{n^{4}}=\frac{q}{c_{0}n^{1-\alpha
}}\frac{c_{0}}{n^{3+\alpha }}\leq \frac{q}{c_{0}n^{1-\alpha
}}|d(m,n)|,
\end{equation*}%
which along with \eqref{32} yields the inequalities
\begin{equation}
\frac{1}{2}|d(m,n)|\leq |d(m,n)\pm \frac{q}{n^{4}}|\leq \frac{3}{2}%
|d(m,n|,\quad |d(m,n)|\geq 2\frac{q}{n^{4}}\text{~~for~~}n\geq
N_{3}. \label{44}
\end{equation}%
From this and the equality $D(0,m,n)=d(m,n)$ we conclude that case
\eqref{53} is impossible and $\beta (m,n)>0$ for all intervals
$J(m,n)$. Let us consider case \eqref{54}. Since the function
$-D(\lambda ,m,n)$ increases in $\lambda $ on the segment
$[0,2r_{2}]$, there is a maximal $\gamma \ast $ in $(0,2r_{2}]$
such that
\begin{eqnarray*}
(\beta (m,n),r) &\subset &(\beta (m,n),\gamma ^{\ast })\subset
(\beta
(m,n),2r_{2})\text{,} \\
|D(\lambda ,m,n)| &\leq &qn^{-4}\text{~in~}(\beta (m,n),\gamma
^{\ast }).
\end{eqnarray*}%
Let us show that $\gamma ^{\ast }<2r_{2}$. If the assertion is false, then $%
r_{2},2r_{2}\in (\beta (m,n),\gamma \ast ]$ which yields
\begin{equation*}
-qn^{-4}\leq -D(r_{2},m,n)\leq -D(2r_{2},m,n)\leq qn^{-4}.
\end{equation*}%
Thus we get
\begin{equation*}
r_{2}\min\limits_{[r_{2},2r_{2}]}\{-\partial _{\lambda }D(\lambda
,m,n)\}\leq D(r_{2},m,n)-D(2r_{2},m,n)\leq 2qn^{-4}\leq
2qN_{3}^{-4}.
\end{equation*}%
From this and \eqref{48} we obtain the inequality
\begin{equation*}
\omega _{1}\frac{r_{2}}{2}\leq 2qN_{3}^{-4},
\end{equation*}%
which contradicts \eqref{32}, and the assertion follows. In
particular, we have $-D(\gamma ^{\ast },m,n)=qn^{-4}$. Since the
equations $-D(\lambda
_{\pm },m,n)=\pm qn^{-4}$ are equivalent to \eqref{42}, in the cases %
\eqref{51}, \eqref{52} we have in the cases \eqref{51}, \eqref{52}
that
\begin{equation*}
\beta (m,n)=\lambda ^{-}(m,n),\quad \gamma (m,n)=\lambda
^{+}(m,n),\quad J(m,n)=I(m,n),
\end{equation*}%
in the case \eqref{54} we have
\begin{equation*}
\beta (m,n)=\lambda ^{-}(m,n),\quad \gamma ^{\ast }=\lambda
^{+}(m,n),\quad J(m,n)\subset I(m,n).
\end{equation*}%
Hence the set $[0,r]\setminus \mathcal{H}(N,q,r)$ is covered by
the intervals $I(m,n)$ with the ends satisfying \eqref{42},
\eqref{50}. It follows from \eqref{32} and (\ref{HYPSA}) that for
$n>N_{3}(q)$, $\lambda <r_{2}$,
\begin{equation}
\frac{3\omega _{1}}{4}\lambda \leq \lambda \Bigl(\omega _{1}-\frac{\varphi }{%
n^{2}}\Bigr)-\lambda ^{3/2}\Omega (\lambda )\leq \frac{5}{4}\omega
_{1}\lambda .  \label{43}
\end{equation}%
In particular, the left side of equation \eqref{42} is positive,
which along
with \eqref{44} yields positivity of $d(m,n)$ in the right side of equations %
\eqref{42}. Combining \eqref{42}, \eqref{43}, \eqref{44} gives
\eqref{41}. Next equation \eqref{42} implies
\begin{equation*}
D(\lambda ^{-}(m,n),m,n)-D(\lambda ^{+}(m,n),m,n)=2qn^{-4}.
\end{equation*}%
From this and \eqref{48} we conclude that
\begin{gather*}
\omega _{1}(\lambda ^{+}(m,n)-\lambda ^{-}(m,n))\leq 4qn^{-4}, \\
\omega _{1}(\lambda ^{+}(m,n)-\lambda ^{-}(m,n))\geq
\frac{4}{3}qn^{-4},
\end{gather*}%
which yields
\begin{equation*}
\frac{4q}{3\omega _{1}}\frac{1}{n^{4}}\leq \lambda
^{+}(m,n)-\lambda _{-}(m,n)\leq \frac{4q}{\omega
_{1}}\frac{1}{n^{4}}
\end{equation*}%
and the assertion follows.

Let us turn to the proof of (\textit{ii}). By Lemma \ref{LO}(i),
$I(m,n)\cap \lbrack 0,r]\neq \emptyset $ if and only if $\lambda
^{-}(m,n)\leq r$. In particular, $0<\lambda
^{-}(m_{n}(r),n),\lambda ^{-}(M_{n}(r),n)\leq r$. On the other
hand, $\lambda ^{-}(m,n)$ is a solution to the equation
\begin{equation}
-D(\lambda ^{-}(m,n),m,n)=-qn^{-4}.  \label{481}
\end{equation}%
By \eqref{48} we have
\begin{equation*}
-\partial _{m}D(\lambda ,m,n)=\frac{1}{n^{2}},\quad \frac{\omega _{1}}{2}%
\leq -\partial _{\lambda }D(\lambda ,m,n)\leq \frac{3\omega _{1}}{2}\text{%
~~for~~}\lambda \in \lbrack 0,r].
\end{equation*}%
Hence for fixed $n\geq N_{3}$, the implicit function $\lambda
^{-}(m,n)$ satisfying the equation \eqref{481} is defined and
strongly decreases on the
interval $[m_{n}(r),M_{n}(r)]$ which yields \eqref{451}. The relation $%
I(k,n)\cap \lbrack 0,r]\neq \emptyset $ for each integer $k\in
\lbrack m_{n}(r),M_{n}(r)]$ follows from \eqref{451}. Next if
$m_{n}(r)\leq k-1\leq k\leq M_{n}(r)$, then we have
\begin{gather*}
0=-D(\lambda ^{-}(k-1,n),k-1,n)+D(\lambda ^{-}(k,n),k,n)= \\
=-D(\lambda ^{-}(k-1,n),k-1,n)+D(\lambda ^{-}(k,n),k-1,n)+ \\
-D(\lambda ^{-}(k,n),k-1,n)+D(\lambda ^{-}(k,n),k,n) \\
=-D(\lambda ^{-}(k-1,n),k-1,n)+D(\lambda ^{-}(k,n),k-1,n)-\frac{1}{n^{2}} \\
\leq \max\limits_{\lambda \in \lbrack 0,r]}\{-\partial _{\lambda
}D(\lambda
,k-1,n)\}(\lambda ^{-}(k-1,n)-\lambda ^{-}(k,n))-\frac{1}{n^{2}} \\
\leq \frac{3\omega _{1}}{2}(\lambda ^{-}(k-1,n)-\lambda ^{-}(k,n))-\frac{1}{%
n^{2}},
\end{gather*}%
and hence
\begin{equation*}
\frac{3\omega _{1}}{2}(\lambda ^{-}(k-1,n)-\lambda ^{-}(k,n))\geq \frac{1}{%
n^{2}},
\end{equation*}%
which completes the proof of (\textit{ii}).

Next note that \eqref{45} is a consequence of the inequalities
\begin{equation*}
\frac{2c_{0}}{5\omega _{1}}\frac{1}{n^{3+\alpha }}\leq \frac{2}{5\omega _{1}}%
d(m,n)\leq \lambda^{-}(m,n)\leq r,
\end{equation*}
and (\textit{iii}) follows.

It remains to prove (\textit{iv}). We have
\begin{equation*}
D_{\infty }(\lambda _{\infty }^{\pm },m,n)-D_{\infty }(\lambda
_{j}^{\pm },m,n)=D_{j}(\lambda _{j}^{\pm },m,n)-D_{\infty
}(\lambda _{j}^{\pm },m,n).
\end{equation*}%
On the other hand,
\begin{gather*}
|D_{j}(\lambda _{j}^{\pm },m,n)-D_{\infty }(\lambda _{j}^{\pm
},m,n)|\leq \frac{\lambda _{j}^{\pm }}{N^{2}}|\varphi _{\infty
}-\varphi _{j}|+(\lambda
_{j}^{\pm })^{3/2}|\Omega _{\infty }-\Omega _{j}|\leq \\
2r_{2}R2^{-j}(N^{-2}+1)\leq 4r_{2}R2^{-j},
\end{gather*}%
which along with \eqref{48} and \eqref{32} implies
\begin{equation*}
|\lambda _{\infty }^{\pm }(m,n)-\lambda _{j}^{\pm }(m,n)|\leq \frac{8R}{%
\omega _{1}}r_{2}2^{-j}\leq \frac{\omega _{1}}{16}2^{-j},
\end{equation*}%
and the lemma follows.
\end{proof}

\paragraph{ Cardinality of the set $\{I_{j}(m,n)\}$.}

Set
\begin{equation*}
\theta _{n}=\{\omega _{0}n^{2}-C\}=\text{ decimal part of }\omega
_{0}n^{2}-C,
\end{equation*}%
and introduce the sequence of numbers $\delta _{n}(r)$ defined by
\begin{equation*}
\delta _{n}(r)=1\text{~~ when~~}\theta _{n}\leq \frac{5}{2}\omega _{1}rn^{2},%
\text{~~and~~}\delta _{n}(r)=0\text{~~ otherwise}.
\end{equation*}

\begin{lemma}
\label{CARD} For each $n\geq N_3$ and $0<r<r_2$, $1\leq j\leq
\infty$,
\begin{equation}  \label{card1}
\text{card~~}\{m:I_j(m,n)\cap [0,r]\neq\emptyset\}\leq
\delta_n+crn^2.
\end{equation}
\end{lemma}

\begin{proof}
For simplicity we omit the index $j$. Denote by $\mathcal{I}$ the
totality of all intervals $I(m,n)$ given by Lemma \ref{LO} such
that $I(m,n)\cap
\lbrack 0,r_{2}]\neq \emptyset .$ By assertion (\textit{ii}) of Lemma \ref%
{LO} there is one-to-one correspondence between the intervals
$I(m,n)\in \mathcal{I}$, having nonempty intersection with
$[0,r]$, and the sequence of $\lambda^{-}(k,n)$ given by
\eqref{451}. In particular, we have
\begin{equation}
\text{card~}\{m:I(m,n)\cap \lbrack 0,r]\neq \emptyset
\}=M_{n}(r)-m_{n}(r)+1. \label{card11}
\end{equation}%
On the other hand, inequality \eqref{452} yields
\begin{equation*}
\lambda^{-}(m_{n}(r),n)-\lambda^{-}(M_{n}(r),n)\geq (M_{n}(r)-m_{n}(r))\frac{%
2}{3\omega _{1}n^{2}},
\end{equation*}%
Since $m_{n}(r),M_{n}(r)\in \lbrack 0,r]$, we conclude from this
that
\begin{equation*}
M_{n}(r)-m_{n}(r)\leq cn^{2}r.
\end{equation*}%
Combining this inequality with (\ref{card11}) we obtain
\begin{equation}
\text{card~}\{m:I(m,n)\cap \lbrack 0,r]\neq \emptyset \}\leq
1+crn^{2} \label{94}
\end{equation}%
On the other hand, for given $n>N_{3}$,
\begin{equation*}
\Big( \bigcup\limits_{m:I(m,n)\in \mathcal{I}}I(m,n)\Big) \cap
\lbrack 0,r]=\emptyset \text{~~if~~}\min\limits_{m:I(m,n)\in
\mathcal{I}}\lambda _{-}(m,n)>r.
\end{equation*}%
Since, by Lemma \ref{LO}, $0<\frac{2}{5\omega _{1}}d(m,n)\leq
\lambda _{-}(m,n)$, we conclude from this that
\begin{equation*}
\Big( \bigcup\limits_{m:I(m,n)\in \mathcal{I}}I(m,n)\Big) \cap
\lbrack 0,r]=\emptyset
\text{~~if~~}\min\limits_{m:d(m,n)>0}d(m,n)>\frac{5}{2}\omega
_{1}r
\end{equation*}%
Obviously
\begin{equation*}
\min\limits_{m:d(m,n)>0}d(m,n)=\frac{1}{n^{2}}\min\limits_{m:\omega
_{0}n^{2}-C-m>0}(\omega _{0}n^{2}-C-m)=\frac{\theta _{n}}{n^{2}},
\end{equation*}%
which yields
\begin{equation}
\text{card~}\{m:I(m,n)\cap \lbrack 0,r]\neq \emptyset \}=0\text{~~for~~}%
\frac{\theta _{n}}{n^{2}}>\frac{5}{2}\omega _{1}r  \label{95}
\end{equation}%
Combining \eqref{94} and \eqref{95} we obtain
\begin{equation*}
\sum\limits_{m:I(m,n)\cap \lbrack 0,r]\neq \emptyset }1=\text{card~}%
\{m:I(m,n)\cap \lbrack 0,r]\neq \emptyset \}\leq \delta
_{n}(r)+crn^{2},
\end{equation*}%
and the lemma follows.
\end{proof}

\paragraph{Proof of Theorem \protect\ref{thm1}.}

First note that
\begin{equation*}
\text{\textrm{~meas~}}\Big([0,r]\setminus \bigcap\limits_{j=1}^{\infty }%
\mathcal{H}_{j}(N,q,r)\Big)\leq \sum\limits_{{\tiny
\begin{array}{c}
(m,n):n\geq N_{3} \\
I_{j}(m,n)\cap \lbrack 0,r]\neq \emptyset%
\end{array}%
}}\text{\textrm{~meas~}}\Big(\bigcup\limits_{j=1}^{\infty }I_{j}(m,n)\Big)%
:=\wp (r).
\end{equation*}%
The following lemma gives an estimate for $\wp $ in terms of the sequence $%
\delta _{n}(r)$.

\begin{lemma}
\label{X} For each $N\geq N_3$, $0<r<r_2$, and $\sigma\in
(0,1/(3+\alpha))$,
\begin{equation}  \label{92}
\wp(r)\leq \sum\limits_{n\geq c_3r^{-1/(3+\alpha)}}c\frac{\delta_n(r)}{n^4}%
(\ln n +1)+c r^{1+\sigma},
\end{equation}
where $c$ depends on $\omega_1$, $\alpha$, $\sigma$ and $R$ only, $c_3=\left(%
\frac{2c_0}{5\omega_1}\right)^{1/(3+\alpha)}$, $c_0$ is the
constant from Condition $M$.
\end{lemma}

\begin{proof}
Introduce the sequence
\begin{equation*}
\varsigma (n)=\frac{4}{\ln 2}\ln n-\frac{1}{\ln 2}\ln
\big[\frac{64}{3\omega _{1}^{2}}q\big],\quad n\geq 1.
\end{equation*}%
It is easy to see that $\omega _{1}2^{-j-4}<4q/(3\omega _{1}n^{4})$ for all $%
j>\varsigma (n)$. It follows from this and \eqref{LINF} that for all such $j$%
, the intervals $I_{j}(m,n)$ have nonempty intersections with
$I_{\infty }(m,n)$ and $\cup _{j>\varsigma (n)}I_{j}(m,n)\subset
\tilde{I}_{\infty }(m,n)$, which yields
\begin{equation*}
\text{\textrm{~meas~}}\Big(\bigcup\limits_{j>\varsigma (n)}I_{j}(m,n)\Big)%
\leq \text{\textrm{~meas~}}\tilde{I}_{\infty }(m,n)\leq
\frac{20q}{3\omega _{1}}\frac{1}{n^{4}},
\end{equation*}%
where $\widetilde{I}_{\infty }(m,n)$ is the $4q/(3\omega
_{1}n^{4})$- neighborhood of $I_{\infty }(m,n)$. Thus we get
\begin{gather*}
\wp (r)\leq \sum\limits_{{\tiny
\begin{array}{c}
n\geq N_{3} \\
I_{j}(m,n)\cap \lbrack 0,r]\neq \emptyset%
\end{array}%
}}\sum\limits_{j=1}^{\varsigma (n)}\sum\limits_{(m:I_{j}(m,n)\cap
\lbrack
0,r]\neq \emptyset }\text{\textrm{~meas~}}I_{j}(m,n)+ \\
\sum\limits_{{\tiny
\begin{array}{c}
n\geq N_{3} \\
I_{j}(m,n)\cap \lbrack 0,r]\neq \emptyset%
\end{array}%
}}\sum\limits_{(m:I_{\infty }(m,n)\cap \lbrack 0,r]\neq \emptyset }\text{%
\textrm{~meas~}}\Big(\bigcup\limits_{j>\varsigma
(n)}I_{j}(m,n)\Big),
\end{gather*}%
which leads to
\begin{gather*}
\wp (r)\leq \sum\limits_{{\tiny
\begin{array}{c}
n\geq N_{3} \\
I_{j}(m,n)\cap \lbrack 0,r]\neq \emptyset%
\end{array}%
}}\frac{c}{n^{4}}\sum\limits_{j=1}^{\varsigma
(n)}\sum\limits_{(m:I_{j}(m,n)\cap \lbrack 0,r]\neq \emptyset }1+ \\
\sum\limits_{n\geq N_{3}}\frac{c}{n^{4}}\sum\limits_{(m:I_{\infty
}(m,n)\cap \lbrack 0,r]\neq \emptyset }1,
\end{gather*}%
Since, by Lemma \ref{LO}, the intersection $I_{j}(m,n)\cap \lbrack
0,r]$ is empty for $n<c_{3}r^{-1/(3+\alpha )}$, we have
\begin{gather*}
\wp (r)\leq \sum\limits_{n\geq c_{3}r^{-1/(3+\alpha )}}\frac{c}{n^{4}}%
\sum\limits_{j=1}^{\varsigma (n)}\text{card~}\{m:I_{j}(m,n)\cap
\lbrack
0,r]\neq \emptyset \}+ \\
\sum\limits_{n\geq c_{3}r^{-1/(3+\alpha )}}\frac{c}{n^{4}}\text{card~}%
\{m:I_{\infty }\cap \lbrack 0,r]\neq \emptyset \},
\end{gather*}%
which along with \eqref{card1} yields
\begin{gather*}
\wp (r)\leq \sum\limits_{n\geq c_{3}r^{-1/(3+\alpha )}}\frac{c}{n^{4}}%
\varsigma (n)(\delta _{n}(r)+crn^{2})+\sum\limits_{n\geq
c_{3}r^{-1/(3+\alpha )}}\frac{c}{n^{4}}(\delta _{n}+crn^{2}), \\
\leq c\sum\limits_{n\geq c_{3}r^{-1/(3+\alpha )}}\frac{\delta _{n}(r)}{n^{4}}%
(1+\ln n)+cr\sum\limits_{n\geq c_{3}r^{-1/(3+\alpha )}}\frac{1+\ln
n}{n^{2}}
\\
\leq \sum\limits_{n\geq c_{3}r^{-1/(3+\alpha )}}c\frac{\delta _{n}(r)}{n^{4}}%
(1+\ln n)+cr^{1+\sigma },
\end{gather*}%
which completes the proof.
\end{proof}

Introduce the sequence
\begin{equation*}
\nu _{n}(r)=1\text{~~when~~}\theta _{n}\leq \frac{5}{2}\omega _{1}r^{\frac{%
1-\alpha }{3-\alpha }},\text{~~and~~}\nu
_{n}(r)=0\text{~~otherwise}.
\end{equation*}

\begin{lemma}
\label{LT} Let $0<r\leq r_{2}$and $0<\beta <\alpha /(3-\alpha )$.
\newline If $c_{3}r^{\frac{-1}{3+\alpha }}\geq
r^{\frac{-1}{3-\alpha }}$, then
\begin{equation}
\wp (r)\leq cr^{1+\beta }.  \label{97}
\end{equation}%
If $c_{3}r^{\frac{-1}{3+\alpha }}<r^{\frac{-1}{3-\alpha }}$, then
\begin{equation}
\wp (r)\leq c\sum\limits_{c_{3}r^{-1/(3+\alpha )}\leq n\leq
r^{-1/(3-\alpha )}}\frac{\nu _{n}(r)}{n^{4}}(1+\ln n)+cr^{1+\beta
},  \label{98}
\end{equation}%
where $c$ does not depend on $r$.
\end{lemma}

\begin{proof}
Since $\beta <1/(3+\alpha )$, we have from inequality \eqref{92}
\begin{equation}
\wp (r)\leq c\Pi (r)+\sum\limits_{n\geq r^{-1/(3-\alpha
)}}c\frac{\delta _{n}(r)}{n^{4}}(1+\ln n)+cr^{1+\beta },
\label{99}
\end{equation}%
where
\begin{equation*}
\Pi (r)=\sum\limits_{c_{3}r^{-1/(3+\alpha )}\leq n\leq r^{-1/(3-\alpha )}}%
\frac{\delta _{n}(r)}{n^{4}}(1+\ln n)\text{~~for~~}c_{3}r^{\frac{-1}{%
3+\alpha }}<r^{\frac{-1}{3-\alpha }},
\end{equation*}%
and $\Pi (r)=0$ otherwise. It is easy to see that
\begin{equation}
\sum\limits_{n\geq r^{-1/(3-\alpha )}}\frac{\delta
_{n}(r)}{n^{4}}(1+\ln n)\leq \sum\limits_{n\geq r^{-1/(3-\alpha
)}}\frac{1}{n^{4}}(1+\ln n)\leq cr^{1+\beta }  \label{992}
\end{equation}%
If $n\leq r^{-1/(3-\alpha )}$ then $\theta _{n}>\frac{5}{2}\omega
_{1}r^{(1-\alpha )/(3-\alpha )}$ yields $\theta
_{n}>\frac{5}{2}\omega _{1}n^{2}r$. In other words, equality $\nu
_{n}(r)=0$ yields $\delta _{n}(r)=0$. Hence
\begin{equation}
\Pi (r)\leq \sum\limits_{c_{3}r^{-1/(3+\alpha )}\leq n\leq
r^{-1/(3-\alpha
)}}\frac{\nu _{n}(r)}{n^{4}}(1+\ln n)\text{~~for~~}c_{3}r^{\frac{-1}{%
3+\alpha }}<r^{\frac{-1}{3-\alpha }},  \label{991}
\end{equation}%
and $\Pi (r)=0$ otherwise. Substituting \eqref{992} and \eqref{991} into %
\eqref{99} and noting that we obtain needed inequalities \eqref{97}, %
\eqref{98}.
\end{proof}

Now set
\begin{equation*}
r_{3}=\min \{((\frac{5}{2}\omega
_{1})^{-1}c_{3}^{-\frac{1}{78}})^{1/\iota
},c_{3}^{3+\alpha }\left( \frac{1}{4}\right) ^{78(3+\alpha )}\},\quad \iota =%
\frac{1-\alpha }{3-\alpha }-\frac{1}{78(3+\alpha )}.
\end{equation*}%
\textbf{Remark} This complicated formulae simply express the fact
that the number $\rho =\left( c_{3}^{-1}r^{1/(3+\alpha )}\right)
^{1/78}$ satisfies the inequalities
\begin{equation*}
\frac{5}{2}\omega _{1}r^{(1-\alpha )/(3-\alpha )}\leq \rho \leq 1/4\text{%
~~for~~}r\leq r_{3}.
\end{equation*}

\begin{lemma}
\label{LY} Let , $0<r\leq r_{3}$, $c_{3}r^{-1/(3+\alpha
)}<r^{-1/(3-\alpha )} $, and $0<\gamma <(1/78-\alpha )/(3+\alpha
).$ Then
\begin{equation}
\sum\limits_{c_{3}r^{-1/(3+\alpha )}\leq n\leq r^{-1/(3-\alpha
)}}\frac{\nu _{n}(r)}{n^{4}}(1+\ln n)\leq cr^{1+\gamma }.
\label{993}
\end{equation}
\end{lemma}

\begin{proof}
Let $p$ and $q$ be the minimal and maximal integers from the interval $%
(c_{3}r^{-1/(3+\alpha )},r^{-1/(3-\alpha )}).$ Introduce the
average
\begin{equation*}
S_{n}(r)=\frac{1}{n}\sum\limits_{k=1}^{n}\nu _{k}(r),
\end{equation*}%
Noting that $\nu _{n}(r)=nS_{n}-(n-1)S_{n-1}$, we obtain
\begin{gather*}
\sum\limits_{c_{3}r^{-1/(3+\alpha )}\leq n\leq r^{-1/(3-\alpha
)}}\frac{\nu
_{n}(r)}{n^{4}}=\sum\limits_{p\leq n\leq q}\frac{\nu _{n}(r)}{n^{4}} \\
\leq \sum\limits_{p\leq n\leq q-1}\left( \frac{1+\ln
n}{n^{4}}-\frac{1+\ln (n+1)}{(n+1)^{4}}\right)
nS_{n}+\frac{1}{q^{3}}S_{q}.
\end{gather*}%
Since $q\geq c_{3}r^{-1/(3+\alpha )}$ and
\begin{equation*}
\left( \frac{1+\ln n}{n^{4}}-\frac{1+\ln (n+1)}{(n+1)^{4}}\right)
n\leq \frac{c}{n^{4}}(1+\ln n),
\end{equation*}%
we conclude from this that
\begin{gather}
\sum\limits_{c_{3}r^{-1/(3+\alpha )}\leq n\leq r^{-1/(3-\alpha
)}}\frac{\nu
_{n}(r)}{n^{4}}(1+\ln n)  \notag \\
\leq c\big\{\sum\limits_{c_{3}r^{-1/(3+\alpha )}\leq n\leq
r^{-1/(3-\alpha )}}\frac{1}{n^{4}}(1+\ln n)+r^{3/(3+\alpha
)}\bigr\}\sup\limits_{n\geq
c_{3}r^{-1/(3+\alpha )}}S_{n}\leq  \notag \\
cr^{(3-\vartheta )/(3+\alpha )}\sup\limits_{n\geq
c_{3}r^{-1/(3+\alpha )}}S_{n},  \label{9941}
\end{gather}%
where $\vartheta $ is an arbitrary positive number. Now set
\begin{equation*}
\rho =\left( c_{3}^{-1}r^{1/(3+\alpha )}\right) ^{1/78}.
\end{equation*}%
It follows from the choice of $r_{3}$ and the inequality $r<r_{3}$
that
\begin{equation}
\frac{5}{2}\omega _{1}r^{(1-\alpha )/(3-\alpha )}\leq \rho
<1/4,\quad n\geq c_{3}r^{-1/(3+\alpha )}\Rightarrow n\geq \rho
^{-78}.  \label{995}
\end{equation}%
It follows from this that
\begin{equation*}
S_{n}\equiv \frac{1}{n}\sum\limits_{{\tiny
\begin{array}{c}
1\leq k\leq n \\
\theta _{n}\leq \frac{5}{2}\omega _{1}r^{(1-\alpha )/(3-\alpha )}%
\end{array}%
}}1\leq \frac{1}{n}\sum\limits_{{\tiny
\begin{array}{c}
1\leq k\leq n, \\
\theta _{n}\in \lbrack 0,\rho ]%
\end{array}%
}}1.
\end{equation*}%
Applying Proposition \ref{distribution} we obtain that for any
$n\geq c_{3}r^{-1/(3+\alpha )}=\rho ^{-78}$
\begin{equation*}
S_{n}\leq c\rho =cr^{1/78(3+\alpha )},
\end{equation*}%
Combining this inequality with (\ref{9941}) we finally obtain
\begin{equation*}
\sum\limits_{c_{3}r^{-1/(3+\alpha )}\leq n\leq r^{-1/(3-\alpha
)}}\frac{\nu _{n}(r)}{n^{4}}(1+\ln n)\leq cr^{1/78(3+\alpha
)}r^{(3-\vartheta )/(3+\alpha )}=cr^{1+\gamma },
\end{equation*}%
and the lemma follows.
\end{proof}

\noindent Finally, combining inequalities \eqref{97}, \eqref{98} and %
\eqref{993} we conclude that for $0<r<r_{3}$,
\begin{equation*}
\wp (r)\leq cr(r^{\beta }+r^{\gamma })\leq cr^{1+\varpi },
\end{equation*}%
which completes the proof of Theorem \ref{thm1}.

\section{Descent method-Inversion of the linearized
operator}\label{resolvent}

In this section we give a general method of reduction, the descent
method, allowing to transform the original linear operator of
order 2 into the sum of a main operator with constant coefficients
and a smoothing perturbation operator.

Let us consider the basic operator equation
\begin{equation}
\mathfrak{L}u+\mathfrak{AD}_{1}u+\mathfrak{B}u+\mathfrak{L}_{-1}u=f,
\label{121}
\end{equation}%
with zero-order pseudodifferential operators $\mathfrak{A}$,
$\mathfrak{B}$ and an integro-differential operator
$\mathfrak{L}_{-1}$ of order $-1$.
Assume that they satisfy the following conditions.\\[0.3ex]
\textbf{Symmetry condition:}

(\textit{i}) For all $Y\in \mathbb{T}^{2}$ and $|\xi |\leq 1$,
\begin{equation}
A(Y,-\xi )=\overline{A(Y,\xi )},  \label{104}
\end{equation}%
which means in particular that $\mathfrak{A}u$ is real for
real-valued
functions $u$. 

(\textit{ii}) Equation \eqref{121} is invariant with respect to
the symmetry $Y\rightarrow Y^{\ast }=(-y_{1},y_{2})$. This is
equivalent to the equivariant property
\begin{equation}
\begin{split}
\mathfrak{AD}_{1}u(Y^{\ast })& =\mathfrak{AD}_{1}u^{\ast
}(Y),\quad
\mathfrak{B}u(Y^{\ast })=\mathfrak{B}u^{\ast }(Y), \\
\mathfrak{L}_{-1}u(Y^{\ast })& =\mathfrak{L}_{-1}u^{\ast
}(Y),\quad u^{\ast }(Y)=u(Y^{\ast }),
\end{split}
\label{irr}
\end{equation}%
which can be also written in the form
\begin{equation}
A\big(Y^{\ast },\xi ^{\ast }\big)=-A\big(Y,\xi \big),\quad
B\big(Y^{\ast },\xi ^{\ast }\big)=B\big(Y,\xi \big).
\label{irreversibility}
\end{equation}%
%
%
%
%
%
%
%
%
%
%
%
%
%
%
%
%
%
%
%
%
%
%
%
%
%
%
%
%
%
%
%
%
%
%

(\textit{iii}) For each $s\in \lbrack 1,l-3]$, $H_{o,e}^{s}(%
\mathbb{R}
^{2}/\Gamma )$ is invariant subspace of operators $\mathfrak{A}\mathfrak{D}%
_{1}$, $\mathfrak{B}$ and $\mathfrak{L}_{-1}$.

\noindent \textbf{Metric condition}: 

(\textit{iv}) There are exponents $r$, $s$, and $l$, satisfying
inequalities
\begin{equation*}
1\leq r\leq s\leq l-10,
\end{equation*}%
and $\varepsilon \in (0,1]$ so that
\begin{equation}
|\mathfrak{A}|_{4,r}+|\mathfrak{B}|_{4,r}\leq \varepsilon ,\quad |\mathfrak{A%
}|_{4,l}+|\mathfrak{B}|_{4,l}<\infty .  \label{124}
\end{equation}%
For each $s\in \lbrack r,l-10]$ there is a constant $C_{s}$ so
that
\begin{gather}
\Vert \mathfrak{L}_{-1}u\Vert _{r}\leq c\varepsilon \Vert u\Vert
_{r-1},
\label{126aa} \\
\Vert \mathfrak{L}_{-1}u\Vert _{s}\leq c\big(\varepsilon \Vert
u\Vert _{s-1}+C_{s}\Vert u\Vert _{0})\big).  \notag
\end{gather}%
\noindent \textbf{Restrictions on the spectrum and resolvent of
$\mathfrak{L}
$}: 

(\textit{v}) The selfadjoint operator $\mathfrak{L}:H_{o,e}^{s}(%
\mathbb{R}
^{2}/\Gamma )\mapsto H_{o,e}^{s}(%
\mathbb{R}
^{2}/\Gamma )$ has a simple eigenvalue $\nu _{1}\varepsilon
^{2}+O(\varepsilon ^{3})$ with corresponding eigenfunction $\psi
^{(0)}$;
the space $H_{o,e}^{s}(%
\mathbb{R}
^{2}/\Gamma )$ is the sum of orthogonal subspaces
\begin{equation*}
H_{o,e}^{s}(%
\mathbb{R}
^{2}/\Gamma )=\text{\textrm{~~span~}}\{\psi ^{(0)}\}\oplus
H_{o,e}^{s,\perp };
\end{equation*}%
$H_{o,e}^{s,\perp }$ are invariant subspaces of the operator
$\mathfrak{L}$.
Denote by $\mathcal{Q}$ the orthogonal projector of $H_{o,e}^{s}(%
\mathbb{R}
^{2}/\Gamma )$ onto $H_{o,e}^{s,\perp }.$

(\textit{vi}) For
\begin{equation}
\varkappa =\frac{1}{16\nu \pi ^{2}}\int\limits_{\mathbb{T}^{2}}\Big(%
A(Y,0,1)^{2}-4\nu B(Y,0,1)\Big)dY,  \label{descent4}
\end{equation}%
and all $s\geq 1$, the inverse $(\mathfrak{L}-\varkappa
)^{-1}:H_{o,e}^{s,\perp }\mapsto H_{o,e}^{s-1,\perp }$ is
continuous and
\begin{equation}
\Vert (\mathfrak{L}-\varkappa )^{-1}u\Vert _{{s-1}}\leq c(s)\Vert
u\Vert _{s}.  \label{123}
\end{equation}

\noindent \textbf{Nondegeneracy condition}: 

(\textit{vii})
\begin{equation}
\varkappa =O(\varepsilon ^{2}),\quad \int\limits_{\mathbb{T}^{2}}\Big((%
\mathfrak{L}+\mathfrak{H})\psi ^{(0)}-\mathfrak{H}\mathcal{Q}(\mathfrak{L}%
-\varkappa )^{-1}\mathcal{Q}\mathfrak{H}\psi ^{(0)}\Big)\psi
^{(0)}\,dY=@\varepsilon ^{2}+O(\varepsilon ^{3}),  \label{125aa}
\end{equation}%
where $@$ is a non-zero absolute constant, and the operator $\mathfrak{H}=%
\mathfrak{A}\mathfrak{D}_{1}+\mathfrak{B}+\mathfrak{L}_{-1}$.\\[1ex]

\textsl{Here and below we denote by $c$ generic constants
depending on $r$ and $l$ only, and use the standard notation
$O(\varepsilon^n)$ for quantities which absolute value does not
exceed $c\varepsilon^n$.} The following theorem is the main result
of this section.

\begin{theorem}
\label{t121} Under the above assumptions, there is a positive constant $%
\varepsilon _{0}$ depending on $l,r$ only so that for any $f\in H_{o,e}^{s}(%
\mathbb{R}
^{2}/\Gamma )$ and $\varepsilon <\varepsilon _{0}$ , equation
\eqref{121}
has a unique solution $u\in H_{o,e}^{s-1}(%
\mathbb{R}
^{2}/\Gamma )$ satisfying the inequalities
\begin{gather}
\Vert u\Vert _{{r-1}}\leq \varepsilon ^{-2}c\Vert f\Vert _{r},
\label{1250}
\\
\Vert u\Vert _{s-1}\leq \varepsilon ^{-2}c\Vert f\Vert _{r}\big(1+C_{s}+|%
\mathfrak{A}|_{4,s+10}+|\mathfrak{B}|_{4,s+10}\big)+c\Vert f\Vert
_{s}. \label{126}
\end{gather}
\end{theorem}

\begin{proof}
The proof constitutes the next two subsections.
\end{proof}

\subsection{Descent method}\label{descentmethod}

In this subsection we develop an algebraic method which allows us to reduce %
\eqref{121} to a Fredholm -type equation. The main result in this
direction is the following

\begin{theorem}
\label{descent1} Let zero-order pseudodifferential operators
$\mathfrak{A}$ and $\mathfrak{B}$ satisfy symmetry conditions
(\textit{i}), ({ii}), metric condition $(\mathit{iv})$, and
$\varkappa $ is given by \eqref{descent4}. Then there exist
integro-differential operators $\mathfrak{C}$, $\mathfrak{E}
$, $\mathfrak{F}$ satisfying symmetry condition \eqref{irr} with $\mathfrak{B%
}$ replaced by $\mathfrak{C}$, $\mathfrak{E}$, $\mathfrak{F}$ so
that
\begin{equation}
\big(\mathfrak{L}+\mathfrak{AD}_{1}+\mathfrak{B}\big)(1+\mathfrak{C})\Pi
_{1}=(1+\mathfrak{E})\Pi _{1}(\mathfrak{L}-\varkappa
)+\mathfrak{F}, \label{descent2}
\end{equation}%
and for $1\leq s\leq l-10$,
\begin{equation}
\begin{split}
\Vert \mathfrak{C}u\Vert _{s}+\Vert \mathfrak{E}u\Vert _{s}& \leq
c\varepsilon \Vert u\Vert _{s}+c\big(|\mathfrak{A}|_{4,s+9}+|\mathfrak{B}%
|_{4,s+9}\big)\Vert u\Vert _{0}, \\
\Vert \mathfrak{F}u\Vert _{s}& \leq c\varepsilon \Vert u\Vert _{s-1}+c\big(|%
\mathfrak{A}|_{4,s+10}+|\mathfrak{B}|_{4,s+10}\big)\Vert u\Vert
_{0}.
\end{split}
\label{descent3}
\end{equation}
\end{theorem}

The proof falls into three steps and is based on the following
proposition, which gives the special decomposition of zero-order
operators and plays the key role in our analysis. In order to
formulate it we introduce the important notion of an elementary
operator.

\begin{definition}
\label{elementary} Let $\tilde{W}:\mathbb{T}^{2}\mapsto
\mathbb{C}$ be a
function of class $C^{l}(%
\mathbb{R}
^{2}/\Gamma ).$ We say that $\mathfrak{W}$ is the elementary
operator associated with $\tilde{W}$, if $\mathfrak{W}$ is a
zero-order
pseudodifferential operator with the symbol $W(Y,\xi _{2})=\text{Re~}\tilde{W%
}(Y)+i\xi _{2}\text{Im~}\tilde{W}(Y)$.
\end{definition}

\begin{proposition}
\label{101l} Let a zero-order pseudodifferential operators
$\mathfrak{A}$ with symbol $A(Y,\xi )$, and an elementary operator
$\mathfrak{W,}$ satisfy conditions \eqref{103}, \eqref{104}, and
$\mathfrak{S}$ be a pseudodifferential operator with the symbol
$S(Y,\xi )=W(Y,\xi _{2})A(Y,\xi
) $. Then there exist pseudodifferential operators $\mathfrak{M}_{A}$, $%
\mathfrak{P}_{A}$, $\mathfrak{U}_{S}$ and $\mathfrak{N}_{A}$, $\mathfrak{Q}%
_{A}$, $\mathfrak{V}_{S}$ so that
\begin{gather}
\mathfrak{A}\Pi _{1}=\sum\limits_{j=0}^{1}\mathfrak{A}_{j}\mathfrak{D}%
_{1}^{-j}+\mathfrak{M}_{A}\mathfrak{L}+\mathfrak{N}_{A},  \label{10190} \\
\mathfrak{AD}_{1}=\sum\limits_{j=0}^{2}\mathfrak{A}_{j}\mathfrak{D}%
_{1}^{1-j}+\mathfrak{P}_{A}\mathfrak{L}+\mathfrak{Q}_{A},  \label{10200} \\
\mathfrak{SD}_{1}=\sum\limits_{j=0}^{2}\mathfrak{S}_{j}\mathfrak{D}%
_{1}^{1-j}+\mathfrak{U}_{S}\mathfrak{L}+\mathfrak{V}_{S}.
\label{1021}
\end{gather}%
Here $\mathfrak{A}_{j}$, $\mathfrak{S}_{j}$ are elementary
pseudodifferential operators associated with the complex-valued
functions
\begin{gather}
\tilde{A}_{0}(Y)=A(Y,0,1),\quad \tilde{A}_{1}(Y)=\frac{1}{\nu
}[\partial
_{\xi _{1}}A](Y,0,1),  \label{10120} \\
\tilde{A}_{2}(Y)=\frac{1}{2\nu ^{2}}[(\partial _{\xi
_{1}}^{2}-\partial
_{\xi _{2}})A+i\text{\textrm{Im~}}A](Y,0,1),  \label{1012a0} \\
\tilde{S}_{j}=\tilde{A}_{j}\tilde{W}\text{~~for~~}j=0,1,\quad \tilde{S_{2}}=%
\tilde{A_{2}}\tilde{W}-\Big(\frac{1}{\nu }\Big)^{2}\text{\textrm{Im~}}\tilde{%
A_{0}}\text{\textrm{Im~}}\tilde{W}.  \label{108}
\end{gather}%
For any $u\in H^{s}(%
\mathbb{R}
^{2}/\Gamma )$ with $1\leq s<l-4$,
\begin{gather}
\Vert \mathfrak{M}_{A}u\Vert _{s}+\Vert \mathfrak{N}_{A}u\Vert
_{s}+\Vert \mathfrak{P}_{A}u\Vert _{s}+\Vert
\mathfrak{Q}_{A}u\Vert _{s}\leq
\label{1020e0} \\
c(s)\Big(|\mathfrak{A}|_{3,s}\Vert u\Vert
_{0}+|\mathfrak{A}|_{3,3}\Vert
u\Vert _{s-1}\Big),  \notag \\
\Vert \mathfrak{U}_{S}u\Vert _{s-1}+\Vert \mathfrak{V}_{S}u\Vert
_{s}\leq
\label{1020ee} \\
c(s)\Big(\big(|\mathfrak{A}|_{3,3}|\mathfrak{W}|_{3,s+3}+|\mathfrak{A}%
|_{3,s+3}||\mathfrak{W}|_{3,3}\big)\Vert u\Vert _{0}+|\mathfrak{A}|_{3,3}|%
\mathfrak{W}|_{3,3}\Vert u\Vert _{s-1}\Big)  \notag
\end{gather}
\end{proposition}

\textbf{Proof:} The proof is given in Appendix
\ref{pseudodifferential}. Let us turn to the proof of Theorem
\ref{descent1}.

\textbf{First step}

We begin with the calculation of a commutator of an elementary
operator and the left side of \eqref{descent2}.

\begin{lemma}
\label{102l}Let $\mathfrak{W}$ be the elementary operator
associated with a
function $\tilde{W}\in C^{l}(%
\mathbb{R}
^{2}/\Gamma )$. Then
\begin{equation}
\Big(\mathfrak{L}+\mathfrak{AD}_{1}+\mathfrak{B}\Big)\mathfrak{W}=\mathfrak{W%
}\mathfrak{L}+\Big(\mathfrak{W}_{1}+\mathfrak{S}\Big)\mathfrak{D}_{1}+%
\mathfrak{T}+\mathfrak{R},  \label{1022}
\end{equation}%
where symbols of an elementary operator $\mathfrak{W}_{1}$ and
zero-order pseudodifferential operators $\mathfrak{S}$,
$\mathfrak{T}$ are given by
\begin{gather}
W_{1}=2\nu \partial _{y_{1}}W,\quad S=AW,  \notag \\
T=\nu \partial _{y_{1}}^{2}W-i\xi _{1}\partial _{y_{1}}W-i\tau \xi
_{2}\partial _{y_{2}}W+A\partial _{y_{1}}W+BW+  \label{1024} \\
\xi _{1}\big(\xi ^{\perp }\cdot \partial _{\xi }A\big)\big(\xi
_{2}\partial _{y_{1}}-\tau \xi _{1}\partial _{y_{2}}\big)W, \notag
\end{gather}%
where $\xi ^{\perp }=(\xi _{2},-\xi _{1})$. The remainder has the
estimate
\begin{equation}
\begin{split}
\Vert \mathfrak{R}u\Vert _{s}& \leq
c(s)\Big(|1-\mathfrak{W}|_{0,6}\Vert
u\Vert _{s-1}+|1-\mathfrak{W}|_{0,s+6}\Vert u\Vert _{0}\Big)+ \\
& c(s)\Big(|\mathfrak{W}|_{0,6}|\mathfrak{A}|_{2,3}\Vert u\Vert _{s-1}+(|%
\mathfrak{W}|_{0,s+6}|\mathfrak{A}|_{2,3}+|\mathfrak{W}|_{0,6}|\mathfrak{A}%
|_{2,s})\Vert u\Vert _{0}\Big)+ \\
& c(s)\Big(|\mathfrak{W}|_{0,4}|\mathfrak{B}|_{2,3}\Vert u\Vert _{s-1}+(|%
\mathfrak{W}|_{0,s+4}|\mathfrak{B}|_{2,3}+|\mathfrak{W}|_{0,5}|\mathfrak{B}%
|_{2,s})\Vert u\Vert _{0}\Big).
\end{split}
\label{1025}
\end{equation}
\end{lemma}

\begin{proof}
It is easy to see that
\begin{equation}
\nu \mathfrak{D}_{1}^{2}\mathfrak{W}=\nu \mathfrak{WD}_{1}^{2}+\mathfrak{W}%
_{1}\mathfrak{D}_{1}+\mathfrak{W}_{2},  \label{1026}
\end{equation}%
where $\mathfrak{W}_{2}$ is the elementary operator associated
with the functions $\nu \partial _{y^{1}}^{2}\tilde{W}$. Next
since $(-\Delta )^{-1/2} $ and $\mathfrak{W}$ are first and zero
order pseudodifferential operators, formulae for commutators
\eqref{500}, \eqref{490} from Proposition \ref{p48} imply
\begin{equation}
(-\Delta )^{1/2}\mathfrak{W}=\mathfrak{W}(-\Delta )^{1/2}+[(-\Delta )^{1/2},%
\mathfrak{W}]_{1}+\mathfrak{D}^{[N,W]},  \label{1027}
\end{equation}%
where the symbol of the operator $[(-\Delta
)^{1/2},\mathfrak{W}]_{1}$ is equal to
\begin{equation}
-i\big(\xi _{1}\partial _{y_{1}}+\tau \xi _{2}\partial
_{y_{2}}\big)W. \label{1028}
\end{equation}%
Moreover, since $|1-\mathfrak{W}|_{m,l}\leq \Vert 1-\tilde{W}\Vert
_{C^{l}}$ and the symbol of the operator $(-\Delta )^{1/2}$ does
not depend on $Y$, inequality \eqref{510} yields the estimate
\begin{equation}
\Vert \mathfrak{D}^{[N,W]}u\Vert _{s}\leq c(s)\big(\Vert
1-\tilde{W}\Vert _{C^{6+s}}\Vert u\Vert _{0}+\Vert
1-\tilde{W}\Vert _{C^{6}}\Vert u\Vert _{s-1}\big).  \label{1029}
\end{equation}%
Next applying formulae \eqref{50a},\eqref{49} to the composition
of the pseudodifferential operators $\mathfrak{A}\mathfrak{D}_{1}$
and $\mathfrak{W} $ we arrive at the equality
\begin{equation}
\mathfrak{AD}_{1}\mathfrak{W}=\mathfrak{SD}_{1}+(\mathfrak{AD}_{1}\mathfrak{W%
})_{1}+\mathfrak{D}^{(AD_{1}W)},  \label{1031}
\end{equation}%
where $(\mathfrak{A}\mathfrak{D}_{1}\mathfrak{W})_{1}$ is a zero
order pseudodifferential operator which symbol is given by formula
\eqref{49} with
$A$ replaced by $ik_{1}A$ and $B$ replaced by $W$. Noting that for $k\neq 0$%
,
\begin{equation*}
k_{1}\partial _{k_{1}}\xi =\xi _{1}\xi _{2}\xi ^{\perp },\quad
k_{1}\partial _{k_{2}}\xi =-\tau \xi _{1}^{2}\xi ^{\perp
}\text{~~with~~}\xi ^{\perp }=(\xi _{2},-\xi _{1})
\end{equation*}%
we can rewrite expression \eqref{49} for the symbol $(\mathfrak{A}\mathfrak{D%
}_{1}\mathfrak{W})_{1}$ in the form
\begin{equation}
A\partial _{y^{1}}W+\xi _{1}\Big[\xi ^{\perp }\cdot \partial _{\xi }A\Big]%
\Big [\xi _{2}\partial _{y_{1}}W-\tau \xi _{1}\partial
_{y_{2}}W\Big]. \label{1032}
\end{equation}%
Recall that it vanishes for $k=0$. Since
$\mathfrak{A}\mathfrak{D}_{1}$ is a
first order operator with $|\mathfrak{A}\mathfrak{D}_{1}|_{1,l}^{1}\leq c|%
\mathfrak{A}|_{1,l}$ inequality \eqref{51a} from Proposition
\ref{p48} implies the estimate
\begin{equation}
\Vert \mathfrak{D}^{AD_{1}W}u\Vert _{s}\leq c\Big(|\mathfrak{A}|_{2,s}|%
\mathfrak{W}|_{0,6}+|\mathfrak{A}|_{2,3}|\mathfrak{W}|_{0,6+s}\Big)\Vert
u\Vert _{0}+c|\mathfrak{A}|_{2,3}|W|_{0,6}\Vert u\Vert _{s-1}.
\label{1033}
\end{equation}%
Setting
\begin{equation*}
\mathfrak{T}=\mathfrak{W}_{2}+[(-\Delta )^{1/2},\mathfrak{W}]_{1}+(\mathfrak{%
AD}_{1}\mathfrak{W})_{1}+(\mathfrak{B}\mathfrak{W})_{0}
\end{equation*}%
we obtain the needed representation. with the remainder
\begin{equation*}
\mathfrak{D}^{[N,W]}+\mathfrak{D}^{(AD_{1}W)}+\mathfrak{D}^{(BW)},
\end{equation*}%
and the lemma follows.
\end{proof}

\textbf{Second step}

Now we give a formal construction of the operators $\mathfrak{C}$, $%
\mathfrak{E}$ and $\mathfrak{F}$. We take the operator
$\mathfrak{C}$ in the form
\begin{equation}
\mathfrak{C}=\sum\limits_{p=0}^{2}\mathfrak{W}^{(p)}\mathfrak{D}%
_{1}^{-p}-\Pi _{1},  \label{111}
\end{equation}%
where elementary operators $\mathfrak{W}^{(p)}$ will be specified
below.
Applying Lemma \ref{102l} and noting that $\mathfrak{L}$ commutes with $%
\mathfrak{D}_{1}^{j}$ we obtain
\begin{eqnarray*}
\big(\mathfrak{L}+\mathfrak{AD}_{1}+\mathfrak{B}\big)\mathfrak{W}^{(p)}%
\mathfrak{D}_{1}^{-p} &=&\mathfrak{W}^{(p)}\mathfrak{D}_{1}^{-p}\mathfrak{L}+%
\mathfrak{W}_{1}^{(p)}\mathfrak{D}_{1}^{1-p}+ \\
&&+\mathfrak{S}^{(p)}\mathfrak{D}_{1}^{1-p}+\mathfrak{T}^{(p)}\mathfrak{D}%
_{1}^{-p}+\mathfrak{R}^{(p)}\mathfrak{D}_{1}^{-p}.
\end{eqnarray*}%
Here $\mathfrak{S}^{(p)}$, $\mathfrak{T}^{(p)}$,
$\mathfrak{R}^{(p)}$ are
given by Lemma \ref{102l} with $\mathfrak{W}$ replaced by $\mathfrak{W}%
^{(p)} $. Combining these identities we arrive at
\begin{equation}
\begin{split}
\big(\mathfrak{L}+\mathfrak{AD}_{1}+\mathfrak{B}\big)(\Pi _{1}+\mathfrak{C}%
)& =(\Pi _{1}+\mathfrak{C})\mathfrak{L}+ \\
& +\sum\limits_{p=0}^{2}\Big(\mathfrak{W}_{1}^{(p)}\mathfrak{D}_{1}^{1-p}+%
\mathfrak{S}^{(p)}\mathfrak{D}_{1}^{1-p}+\mathfrak{T}^{(p)}\mathfrak{D}%
_{1}^{-p}+\mathfrak{R}^{(p)}\mathfrak{D}_{1}^{-p}\Big).
\end{split}
\label{112}
\end{equation}%
Recall that $\mathfrak{S}^{(p)}$ and $\mathfrak{T}^{(p)}$ are
zero-order pseudodifferential operators, hence by Proposition
\ref{101l}, they have the decomposition
\begin{gather*}
\mathfrak{S}^{(p)}\mathfrak{D}_{1}=\sum\limits_{j=0}^{2}\mathfrak{S}%
_{j}^{(p)}\mathfrak{D}_{1}^{1-j}+\mathfrak{U}_{S^{(p)}}\mathfrak{L}+%
\mathfrak{V}_{S^{(p)}}, \\
\mathfrak{T}^{(p)}\Pi _{1}=\sum\limits_{j=0}^{1}\mathfrak{T}_{j}^{(p)}%
\mathfrak{D}_{1}^{1-j}+\mathfrak{M}_{T^{(p)}}\mathfrak{L}+\mathfrak{N}%
_{T^{(p)}},
\end{gather*}%
in which the symbols of elementary operators
$\mathfrak{S}_{j}^{(p)}$ are given by the formulae \eqref{108}
with $W$ replaced by $W^{(p)}$, and the symbols of elementary
operators $\mathfrak{T}_{j}^{(p)}$ are given by formulae
\eqref{10120}, \eqref{1012a0} with $A$ replaced by $T$.
Substituting these relations into \eqref{112} we obtain the
identity
\begin{equation}
\begin{split}
\big(\mathfrak{L}+\mathfrak{AD}_{1}+\mathfrak{B}\big)(\Pi _{1}+\mathfrak{C}%
)& =(\Pi _{1}+\mathfrak{C})(\mathfrak{L}-\varkappa )+ \\
& \sum\limits_{p=0}^{3}\mathfrak{I}_{p}\mathfrak{D}_{1}^{1-p}+\mathfrak{R}%
^{\prime }\mathfrak{L}+\mathfrak{R}^{\prime \prime },
\end{split}
\label{113}
\end{equation}%
where the elementary operators $\mathfrak{I}_{j}$ and the reminders $%
\mathfrak{R}^{\prime }$, $\mathfrak{R}^{\prime \prime }$ are given
by
\begin{equation}
\begin{split}
\mathfrak{I}_{0}& =\mathfrak{W}_{1}^{(0)}+\mathfrak{S}_{0}^{(0)}, \\
\mathfrak{I}_{1}& =\mathfrak{W}_{1}^{(1)}+\mathfrak{S}_{0}^{(1)}+\mathfrak{S}%
_{1}^{(0)}+\mathfrak{T}_{0}^{(0)}+\varkappa \mathfrak{W}^{(0)}, \\
\mathfrak{I}_{2}& =\mathfrak{W}_{1}^{(2)}+\mathfrak{S}_{0}^{(2)}+\mathfrak{S}%
_{2}^{(0)}+\mathfrak{S}_{1}^{(1)}+\mathfrak{T}_{1}^{(0)}+\mathfrak{T}%
_{0}^{(1)}+\varkappa \mathfrak{W}^{(1)}, \\
\mathfrak{I}_{3}& =\mathfrak{S}_{2}^{(1)}+\mathfrak{S}_{1}^{(2)}+\mathfrak{S}%
_{2}^{(2)}\mathfrak{D}_{1}^{-1}+\mathfrak{T}_{1}^{(1)}+\mathfrak{T}%
_{0}^{(2)}+\mathfrak{T}_{1}^{(2)}\mathfrak{D}_{1}^{-1}+\varkappa \mathfrak{W}%
^{(2)},
\end{split}
\label{114}
\end{equation}%
\begin{equation}
\begin{split}
\mathfrak{R}^{\prime }& =\sum\limits_{p=0}^{2}\mathfrak{U}_{S^{(p)}}%
\mathfrak{D}_{1}^{-p}+\sum\limits_{p=0}^{1}\mathfrak{M}_{T^{(p)}}\mathfrak{D}%
_{1}^{-p}, \\
\mathfrak{R}^{\prime \prime }& =\sum\limits_{p=0}^{2}\mathfrak{V}_{S^{(p)}}%
\mathfrak{D}_{1}^{-p}+\sum\limits_{p=0}^{1}\mathfrak{N}_{T^{(p)}}\mathfrak{D}%
_{1}^{-p}+\sum\limits_{p=0}^{2}\mathfrak{R}^{(p)}\mathfrak{D}_{1}^{-p}.
\end{split}
\label{115}
\end{equation}%
Noting that
\begin{equation*}
\mathfrak{D}_{1}^{-2}=(-\Delta
)^{-1/2}\mathfrak{D}_{1}^{-2}\mathfrak{L}-\nu (-\Delta )^{-1/2}.
\end{equation*}%
we can rewrite \eqref{113} in the form
\begin{equation}
\big(\mathfrak{L}+\mathfrak{AD}_{1}+\mathfrak{B}\big)(\Pi _{1}+\mathfrak{C}%
)=(\Pi _{1}+\mathfrak{E})(\mathfrak{L}-\varkappa )+\sum\limits_{p=0}^{2}%
\mathfrak{I}_{p}\mathfrak{D}_{1}^{1-p}+\mathfrak{F}  \label{115a}
\end{equation}%
with
\begin{equation}
\begin{split}
\mathfrak{E}& =\mathfrak{C}+\mathfrak{R}^{\prime }+\mathbf{I}_{3}(-{\Delta }%
)^{-1/2}\mathfrak{D}_{1}^{-2}, \\
\mathfrak{F}& =\mathfrak{R}^{\prime \prime }-\nu \mathfrak{I}_{3}(-{\Delta }%
)^{-1/2}+\varkappa \big(\mathfrak{R}^{\prime
}+\mathfrak{I}_{3}(-\Delta )^{-1/2}\mathfrak{D}_{1}^{-2}.
\end{split}
\label{116}
\end{equation}

Now our task is to show that $\mathfrak{I}_{j}$, $j\leq 2$, vanish
for an appropriate choice of elementary operators
$\mathfrak{W}^{(p)}$. Note that the operator equations
$\mathfrak{I}_{j}=0$ are equivalent to the scalar
equations $\widetilde{I}_{j}(Y)=0$, in which the complex-valued functions $%
\widetilde{I}_{j}$ are associated with the operators
$\mathfrak{I}_{j}$. This observation along with formulae
\eqref{114} leads to the equations
\begin{gather}
\tilde{W}_{1}^{(0)}+\tilde{S}_{0}^{(0)}=0,  \notag \\
\tilde{W}_{1}^{(1)}+\tilde{S}_{0}^{(1)}+\tilde{S}_{1}^{(0)}+\tilde{T}%
_{0}^{(0)}+\kappa \tilde{W}^{(0)}=0,  \label{118} \\
\tilde{W}_{1}^{(2)}+\tilde{S}_{0}^{(2)}+\tilde{S}_{2}^{(0)}+\tilde{S}%
_{1}^{(1)}+\tilde{T}_{1}^{(0)}+\tilde{T}_{0}^{(1)}+\varkappa \tilde{W}%
^{(1)}=0  \notag
\end{gather}%
By Proposition \ref{101l}, we have $\tilde{T}_{0}^{(p)}=T^{(p)}(Y,0,1)$ and $%
\tilde{T}_{1}^{(p)}=\nu ^{-1}[\partial _{\xi
_{1}}T^{(p)}](Y,0,1)$, which along with equality \eqref{1024}
yields
\begin{gather*}
\tilde{T}_{0}^{(p)}=\nu \partial _{y_{1}}^{2}\tilde{W}^{(p)}-i\tau
\partial
_{y_{2}}\tilde{W}^{(p)}+\tilde{A}_{0}\partial _{y_{1}}\tilde{W}^{(p)}+\tilde{%
B}_{0}\tilde{W}^{(p)}, \\
\tilde{T}_{1}^{(p)}=-i\frac{1}{\nu }\partial _{y_{1}}\tilde{W}^{(p)}+2\tilde{%
A}_{1}\partial _{y_{1}}\tilde{W}^{(p)}+\tilde{B}_{1}\partial _{y_{1}}\tilde{W%
}^{(p)}.
\end{gather*}%
Substituting these identities into \eqref{118} and using $\tilde{S}%
_{j}^{(p)}=\tilde{A}_{j}\tilde{W}^{(p)}$ we obtain the recurrent
system of ordinary differential equations for functions
$\tilde{W}^{(p)}$, $p=0,1,2$,
\begin{gather}
(2\nu \mathfrak{\partial
}_{y_{1}}+\tilde{A}_{0})\tilde{W}^{(0)}=0,
\label{1190} \\
(2\nu \mathfrak{\partial }_{y_{1}}+\tilde{A}_{0})\tilde{W}%
^{(j)}+g_{j}=0,\quad j=1,2,  \label{119a}
\end{gather}%
where
\begin{gather}
g_{1}=\Big(\nu \partial
_{y_{1}}^{2}\tilde{W}^{(0)}+\tilde{A}_{0}\partial
_{y_{1}}\tilde{W}^{(0)}+\tilde{B}_{0}\tilde{W}^{(0)}-i\tau \partial _{y_{2}}%
\tilde{W}^{(0)}\Big)+\varkappa \tilde{W}^{(0)},  \label{1110} \\
g_{2}=\Big(\tilde{A}_{2}\tilde{W}^{(0)}+\tilde{A}_{1}\tilde{W}^{(1)}\Big)+%
\Big(2\tilde{A}_{1}\partial
_{y_{1}}\tilde{W}^{(0)}+\tilde{B}_{1}\partial
_{y_{1}}\tilde{W}^{(0)}-i\frac{1}{\nu }\partial
_{y_{1}}\tilde{W}^{(0)}\Big)
\notag \\
+\Big(\nu \partial _{y_{1}}^{2}\tilde{W}^{(1)}+\tilde{A}_{0}\partial _{y_{1}}%
\tilde{W}^{(1)}+\tilde{B}_{0}\tilde{W}^{(1)}-i\tau \partial _{y_{2}}\tilde{W}%
^{(1)}\Big)+\varkappa \tilde{W}^{(1)}.  \label{11110}
\end{gather}%
By the equivariant property, the function $\tilde{A}_{0}$ is odd
in $y_{1}$ and $\Pi _{1}\tilde{A}_{0}=\tilde{A}_{0}$. Therefore,
the general solutions of homogeneous equation \eqref{1190} has the
form
\begin{equation}
\tilde{W}^{(0)}(Y)=C(y_{2})a^{-}(Y),\text{~~where~~}a^{\pm
}(Y)=\exp (\pm \frac{1}{2\nu }\mathfrak{D}_{1}^{-1}\tilde{A}_{0})
\label{1112}
\end{equation}%
and $C(y_{2})$ is an arbitrary function. On the other hand,
inhomogeneous equation \eqref{119a} has a periodic solution if and
only if
\begin{equation}
\int\limits_{-\pi }^{\pi }g_{j}(y_{1},y_{2})a^{+}(y_{1},y_{2})\,dy_{1}=0%
\text{~~for all ~~}y_{2}.  \label{1113}
\end{equation}%
Now we aim to show that solvability condition \eqref{1113} is
fulfilled for an appropriate choice of $C$. Substituting
\eqref{1112} into the expression for $g_{1}$ and next to
\eqref{1113} we obtain the ordinary differential equation for $C$,
\begin{gather*}
-i2\pi C^{\prime }+2\pi C\varkappa -iC\tau \int\limits_{-\pi
}^{\pi }a^{+}\partial _{y_{2}}a^{-}\,dy_{1}+C\nu \int\limits_{-\pi
}^{\pi
}a^{+}\partial _{y_{1}}^{2}a^{-}dy_{1}+ \\
+C\int\limits_{-\pi }^{\pi }\tilde{A}_{0}a^{+}\partial
_{y_{1}}a^{-}\,dy_{1}+C\int\limits_{-\pi }^{\pi
}\tilde{B}_{0}\,dy_{1}=0.
\end{gather*}%
Noting that
\begin{gather*}
\int\limits_{-\pi }^{\pi }a^{+}\partial _{y_{2}}a^{-}dy_{1}=-\frac{1}{2\nu }%
\int\limits_{-\pi }^{\pi }\partial _{y_{2}}\mathfrak{D}_{1}^{-1}\tilde{A}%
_{0}\,dy_{1}=0, \\
\nu \int\limits_{-\pi }^{\pi }a^{+}\partial
_{y_{1}}^{2}a^{-}dy_{1}+\int\limits_{-\pi }^{\pi
}\tilde{A}_{0}a^{+}\partial
_{y_{1}}a^{-}\,dy_{1}=-\frac{1}{4\nu }\int\limits_{-\pi }^{\pi }\tilde{A}%
_{0}^{2}\,dy_{1},
\end{gather*}%
we can rewrite it in the form
\begin{equation}
C^{\prime }(y_{2})+b(y_{2})C(y_{2})=0,  \label{1115}
\end{equation}%
with the coefficient
\begin{equation}
b\left( y_{2}\right) =\frac{i}{2\pi }\int\limits_{-\pi }^{\pi }\big(\tilde{B}%
_{0}-\frac{1}{4\nu }\tilde{A}_{0}^{2}\big)\,dy_{1}+i\varkappa .
\label{1116}
\end{equation}%
It follows from formula \eqref{descent4} for $\varkappa $, that
the mean value of $b$ over a period is zero. Therefore, equation
\eqref{1115} has a periodic solution
\begin{equation}
C=\exp (-\mathfrak{D}_{2}^{-1}b).  \label{1116a}
\end{equation}%
In this case the particular periodic solution to \eqref{1113} is
\begin{equation}
\tilde{W}^{(1)}=-\frac{1}{2\nu }a^{-}\mathfrak{D}_{1}^{-1}\big(g_{1}a^{+}%
\big).  \label{1114}
\end{equation}%
Next note that, by the equivariant property, the functions $\tilde{A}_{0}$, $%
\tilde{A}_{2}$, and $\tilde{B}_{1}$ are odd and the functions $\tilde{A}_{1}$%
, $\tilde{B}_{0}$ are even in $y_{1}$. Hence $a^{\pm }$,
$\tilde{W}^{(0)}$ are even, and $\tilde{W}^{(1)}$ is odd in
$y_{1}$. In particular, $g_{2}$ is
odd in $y_{1}$ and automatically satisfies solvability condition \eqref{1113}%
. From this we conclude that the operators $\mathfrak{I}_{j}$,
$j\leq 2$ vanish when
\begin{equation}
\tilde{W}^{(0)}=\exp \Big(-\mathfrak{D}_{2}^{-1}b\Big)a^{-},\quad \tilde{W}%
^{(j)}=-\frac{1}{2\nu
}a^{-}\mathfrak{D}_{1}^{-1}\Big(g_{j}a^{+}\Big),\quad j=1,2,
\label{1117}
\end{equation}%
which along with \eqref{115a} leads to identity \eqref{descent2}.

\textbf{Third step}

It remains to show that operators introduced above are
well-defined and satisfy inequalities \eqref{descent3}. We begin
with the estimating of the functions $\tilde{W}^{(p)}$. Recall
that for all smooth functions $a,b\in
C^{s}(%
\mathbb{R}
^{2}/\Gamma )$,
\begin{equation*}
\Vert 1-\exp a\Vert _{C^{s}}\leq c(\Vert a\Vert _{C^{0}})\Vert
a\Vert _{C^{s}},\quad \Vert ab\Vert _{C^{s}}\leq c(s)\Vert a\Vert
_{C^{0}}\Vert b\Vert _{C^{s}}+\Vert a\Vert _{C^{s}}\Vert b\Vert
_{C^{0}}.
\end{equation*}%
It follows from this and \eqref{descent4}, \eqref{1112}, \eqref{116}, and %
\eqref{1116a} that
\begin{equation}
\Vert 1-\tilde{W}^{(0)}\Vert _{C^{s}}\leq c(\Vert
\tilde{A}_{0}\Vert _{C^{0}},\Vert \tilde{B}_{0}\Vert
_{C^{0}})\Big(\Vert \tilde{A}_{0}\Vert _{C^{s}}+\Vert
\tilde{B}_{0}\Vert _{C^{s}}\Big).  \label{1118}
\end{equation}%
On the other hand, \eqref{1110} implies the estimate
\begin{gather*}
\Vert g_{1}\Vert _{C^{s}}\leq c\Big(\Vert 1-\tilde{W}^{(0)}\Vert
_{C^{s+2}}+\Vert \tilde{A}_{0}\Vert _{C^{0}}\Vert
\tilde{W}^{(0)}\Vert
_{C^{s+1}}+ \\
\Vert \tilde{A}_{0}\Vert _{C^{s}}\Vert \tilde{W}^{(0)}\Vert
_{C^{1}}+\Vert \tilde{A}_{0}\Vert _{C^{s}}\Big)+\Big(\Vert
\tilde{B}_{0}\Vert _{C^{s}}\Vert
\tilde{W}^{(0)}\Vert _{C^{0}}+\Vert \tilde{B}_{0}\Vert _{C^{0}}\Vert \tilde{W%
}^{(0)}\Vert _{C^{s}}\Big).
\end{gather*}%
Combining it with \eqref{1118} and noting that
\begin{equation*}
\Vert \tilde{W}^{(0)}\Vert _{C^{1}}\leq c\Vert
\tilde{W}^{(0)}\Vert _{C^{0}}+c\Vert \tilde{W}^{(0)}\Vert _{C^{s}}
\end{equation*}%
we arrive at
\begin{equation}
\Vert \tilde{W}^{1}\Vert _{C^{s}}\leq c(\Vert \tilde{A}_{0}\Vert
_{C^{0}},\Vert \tilde{B}_{0}\Vert _{C^{0}})\Big(\Vert
\tilde{A}_{0}\Vert _{C^{s+2}}+\Vert \tilde{B}_{0}\Vert
_{C^{s+2}}\Big).  \label{1119}
\end{equation}%
Arguing as before and using \eqref{11110} we obtain
\begin{multline*}
\Vert g_{2}\Vert _{C^{s}}\leq c(\Vert \tilde{A}_{0}\Vert
_{C^{0}},\Vert
\tilde{B}_{0}\Vert _{C^{0}},\Vert \tilde{A}_{1}\Vert _{C^{0}},\Vert \tilde{A}%
_{2}\Vert _{C^{0}},\Vert \tilde{B}_{1}\Vert _{C^{0}})\Big(\Vert \tilde{A}%
_{0}\Vert _{C^{s+4}}+ \\
\Vert \tilde{B}_{0}\Vert _{C^{s+4}}+\Vert \tilde{A}_{1}\Vert
_{C^{s+3}}+\Vert \tilde{B}_{1}\Vert _{C^{s+2}}+\Vert
\tilde{A}_{2}\Vert _{C^{s}}\Big),
\end{multline*}%
which along with \eqref{1117} leads to the estimate
\begin{equation}
\Vert \tilde{W}^{(2)}\Vert _{C^{s}}\leq c(\Vert \tilde{A}_{j}\Vert
_{C^{0}},\Vert \tilde{B}_{j}\Vert
_{C^{0}})\Big(\sum\limits_{j=0}^{2}\Vert \tilde{A}_{j}\Vert
_{C^{s+4}}+\sum\limits_{j=0}^{1}\Vert \tilde{B}_{j}\Vert
_{C^{s+4}}\Big).  \label{1122}
\end{equation}%
Noting that $\Vert \tilde{A}_{j}\Vert _{C^{s}}\leq
|\mathfrak{A}|_{3,s}$, we conclude from \eqref{1118},
\eqref{1119}, \eqref{1122} that
\begin{equation*}
\Vert 1-\tilde{W}^{(0)}\Vert _{C^{s}}+\Vert \tilde{W}^{(1)}\Vert
_{C^{s}}+\Vert \tilde{W}^{(2)}\Vert _{C^{s}}\leq c(|\mathfrak{A}|_{3,0},|%
\mathfrak{B}|_{3,0})\Big(|\mathfrak{A}|_{3,s+4}+|\mathfrak{B}|_{3,s+4}\Big),
\end{equation*}%
which leads to
\begin{equation}
|1-\mathfrak{W}^{(0)}|_{m,s}+|\mathfrak{W}^{(1)}|_{m,s}+|\mathfrak{W}%
^{(2)}|_{m,s}\leq c(|\mathfrak{A}|_{3,0},|\mathfrak{B}|_{3,0})\Big(|%
\mathfrak{A}|_{3,s+4}+|\mathfrak{B}|_{3,s+4}\Big)  \label{11230}
\end{equation}%
for all $m\geq 0$. Now we can estimate the norm of the operator
$\mathfrak{C} $. Since $\mathfrak{C}\Pi _{1}=\mathfrak{C}$,
Proposition \ref{boundedness} along with \eqref{11230} implies
\begin{equation}
\begin{split}
\Vert \mathfrak{C}u\Vert _{s}& \leq \Vert
(1-\mathfrak{W}^{(0)})u\Vert _{s}+\Vert \mathfrak{W}^{(1)}u\Vert
_{s}+\Vert \mathfrak{W}^{(2)}u\Vert
_{s}\leq \\
& c(|\mathfrak{A}|_{3,0},|\mathfrak{B}|_{3,0})\Big(\big(|\mathfrak{A}%
|_{3,s+4}+|\mathfrak{B}|_{3,s+4}\big)\Vert u\Vert _{0}+\big(|\mathfrak{A}%
|_{3,4}+|\mathfrak{B}|_{3,4}\big)\Vert u\Vert _{s}\Big).
\end{split}
\label{1123a}
\end{equation}%
Let us estimate the operators $\mathfrak{E}$, $\mathfrak{F}$.
First note that inequalities \eqref{1020e0} \eqref{1020ee} from
Proposition \ref{101l} imply the estimate
\begin{gather*}
\Vert \mathfrak{U}_{S^{(p)}}u\Vert _{s-1}+\Vert
\mathfrak{V}_{S^{(p)}}u\Vert
_{s}\leq c(|\mathfrak{A}|_{3,0},|\mathfrak{B}|_{3,0})\Big(\big(|\mathfrak{A}%
|_{3,3}|\mathfrak{W}^{(p)}|_{3,s+3}+ \\
+|\mathfrak{A}|_{3,s+3}|\mathfrak{W}^{(p)}|_{3,3}\big)\Vert u\Vert _{0}+|%
\mathfrak{A}|_{3,3}|\mathfrak{W}^{(p)}|_{3,3}\Vert u\Vert
_{s-1}\Big).
\end{gather*}%
From this and \eqref{11230} we conclude that
\begin{gather}
\Vert \mathfrak{U}_{S^{(p)}}u\Vert _{s-1}+\Vert
\mathfrak{V}_{S^{(p)}}u\Vert
_{s}\leq  \label{1129} \\
c(|\mathfrak{A}|_{3,7},|\mathfrak{B}|_{3,7})\Big(\big(|\mathfrak{A}%
|_{3,s+7}+|\mathfrak{B}|_{3,s+7}\Big)\Vert u\Vert _{0}+\big(|\mathfrak{A}%
|_{3,7}+|\mathfrak{B}|_{3,7}\Big)\Vert u\Vert _{s-1}.  \notag
\end{gather}%
Next \eqref{11230} along with inequality \eqref{1025} from Lemma
\ref{102l} implies the inequality
\begin{eqnarray}
\Vert \mathfrak{R}^{(p)}u\Vert _{s} &\leq &c(|\mathfrak{A}|_{3,10},|%
\mathfrak{B}|_{3,10})\big((|\mathfrak{A}|_{3,s+10}+|\mathfrak{B}%
|_{3,s+10})\Vert u\Vert _{0}+  \notag \\
&&+(|\mathfrak{A}|_{3,10}+|\mathfrak{B}|_{3,10})\Vert u\Vert
_{s-1}\big). \label{11300}
\end{eqnarray}%
Recalling that for arbitrary zero-order operators $\mathfrak{A}$,$\mathfrak{W%
}$ and the operator $\mathfrak{S}$ with a symbol $S=AW$,
\begin{equation*}
|\mathfrak{S}|_{m,s}\leq c(s,m)(|\mathfrak{A}|_{m,0}|\mathfrak{W}|_{m,s}+|%
\mathfrak{A}|_{m,s}||\mathfrak{W}|_{m,0}),
\end{equation*}%
and using formula \eqref{1024} for symbol $\mathfrak{T}$ we obtain
\begin{equation*}
|\mathfrak{T}^{(p)}|_{3,s}\leq |1-\mathfrak{W}^{(p)}|_{3,s+2}+(|\mathfrak{A}%
|_{4,0}+|\mathfrak{B}|_{3,0})|\mathfrak{W}^{(p)}|_{0,s+1}+(|\mathfrak{A}%
|_{4,s}+|\mathfrak{B}|_{3,s})|\mathfrak{W}^{(p)}|_{0,1},
\end{equation*}%
which being combined with \eqref{11230} gives
\begin{equation}
|\mathfrak{T}^{(p)}|_{3,s}\leq c(|\mathfrak{A}|_{4,0},|\mathfrak{B}|_{3,0})%
\Big(|\mathfrak{A}|_{4,s+6}+|\mathfrak{B}|_{4,s+6}\Big).
\label{1126}
\end{equation}%
In particular, inequality \eqref{1126} along with \eqref{1020e0}
yields the estimate
\begin{gather}
\Vert \mathfrak{M}_{T^{(p)}}u\Vert _{s}+\Vert
\mathfrak{N}_{T^{(p)}}u\Vert
_{s}\leq  \label{1129a} \\
c(|\mathfrak{A}|_{4,0},|\mathfrak{B}|_{4,0})\Big(\big(|\mathfrak{A}%
|_{4,s+6}+|\mathfrak{B}|_{4,s+6}|\big)\Vert u\Vert _{0}+\big(|\mathfrak{A}%
|_{4,9}+|\mathfrak{B}|_{4,9}\big)\Vert u\Vert _{s-1}\Big).  \notag
\end{gather}%
Similar arguments applying to the formula \eqref{1024} for the
symbol of the operator $\mathfrak{S}$ give the estimate
\begin{equation}
|\mathfrak{S}^{(p)}|_{3,s}\leq c(|\mathfrak{A}|_{3,0}|\mathfrak{B}|_{3,0})%
\Big(|\mathfrak{A}|_{3,s+4}+|\mathfrak{B}|_{3,s+4}\Big).
\label{1127}
\end{equation}%
Finally estimate $\mathfrak{I}_{3}(-\Delta )^{-1/2}$. Applying \eqref{11230}%
, \eqref{1126}, \eqref{1127} to \eqref{114} we arrive at
\begin{equation*}
|\mathfrak{I}_{3}|_{3,s}\leq c(|\mathfrak{A}|_{4,0},|\mathfrak{B}|_{3,0})(|%
\mathfrak{A}|_{4,s+6}+|\mathfrak{B}|_{4,s+6}),
\end{equation*}%
which along with inequality \eqref{041} from Proposition
\ref{boundedness} implies
\begin{equation}
\Vert \mathfrak{I}_{3}(-\Delta )^{-1/2}u\Vert _{s}\leq c(|\mathfrak{A}%
|_{4,0},|\mathfrak{B}|_{3,0})\Big((|\mathfrak{A}|_{4,s+9}+|\mathfrak{B}%
|_{4,s+9})\Vert u\Vert
_{0}+(|\mathfrak{A}|_{4,9}+|\mathfrak{B}|_{4,9})\Vert u\Vert
_{s-1}\big).  \label{IDELT}
\end{equation}%
Next note that
\begin{gather*}
\Vert \mathfrak{E}u\Vert _{s}\leq \Vert \mathfrak{C}u\Vert
_{s}+\Vert \mathfrak{I}_{3}(-\Delta )^{-1/2}u\Vert _{s}+\Vert
\mathfrak{R}^{\prime
}u\Vert _{s}, \\
\Vert \mathfrak{F}u\Vert _{s}\leq \Vert \mathfrak{I}_{3}(-\Delta
)^{-1/2}u\Vert _{s}+\Vert \mathfrak{R}^{\prime \prime }u\Vert
_{s}.
\end{gather*}%
It follows from \eqref{1129}, \eqref{11300}, and \eqref{1129a}
that
\begin{gather*}
\Vert \mathfrak{R}^{\prime }u\Vert _{s}\leq
\sum\limits_{p=0}^{2}\Vert
\mathfrak{U}_{S^{(p)}}u\Vert _{s}+\sum\limits_{p=0}^{1}\Vert \mathfrak{M}%
_{T^{(p)}}u\Vert _{s}\leq \\
c(|\mathfrak{A}|_{4,9},|\mathfrak{B}|_{4,9})\Big(\big(|\mathfrak{A}%
|_{4,9+s}+|\mathfrak{B}|_{4,9+s}\big)\Vert u\Vert _{0}+\big(|\mathfrak{A}%
|_{4,9}+|\mathfrak{B}|_{4,9}\big)\Vert u\Vert _{s}\Big), \\
\Vert \mathfrak{R}^{\prime \prime }u\Vert _{s}\leq
\sum\limits_{p=0}^{2}\Vert \mathfrak{V}_{S^{(p)}}u\Vert
_{s}+\sum\limits_{p=0}^{1}\Vert \mathfrak{N}_{T^{(p)}}u\Vert
_{s}+\Vert
\mathfrak{R}^{(p)}u\Vert _{s}\leq \\
c(|\mathfrak{A}|_{4,10},|\mathfrak{B}|_{4,10})\Big(\big(|\mathfrak{A}%
|_{4,10+s}+|\mathfrak{B}|_{4,10+s}\big)\Vert u\Vert _{0}+\big(|\mathfrak{A}%
|_{4,10}+|\mathfrak{B}|_{4,10}\big)\Vert u\Vert _{s}\Big).
\end{gather*}%
Combining these results with \eqref{IDELT} and \eqref{descent3} we
finally obtain
\begin{equation*}
\begin{split}
\Vert \mathfrak{C}u\Vert _{s}+\Vert \mathfrak{E}u\Vert _{s}& \leq c(|%
\mathfrak{A}|_{4,9},|\mathfrak{B}|_{4,9})\Big(\big(|\mathfrak{A}|_{4,9}+|%
\mathfrak{B}|_{4,9}\big)\Vert u\Vert _{s}+ \\
& \big(|\mathfrak{A}|_{4,s+9}+|\mathfrak{B}|_{4,s+9}\big)\Vert u\Vert _{0}%
\Big), \\
\Vert \mathfrak{F}u\Vert _{s}& \leq c(|\mathfrak{A}|_{4,10},|\mathfrak{B}%
|_{4,10})\Big(\big(|\mathfrak{A}|_{4,10}+|\mathfrak{B}|_{4,10}\big)\Vert
u\Vert _{s-1}+ \\
& \big(|\mathfrak{A}|_{4,s+10}+|\mathfrak{B}|_{4,s+10}\big)\Vert u\Vert _{0}%
\Big),
\end{split}%
.
\end{equation*}%
which gives \eqref{descent3} and the theorem \ref{descent1}
follows.

\subsection{Proof of Theorem \protect\ref{t121}}

The proof is based on the following lemma on the invertibility of
the operator $1+\mathfrak{E}$ from Theorem \ref{descent1}.

\begin{lemma}
\label{l122} Under the assumptions of Theorem \ref{t121}, there is
a
positive $\varepsilon _{1}$ depending on $r$ and $l$ only, so that for all $%
\varepsilon \in (0,\varepsilon _{1})$, the operator
$1+\mathfrak{E}$ has the
bounded inverse $(1+\mathfrak{E})^{-1}:H_{o,e}^{s}(%
\mathbb{R}
^{2}/\Gamma )\mapsto H_{o,e}^{s}(%
\mathbb{R}
^{2}/\Gamma )$ and, for all $u\in H_{o,e}^{s}(%
\mathbb{R}
^{2}/\Gamma )$,
\begin{gather}
\Vert (1+\mathfrak{E})^{-1}u\Vert _{r}\leq c\Vert u\Vert _{r},
\label{1270}
\\
\Vert (1+\mathfrak{E})^{-1}u\Vert _{s}\leq c\Vert u\Vert _{r}\big(C_{s}+|%
\mathfrak{A}|_{4,s+10}+|\mathfrak{B}|_{4,s+10}\big)+c\Vert u\Vert
_{s}. \label{128}
\end{gather}
\end{lemma}

\begin{proof}
Formally we have
\begin{equation*}
(1+\mathfrak{E})^{-1}=\sum\limits_{n=0}^{\infty
}(-1)^{n}\mathfrak{E}^{n}.
\end{equation*}%
By Theorem \ref{descent1}, the operator $\mathfrak{E}:H_{o,e}^{r}(%
\mathbb{R}
^{2}/\Gamma )\mapsto H_{o,e}^{r}(%
\mathbb{R}
^{2}/\Gamma )$ is bounded and its norm does not exceed
$c\varepsilon $. It follows from this and inequality
\eqref{descent3} that for $c\varepsilon \leq 1$,
\begin{gather*}
\Vert \mathfrak{E}^{n}u\Vert _{s}\leq c\varepsilon \Vert \mathfrak{E}%
^{n-1}u\Vert
_{s}+c(|\mathfrak{A}|_{4,s+10}+|\mathfrak{B}|_{4,s+10})\Vert
\mathfrak{E}^{n-1}u\Vert _{r}\leq \\
c\varepsilon \Vert \mathfrak{E}^{n-1}u\Vert _{s}+(c\varepsilon )^{n-1}(|%
\mathfrak{A}|_{4,s+10}+|\mathfrak{B}|_{4,s+10})\Vert u\Vert _{r}\leq \\
(c\varepsilon )^{2}\Vert \mathfrak{E}^{n-2}u\Vert
_{s}+2(c\varepsilon
)^{n-1}(|\mathfrak{A}|_{4,s+10}+|\mathfrak{B}|_{4,s+10})\Vert
u\Vert _{r}\leq
\\
...\;(c\varepsilon )^{n-1}\Big(\Vert u\Vert _{s}+n(|\mathfrak{A}|_{4,s+10}+|%
\mathfrak{B}|_{4,s+10})\Vert u\Vert _{r}\Big),
\end{gather*}%
which leads to the estimate
\begin{equation*}
\left\Vert \sum\limits_{n=0}^{\infty
}(-1)^{n}\mathfrak{E}^{n}u\right\Vert
_{s}\leq c\Big(\Vert u\Vert _{s}+(|\mathfrak{A}|_{4,s+10}+|\mathfrak{B}%
|_{4,s+10})\Vert u\Vert _{r}\Big)\Big(1+\sum\limits_{n=1}^{\infty
}n(c\varepsilon )^{n-1}\Big).
\end{equation*}%
It remains to note that for $\varepsilon <1/c$, the series in the
right hand side converges absolutely and the lemma follows.
\end{proof}


Let us turn to the proof of Theorem \ref{t121}. Since $\mathfrak{A}$ and $%
\mathfrak{B}$ satisfy all hypotheses of Theorem \ref{descent1},
the corresponding operators $\mathfrak{C}$, $\mathfrak{E}$,
$\mathfrak{F}$ are
well defined and meet all requirements of this theorem. Moreover, condition (%
\textit{iv}) along with inequalities \eqref{descent3} yields the
estimates
\begin{equation}
\begin{split}
\Vert \mathfrak{C}u\Vert _{s}+\Vert \mathfrak{E}u\Vert _{s}& \leq
c\varepsilon \Vert u\Vert _{s}+c\big(|\mathfrak{A}|_{4,s+10}+|\mathfrak{B}%
|_{4,s+10}\big)\Vert u\Vert _{0}, \\
\Vert \mathfrak{F}u\Vert _{s}& \leq c\varepsilon \Vert u\Vert _{s-1}+c\big(|%
\mathfrak{A}|_{4,s+10}+|\mathfrak{B}|_{4,s+10}\big)\Vert u\Vert
_{0},
\end{split}
\label{descent3a}
\end{equation}%
\begin{equation}
\Vert \mathfrak{C}u\Vert _{t}+\Vert \mathfrak{E}u\Vert _{t}\leq
c\varepsilon \Vert u\Vert _{t},\quad \text{~~for ~~}t\in \lbrack
2,r].  \label{descent3b}
\end{equation}%
We look for a solution to basic equation \eqref{121} in the form
\begin{equation}
u=\lambda \psi ^{(0)}+(1+\mathfrak{C})v\text{~~with~~}v\in {H_{o.e}^{s-1,}}%
^{\perp }.  \label{129}
\end{equation}%
Substituting this representation into \eqref{121} we obtain the
equation
\begin{equation}
\lambda (\mathfrak{L}+\mathfrak{H})\psi ^{(0)}+(\mathfrak{L}+\mathfrak{H})(1+%
\mathfrak{C})v=f.  \label{129a}
\end{equation}%
Since $\Pi _{1}$ coincides with the identity mapping on $H_{o,e}^{r}(%
\mathbb{R}
^{2}/\Gamma )$, it follows from Theorem \ref{descent1} that
\begin{equation}
(\mathfrak{L}+\mathfrak{H})(1+\mathfrak{C})=(1+\mathfrak{E})(\mathfrak{L}%
-\varkappa )+\mathfrak{Z}\text{~~on~~}H_{o,e}^{r}(%
\mathbb{R}
^{2}/\Gamma ),  \label{1210}
\end{equation}%
where
\begin{equation*}
\mathfrak{Z}=\mathfrak{F}+\mathfrak{L}_{-1}(1+\mathfrak{C}).
\end{equation*}%
Hence we can rewrite \eqref{129a} in the equivalent form
\begin{equation*}
\lambda (\mathfrak{L}+\mathfrak{H})\psi ^{(0)}+(1+\mathfrak{E})(\mathfrak{L}%
-\varkappa )v+\mathfrak{Z}v=f.
\end{equation*}%
Applying to both sides the operators $\mathcal{Q}$ and
$1-\mathcal{Q}$ and noting that $\mathcal{Q}\mathfrak{L}\psi
^{(0)}=0$ we obtain the system of
operator equations for a scalar $\lambda $ and unknown function $v\in {%
H_{o,e}^{s-1,}}^{\perp }$,
\begin{gather}
\mathcal{Q}(1+\mathfrak{E})(\mathfrak{L}-\varkappa )v+\mathcal{Q}\mathfrak{Z}%
v=\mathcal{Q}\big(f-\lambda \mathfrak{H}\psi ^{(0)}\big),  \label{12110} \\
\lambda (1-\mathcal{Q})(L+\mathfrak{H})\psi ^{(0)}+(1-\mathcal{Q})(\mathfrak{%
L}+\mathfrak{H})(1+\mathfrak{C})v=(1-\mathcal{Q})f.  \label{1212}
\end{gather}%
Our aim is to resolve the first equation with respect to $v$.
Hence the task now is to prove the solvability of the equation
\begin{equation}
\mathcal{Q}(1+\mathfrak{E})(\mathfrak{L}-\varkappa )v+\mathcal{Q}\mathfrak{Z}%
v=\mathcal{Q}g.  \label{1213}
\end{equation}%
The corresponding result is given by

\begin{lemma}
\label{l123} Under the above assumptions, there is a positive
$\varepsilon
_{2}$ depending on $r,l$ only, so that for $g\in H_{o,e}^{s}(%
\mathbb{R}
^{2}/\Gamma )$ and $0<\varepsilon <\varepsilon _{2}$, equation
\eqref{1213} has the unique solution satisfying the inequalities
\begin{gather}
\Vert v\Vert _{r-1}\leq c\Vert g\Vert _{r},  \label{12140} \\
\Vert v\Vert _{s-1}\leq c\Vert g\Vert _{r}\big(C_{s}+|\mathfrak{A}%
|_{4,s+10}+|\mathfrak{B}|_{4,s+10}\big)+c\Vert g\Vert _{s}.
\label{1215}
\end{gather}
\end{lemma}

\begin{proof}
We look for a solution to equation \eqref{1213} in the form
\begin{equation}
v=\mathcal{Q}(\mathfrak{L}-\varkappa )^{-1}(1+\mathfrak{E})^{-1}\mathcal{Q}%
\varphi  \label{1217}
\end{equation}%
with the new unknown function $\varphi \in {H_{o,e}^{s}}^{\perp
}$. Since the operator $\mathfrak{L}$ commutes with the projector
$\mathcal{Q}$, we
have%
\begin{equation}
\mathcal{Q}(1+\mathfrak{E})(\mathfrak{L}-\varkappa )\mathcal{Q}(\mathfrak{L}%
-\varkappa )^{-1}(1+\mathfrak{E})^{-1}\mathcal{Q}=\mathcal{Q}\mathfrak{X}_{0}%
\mathcal{Q}+\mathcal{Q}  \label{1218}
\end{equation}%
with%
\begin{equation*}
\mathfrak{X}_{0}=-(\mathfrak{E}\mathcal{Q})^{2}+(1+\mathfrak{E})\mathcal{Q}%
\Big((1+\mathfrak{E})^{-1}-1+\mathfrak{E}\Big).
\end{equation*}%
Hence \eqref{1213} is equivalent to the equation
\begin{equation}
\varphi +\mathcal{Q}\mathfrak{X}\mathcal{Q}\varphi =\mathcal{Q}g,
\label{1219}
\end{equation}%
in which the operator $\mathfrak{X}$ is defined by
\begin{equation*}
\mathfrak{X}=\mathfrak{X}_{0}+\mathfrak{Z}\mathcal{Q}(\mathfrak{L}-\varkappa
)^{-1}(1+\mathfrak{E})^{-1}.
\end{equation*}%
We use the regularization method to prove the existence and
uniqueness of solutions to \eqref{1219}, and consider a family of
regularising equations depending on small positive parameter
$\delta $,
\begin{equation}
\varphi +\mathcal{Q}\mathfrak{X}\mathcal{Q}(-\Delta )^{-\delta }\varphi =%
\mathcal{Q}g,\quad \delta >0.  \label{1219a}
\end{equation}%
Let us estimate the norm of the operator $\mathfrak{X}$. It is
easy to see that Lemma \ref{l122} implies the estimates
\begin{gather}
\Vert \mathfrak{X}_{0}\varphi \Vert _{r}\leq c\varepsilon
^{2}\Vert \varphi
\Vert _{r},  \label{1221} \\
\Vert \mathfrak{X}_{0}\varphi \Vert _{s}\leq c\Vert \varphi \Vert _{r}\big(|%
\mathfrak{A}|_{4,s+10}+|\mathfrak{B}|_{4,s+10}\big)+c\varepsilon
^{2}\Vert \varphi \Vert _{s}.  \notag
\end{gather}%
On the other hand, inequality \eqref{123} along with Lemma
\ref{l122} leads to the estimate
\begin{gather}
\Vert \mathcal{Q}(\mathfrak{L}-\varkappa )^{-1}(1+\mathfrak{E})^{-1}\mathcal{%
Q}\varphi \Vert _{r-1}\leq c\Vert \varphi \Vert _{r},  \label{1225} \\
\Vert \mathcal{Q}(\mathfrak{L}-\varkappa )^{-1}(1+\mathfrak{E})^{-1}\mathcal{%
Q}\varphi \Vert _{s-1}\leq c\Vert \varphi \Vert _{r}\big(|\mathfrak{A}%
|_{4,s+10}+|\mathfrak{B}|_{4,s+10}\big)+c\Vert \varphi \Vert _{s}.
\notag
\end{gather}%
Next estimates \eqref{descent3a} for $\mathfrak{C}$ and
$\mathfrak{F}$ along with inequalities \eqref{126aa} for
$\mathfrak{L}_{-1}$ imply
\begin{gather}
\Vert \mathfrak{Z}\varphi \Vert _{r}\leq c\varepsilon \Vert
\varphi \Vert
_{r-1},  \label{1223} \\
\Vert \mathfrak{Z}\varphi \Vert _{s}\leq c\Vert \varphi \Vert _{r-1}\big(%
C_{s}+|\mathfrak{A}|_{4,s+10}+|\mathfrak{B}|_{4,s+10}\big)+c\varepsilon
\Vert \varphi \Vert _{s-1}.  \notag
\end{gather}%
Combining \eqref{1221}-\eqref{1223} we finally arrive at
\begin{gather*}
\Vert \mathfrak{X}\varphi \Vert _{r}\leq c\varepsilon \Vert
\varphi \Vert
_{r}, \\
\Vert \mathfrak{X}\varphi \Vert _{s}\leq c\Vert \varphi \Vert _{r}\big(%
C_{s}+|\mathfrak{A}|_{4,s+10}+|\mathfrak{B}|_{4,s+10}\big)+c\varepsilon
\Vert \varphi \Vert _{s}.
\end{gather*}%
Since the operators $(-\Delta )^{-\delta }$ are uniformly bounded in $%
H_{o,e}^{r}(%
\mathbb{R}
^{2}/\Gamma )$, it follows that each solution to equation
\eqref{1219a} satisfies the inequalities
\begin{gather*}
\Vert \varphi \Vert _{r}(1-c\varepsilon )\leq c\Vert g\Vert _{r}, \\
\Vert \varphi \Vert _{s}(1-c\varepsilon )\leq c\Vert \varphi \Vert _{r}\big(%
C_{s}+|\mathfrak{A}|_{4,s+10}+|\mathfrak{B}|_{4,s+10}\big)+ \\
c\Vert g\Vert _{s}+c\Vert g\Vert _{r}\big(|\mathfrak{A}|_{4,s+10}+|\mathfrak{%
B}|_{4,s+10}\big).
\end{gather*}%
Hence for $\varepsilon <1/2c$,
\begin{gather}
\Vert \varphi \Vert _{r}\leq c\Vert g\Vert _{r},  \label{1228} \\
\Vert \varphi \Vert _{s}\leq c\Vert g\Vert _{s}+c\Vert g\Vert _{r}\big(%
C_{s}+|\mathfrak{A}|_{4,s+10}+|\mathfrak{B}|_{4,s+10}\big), \notag
\end{gather}%
In particular, the solution is unique. Since the operator
\begin{equation*}
\mathcal{Q}\mathfrak{X}\mathcal{Q}(-\Delta )^{-\delta }:H_{o,e}^{r}(%
\mathbb{R}
^{2}/\Gamma )\mapsto H_{o,e}^{r}(%
\mathbb{R}
^{2}/\Gamma )
\end{equation*}%
is compact, uniqueness along with the Fredholm Theorem implies the
solvability of equations \eqref{1219a} for all $\delta >0$ . Hence
they have a family of solutions $\varphi _{\delta }$ satisfying
\eqref{1228}. After passing to a subsequence we can assume that
\begin{equation*}
\varphi _{\delta }\rightarrow \varphi ,\quad (-\Delta )^{-\delta
}\varphi
_{\delta }\rightarrow \varphi \text{~~weakly in~~}H^{s}(%
\mathbb{R}
^{2}/\Gamma )\text{~~as~~}\delta \searrow 0.
\end{equation*}%
Obviously $\varphi $ serves as a solution to equation
\eqref{1219}.
Recalling \eqref{1217} we obtain the solvability and uniqueness result for %
\eqref{1213}. It remains to note that estimates \eqref{12140},
\eqref{1215} follow from formula \eqref{1217} and inequalities
\eqref{1225},\eqref{1228}.
\end{proof}


We are now in a position to complete the proof of Theorem
\ref{t121}. Applying Lemma \ref{l123} to equation \eqref{12110} we
obtain the representation
\begin{equation}
v=\lambda v_{0}+v_{f},  \label{1230}
\end{equation}%
in which $v_{0}$ and $v_{f}$ are solutions to equations \eqref{1213} with $%
g=-\mathcal{Q}\mathfrak{H}\psi ^{(0)}$ and $g=\mathcal{Q}f$. It follows from %
\eqref{12140}, \eqref{1215} that
\begin{gather}
\Vert v_{f}\Vert _{r-1}\leq c\Vert f\Vert _{r},  \notag \\
\Vert v_{f}\Vert _{s-1}\leq c\Vert f\Vert _{r}\big(C_{s}+|\mathfrak{A}%
|_{4,s+10}+|\mathfrak{B}|_{4,s+10}\big)+c\Vert f\Vert _{s},  \label{1215b} \\
\Vert v_{0}\Vert _{r-1}\leq c\varepsilon ,  \notag \\
\Vert v_{0}\Vert _{s-1}\leq c\big(C_{s}+|\mathfrak{A}|_{4,s+10}+|\mathfrak{B}%
|_{4,s+10}\big).  \notag
\end{gather}%
The functions $v_{0}$ and $v_{f}$ are completely defined by the
eigenfunction $\psi _{0}$ and the right hand side $f$ of equation \eqref{121}%
, but the scalar $\lambda $ still remains unknown. In order to
find it we substitute representation \eqref{1230} into
\eqref{1212} to obtain
\begin{equation}
\lambda (1-\mathcal{Q})(\mathfrak{L}+\mathfrak{H})\Big(\psi ^{(0)}+(1+%
\mathfrak{C})v_{0}\Big)=(1-\mathcal{Q})\Big(f-\big(\mathfrak{L}+\mathfrak{H}%
)(1+\mathfrak{C})v_{f}\Big).  \label{1231}
\end{equation}%
It follows from identity \eqref{1210} that equation \eqref{1213}
for $v_{0}$ is equivalent to
\begin{equation*}
\mathcal{Q}(\mathfrak{L}-\varkappa )\mathcal{Q}v_{0}+\mathcal{Q}\big((%
\mathfrak{L}+\mathfrak{H})\mathfrak{C}+\mathfrak{H}+\kappa )\mathcal{Q}%
v_{0}=-\mathcal{Q}\mathfrak{H}\psi ^{(0)}.
\end{equation*}%
Since the operator $\mathcal{Q}(\mathfrak{L}-\varkappa )$ has the
bounded
inverse $\mathcal{Q}(\mathfrak{L}-\varkappa )^{-1}:$ $\ {H_{o,e}^{s,\bot }}%
\mapsto {H_{o,e}^{s-1,\bot }}$, we have
\begin{equation}
v_{0}+\mathcal{Q}(\mathfrak{L}-\varkappa )^{-1}\mathcal{Q}\big((\mathfrak{L}+%
\mathfrak{H})\mathfrak{C}+\mathfrak{H}+\kappa \big)\mathcal{Q}v_{0}=-%
\mathcal{Q}(\mathfrak{L}-\varkappa
)^{-1}\mathcal{Q}\mathfrak{H}\psi ^{(0)}. \label{1232}
\end{equation}%
Next, using inequalities \eqref{1215b} along with \eqref{descent3b} and %
\eqref{12140} and noting that for $2\leq t\leq r$,
\begin{equation*}
\Vert \mathfrak{H}u\Vert _{t-1}\leq c\varepsilon \Vert u\Vert _{t}
\end{equation*}%
we arrive at
\begin{equation*}
\Vert \mathcal{Q}(\mathfrak{L}-\varkappa
)^{-1}\mathcal{Q}\mathfrak{H}\psi
^{(0)}\Vert _{r-2}\leq c\varepsilon ,\quad \Vert \mathcal{Q}(\mathfrak{L}%
-\varkappa )^{-1}\big((\mathfrak{L}+\mathfrak{H})\mathfrak{C}+\mathfrak{H}%
+\varkappa )v_{0}\Vert _{r-3}\leq c\varepsilon ^{2},
\end{equation*}%
which leads to
\begin{equation*}
v_{0}=-\mathcal{Q}(\mathfrak{L}-\varkappa
)^{-1}\mathcal{Q}\mathfrak{H}\psi ^{(0)}+v^{\prime
}\text{~~with~~}\Vert v^{\prime }\Vert _{r-3}\leq c\varepsilon
^{2}.
\end{equation*}%
Substituting this expression into \eqref{1231} gives
\begin{equation}
\lambda (1-\mathcal{Q})(\mathfrak{L}+\mathfrak{H})\Big(\psi ^{(0)}-\mathcal{Q%
}(\mathfrak{L}-\varkappa )^{-1}\mathcal{Q}\mathfrak{H}\psi ^{(0)}\Big)%
+\lambda v^{\prime \prime }=(1-\mathcal{Q})\Big(f-\big(\mathfrak{L}+%
\mathfrak{H})(1+\mathfrak{C})v_{f}\Big)  \label{1231a}
\end{equation}%
where
\begin{equation*}
v^{\prime \prime
}=(1-\mathcal{Q})(\mathfrak{L}+\mathfrak{H})\Big(v^{\prime
}+\mathfrak{C}v_{0}\Big).
\end{equation*}%
Next noting that
\begin{equation*}
\Vert (1-\mathcal{Q})\mathfrak{L}u\Vert _{s}\leq c\varepsilon
\Vert u\Vert _{0}
\end{equation*}%
and using estimate \eqref{descent3b} we obtain
\begin{equation*}
\Vert v^{\prime \prime }\Vert _{r-4}\leq c\varepsilon \Vert
v^{\prime }\Vert _{r-3}+c\varepsilon ^{2}\Vert v_{0}\Vert
_{r-2}\leq c\varepsilon ^{3}.
\end{equation*}%
Multiplying both sides of \eqref{1231} by $\psi ^{(0)}$ and
integrating the result over $\mathbb{T}^{2}$ leads to equality
\begin{equation}
\lambda K=\int\limits_{\mathbb{T}^{2}}\Big(f-\big(\mathfrak{L}+\mathfrak{H}%
)(1+\mathfrak{C})v_{f}\Big)\psi ^{(0)}\,dY,  \label{1233}
\end{equation}%
in which
\begin{equation*}
K=\int\limits_{\mathbb{T}^{2}}(\mathfrak{L}+\mathfrak{H})\Big(\psi ^{(0)}-%
\mathcal{Q}(\mathfrak{L}-\varkappa )^{-1}\mathcal{Q}\mathfrak{H}\psi ^{(0)}%
\Big)\psi ^{(0)}\,dY+\int\limits_{\mathbb{T}^{2}}v^{\prime \prime
}\psi ^{(0)}\,dY=@\varepsilon ^{2}+O(\varepsilon ^{3}).
\end{equation*}%
Hence there is a positive $\varepsilon _{1}<\varepsilon _{2}$ depending on $%
r,l$ and $\psi _{0}$ only so that $2|K|>|@|\varepsilon ^{2}$ for all $%
0<\varepsilon <\varepsilon _{1}$. From this, relation \eqref{1233}
and inequalities \eqref{1215b} we conclude that for such
$\varepsilon $ the unknown $\lambda $ is well defined and
\begin{equation}
\varepsilon ^{2}|\lambda |<\frac{c}{|@|}\Vert f\Vert _{r}.
\label{1234}
\end{equation}%
Hence for $0<\varepsilon <\varepsilon _{1}$, equation \eqref{121}
has the
unique solution $u\in H_{o,e}^{r}$. It remains to note that inequalities %
\eqref{1215b} and \eqref{1234} along with the identity
\begin{equation*}
u=\lambda \psi ^{(0)}+\lambda
(1+\mathfrak{C})v_{0}+(1+\mathfrak{C})v_{f}
\end{equation*}%
imply estimates \eqref{1250}, \eqref{126} and this ends the proof of Theorem %
\ref{t121}.

\subsection{Verification of assumptions of Theorem \protect\ref{t121}}

In this subsection we check all abstract conditions required for
solving the linear equation (\ref{121}), in the case when the
operators involved in this equation are defined by Theorem
\ref{thmChangeCoord}. Note that the symmetry
conditions (\ref{irr}) and (\textit{iii}) follows from assertion (\textit{iii%
}) of this theorem. The same conclusion can be drawn for the
Metric
conditions (\ref{124}), (\ref{126aa}) which are realized, as soon as $%
||U||_{\rho }\leq \varepsilon $ and $\rho =l+3.$

Let us consider the restrictions on the spectrum and resolvent of $\mathfrak{%
L}.$ Part (v) results immediately from Theorem \ref{1411},
restricted to the space of functions odd in $y_{1}$ and even in
$y_{2},$ and once we notice that the eigenvalue corresponding to
the eigenvector $\psi ^{(0)}=\sin y_{1}\cos y_{2}$ satisfies
\begin{equation}
-\nu +\sqrt{1+\tau ^{2}}=\varepsilon ^{2}\nu _{1}(\tau
)+O(\varepsilon ^{3}). \label{eigenvalue}
\end{equation}%
This results immediately from the fact that the linear operator
$\mathcal{L}$ (see (\ref{basiclinop}) corresponds to the
differentiation of the basic system at the point
\begin{equation}
U=U_{\varepsilon }^{(N)}+\varepsilon ^{N}W,\text{ \ \ }N\geq 3
\label{W}
\end{equation}%
where $U_{\varepsilon }^{(N)}$ is the approximate solution at order $%
\varepsilon ^{N},$ and from Lemma \ref{Lemcoef-nu}. Moreover the
coefficient
$\nu _{1}$ is positive (larger than a positive number) for any value of $%
\tau .$

Now, from Lemma \ref{Lemmcoordchange} we see that the Fourier
expansion of the diffeomorphism of the torus has, at order
$\varepsilon ,$ only harmonic
1 terms in $Y$, and at order $\varepsilon ^{2}$ only harmonics 2 and 0 in $%
Y=(y_{1},y_{2}),$ all terms being invariant under the shift $\mathcal{T}_{%
\mathbf{v}_{0}}:Y\mapsto Y+(\pi ,\pi ).$ Let us now consider the
expression of the linear operator $\mathfrak{L}+\mathfrak{H}$
which is obtained after applying the above diffeomorphism on a
linear operator which may be formally expanded in powers of
$\varepsilon ,$ having coefficients with the same property as the
diffeomorphism in terms of harmonics in $Y$. The result is
that the formal expansion in powers of $\varepsilon $ of $\mathfrak{L}+%
\mathfrak{H}$ has the same property as above i.e. the order 0 is
independent of $Y,$ while orders $\varepsilon $ and $\varepsilon
^{2}$ only contain respectively harmonic 1 and harmonics 2 and 0
in $Y.$

Let us consider the formula (\ref{descent4}) giving the coefficient $%
\varkappa .$ In the integral over $\mathbb{T}^{2}$ the functions
$A$ and $B$ are $O(\varepsilon )$ and can be expanded in powers of
$\varepsilon ,$ their principal part containing only harmonic
1-terms in $Y.$ It results
immediately that%
\begin{equation}
\varkappa =O(\varepsilon ^{2}).  \label{kappa}
\end{equation}%
It is not useful in our proof to give more precision on this
coefficient, however the interested reader might check (after few
days of computations) that
\begin{equation}
\varkappa =\varepsilon ^{2}\frac{\mu _{c}(\tau )}{16}\{8\tau ^{2}-5\tau ^{2}-%
\frac{7}{4}\}+O(\varepsilon ^{3}).  \label{varkappa}
\end{equation}%
Now, from the estimate on operators $\mathfrak{A}$ and
$\mathfrak{B}$ in
theorem \ref{thmChangeCoord} and from (\ref{descent4}) we have%
\begin{equation*}
|\varkappa |\leq c||U||_{14},
\end{equation*}%
hence from (\ref{kappa}) and (\ref{W}) we deduce that
\begin{equation}
\varkappa =\varepsilon ^{2}\widetilde{\varkappa },\text{ \ \ }|\widetilde{%
\varkappa }|\leq c(1+\varepsilon ||W||_{14}).
\label{estimvarkappa}
\end{equation}%
It should be noticed that in formulae (\ref{eigenvalue}) and (\ref{varkappa}%
) the terms of order $O(\varepsilon ^{3})$ depend on the unknown $W$ of (\ref%
{W}) which is assumed to be bounded in $\mathbb{H}_{(S)}^{14}$
(because of
the smoothness required for computing operators $\mathfrak{B}$ and $%
\mathfrak{L}_{-1}$ above). Hence, if we want to apply the result of Theorem %
\ref{142t} \ for obtaining the required estimate (\ref{123}) for $(\mathfrak{%
L}-\varkappa )^{-1}$ we need to show that properties
(\ref{perturbationb}) hold for the functions of $\varepsilon $
\begin{equation*}
\nu _{j}(\varepsilon )=\nu _{0}-\varepsilon ^{2}\nu _{1}+\varepsilon ^{3}%
\widetilde{\nu }_{j}(\varepsilon ),\text{ \ }\varkappa =\varepsilon ^{2}%
\widetilde{\varkappa }_{j}(\varepsilon )
\end{equation*}%
where $\nu _{0}=(1+\tau ^{2})^{1/2},$ $\nu _{1}=\nu _{1}(\tau ),$ $%
\widetilde{\nu }_{j}(\varepsilon )$ and $\widetilde{\varkappa }%
_{j}(\varepsilon )$ are obtained with some $W_{j}$ in
$\mathbb{H}_{(S)}^{14}$ through the iteration of the Newton method
of Nash Moser theorem. Because of the smooth dependence of these
coefficients in function of $(\varepsilon ,U)$
(see in particular Theorem \ref{thmChangeCoord} for $\nu ,$ and (\ref%
{estimvarkappa}) for $\varkappa $), the properties
(\ref{perturbationb}) for
$\widetilde{\nu }_{j}(\varepsilon )$ and $\widetilde{\varkappa }%
_{j}(\varepsilon )$ are verified as soon as there exists $R>0$ such that%
\begin{eqnarray}
||W_{j}(\varepsilon ^{\prime })-W_{j}(\varepsilon ^{\prime \prime
})||_{14} &\leq &R|\varepsilon ^{\prime }-\varepsilon ^{\prime
\prime }|,
\label{propW} \\
||W_{j+1}(\varepsilon )-W_{j}(\varepsilon )||_{14} &\leq &2^{-j}R
\notag
\end{eqnarray}%
holds. This property will be checked in section
\ref{sectionnonlin}. Assuming that this is true, we have all
conditions of Theorem \ref{142t} realized. This completes the
verification of restrictions on the spectrum and resolvent of
$\mathfrak{L}.$

The nondegeneracy condition (\ref{125aa}) is fundamental for
obtaining the estimates (\ref{1250}), (\ref{126}) for the solution
of the linearized
system. So, let us now consider the coefficient $@$ defined in (\ref{125aa}%
). We notice that the eigenfunction $\psi ^{(0)}=\sin y_{1}\cos
y_{2}$ only contains harmonic 1-terms in $Y,$ and assuming that
$W$ in formula (\ref{W}) is bounded in $\mathbb{H}_{(S)}^{14},$
and denoting the coefficient of order $\varepsilon $ in the
operator $\mathfrak{H}$ by $\mathfrak{H}^{(1)}$ we have
\begin{equation*}
\mathcal{Q}\mathfrak{H}\psi ^{(0)}=\varepsilon
\mathfrak{H}^{(1)}\psi ^{(0)}+O(\varepsilon ^{2}),
\end{equation*}%
and, since $\mathfrak{L}_{0}$ is invertible on finite Fourier
series
orthogonal to $\psi ^{(0)},$ one obtains%
\begin{equation*}
\mathfrak{H}\mathcal{Q}(\mathfrak{L}-\varkappa )^{-1}\mathcal{Q}\mathfrak{H}%
\psi ^{(0)}=\varepsilon ^{2}\mathfrak{H}^{(1)}\mathfrak{L}_{0}^{-1}\mathfrak{%
H}^{(1)}\psi ^{(0)}+O(\varepsilon ^{3}).
\end{equation*}%
Now defining the coefficient of order $\varepsilon ^{2}$ in the operator $%
\mathfrak{L}+\mathfrak{H}$ by $\mathfrak{H}^{(2)}$ we have%
\begin{equation}
@=\int\limits_{\mathbb{T}^{2}}\Big(\mathfrak{H}^{(2)}\psi ^{(0)}-\mathfrak{H}%
^{(1)}\mathfrak{L}_{0}^{-1}\mathfrak{H}^{(1)}\psi ^{(0)}\Big)\psi
^{(0)}\,dY. \label{@}
\end{equation}%
Then, we prove the following

\begin{lemma}
\label{coef@} The coefficient $@$ of nondegeneracy condition
(\ref{125aa})
is given by%
\begin{equation*}
@=-2\pi ^{2}\frac{\mu _{1}}{\mu _{0}^{2}}=\frac{\pi
^{2}}{8}(\alpha _{0}+\beta _{0}),
\end{equation*}%
where $\mu _{1}(\tau )$ is given by (\ref{mu_1}) and $(\alpha
_{0}+\beta _{0})(\tau )$ is given at Theorem \ref{Lembifurc}. This
coefficient is non zero for $\tau \neq \tau _{c}.$
\end{lemma}

\textbf{Proof:} the proof is made at Appendix \ \ref{app@}.

This completes the verification of Nondegeneracy condition.

\subsection{Inversion of $\mathcal{L}$}

Let us now consider the inversion of the operator $\mathcal{L}$ defined in (%
\ref{basiclinop})%
\begin{equation}
\mathcal{L}(U,\mu )V=F,  \label{linequV}
\end{equation}%
where%
\begin{equation*}
\mathcal{L}(U,\mu )=\left(
\begin{array}{cc}
\mathcal{G}_{\eta } & \mathcal{J}^{\ast } \\
\mathcal{J} & \mathfrak{a}%
\end{array}%
\right) ,
\end{equation*}%
\begin{equation*}
F=(f,g)\in \mathbb{H}_{(S)}^{s},\text{ \ \ }V=(\delta \phi ,\delta
\eta )
\end{equation*}%
and where we look for
\begin{equation*}
\delta U=(\delta \psi ,\delta \eta ),\text{ \ \ }\delta \psi =\delta \phi +%
\mathfrak{b}\delta \eta ,
\end{equation*}%
where $\mathfrak{b}$ is defined in (\ref{Z}). In this subsection
we prove the following Theorem which collects the results of
previous sections on the operator $\mathcal{L}$

\begin{theorem}
\label{inv L}Consider $M>0$, $\imath \geq 14,$ and $\pi /4>\delta
>0,$ and
set%
\begin{eqnarray}
U &=&U_{\varepsilon }^{(N)}+\varepsilon ^{N}W,\text{ \ }N\geq 3,
\label{decompU} \\
\mu &=&\mu _{\varepsilon }^{(N)}=\mu _{c}+\varepsilon ^{2}\mu
_{1}(\tau )+O(\varepsilon ^{3}),  \notag
\end{eqnarray}%
where $||W||_{\imath }\leq M,$ and $(U_{\varepsilon }^{(N)},\mu
_{\varepsilon }^{(N)})$ is the approximate solution at order
$\varepsilon
^{N}$ obtained at Theorem \ref{Lembifurc} in the case of diamond waves $%
(\varepsilon _{1}=\varepsilon _{2}=\varepsilon /2),$ and $\tau
=\tan \theta _{0},$ $\mu _{c}=\cos \theta _{0},$ $\delta <\theta
_{0}<\pi /2-\delta .$ Assume moreover in (\ref{decompU}) that
$W=W_{j}(\varepsilon ),$ $j\in
\mathbb{N}
$ satisfies the property (\ref{propW}), then for $\mu _{c}^{-1}\in \mathfrak{%
N}_{\alpha },$ $\alpha \in (0,1/78),$ there exists $\varepsilon
_{0}>0$ and a subset $\mathcal{E}$ of $[0,\varepsilon _{0}]$ such
that for any $F\in \mathbb{H}_{(S)}^{s},$ $s\geq 5$, and for
$\varepsilon \in \mathcal{E},$ the
linear equation (\ref{linequV}) has a unique solution $V$ corresponding to $%
\delta U\in \mathbb{H}_{(S)}^{s-3},$ such that the following estimates hold%
\begin{eqnarray}
||\delta U||_{r-3} &\leq &\frac{c}{\varepsilon ^{2}}||F||_{r},\text{ \ \ }%
5\leq r\leq \imath -8,  \label{estim inv L} \\
||\delta U||_{s-3} &\leq &\frac{c(s)}{\varepsilon ^{2}}||F||_{r}(1+%
\varepsilon ^{N}||W_{j}||_{s+18})+c(s)||F||_{s}.  \notag
\end{eqnarray}%
Moreover the following property holds for the "good" set $\mathcal{E}:$%
\begin{equation}
\frac{1}{|\mu_{r}^{(N)} -\mu _{c}|}meas\{\mu =\mu _{\varepsilon
}^{(N)};\varepsilon \in \mathcal{E\cap }(0,r)\}\rightarrow
1,\text{ \ as }r\rightarrow 0. \label{measure}
\end{equation}
\end{theorem}

\textbf{Proof.} By construction, we have
\begin{equation*}
\mathcal{G}_{\eta }(\delta \phi )-\mathcal{J}^{\ast }(\frac{1}{\mathfrak{a}}%
\mathcal{J}(\delta \phi ))=h,
\end{equation*}%
where%
\begin{equation*}
h=f-\mathcal{J}^{\ast }(g/\mathfrak{a}),
\end{equation*}%
and which, after the change of coordinates of Theorem
\ref{thmChangeCoord},
becomes%
\begin{equation*}
(\mathfrak{L}+\mathfrak{AD}_{1}+\mathfrak{B}+\mathfrak{L}_{-1})\widetilde{%
\delta \phi }=\frac{\widetilde{h}}{\kappa },
\end{equation*}%
which we learnt to invert at Theorem \ref{t121}.

First, for $U\in \mathbb{H}_{(S)}^{m},$ and $||U||_{r}\leq M_{r},$
we have
the following estimates, due to Lemma \ref{tame coef}%
\begin{equation*}
||h||_{s}\leq c(M_{4})\{||F||_{s+1}+||F||_{2}||U||_{s+2}),
\end{equation*}%
and with Theorem \ref{thmChangeCoord} we obtain
\begin{equation*}
||\widetilde{h}||_{s}\leq
c(M_{4})(||h||_{s}+||U||_{s+4}||h||_{0}).
\end{equation*}%
Taking into account of the estimate on $\kappa (Y)$ in Theorem \ref%
{thmChangeCoord}, we then arrive to%
\begin{equation}
||\frac{\widetilde{h}}{\kappa }||_{s}\leq
c(M_{5})\{||F||_{s+1}+||F||_{2}||U||_{s+5}\}.  \label{estimh/k}
\end{equation}%
In the same way, thanks to the Theorem \ref{thmChangeCoord}, we also have%
\begin{equation*}
||\delta \phi ||_{s}\leq c(M_{4})\{||\widetilde{\delta \phi }%
||_{s}+||U||_{s+4}||\widetilde{\delta \phi }||_{0}\},
\end{equation*}%
and since%
\begin{equation*}
\delta \eta =-\frac{1}{\mathfrak{a}}\mathcal{J}(\delta \phi )+\frac{1}{%
\mathfrak{a}}g,
\end{equation*}%
we obtain%
\begin{eqnarray*}
||\delta \eta ||_{s} &\leq &c(M_{4})\{||\widetilde{\delta \phi }%
||_{s+1}+||U||_{s+5}||\widetilde{\delta \phi }||_{0}\}+ \\
&&+c(M_{4})\{||F||_{s}+||U||_{s+2}||F||_{2}\}.
\end{eqnarray*}%
Hence%
\begin{eqnarray*}
||V||_{s} &\leq &c(M_{4})\{||\widetilde{\delta \phi }||_{s+1}+||U||_{s+5}||%
\widetilde{\delta \phi }||_{0}\}+ \\
&&+c(M_{4})\{||F||_{s}+||U||_{s+2}||F||_{2}\},
\end{eqnarray*}%
and thanks to Lemma \ref{tame coef}%
\begin{eqnarray}
||\delta U||_{s} &\leq
&c(M_{3})\{||V||_{s}+||U||_{s+1}||V||_{2}\},  \notag
\\
&\leq &c(M_{7})\{||\widetilde{\delta \phi }||_{s+1}+||U||_{s+5}||\widetilde{%
\delta \phi }||_{0}+||U||_{s+1}||\widetilde{\delta \phi }||_{3}\}+  \notag \\
&&+c(M_{4})\{||F||_{s}+||U||_{s+2}||F||_{2}\}.
\label{estim1deltaU}
\end{eqnarray}%
It remains to use Theorem \ref{t121} \ which connects
$\widetilde{\delta \phi }$ and $\frac{\widetilde{h}}{\kappa },$
taking into account of
estimates ii) of Theorem \ref{thmChangeCoord}. Assuming that properties (\ref%
{perturbationb}) hold for $\widetilde{\nu }_{j}(\varepsilon )$ and $%
\widetilde{\varkappa }_{j}(\varepsilon )$ we have (see Theorem
\ref{142t})
for $\varepsilon \leq \varepsilon _{0}$ and $\varepsilon \in \mathfrak{N}%
_{\alpha },$ $\alpha \in (0,1/78)$ the following estimates for
$||U||_{\rho
}\leq \varepsilon ,$ $\ \imath \geq 14,$ $1\leq r\leq s$%
\begin{eqnarray*}
||\widetilde{\delta \phi }||_{r-1} &\leq &\frac{c(M_{14)}}{\varepsilon ^{2}}%
||\frac{\widetilde{h}}{\kappa }||_{r},\text{ \ }1\leq r\leq \imath -9 \\
||\widetilde{\delta \phi }||_{s-1} &\leq &\frac{c(M_{14})}{\varepsilon ^{2}}%
||\frac{\widetilde{h}}{\kappa }||_{r}(1+||U||_{s+19})+c(M_{14})||\frac{%
\widetilde{h}}{\kappa }||_{s}.
\end{eqnarray*}%
From this and (\ref{estimh/k}) we deduce that%
\begin{eqnarray*}
||\widetilde{\delta \phi }||_{r-2} &\leq &\frac{c(M_{14)}}{\varepsilon ^{2}}%
||F||_{r},\text{ \ }2\leq r\leq \imath -8, \\
||\widetilde{\delta \phi }||_{s-2} &\leq &\frac{c(M_{14})}{\varepsilon ^{2}}%
||F||_{r}(1+||U||_{s+18})+c(M_{14})||F||_{s},
\end{eqnarray*}%
and from (\ref{estim1deltaU}) we obtain the estimates for $\delta U$ (\ref%
{estim inv L}). The rest of Theorem \ref{inv L} follows directly
from the results of previous subsection and from Theorem
\ref{142t}.

\section{Nonlinear problem. Proof of Theorem \protect\ref{Thmexistence0}}

\label{sectionnonlin}

In this section we complete the proof of the main Theorem
\ref{Thmexistence0} on existence of diamond nonlinear waves of
finite amplitude. To this end we exploit the general version of
the Nash-Moser Implicit Function Theorem proved in Appendix N of
\cite{IPT}. This result concerns the solvability of the operator
equation
\begin{equation}
\Phi (W,\varepsilon )=0  \label{eqmain}
\end{equation}%
in scales of Banach spaces $E_{s}$ and $F_{s}$ parametrized by $s\in \mathbb{%
N}_{0}=\mathbb{N}\cup \{0\}$, and supplemented with the norms
$||\cdot ||_{s} $ and $|{\cdot }|_{s}$. It is supposed that they
satisfy the following conditions.

\begin{itemize}
\item[(A1)] For $t<s$ there exists $c(t,s)$ such that \newline
$\Vert {\cdot }\Vert _{t}\leq c(t,s)\Vert {\cdot }\Vert
_{s},\qquad |\cdot |_{t}\leq c(t,s)|{\cdot }|_{s}.$

\item[(A2)] For $\lambda \in \lbrack 0,1]$ with $\lambda
t+(1-\lambda )s\in \mathbb{N},$
\begin{equation*}
||{\cdot }||_{\lambda t+(1-\lambda )s}\leq c(t,s)\Vert {\cdot
}\Vert
_{t}^{\lambda }\,\Vert {\cdot }\Vert _{s}^{1-\lambda },\quad |{\cdot }%
|_{\lambda t+(1-\lambda )s}\leq c(t,s)|{\cdot }|{t}^{\lambda }\,|{\cdot }%
|_{s}^{1-\lambda }.
\end{equation*}

\item[(A3)] There exists a family of smoothing operators $S_\wp$
defined over the first scale such that for $\wp>0$ and $t <s$,
\begin{gather*}
\|{S_\wp W}\|_{t}\leq c(t,s)\|{W}\|_{s},\qquad \|{\ S_\wp
W}\|_{s}\leq
c(t,s)\wp^{t-s}\|{W}\|_{t}, \\
\|{S_\wp W-W}\|_{t}\leq c(t,s)\wp^{s-t}\|{W}\|_{s},
\end{gather*}
and, if $\varepsilon\mapsto \wp(\varepsilon)$ is a smooth,
increasing,
convex function on $[0,\infty)$ with $\wp (0) =0$, then, for $%
0<\varepsilon_1<\varepsilon_2$,
\begin{equation*}
\|{\ (S_{\wp(\varepsilon_1)}- S_{\wp(\varepsilon_2)})W}\|_{s}\leq
c(t,s)|\varepsilon_1-\varepsilon_2|\wp^{\prime}(\varepsilon_2)
\wp(\varepsilon_1)^{t-s-1}\|{W}\|_{t}.
\end{equation*}
\end{itemize}

\begin{itemize}
\item[(B1)] Operators $\Phi(\cdot, \varepsilon)$, depend on a
small parameter $\varepsilon\in [0, \varepsilon_0]$, and map a
neighborhood of $0$ in $E_r$ into $F_\rho$. Suppose that there
exist
\begin{equation*}
\sigma \leq\rho\leq r-1,~~ \sigma,\,\rho,\,r\in \mathbb{N}_0,
\end{equation*}
and, for all $l\in \mathbb{N}_0$, numbers $c(l)>0$ and
$\varepsilon(l)\in (0,\varepsilon_0]$ with the following
properties for all $W,\,U,\,W_i,\,U_i \in B$ and
$\varepsilon,\,\varepsilon_i\in [0,\varepsilon_0]$, $i = 1,\,2$,
where $B=\{W\in E_r:|{W}|_{r}\leq R_0\}$ for some $R_0 >0$:

\item[(B2)] The operator $\Phi :B\times \lbrack 0,\varepsilon
_{0}]\rightarrow F_{\rho }$ is twice continuously differentiable,
\begin{equation}
|{\Phi (W,\varepsilon )}|_{\rho +l}\leq c(l)(1+\Vert {W}\Vert
_{r+l}) \label{e0}
\end{equation}%
and, for $W,\,U\in E_{r+l}$, $\varepsilon \in \lbrack
0,\varepsilon (l)]$,
\begin{eqnarray}
|{D(W,U,\varepsilon )}|_{\rho +l} &\leq &c(l)(1+\Vert {W}\Vert _{r+l}+\Vert {%
U}\Vert _{r+l})\Vert {W-U}\Vert _{r}^{2}+  \notag \\
&&+c(l)\Vert {W-U}\Vert _{r}\Vert {W-U}\Vert _{r+l},  \label{e5}
\end{eqnarray}%
where
\begin{equation*}
D(W,U,\varepsilon )=\Phi (W,\varepsilon )-\Phi (U,\varepsilon
)-\Phi _{W}^{\prime }(U,\varepsilon )(W-U).
\end{equation*}%
Moreover,
\begin{eqnarray}
|{D(W_{1},U_{1},\varepsilon _{1})-D(W_{2},U_{2},\varepsilon
_{2})}|_{\rho } &\leq &c\,\big(|\varepsilon _{1}-\varepsilon
_{2}|+\Vert {W_{1}-W_{2}}\Vert
_{r}+\Vert {W_{1}-W_{2}}\Vert _{r}\big)  \notag \\
&&\cdot \big(\Vert {W_{1}-U_{1}}\Vert _{r}+\Vert {W_{2}-U_{2}}\Vert _{r}\big)%
.  \label{e500}
\end{eqnarray}

\item[(B3)] There exists a family of bounded linear operators
$\Lambda (W,\varepsilon ):E_{r}\rightarrow F_{\rho }$, depending
on $(W,\varepsilon )\in B\times \lbrack 0,\varepsilon _{0}]$, with
\begin{equation}
|{\Lambda (W,\varepsilon )U}|_{\rho }\leq c(0)\Vert U\Vert
_{r},~~~U\in E_{r},  \label{e51}
\end{equation}%
that approximates the Fr\'{e}chet derivative $\Phi _{W}^{\prime }$
as follows. For $W\in E_{r+l}\cap B$, $\varepsilon \in \lbrack
0,\varepsilon (l)]$ and $U\in E_{r+l}$,
\begin{eqnarray}
|{\Lambda (W,\varepsilon )U-\Phi _{W}^{\prime }(W,\varepsilon
)U}|_{\rho +l}
&\leq &c(l)(1+\Vert {W}\Vert _{r+l})|{\Phi (W,\varepsilon )}|_{r}\Vert {U}%
\Vert _{r}+  \label{e520} \\
&&+c(l)|{\Phi (W,\varepsilon )}|_{r+l}\Vert {U}\Vert
_{r}+c(l)|{\Phi (W,\varepsilon )}|_{r}\Vert {U}\Vert _{r+l}.
\notag
\end{eqnarray}

\item[(B4)] When $W_{i}\in B\cap E_{r+l}$, $\varepsilon _{i}\in
\lbrack 0,\varepsilon (l)]$, $i=1,2$,
\begin{eqnarray}
|{\Phi (W_{1},\varepsilon _{1})-\Phi (W_{2},\varepsilon
_{2})}|_{\rho +l} &\leq &c(l)\big(1+\Vert {W_{1}}\Vert
_{r+l}+\Vert {W_{2}}\Vert _{r+l}\big)
\label{e10} \\
&&\cdot \big(|\varepsilon _{1}-\varepsilon _{2}|+\Vert
{W_{1}-W_{2}}\Vert _{r}\big)+c(l)\Vert {W_{1}-W_{2}}\Vert _{r+l},
\notag
\end{eqnarray}%
\begin{multline}  \label{e3}
|{(\Phi^{\prime}_W(W_1,\varepsilon_1)-\Phi^{\prime}_W(W_2,\varepsilon_2))U}%
|_{\rho+l}+ |{(\Lambda(W_1,\varepsilon_1)-\Lambda(W_2,\varepsilon_2))U}%
|_{\rho+l}\leq \\
\leq c(l)\Big(\,\|{W_1-W_2}\|_{r+l}+(|\varepsilon_1-\varepsilon_2|+\|{W_1-W_2%
}\|_{r}) (\|{W_1}\|_{r+l}+\|{W_2}\|_{r+l})\Big)\|{U}\|_{r} \\
+\big(|\varepsilon_1-\varepsilon_2|+\|{W_1-W_2}\|_{r}\big)\|{U}\|_{r+l},
\end{multline}
\end{itemize}

\begin{itemize}
\item A set $\mathcal{E}\subset \lbrack 0,\infty )$ is dense at 0 if $%
\displaystyle\lim_{r\searrow
0}\frac{2}{r^{2}}\int\limits_{\mathcal{E}\cap \lbrack
0,r]}\varepsilon \,d\varepsilon =1.$
\end{itemize}

\begin{itemize}
\item[(B5)] If a set $\mathcal{E}\subset \lbrack 0,\varepsilon
(l)]$ is dense at $0$ and a mapping $\vartheta
:\mathcal{E}\rightarrow B\cap E_{r+l}$ is Lipschitz in the sense
that for $\varepsilon _{1},\,\varepsilon _{2}\in \mathcal{E}$,
\begin{equation*}
\Vert {\vartheta (\varepsilon _{1})-\vartheta (\varepsilon
_{2})}\Vert
_{r}\leq C|\varepsilon _{1}-\varepsilon _{2}|\text{ where $C=C(\vartheta )$%
,~ constant,}
\end{equation*}%
then there is a set $\mathcal{E}(\vartheta )\subset \mathcal{E}$,
which is also dense at 0, such that, for any $\varepsilon \in
\mathcal{E}(\vartheta )$ and $f\in F_{\rho +l}$, the equation
$\Lambda (\vartheta (\varepsilon ),\varepsilon )\delta W=f$ has a
unique solution satisfying
\begin{equation}
\Vert {\delta W}\Vert _{\rho -\sigma +l}\leq \varepsilon ^{-\varrho }c(l)(|{f%
}|_{{\rho +l}}+\Vert {\vartheta (\varepsilon )}\Vert
_{r+l}|{f}|_{\rho }). \label{dDd}
\end{equation}

\item[(B6)] Suppose that $\vartheta
_{0}:\mathcal{E}_{0}\rightarrow B\cap E_{r+l}$ and mappings
$\vartheta _{k}:\cap _{i=0}^{k-1}\mathcal{E}(\vartheta
_{i})\rightarrow B\cap E_{r+l}$ satisfy, for a constant $C$ independent of $%
k\in \mathbb{N}$ sufficiently large,
\begin{gather*}
\Vert {\vartheta _{k}(\varepsilon _{1})-\vartheta _{k}(\varepsilon _{2})}%
\Vert _{r}\leq C|\varepsilon _{1}-\varepsilon _{2}|,\quad
\varepsilon
_{1},\,\varepsilon _{2}\in \cap _{j=0}^{k-1}\mathcal{E}(\vartheta _{j}), \\
\Vert {\vartheta _{k+1}(\varepsilon )-\vartheta _{k}(\varepsilon
)}\Vert
_{r}\leq \frac{1}{2^{k}},\quad \varepsilon \in \cap _{j=0}^{k}\mathcal{E}%
(\vartheta _{j}).
\end{gather*}%
Then $\cap _{j=0}^{\infty }\mathcal{E}(\vartheta _{j})$ is dense
at $0$, where the sets $\mathcal{E}(\vartheta _{j})$ are defined
in (B5).
\end{itemize}

\begin{theorem}
\label{implicit}Suppose (A1)--(B6) hold and, for $N\in \mathbb{N}$ with $%
N\geq 2$, equation \eqref{eqmain} has approximate solution
$W=W_{\varepsilon }^{(N)}\in \cap _{s\in \mathbb{N}_{0}}E_{s},$
with, for a constant $k(N,s)$,
\begin{equation}
\Vert W_{\varepsilon }^{(N)}\Vert _{s}\leq k(N,s)\varepsilon
,~~~|\Phi (\varepsilon ,W_{\varepsilon }^{(N)})|_{s}\leq
k(N,s)|\varepsilon |^{N+1} \label{zero}
\end{equation}%
and
\begin{equation}
\Vert W_{\varepsilon _{1}}^{(N)}-W_{\varepsilon _{2}}^{(N)}\Vert
_{s}\leq k(N,s)|\varepsilon _{1}-\varepsilon _{2}|.  \label{orez}
\end{equation}%
Then there is a set $\mathcal{E}$, which is dense at $0$, and a
family
\begin{equation*}
\{W=\vartheta (\varepsilon ):\varepsilon \in \mathcal{E}\}\subset
E_{r}
\end{equation*}%
of solutions to \eqref{eqmain} with $\Vert {\vartheta (\varepsilon
_{1})-\vartheta (\varepsilon _{2})}\Vert _{r}\leq c|\varepsilon
_{1}-\varepsilon _{2}|$ for some constant $c$.
\end{theorem}

In order to apply Theorem \ref{implicit} to the 3D wave problem,
let us
introduce some notations. Fix an arbitrary $\alpha \in (0,1/78)$ and $%
0<\delta <1$. Denote by $\mathcal{N}$ the set of all $\mu _{c}$ so
that $\nu _{0}=\mu _{c}^{-1}$ belongs to the set
$\mathfrak{N}_{\alpha }$ given by Theorem \ref{142t}. Since
$\mathfrak{N}_{\alpha }\subset \lbrack 1,\infty )$
is a set of full measure, $\mathfrak{N}$ is the set of full measure in $%
(0,1) $. Choose an arbitrary $\mu _{c}=(1+\tau ^{2})^{-1/2}\in
\mathfrak{N}$ so that $\tau \in (\delta ,1/\delta )$ and set $\tau
=\tan \theta $. It is
clear that $\mu _{c}$ and $\tau $ meet all requirements of Theorem \ref%
{Thmexistence0}.

Next we fix the lattice $\Gamma $ such that the dual lattice is
spanned by the wave vectors $(1,\pm \tau )$, and set
\begin{equation*}
E_{s}=F_{s}=H_{o,e}^{s}(\mathbb{R}^{2}/\Gamma )\times H_{e,e}^{s}(\mathbb{R}%
^{2}/\Gamma ).
\end{equation*}%
Since the scaling mapping $u\rightarrow u\circ \mathbb{T}^{-1}$
establishes an isomorphism between $H^{s}(\mathbb{R}^{2}/\Gamma )$
and a closed subspace of the Sobolev space $H^{s}$ of doubly $2\pi
-$ periodic functions, the properties $(A1)$ and $(A2)$ are clear.
A smoothing operator with the required properties can be defined
by
\begin{equation*}
S_{\wp }u=\frac{\sqrt{\tau }}{2\pi }\sum\limits_{K\in \Gamma
^{\prime }}\varsigma (\wp |K|)\hat{u}_{K}e^{iK\cdot X},
\end{equation*}%
where $\varsigma :\mathbb{R}^{+}\mapsto \mathbb{R}^{+}$ is a
smooth function which equals $1$ on $[0,1]$ and $0$ on $[2,\infty
)$.

Fix an arbitrary $N>3$ and define the operator $\Phi =(\Phi
_{1},\Phi _{2})$ by the equalities
\begin{equation}
\Phi (W,\varepsilon )=\varepsilon
^{-N}\mathcal{F}\big(U_{\varepsilon }^{(2N)}+\varepsilon ^{N}W,\mu
_{\varepsilon }^{(2N)}\big),  \label{PHI}
\end{equation}%
where $(U_{\varepsilon }^{(2N)},\mu _{\varepsilon }^{(2N)})$ is
the
approximate solution at order $\varepsilon ^{2N}$ obtained at Theorem \ref%
{Lembifurc} in the case of diamond waves $(\varepsilon
_{1}=\varepsilon
_{2}=\varepsilon /2),$ and $\tau =\tan \theta ,$ $\mu _{c}=\cos \theta ,$ $%
\delta <\tau <1/\delta .$ By construction $W_{\varepsilon }^{(N)}=0,$ hence (%
\ref{zero}) and (\ref{orez}) are satisfied.

It follows from Lemma \ref{F(U,mu)} that the continuous mapping $%
(W,\varepsilon )\mapsto \Phi (W,\varepsilon )$ from $E_{s+1}\times
\mathbb{R}
$ into $F_{s}$ for $s\geq 2,$ is of class of $C^{\infty }$.
Applying the
same arguments as in section 9 of \cite{IPT} we conclude from this that for $%
5\leq \rho \leq r-18$, the operator $\Phi $ satisfies Conditions
(B1) and
(B2 ), and inequality \eqref{e10} from Condition (B4). Let us denote by $%
\Lambda $ the approximate differential, defined as
\begin{equation*}
\Lambda (W,\varepsilon )=\mathcal{L}(U_{\varepsilon
}^{(2N)}+\varepsilon ^{N}W,\mu _{\varepsilon }^{(2N)})\left(
\begin{array}{cc}
1 & -\mathfrak{b} \\
0 & 1%
\end{array}%
\right)
\end{equation*}%
where the operator $\mathfrak{L}$ and coefficient $\mathfrak{b}$
are defined
by formulae \eqref{Z} and \eqref{basiclinop}. It follows from this and (\ref%
{defR}) that the linear operator
\begin{equation*}
\mathcal{R}=\Lambda (W,\varepsilon )-\partial _{W}\Phi
(W,\varepsilon ),
\end{equation*}%
is defined by $\mathcal{R}\delta W=(R_{1}\delta W,0)$ with
\begin{equation*}
R_{1}\delta W=\mathcal{G}_{\eta }\left( \frac{\Phi _{1}\delta \eta }{%
1+(\nabla \eta )^{2}}\right) -\nabla \cdot \left( \frac{\Phi
_{1}\delta \eta }{1+(\nabla \eta )^{2}}\nabla \eta \right) ,
\end{equation*}%
where $\delta W=(\delta \psi ,\delta \eta )$ and $\Phi _{1}=\Phi
_{1}(W,\varepsilon )$. Then, thanks to Lemma \ref{tame coef}, we have for $%
s\geq 2,$ and $||U||_{3}\leq M_{3}$%
\begin{equation*}
||\Lambda (W,\varepsilon )u||_{s}\leq
c(M_{3})(1+||W||_{s+2})||u||_{s+1},
\end{equation*}%
\begin{equation*}
||\mathcal{R}u||_{s}\leq c_{s}(M_{3})\{||\Phi
||_{2}(1+||W||_{s+2})||u||_{2}+||u||_{s+1})+||\Phi
||_{s+1}||u||_{2}.
\end{equation*}%
Therefore, the operators $\Phi $ and $\Lambda $ satisfy Conditions
(B3) and (B4). Now taking into account the result of Theorem
\ref{inv L}, it appears that Condition (B5) is satisfied for
\begin{eqnarray*}
5 &\leq &\rho \leq r-18 \\
\sigma &=&3,\text{ \ \ }\varrho =2,
\end{eqnarray*}%
Finally, since $N\geq 3,$ condition (B6) is also satisfied thanks to Theorem %
\ref{142t}, so we can apply Theorem \ref{implicit}, and Theorem \ref%
{Thmexistence0} is proved.

\appendix

\section{Analytical study of $\mathcal{G}_{\protect\eta }$}

\subsection{Computation of the differential of $\mathcal{G}_{\protect\eta }$%
\label{a0}}

Let us make the following change of coordinate for $x_{3}$%
\begin{equation*}
s=x_{3}-\eta (X)
\end{equation*}%
and denote by $\theta $ the function defined by%
\begin{equation*}
\theta (X,s)=\varphi (X,s+\eta (X)).
\end{equation*}%
Then the Dirichlet-Neumann operator is defined by%
\begin{equation}
\Delta \theta -2\nabla _{X}\eta \cdot \nabla _{X}(\frac{\partial \theta }{%
\partial s})-\frac{\partial \theta }{\partial s}\Delta _{X}\eta +\frac{%
\partial ^{2}\theta }{\partial s^{2}}(\nabla _{X}\eta )^{2}=0,\text{ \ }s<0
\label{equtheta}
\end{equation}%
\begin{eqnarray*}
\theta |_{s=0} &=&\psi \\
\nabla \theta &\rightarrow &0\text{ \ as \ }s\rightarrow -\infty ,
\end{eqnarray*}%
and%
\begin{equation}
\mathcal{G}_{\eta }\psi =(1+(\nabla _{X}\eta )^{2})\frac{\partial \theta }{%
\partial s}|_{s=0}-\nabla _{X}\eta \cdot \nabla _{X}\psi .
\label{Gdemi-plan}
\end{equation}%
Notice that the above equations (\ref{equtheta}),
(\ref{Gdemi-plan}) may be
written into the form%
\begin{eqnarray*}
\nabla \cdot (P_{\eta }\nabla \theta ) &=&0,\text{ \ }s<0, \\
\mathcal{G}_{\eta }\psi &=&(P_{\eta }\nabla \theta )\cdot
\mathbf{e}_{3}
\end{eqnarray*}%
where $\mathbf{e}_{3}$ is the unit vertical vector, and $P_{\eta
}$ is the
following symmetric matrix%
\begin{equation*}
P_{\eta }=\left(
\begin{array}{cc}
\mathbb{I} & -\nabla _{X}\eta \\
-(\nabla _{X}\eta )^{t} & 1+(\nabla _{X}\eta )^{2}%
\end{array}%
\right) .
\end{equation*}%
We can see easily that the operator $\mathcal{G}_{\eta }$\textbf{\
}\emph{is
symmetric and non negative in} $L^{2}(%
\mathbb{R}
^{2}/\Gamma ):$ for $\eta ,$ $\psi _{1},$ $\psi _{2}$ smooth
enough
bi-periodic functions%
\begin{eqnarray*}
\langle \mathcal{G}_{\eta }\psi _{1},\psi _{2}\rangle &=&\langle
(P_{\eta
}\nabla \theta _{1})\cdot \mathbf{e}_{3},\psi _{2}\rangle \\
&=&\int_{-\infty }^{0}\int_{%
\mathbb{R}
^{2}/\Gamma }(\nabla \theta _{2}\cdot P_{\eta }\nabla \theta
_{1})dXds
\end{eqnarray*}%
which is symmetric. Moreover, we have%
\begin{equation*}
\nabla \theta \cdot P_{\eta }\nabla \theta =(\nabla _{X}\theta
-\nabla
_{X}\eta \frac{\partial \theta }{\partial s})^{2}+(\frac{\partial \theta }{%
\partial s})^{2}\geq 0,
\end{equation*}%
hence, for $\eta $ and $\psi $ smooth enough bi-periodic functions%
\begin{equation*}
\langle \mathcal{G}_{\eta }\psi ,\psi \rangle \geq 0.
\end{equation*}%
Let us now compute formally the differential of $\mathcal{G}_{\eta }.$ >From (%
\ref{Gdemi-plan}), we obtain%
\begin{equation}
\partial _{\eta }\mathcal{G}_{\eta }[h]\psi =(P_{\eta }\nabla \theta
_{1})\cdot \mathbf{e}_{3}+(Q_{\eta }[h]\nabla \theta )\cdot
\mathbf{e}_{3} \label{diffG_1}
\end{equation}%
with $\theta $ as above, and where $\theta _{1}$ satisfies the system%
\begin{eqnarray*}
\nabla \cdot (P_{\eta }\nabla \theta _{1}) &=&-\nabla \cdot
(Q_{\eta
}[h]\nabla \theta ),\text{ \ }s<0 \\
\theta _{1}|_{s=0} &=&0, \\
\nabla \theta _{1} &\rightarrow &0,\text{ \ \ }s\rightarrow
-\infty ,
\end{eqnarray*}%
\begin{equation*}
Q_{\eta }[h]=\left(
\begin{array}{cc}
0 & -\nabla _{X}h \\
-(\nabla _{X}h)^{t} & 2\nabla _{X}\eta \cdot \nabla _{X}h%
\end{array}%
\right) .
\end{equation*}%
Let us notice from (\ref{equtheta}), (\ref{Gdemi-plan}), that we have%
\begin{equation*}
\nabla \cdot \left( P_{\eta }\nabla (h\frac{\partial \theta }{\partial s}%
)\right) +\nabla \cdot (Q_{\eta }[h]\nabla \theta )=0,\text{ \
}s<0.
\end{equation*}%
Indeed, we have%
\begin{eqnarray*}
P_{\eta }\nabla (h\frac{\partial \theta }{\partial s}) &=&hP_{\eta
}(\nabla
\frac{\partial \theta }{\partial s})+\frac{\partial \theta }{\partial s}%
\left(
\begin{array}{c}
\nabla _{X}h \\
-\nabla _{X}\eta \cdot \nabla _{X}h%
\end{array}%
\right) , \\
Q_{\eta }[h]\nabla \theta &=&\left(
\begin{array}{c}
-\nabla _{X}h\frac{\partial \theta }{\partial s} \\
-\nabla _{X}\theta \cdot \nabla _{X}h+2(\nabla _{X}\eta \cdot \nabla _{X}h)%
\frac{\partial \theta }{\partial s}%
\end{array}%
\right) ,
\end{eqnarray*}%
hence%
\begin{equation*}
\nabla \cdot \left( P_{\eta }\nabla (h\frac{\partial \theta }{\partial s}%
)\right) +\nabla \cdot (Q_{\eta }[h]\nabla \theta )=
\end{equation*}%
\begin{equation*}
=h\nabla \cdot (P_{\eta }\nabla \frac{\partial \theta }{\partial
s})+\nabla
_{X}h\cdot \left\{ (P_{\eta }\nabla \frac{\partial \theta }{\partial s}%
)-\nabla _{X}\frac{\partial \theta }{\partial s}+\nabla _{X}\eta \frac{%
\partial ^{2}\theta }{\partial s^{2}}\right\} ,
\end{equation*}%
and this cancels, thanks to $\nabla \cdot (P_{\eta }\nabla \theta
)=0,$ and to the definition of $P_{\eta }.$

It results that%
\begin{eqnarray*}
\theta _{1} &=&h\frac{\partial \theta }{\partial s}+\theta _{2}, \\
\nabla \cdot (P_{\eta }\nabla \theta _{2}) &=&0,\text{ \ }s<0, \\
\theta _{2}|_{s=0} &=&-h\frac{\partial \theta }{\partial s}|_{s=0}, \\
\nabla \theta _{2} &\rightarrow &0,\text{ \ \ }s\rightarrow
-\infty ,
\end{eqnarray*}%
hence looking at relationship (\ref{diffG_1}), we obtain%
\begin{eqnarray*}
\partial _{\eta }\mathcal{G}_{\eta }[h]\psi &=&(P_{\eta }\nabla \theta
_{2})\cdot \mathbf{e}_{3}+(P_{\eta }\nabla (h\frac{\partial \theta }{%
\partial s}))\cdot \mathbf{e}_{3}+(Q_{\eta }[h]\nabla \theta )\cdot \mathbf{e%
}_{3} \\
&=&-\mathcal{G}_{\eta }(h\frac{\partial \theta }{\partial s}%
|_{0})+h\{(1+(\nabla _{X}\eta )^{2})\frac{\partial ^{2}\theta
}{\partial
s^{2}}|_{0}-\nabla _{X}\eta \cdot \nabla _{X}\frac{\partial \theta }{%
\partial s}|_{0}\}+ \\
&&+\nabla _{X}h\cdot \{\nabla _{X}\eta \frac{\partial \theta }{\partial s}%
|_{0}-\nabla _{X}\psi \},
\end{eqnarray*}%
and using $\nabla \cdot (P_{\eta }\nabla \theta )=0,$ we get%
\begin{eqnarray*}
\partial _{\eta }\mathcal{G}_{\eta }[h]\psi &=&-\mathcal{G}_{\eta }(h\frac{%
\partial \theta }{\partial s}|_{0})+h\{-\Delta _{X}\psi +\nabla _{X}\eta
\cdot \nabla _{X}\frac{\partial \theta }{\partial s}|_{0}+\Delta
_{X}\eta
\frac{\partial \theta }{\partial s}|_{0}\}+ \\
&&+\nabla _{X}h\cdot \{\nabla _{X}\eta \frac{\partial \theta }{\partial s}%
|_{0}-\nabla _{X}\psi \} \\
&=&-\mathcal{G}_{\eta }(h\frac{\partial \theta }{\partial
s}|_{0})+\nabla
_{X}\cdot \{h(\nabla _{X}\eta \frac{\partial \theta }{\partial s}%
|_{0}-\nabla _{X}\psi )\}
\end{eqnarray*}%
as required in Lemma \ref{Lemma diff G}, thanks to
(\ref{Gdemi-plan}).

\subsection{Second order Taylor expansion of $\mathcal{G}_{\protect\eta }$
in $\protect\eta =0$}

\label{a01}Let us consider the Taylor expansion of
$\mathcal{G}_{\eta }$ in
"powers" of $\eta :$%
\begin{equation*}
\mathcal{G}_{\eta }\psi =\mathcal{G}^{(0)}\psi
+\mathcal{G}^{(1)}\{\eta \}\psi +\mathcal{G}^{(2)}\{\eta
^{(2)}\}\psi +...
\end{equation*}%
where $\mathcal{G}^{(k)}$ is $k-$ linear symmetric with respect to
$\eta ,$
and linear in $\psi .$ Moreover, for any $k\geq 1$ and $m\geq 2,$ $\mathcal{G%
}^{(k)}$ is bounded from%
\begin{equation*}
\{H^{m+1}(%
\mathbb{R}
^{2}/\Gamma )\}^{k}\text{ \ into \ }\mathcal{L}(H_{0}^{m+1}(%
\mathbb{R}
^{2}/\Gamma ),H_{0}^{m}(%
\mathbb{R}
^{2}/\Gamma )).
\end{equation*}%
From the lemma above and (\ref{diffG_2}), (\ref{diffG_3}), we obtain%
\begin{equation}
\mathcal{G}^{(1)}\{\eta \}\psi =-\mathcal{G}^{(0)}(\eta \mathcal{G}%
^{(0)}\psi )-\nabla \cdot (\eta \nabla \psi ).  \label{G^1}
\end{equation}%
Differentiating once more (\ref{diffG_3}) with respect to $\eta ,$
we now
obtain in 0%
\begin{equation*}
\partial _{\eta }\zeta \lbrack h]=-\mathcal{G}^{(0)}(h\mathcal{G}^{(0)}\psi
)-h\Delta \psi ,
\end{equation*}%
and differentiating (\ref{diffG_2}) with respect to $\eta $, we
then obtain
the (symmetric) second order derivative in 0, and%
\begin{equation}
\mathcal{G}^{(2)}\{\eta ^{(2)}\}\psi =\mathcal{G}^{(0)}(\eta \mathcal{G}%
^{(0)}(\eta \mathcal{G}^{(0)}\psi
))+\frac{1}{2}\mathcal{G}^{(0)}(\eta ^{2}\Delta \psi
)+\frac{1}{2}\Delta (\eta ^{2}\mathcal{G}^{(0)}\psi ). \label{G^2}
\end{equation}%
Now, if we Fourier expand any bi-periodic function as%
\begin{equation*}
\psi =\sum_{K\in \Gamma ^{\prime }}\psi _{K}e^{iK\cdot X}
\end{equation*}%
then we have%
\begin{equation}
\{\mathcal{G}^{(0)}\psi \}_{K}=|K|\psi _{K},  \label{G0}
\end{equation}%
\begin{equation}
\{\mathcal{G}^{(1)}\{\eta \}\psi \}_{K}=\sum_{K_{1}+K_{2}=K,\text{ \ }%
K_{j}\in \Gamma ^{\prime }}\{(K\cdot K_{1})-|K||K_{1}|\}\psi
_{K_{1}}\eta _{K_{2}},  \label{G1}
\end{equation}%
\begin{eqnarray}
\{\mathcal{G}^{(2)}\{\eta ^{(2)}\}\psi \}_{K} &=&\sum_{K_{1}+K_{2}+K_{3}=K,%
\text{ \ }K_{j}\in \Gamma ^{\prime }}  \label{G2} \\
&&\frac{|K||K_{1}|}{2}\{|K_{1}+K_{2}|+|K_{1}+K_{3}|-|K|-|K_{1}|\}\psi
_{K_{1}}\eta _{K_{2}}\eta _{K_{3}},  \notag
\end{eqnarray}%
and we check that, for $m\geq 2,$ the operators are bounded from $\{H^{m+1}(%
\mathbb{R}
^{2}/\Gamma )\}^{k}$ \ into \ $\mathcal{L}(H_{0}^{m+1}(%
\mathbb{R}
^{2}/\Gamma ),H_{0}^{m}(%
\mathbb{R}
^{2}/\Gamma ))$ since we have%
\begin{equation*}
|(K_{1}+K_{2})\cdot K_{1})-|K_{1}+K_{2}||K_{1}||\leq
4|K_{1}||K_{2}|,
\end{equation*}%
and there exists a constant $c$ such that%
\begin{equation*}
|K_{1}+K_{2}+K_{3}|||K_{1}+K_{2}|+|K_{1}+K_{3}|-|K_{1}+K_{2}+K_{3}|-|K_{1}||%
\leq
\end{equation*}%
\begin{equation*}
\leq c|K_{2}||K_{3}|.
\end{equation*}

\section{Formal computation of 3-dimensional waves\label{a02}}

In taking $\mathbf{u}_{0}=(1,0)$, the symmetric linearized
operator for $\mu =\mu _{c}$ $\ $and $\mathbf{u}=\mathbf{u}_{0}$
reads
\begin{equation*}
\mathcal{L}_{0}=\left(
\begin{array}{cc}
\mathcal{G}^{(0)} & -\frac{\partial }{\partial x_{1}} \\
\frac{\partial }{\partial x_{1}} & \mu _{c}%
\end{array}%
\right) ,
\end{equation*}%
and $\mathcal{L}_{0}$ has a four-dimensional kernel, spanned by the vectors%
\begin{eqnarray*}
\zeta _{0} &=&(1,-i/\mu _{c})e^{iK_{1}\cdot X},\text{ }\overline{\zeta }%
_{0}=(1,i/\mu _{c})e^{-iK_{1}\cdot X}, \\
\zeta _{1} &=&(1,-i/\mu _{c})e^{iK_{2}\cdot X},\text{ }\overline{\zeta }%
_{1}=(1,i/\mu _{c})e^{-iK_{2}\cdot X}.
\end{eqnarray*}%
We observe that the action of different symmetries of the system
on
eigenvectors is as follows:%
\begin{eqnarray*}
\mathcal{T}_{\mathbf{v}}\zeta _{0} &=&\zeta _{0}e^{iK_{1}\cdot \mathbf{v}},%
\text{ }T_{\mathbf{v}}\zeta _{1}=\zeta _{1}e^{iK_{2}\cdot \mathbf{v}}, \\
\mathcal{S}_{0}\zeta _{0} &=&-\overline{\zeta }_{0},\text{ \ }\mathcal{S}%
_{0}\zeta _{1}=-\overline{\zeta }_{1}, \\
\mathcal{S}_{1}\zeta _{0} &=&\zeta _{1},\text{ \ \
}\mathcal{S}_{1}\zeta _{1}=\zeta _{0}.
\end{eqnarray*}%
Let us write formally the nonlinear system (\ref{basic1}),
(\ref{basic2})
under the form%
\begin{equation}
\mathcal{L}_{0}U+\tilde{\mu}\mathcal{L}_{1}U+\mathcal{L}_{2}(\mathbf{\omega }%
,U)+\mathcal{N}_{2}(U,U)+\mathcal{N}_{3}(U,U,U)+O(||U||^{4}+|\mathbf{\omega }%
|||U||^{2})=0,  \label{expandsyst}
\end{equation}%
with $U=(\psi ,\eta ),$ $\tilde{\mu}=\mu -\mu _{c},$ $\mathbf{\omega }=%
\mathbf{u}-\mathbf{u}_{0}$%
\begin{eqnarray*}
\mathcal{L}_{1}U &=&(0,\eta ), \\
\mathcal{L}_{2}(\mathbf{\omega },U) &=&(-\mathbf{\omega }\cdot \nabla \eta ,%
\mathbf{\omega }\cdot \nabla \psi ),
\end{eqnarray*}%
\begin{equation*}
\mathcal{N}_{2}(U,U)=\left\{
\begin{array}{c}
\mathcal{G}^{(1)}\{\eta \}\psi \\
\frac{1}{2}\nabla \psi ^{2}-\frac{1}{2}(\frac{\mathcal{\partial \eta }}{%
\partial x_{1}})^{2}%
\end{array}%
\right. ,
\end{equation*}%
\begin{equation*}
\mathcal{N}_{3}(U,U,U)=\left\{
\begin{array}{c}
\mathcal{G}^{(2)}\{\eta ,\eta \}\psi \\
-\frac{\mathcal{\partial \eta }}{\partial x_{1}}(\nabla \eta \cdot
\nabla
\psi )%
\end{array}%
\right. .
\end{equation*}

\subsection{Formal Fredholm alternative}

Let us consider the formal resolution of the linear system%
\begin{equation*}
\mathcal{L}_{0}U=F=(f,g),
\end{equation*}%
with%
\begin{eqnarray*}
U &=&\sum_{n=(n_{1},n_{2})\in
\mathbb{Z}
^{2}}U_{n}e^{i(n_{1}K_{1}\cdot X+n_{2}K_{2}\cdot X)},\text{ \
}U_{n}=(\psi
_{n},\eta _{n}),\text{ \ }\psi _{0}=0, \\
F &=&\sum_{n=(n_{1},n_{2})\in
\mathbb{Z}
^{2}}F_{n}e^{i(n_{1}K_{1}\cdot X+n_{2}K_{2}\cdot X)},\text{ \ }%
F_{n}=(f_{n},g_{n}),\text{ \ }f_{0}=0.
\end{eqnarray*}%
Then, we have%
\begin{eqnarray*}
|n_{1}K_{1}+n_{2}K_{2}|\psi _{n}-i\{(n_{1}K_{1}+n_{2}K_{2})\cdot \mathbf{u}%
_{0}\}\eta _{n} &=&f_{n} \\
i\{(n_{1}K_{1}+n_{2}K_{2})\cdot \mathbf{u}_{0}\}\psi _{n}+\mu
_{c}\eta _{n} &=&g_{n},
\end{eqnarray*}%
hence for
\begin{equation*}
\{(n_{1}K_{1}+n_{2}K_{2})\cdot \mathbf{u}_{0}\}^{2}-\mu
_{c}|n_{1}K_{1}+n_{2}K_{2}|\neq 0
\end{equation*}%
i.e. by assumption for $(n_{1},n_{2})\neq (\pm 1,0),(0,\pm 1),$ this leads to%
\begin{eqnarray}
\psi _{n} &=&-\frac{\mu _{c}f_{n}+i\{(n_{1}K_{1}+n_{2}K_{2})\cdot \mathbf{u}%
_{0}\}g_{n}}{\{(n_{1}K_{1}+n_{2}K_{2})\cdot
\mathbf{u}_{0}\}^{2}-\mu
_{c}|n_{1}K_{1}+n_{2}K_{2}|},  \label{fredholm1} \\
\eta _{n} &=&\frac{i\{(n_{1}K_{1}+n_{2}K_{2})\cdot \mathbf{u}%
_{0}\}f_{n}-|n_{1}K_{1}+n_{2}K_{2}|g_{n}}{\{(n_{1}K_{1}+n_{2}K_{2})\cdot
\mathbf{u}_{0}\}^{2}-\mu _{c}|n_{1}K_{1}+n_{2}K_{2}|},
\label{fredholm2}
\end{eqnarray}%
and for $(n_{1},n_{2})=(0,0)$%
\begin{equation*}
\psi _{0,0}=0,\text{ \ }\eta _{0,0}=\frac{1}{\mu _{c}}g_{0,0},
\end{equation*}%
while for $(n_{1},n_{2})=(\pm 1,0),(0,\pm 1),$ we need to satisfy
the compatibility conditions
\begin{equation*}
(F,\zeta _{0})=(F,\overline{\zeta }_{0})=(F,\zeta _{1})=(F,\overline{\zeta }%
_{1})=0
\end{equation*}%
which gives%
\begin{eqnarray*}
\mu _{c}f_{1,0}+ig_{1,0} &=&0, \\
\mu _{c}f_{-1,0}-ig_{-1,0} &=&0, \\
\mu _{c}f_{0,1}+ig_{0,1} &=&0, \\
\mu _{c}f_{0,-1}-ig_{0,-1} &=&0.
\end{eqnarray*}%
For uniqueness of the definition of the pseudo-inverse $\widetilde{\mathcal{L%
}}_{0}^{-1}$, we fix $U$ such that
\begin{equation*}
(U,\zeta _{0})=(U,\overline{\zeta }_{0})=(U,\zeta _{1})=(U,\overline{\zeta }%
_{1})=0,
\end{equation*}%
hence this leads to%
\begin{eqnarray*}
\psi _{1,0} &=&-i\frac{1+\tau ^{2}}{2+\tau
^{2}}g_{1,0}=\frac{1}{\mu
_{c}(2+\tau ^{2})}f_{1,0}, \\
\eta _{1,0} &=&\frac{1}{\mu _{c}(2+\tau ^{2})}g_{1,0}=\frac{i}{2+\tau ^{2}}%
f_{1,0},
\end{eqnarray*}%
\begin{eqnarray*}
\psi _{0,1} &=&-i\frac{1+\tau ^{2}}{2+\tau
^{2}}g_{0,1}=\frac{1}{\mu
_{c}(2+\tau ^{2})}f_{0,1}, \\
\eta _{0,1} &=&\frac{1}{\mu _{c}(2+\tau ^{2})}g_{0,1}=\frac{i}{2+\tau ^{2}}%
f_{0,1}.
\end{eqnarray*}

\subsection{Bifurcation equation}

Now coming back to (\ref{expandsyst}), we use a \emph{formal}
Lyapunov -
Schmidt method and decompose $U$ as follows%
\begin{eqnarray*}
U &=&W+V \\
W &=&A\zeta _{0}+\overline{A}\overline{\zeta }_{0}+B\zeta _{1}+\overline{B}%
\overline{\zeta }_{1}=\mathcal{P}_{0}U \\
(V,\zeta _{0}) &=&(V,\overline{\zeta }_{0})=(V,\zeta _{1})=(V,\overline{%
\zeta }_{1})=0.
\end{eqnarray*}%
We can solve formally with respect to $V$ the part of equ. (\ref{expandsyst}%
) orthogonal to the 4-dimensional kernel of $\mathcal{L}_{0},$ as
a uniquely
determined formal power series in $\omega ,\tilde{\mu},A,\overline{A},B,%
\overline{B},$ which we write as%
\begin{equation*}
V=\mathcal{V}(\tilde{\mu},\mathbf{\omega
},A,\overline{A},B,\overline{B}).
\end{equation*}%
The uniqueness of the series and the symmetries of our system lead
to the following identities
\begin{eqnarray*}
\mathcal{T}_{\mathbf{v}}\mathcal{V}(\tilde{\mu},\mathbf{\omega },A,\overline{%
A},B,\overline{B}) &=&\mathcal{V}(\tilde{\mu},\mathbf{\omega },Ae^{iK_{1}.%
\mathbf{v}},\overline{A}e^{-iK_{1}.\mathbf{v}},Be^{iK_{2}.\mathbf{v}},%
\overline{B}e^{-iK_{2}.\mathbf{v}}), \\
\mathcal{S}_{0}\mathcal{V}(\tilde{\mu},\mathbf{\omega },A,\overline{A},B,%
\overline{B}) &=&\mathcal{V}(\tilde{\mu},\mathbf{\omega },-\overline{A},-A,-%
\overline{B},-B) \\
\mathcal{SV}(\tilde{\mu},\mathbf{0},A,\overline{A},B,\overline{B}) &=&%
\mathcal{V}(\tilde{\mu},\mathbf{0},B,\overline{B},A,\overline{A}).
\end{eqnarray*}%
The principal part of $V$ is given by%
\begin{eqnarray*}
V &=&-\widetilde{\mathcal{L}}_{0}^{-1}(\mathbb{I}-P_{0})\mathcal{N}%
_{2}(W,W)+O\{(|\tilde{\mu}|+|\mathbf{\omega }|)||W||+||W||^{3}\}, \\
-\widetilde{\mathcal{L}}_{0}^{-1}\mathcal{N}_{2}(W,W)
&=&A^{2}U_{2000}+|A|^{2}U_{1100}+\overline{A}^{2}U_{0200}+ABU_{1010}+%
\overline{A}BU_{0110}+ \\
&&+A\overline{B}U_{1001}+\overline{AB}U_{0101}+B^{2}U_{0020}+|B|^{2}U_{0011}+%
\overline{B}^{2}U_{0002}
\end{eqnarray*}%
where we observe easily that
\begin{equation*}
\mathcal{P}_{0}\mathcal{N}_{2}(W,W)=0.
\end{equation*}%
Using (\ref{G1}), and%
\begin{equation}
2\mathcal{N}_{2}(U_{1},U_{2})=\left\{
\begin{array}{c}
\mathcal{G}^{(1)}\{\eta _{1}\}\psi _{2}+\mathcal{G}^{(1)}\{\eta
_{2}\}\psi
_{1} \\
\nabla \psi _{1}\cdot \nabla \psi _{2}-\frac{\partial \eta
_{1}}{\partial
x_{1}}\frac{\partial \eta _{2}}{\partial x_{1}}%
\end{array}%
\right. ,  \label{bilin}
\end{equation}%
we find%
\begin{eqnarray*}
\mathcal{N}_{2}(W,W) &=&A^{2}V_{2000}+|A|^{2}V_{1100}+\overline{A}%
^{2}V_{0200}+ABV_{1010}+\overline{A}BV_{0110}+ \\
&&+A\overline{B}V_{1001}+\overline{AB}V_{0101}+B^{2}V_{0020}+|B|^{2}V_{0011}+%
\overline{B}^{2}V_{0002}
\end{eqnarray*}%
with%
\begin{eqnarray*}
V_{2000} &=&(0,-\frac{1}{\mu _{c}^{2}})e^{2iK_{1}\cdot X},\text{ \ }%
V_{0020}=(0,-\frac{1}{\mu _{c}^{2}})e^{2iK_{2}\cdot X} \\
V_{1100} &=&0,\text{ \ }V_{0011}=0,\text{ \ }V_{0200}=\overline{V}_{2000},%
\text{ \ \ }V_{0002}=\overline{V}_{0020} \\
V_{1010} &=&(\frac{4i}{\mu _{c}^{2}}(1-\mu
_{c}),-2)e^{i(K_{1}\cdot
X+K_{2}\cdot X)},\text{ \ }V_{0101}=\overline{V}_{1010} \\
V_{0110} &=&(0,-2\tau ^{2})e^{i(-K_{1}\cdot X+K_{2}\cdot X)},\text{ \ \ }%
V_{1001}=\overline{V}_{0110}.
\end{eqnarray*}%
Now, thanks to (\ref{fredholm1}), (\ref{fredholm2}) we obtain $-\widetilde{%
\mathcal{L}}_{0}^{-1}\mathcal{N}_{2}(W,W)$ as follows%
\begin{eqnarray*}
U_{2000} &=&\left( \frac{-i}{\mu _{c}^{2}},-\frac{1}{\mu
_{c}^{3}}\right)
e^{2iK_{1}\cdot X},\text{ \ }U_{0020}=\left( \frac{-i}{\mu _{c}^{2}},-\frac{1%
}{\mu _{c}^{3}}\right) e^{2iK_{2}\cdot X}, \\
U_{1100} &=&0,\text{ \ }U_{0011}=0,\text{ \ }U_{0200}=\overline{U}_{2000},%
\text{ \ \ }U_{0002}=\overline{U}_{0020}, \\
U_{1010} &=&\left( \frac{-2i(2\mu _{c}-1)}{\mu _{c}(2-\mu _{c})},\frac{%
-2(\mu _{c}^{2}+2\mu _{c}-2)}{\mu _{c}^{2}(2-\mu _{c})}\right)
e^{i(K_{1}\cdot X+K_{2}\cdot X)},\text{ \ }U_{0101}=\overline{U}_{1010}, \\
U_{0110} &=&\left( 0,\frac{2\tau ^{2}}{\mu _{c}}\right)
e^{i(-K_{1}\cdot X+K_{2}\cdot X)},\text{ \ \
}U_{1001}=\overline{U}_{0110}.
\end{eqnarray*}%
Replacing $V$ by $\mathcal{V}(\tilde{\mu},\mathbf{\omega },A,\overline{A},B,%
\overline{B})$ in the compatibility conditions, i.e. the components of (\ref%
{expandsyst}) on $\ker \mathcal{L}_{0}$,%
\begin{eqnarray*}
\langle \tilde{\mu}\mathcal{L}_{1}(W+\mathcal{V})+\mathcal{L}_{2}(\mathbf{%
\omega },W+\mathcal{V})+\mathcal{N}_{2}(W+\mathcal{V},W+\mathcal{V}%
)+...,\zeta _{0}\rangle &=&0, \\
\langle \tilde{\mu}\mathcal{L}_{1}(W+\mathcal{V})+\mathcal{L}_{2}(\mathbf{%
\omega },W+\mathcal{V})+\mathcal{N}_{2}(W+\mathcal{V},W+\mathcal{V}%
)+...,\zeta _{1}\rangle &=&0,
\end{eqnarray*}%
noticing that the complex conjugate equations are then
automatically
satisfied, lead to two complex equations of the form%
\begin{eqnarray*}
f(\tilde{\mu},\mathbf{\omega },A,\overline{A},B,\overline{B}) &=&0, \\
g(\tilde{\mu},\mathbf{\omega },A,\overline{A},B,\overline{B})
&=&0,
\end{eqnarray*}%
for which the equivariance of the system (\ref{expandsyst}) with
respect to
various symmetries leads to the following properties%
\begin{eqnarray*}
f(\tilde{\mu},\mathbf{\omega },Ae^{iK_{1}.\mathbf{v}},\overline{A}e^{-iK_{1}.%
\mathbf{v}},Be^{iK_{2}.\mathbf{v}},\overline{B}e^{-iK_{2}.\mathbf{v}})
&=&e^{iK_{1}.\mathbf{v}}f(\tilde{\mu},\mathbf{\omega },A,\overline{A},B,%
\overline{B}) \\
g(\tilde{\mu},\mathbf{\omega },Ae^{iK_{1}.\mathbf{v}},\overline{A}e^{-iK_{1}.%
\mathbf{v}},Be^{iK_{2}.\mathbf{v}},\overline{B}e^{-iK_{2}.\mathbf{v}})
&=&e^{iK_{2}.\mathbf{v}}g(\tilde{\mu},\mathbf{\omega },A,\overline{A},B,%
\overline{B}) \\
f(\tilde{\mu},\mathbf{\omega },-\overline{A},-A,-\overline{B},-B) &=&-%
\overline{f}(\tilde{\mu},\mathbf{\omega },A,\overline{A},B,\overline{B}) \\
g(\tilde{\mu},\mathbf{\omega },-\overline{A},-A,-\overline{B},-B) &=&-%
\overline{g}(\tilde{\mu},\mathbf{\omega },A,\overline{A},B,\overline{B}) \\
f(\tilde{\mu},\mathbf{0},B,\overline{B},A,\overline{A}) &=&g(\tilde{\mu},%
\mathbf{0},A,\overline{A},B,\overline{B}).
\end{eqnarray*}%
Since $K_{1}$ and $K_{2}$ are linearly independent, it results that $f$ and $%
g$ take formally the form%
\begin{eqnarray*}
f(\tilde{\mu},\mathbf{\omega },A,\overline{A},B,\overline{B}) &=&A\phi _{0}(%
\tilde{\mu},\mathbf{\omega },|A|^{2},|B|^{2}) \\
g(\tilde{\mu},\mathbf{\omega },A,\overline{A},B,\overline{B}) &=&B\phi _{1}(%
\tilde{\mu},\mathbf{\omega },|A|^{2},|B|^{2}),
\end{eqnarray*}%
where functions $\phi _{0}$ and $\phi _{1}$ are real valued. Moreover, for $%
\mathbf{\omega }=0$ we have%
\begin{equation*}
\phi _{1}(\tilde{\mu},\mathbf{0},|A|^{2},|B|^{2})=\phi _{0}(\tilde{\mu},%
\mathbf{0},|B|^{2},|A|^{2}).
\end{equation*}%
It results immediately that we have the following (formal)
solutions of our system (in addition to the trivial solution 0):

\begin{proof}
i) $B=0,$ $|A|$ satisfying $\phi _{0}(\tilde{\mu},\mathbf{\omega }%
,|A|^{2},0)=0,$ which is not else that the 2-dimensional
travelling wave with basic wave vector $K_{1},$ and where, with no
loss of generality, we can choose the velocity $\mathbf{c}$ in the
direction of $K_{1}.$

ii) $A=0,$ \ $|B|$ satisfying $\phi _{1}(\tilde{\mu},\mathbf{\omega }%
,0,|B|^{2})=0,$ which is not else that the 2-dimensional
travelling wave with basic wave vector $K_{2},$ and where, with no
loss of generality, we can choose the velocity $\mathbf{c}$ in the
direction of $K_{2}.$

iii) $|A|$ and $|B|$ such that
\begin{eqnarray}
\phi _{0}(\tilde{\mu},\mathbf{\omega },|A|^{2},|B|^{2}) &=&0,
\label{bifurcsyst} \\
\phi _{1}(\tilde{\mu},\mathbf{\omega },|A|^{2},|B|^{2}) &=&0,
\notag
\end{eqnarray}%
which gives a family of 3-dimensional travelling waves. Moreover
we notice that if $\mathbf{\omega }=0$ there is a family of
solutions where $|A|=|B|$
and%
\begin{equation*}
\phi _{0}(\tilde{\mu},\mathbf{0},|A|^{2},|A|^{2})=0,
\end{equation*}%
representing the "diamond waves" of \cite{Reed-Shin}. The leading terms of $%
\phi _{0}$ and $\phi _{1}$ are computed by Bridges, Dias, Menasce in \cite%
{B-D-M}, even in cases with a finite depth and with surface
tension. Since our case has less parameters, our computations may
look simpler. The leading
terms independent of $|A|$ and $|B|$ in $\phi _{0}$ come from%
\begin{equation*}
\langle \tilde{\mu}\mathcal{L}_{1}\zeta _{0}+\mathcal{L}_{2}(\mathbf{\omega }%
,\zeta _{0}),\zeta _{0}\rangle =\frac{4\pi ^{2}}{\tau }\frac{1}{\mu _{c}}%
\left( \frac{\tilde{\mu}}{\mu _{c}}-2\omega \cdot K_{1}\right) ,
\end{equation*}%
and we notice that (in using $\frac{\delta \mu }{\mu }=-2\frac{\delta c}{c})$%
\begin{equation*}
\frac{\tilde{\mu}}{\mu _{c}}-2\omega \cdot K_{1}\sim -2\frac{(\mathbf{c}-%
\mathbf{c}_{0})\cdot K_{1}}{c_{0}}.
\end{equation*}%
We then have%
\begin{eqnarray*}
\phi _{0}(\tilde{\mu},\mathbf{\omega },|A|^{2},|B|^{2}) &=&\frac{1}{\mu _{c}}%
\left( \frac{\tilde{\mu}}{\mu _{c}}-2\omega \cdot K_{1}\right)
+\alpha
_{0}|A|^{2}+\beta _{0}|B|^{2}+O\{(|\tilde{\mu}|+|\mathbf{\omega }%
|+|A|^{2}+|B|^{2})^{2}\} \\
\phi _{1}(\tilde{\mu},\mathbf{\omega },|A|^{2},|B|^{2}) &=&\frac{1}{\mu _{c}}%
\left( \frac{\tilde{\mu}}{\mu _{c}}-2\omega \cdot K_{2}\right)
+\beta
_{0}|A|^{2}+\alpha _{0}|B|^{2}+O\{(|\tilde{\mu}|+|\mathbf{\omega }%
|+|A|^{2}+|B|^{2})^{2}\},
\end{eqnarray*}%
with%
\begin{eqnarray*}
\alpha _{0} &=&\frac{\tau }{4\pi ^{2}}\langle 2\mathcal{N}_{2}(\overline{%
\zeta }_{0},U_{2000})+3\mathcal{N}_{3}(\zeta _{0},\zeta
_{0},\overline{\zeta
}_{0}),\zeta _{0}\rangle , \\
\beta _{0} &=&\frac{\tau }{4\pi ^{2}}\langle
2\mathcal{N}_{2}(\zeta
_{1},U_{1001})+2\mathcal{N}_{2}(\overline{\zeta }_{1},U_{1010})+6\mathcal{N}%
_{3}(\zeta _{0},\zeta _{1},\overline{\zeta }_{1}),\zeta
_{0}\rangle .
\end{eqnarray*}%
We can formally solve the system of equations
\begin{eqnarray*}
\phi _{0}(\tilde{\mu},\mathbf{\omega },|A|^{2},|B|^{2}) &=&0, \\
\phi _{1}(\tilde{\mu},\mathbf{\omega },|A|^{2},|B|^{2}) &=&0,
\end{eqnarray*}%
with respect to $\tilde{\mu},\omega _{2}=\frac{1}{2\tau }\mathbf{\omega }%
\cdot (K_{1}-K_{2}).$ Indeed, we obtain respectively in adding and
subtracting the two equations, a system easy to solve, in taking
into
account of $\omega _{1}=O(\omega _{2}^{2}),$%
\begin{eqnarray*}
\tilde{\mu} &=&-\frac{\mu _{c}^{2}}{2}(\alpha _{0}+\beta
_{0})(|A|^{2}+|B|^{2})+O\{(|A|^{2}+|B|^{2})^{2}\}, \\
\omega _{2} &=&(|A|^{2}-|B|^{2})\left( \frac{\mu _{c}}{4\tau
}(\alpha _{0}-\beta _{0})+O\{|A|^{2}+|B|^{2}\}\right) .
\end{eqnarray*}%
For the computation of coefficients $\alpha _{0}$ and $\beta _{0}$
we use again (\ref{G1}) and (\ref{bilin}) and obtain
\begin{eqnarray*}
2\mathcal{N}_{2}(\overline{\zeta }_{0},U_{2000}) &=&\left(
\frac{2}{\mu
_{c}^{5}},0\right) e^{iK_{1}\cdot X}, \\
2\mathcal{N}_{2}(\zeta _{1},U_{1001}) &=&\left( -\frac{4\tau ^{4}}{\mu _{c}}%
,0\right) e^{iK_{1}\cdot X}, \\
2\mathcal{N}_{2}(\overline{\zeta }_{1},U_{1010}) &=&\left(
\frac{4(\mu _{c}^{3}+4\mu _{c}^{2}-5\mu _{c}+1)}{\mu
_{c}^{3}(2-\mu _{c})},\frac{8i(\mu _{c}-1)(1-\mu _{c}^{2})}{\mu
_{c}^{3}(2-\mu _{c})}\right) e^{iK_{1}\cdot X}.
\end{eqnarray*}%
Now, with (\ref{G2}) and the form of $\mathcal{N}_{3}$ we have%
\begin{eqnarray*}
3\mathcal{N}_{3}(\zeta _{0},\zeta _{0},\overline{\zeta }_{0})
&=&\left(
\frac{1}{\mu _{c}^{5}},\frac{-i}{\mu _{c}^{4}}\right) e^{iK_{1}\cdot X}, \\
6\mathcal{N}_{3}(\zeta _{0},\zeta _{1},\overline{\zeta }_{1})
&=&\left(
\frac{2}{\mu _{c}^{4}}(2-\frac{1}{\mu _{c}}),\frac{-2i}{\mu _{c}^{2}}(2-%
\frac{1}{\mu _{c}^{2}})\right) e^{iK_{1}\cdot X}.
\end{eqnarray*}%
Finally we obtain%
\begin{equation*}
\alpha _{0}=\frac{4}{\mu _{c}^{5}}
\end{equation*}%
\begin{equation*}
\beta _{0}=-\frac{8}{\mu _{c}^{5}}+\frac{8}{\mu
_{c}^{4}}+\frac{12}{\mu _{c}^{3}}-\frac{32}{\mu
_{c}^{2}}-\frac{8}{\mu _{c}}+\frac{36}{\mu _{c}^{2}(2-\mu _{c})}.
\end{equation*}%
If $\alpha _{0}+\beta _{0}$ and $\alpha _{0}-\beta _{0}$ are both
different from 0, we can solve the system (\ref{bifurcsyst}) with
respect to $|A|^{2}$
and $|B|^{2}$ and obtain a formal expansion in powers of $(\tilde{\mu},%
\mathbf{\omega })$ of the form%
\begin{eqnarray*}
|A|^{2} &=&\frac{\tilde{\mu}}{-\mu _{c}^{2}(\alpha _{0}+\beta _{0})}+\frac{%
\mathbf{\omega }\cdot (K_{1}-K_{2})}{\mu _{c}(\alpha _{0}-\beta _{0})}+O(|(%
\tilde{\mu},\mathbf{\omega })|^{2}), \\
|B|^{2} &=&\frac{\tilde{\mu}}{-\mu _{c}^{2}(\alpha _{0}+\beta _{0})}-\frac{%
\mathbf{\omega }\cdot (K_{1}-K_{2})}{\mu _{c}(\alpha _{0}-\beta _{0})}+O(|(%
\tilde{\mu},\mathbf{\omega })|^{2}).
\end{eqnarray*}%
The indeterminacy on the phases of $A$ and $B$ means that we
obtain in fact, for each fixed $(\tilde{\mu},\mathbf{\omega })$
leading to positive expressions for $|A|^{2}$ and $|B|^{2}$, a
torus of solutions, which is generated by acting the operator
$T_{\mathbf{v}}$ on a particular solution, for instance with $A$
and $B$ pure imaginary. This two-parameter family of tori of
3-dimensional waves connects with the 2-dimensional travelling
waves respectively of wave vectors $K_{1}$ and $K_{2}$ (choosing
$\mathbf{\omega }$ orthogonal to $K_{2}$ or to $K_{1}.$ If we
choose $\mathbf{\omega }=0$ which means that we choose the
direction of the waves as $x_{1}$ axis, then we obtain the
"diamond waves" as in \cite{Reed-Shin}, and \cite{B-D-M}, here
without surface tension.

Let us study the sign of $(\alpha _{0}+\beta _{0})$ and $(\alpha
_{0}-\beta
_{0}).$ We have%
\begin{eqnarray*}
\alpha _{0}+\beta _{0} &=&\frac{4}{\mu _{c}^{2}}\left( -\frac{1}{\mu _{c}^{3}%
}+\frac{2}{\mu _{c}^{2}}+\frac{3}{\mu _{c}}-8-2\mu _{c}+\frac{9}{2-\mu _{c}}%
\right) \\
\beta _{0}-\alpha _{0} &=&\frac{4}{\mu _{c}^{2}}\left( -\frac{3}{\mu _{c}^{3}%
}+\frac{2}{\mu _{c}^{2}}+\frac{3}{\mu _{c}}-8-2\mu _{c}+\frac{9}{2-\mu _{c}}%
\right) ,
\end{eqnarray*}%
where we notice that%
\begin{equation*}
\alpha _{0}+\beta _{0}=\frac{4}{\mu _{c}}\left\{ (\frac{1}{\mu _{c}^{2}}-%
\frac{2}{\mu _{c}}-1)(2-\frac{1}{\mu _{c}^{2}})+\frac{1}{2\mu _{c}}+\frac{9}{%
2(2-\mu _{c})}\right\}
\end{equation*}%
and it is easy to show (study of the factor as a function of
$1/\mu _{c})$ that $\alpha _{0}+\beta _{0}>0$ for $\tau \in
(0,\tau _{c}),$ and $\alpha
_{0}+\beta _{0}<0$ for $\tau >\tau _{c}$ where $\tau _{c}\in (\sqrt{6},\sqrt{%
7}),$ i.e. more precisely $\mu _{c,c}\sim 0.374$. We also notice
that the function of $\mu _{c}$ in the factor for $\beta
_{0}-\alpha _{0}$ is strictly increasing for $\mu _{c}\in (0,1),$
and cancels for $\mu _{c}\sim 0.893.$ Now, defining new parameters
$\varepsilon _{1}$ and $\varepsilon
_{2} $ by%
\begin{equation*}
A=-\frac{i\varepsilon _{1}}{2},B=-\frac{i\varepsilon _{2}}{2},
\end{equation*}%
for a particular choice of a solution belonging to the torus, we
get the results of Theorem \ref{Lembifurc}.
\end{proof}

\section{Proof of Lemma \protect\ref{changevar1}\label{a2}}

Let first consider the unique solution of the Cauchy problem (notice that $%
\frac{V_{2}}{V_{1}}\in H_{o,o}^{m-1}(%
\mathbb{R}
^{2}/\Gamma )\subset C_{o,o}^{m-3}(%
\mathbb{R}
^{2}/\Gamma ),$ hence the vector field is Lipschitz for $m\geq 4$)%
\begin{equation*}
\frac{dX}{dx_{1}}=\frac{V_{2}}{V_{1}}(x_{1},X),\text{ \ \
}X|_{x_{1}=0}=y,
\end{equation*}%
which we denote by $X(x_{1},y)\in C^{m-3}.$ We can successively show that $%
X(x_{1},y)$ is even in $x_{1},$ odd in $y,$ and such that%
\begin{equation*}
X(x_{1},y)=X(x_{1}+2\pi ,y)=X(x_{1},y+\frac{2\pi }{\tau
})-\frac{2\pi }{\tau },
\end{equation*}%
in using the uniqueness of the solution of the Cauchy problem, in
looking
for the system satisfied by $X(-x_{1},y),$ $-X(x_{1},-y),$ $X(x_{1},y+\frac{%
2\pi }{\tau })-\frac{2\pi }{\tau },$ which is the same as
$X(x_{1},y),$ and
in comparing the systems satisfied by $X(x_{1}+\pi ,y)$ and $X(x_{1}-\pi ,y)$%
, using the evenness in $x_{1}.$ Notice that%
\begin{equation*}
\frac{\partial X}{\partial y}(x_{1},y)=\exp \left(
\int_{0}^{x_{1}}\partial
_{x_{2}}(\frac{V_{2}}{V_{1}})(t,y)dt\right)
\end{equation*}%
equals 1 for $\frac{V_{2}}{V_{1}}=0,$ i.e.if $U=0.$ Hence, for
$||U||_{4}$
small enough, we can solve (implicit function theorem) with respect to $%
y=y(z_{2})$ the equation%
\begin{equation*}
\Pi _{1}X(\cdot ,y)=z_{2}
\end{equation*}%
where $\Pi _{1}$ represents the average over a period in $x_{1}.$
Then, we set
\begin{equation*}
\mathcal{Z}(Z)=X(z_{1},y(z_{2}))
\end{equation*}%
and, thanks to $\frac{V_{2}}{V_{1}}(x_{1}+\pi ,x_{2}+\frac{\pi }{\tau })=%
\frac{V_{2}}{V_{1}}(x_{1},x_{2})$ we observe that the function
\begin{equation*}
\mathcal{Z}(z_{1}+\pi ,z_{2}+\frac{\pi }{\tau })-\frac{\pi }{\tau
}
\end{equation*}%
satisfies the same differential equation as
$\mathcal{Z}(z_{1},z_{2})$ with
the same average in $z_{1},$ hence by uniqueness it is identical to $%
\mathcal{Z}(z_{1},z_{2}).$ It is then clear that $(z_{1},z_{2})\mapsto z_{2}-%
\mathcal{Z}(Z)=\widetilde{d_{1}}(Z)$ satisfies the properties
indicated at
Lemma \ref{changevar1} hence lies in $C_{e,o}^{m-3}(%
\mathbb{R}
^{2}/\Gamma )$. The proof of the tame estimates for $d_{1}$ and $\widetilde{%
d_{1}}$ is identical to the one made in Appendix G of \cite{IPT}.
Notice in addition, that replacing the initial condition at
$x_{1}=0$ by an average condition, allows to keep the equivariance
under shifts $\mathcal{T}_{\delta }$ parallel to $x_{1}$
direction.

\section{Proofs of Lemmas \protect\ref{Lemcoef-nu} and \protect\ref%
{Lemmcoordchange}\label{a3}}

Let us start with the expressions for $U,$ given in (\ref{Diamond})%
\begin{eqnarray*}
U &=&\varepsilon \xi _{0}+\varepsilon ^{2}U^{(2)}+O(\varepsilon ^{3}), \\
\mu &=&\mu _{c}+\mu _{1}(\tau )\varepsilon ^{2}+O(\varepsilon
^{4}),
\end{eqnarray*}%
with $\xi _{0}=(\psi _{1},\eta _{1})$ and $\mu _{1}(\tau )$ given in (\ref%
{Diamond}) and (\ref{mu_1}), and where
\begin{eqnarray*}
U^{(2)} &=&\left\{
\begin{array}{c}
\frac{1-2\mu _{c}}{4\mu _{c}(2-\mu _{c})}\sin 2x_{1}-\frac{1}{4\mu _{c}^{2}}%
\sin 2x_{1}\cos 2\tau x_{2} \\
\frac{\mu _{c}^{2}+2\mu _{c}-2}{4\mu _{c}^{2}(2-\mu _{c})}\cos 2x_{1}+\frac{%
\tau ^{2}}{4\mu _{c}}\cos 2\tau x_{2}+\frac{1}{4\mu _{c}^{3}}\cos
2x_{1}\cos
2\tau x_{2}%
\end{array}%
\right. \\
&=&(\psi _{2},\eta _{2}).
\end{eqnarray*}%
We successively find from (\ref{Z}) and (\ref{a})%
\begin{eqnarray*}
\mathfrak{b} &=&\frac{\varepsilon }{\mu _{c}}\sin x_{1}\cos \tau
x_{2}+\varepsilon ^{2}\mathfrak{b}^{(2)}+O(\varepsilon ^{3}), \\
\mathfrak{b}^{(2)} &=&\mathcal{\partial }_{x_{1}}\eta _{2}+\nabla
\eta
_{1}\cdot \nabla \psi _{1} \\
&=&\left( \frac{1}{4\mu _{c}}+\frac{1}{2\mu _{c}^{2}}-\frac{1}{4\mu _{c}^{3}}%
-\frac{3}{4(2-\mu _{c})}\right) \sin 2x_{1}+ \\
&&-\frac{1}{4\mu _{c}^{3}}\sin 2x_{1}\cos 2\tau x_{2},
\end{eqnarray*}%
\begin{eqnarray*}
V_{1} &=&1+\varepsilon \cos x_{1}\cos \tau x_{2}+\varepsilon
^{2}V_{1}^{(2)}+O(\varepsilon ^{3}), \\
V_{2} &=&-\varepsilon \tau \sin x_{1}\sin \tau x_{2}+\varepsilon
^{2}V_{2}^{(2)}+O(\varepsilon ^{3}),
\end{eqnarray*}%
\begin{eqnarray*}
V_{1}^{(2)} &=&\partial _{x_{1}}\psi
_{2}-\mathfrak{b}^{(1)}\partial
_{x_{1}}\eta _{1}, \\
V_{1}^{(2)} &=&-\frac{1}{4\mu _{c}^{2}}(1+\cos 2\tau x_{2}+\cos
2x_{1}\cos
2\tau x_{2})+ \\
&&+\left( \frac{1}{4\mu _{c}}+\frac{1}{4\mu _{c}^{2}}-\frac{3}{4(2-\mu _{c})}%
\right) \cos 2x_{1},
\end{eqnarray*}%
\begin{eqnarray*}
V_{2}^{(2)} &=&\partial _{x_{2}}\psi
_{2}-\mathfrak{b}^{(1)}\partial
_{x_{2}}\eta _{1} \\
&=&\frac{\tau }{4\mu _{c}^{2}}\sin 2x_{1}\sin 2\tau x_{2},
\end{eqnarray*}%
\begin{eqnarray*}
\mathfrak{a} &\mathfrak{=}&\mu _{c}+\frac{\varepsilon }{\mu
_{c}}\cos x_{1}\cos \tau x_{2}+\varepsilon
^{2}\mathfrak{a}^{(2)}+O(\varepsilon ^{3}),
\\
\mathfrak{a}^{(2)} &=&\mu _{1}+V^{(1)}\cdot \nabla \mathfrak{b}%
^{(1)}+\partial _{x_{1}}\mathfrak{b}^{(2)},
\end{eqnarray*}%
\begin{eqnarray*}
\mathfrak{a}^{(2)} &=&\mu _{1}+\left( \frac{1}{\mu
_{c}}+\frac{1}{\mu _{c}^{2}}-\frac{3}{4\mu
_{c}^{3}}-\frac{3}{2(2-\mu _{c})}\right) \cos 2x_{1}+
\\
&&+\frac{1}{4\mu _{c}^{3}}(1-\cos 2x_{1}\cos 2\tau
x_{2})+\frac{2\mu _{c}^{2}-1}{4\mu _{c}^{3}}\cos 2\tau x_{2},
\end{eqnarray*}%
\begin{equation*}
\frac{\{V^{2}+(V\cdot \nabla \eta
)^{2}\}^{1/2}}{V_{1}^{3}}=1-2\varepsilon \cos x_{1}\cos \tau
x_{2}+\varepsilon ^{2}Y^{(2)}+O(\varepsilon ^{3})
\end{equation*}%
\begin{eqnarray*}
Y^{(2)} &=&\frac{3}{4\mu _{c}^{2}}+\frac{5}{8}+\left( \frac{7}{8}-\frac{1}{%
2\mu _{c}}-\frac{3}{4\mu _{c}^{2}}+\frac{3}{2(2-\mu _{c})}\right)
\cos
2x_{1}+ \\
&&+(\frac{7}{8}+\frac{1}{2\mu _{c}^{2}})\cos 2\tau x_{2}+(\frac{5}{8}+\frac{1%
}{2\mu _{c}^{2}})\cos 2x_{1}\cos 2\tau x_{2}.
\end{eqnarray*}%
For the determination of the function $\mathcal{Z}(Z)=z_{2}-\widetilde{d_{1}}%
(Z)$ of Lemma \ref{changevar1}, we have%
\begin{equation*}
\partial _{x_{1}^{\prime }}\widetilde{d_{1}}=-\frac{V_{2}}{V_{1}}%
(x_{1}^{\prime },x_{2}^{\prime }-\widetilde{d_{1}}),\text{ \ \ \ \
\ }\Pi _{1}\widetilde{d_{1}}=0,
\end{equation*}%
i.e.%
\begin{equation*}
\widetilde{d_{1}}=\varepsilon \widetilde{d_{1}}^{(1)}+\varepsilon ^{2}%
\widetilde{d_{1}}^{(2)}+O(\varepsilon ^{3})
\end{equation*}%
with%
\begin{equation*}
\partial _{z_{1}}\widetilde{d_{1}}^{(1)}=-V_{2}^{(1)},\text{ \ \ \ \ \ }
\Pi_{1}\widetilde{d_{1}}^{(1)}=0,
\end{equation*}

\begin{equation*}
\partial _{z_{1}}\widetilde{d_{1}}^{(2)}=V_{1}^{(1)}V_{2}^{(1)}-V_{2}^{(2)}+
\partial _{x_{2}}V_{2}^{(1)}
\widetilde{d_{1}}^{(1)},\text{ \ \ \ \ }\Pi _{1}
\widetilde{d_{1}}^{(2)}=0.
\end{equation*}

This leads to%
\begin{eqnarray*}
\widetilde{d_{1}}^{(1)} &=&-\tau \cos z_{1}\sin \tau z_{2}, \\
\widetilde{d_{1}}^{(2)} &=&\frac{\tau }{4}\cos 2z_{1}\sin 2\tau
z_{2},
\end{eqnarray*}%
hence the diffeomorphism $\mathcal{U}_{1}^{-1}$ takes the form%
\begin{equation*}
x_{1}=z_{1},\text{ \ }x_{2}=z_{2}+\varepsilon \tau \cos z_{1}\sin \tau z_{2}-%
\frac{\varepsilon ^{2}\tau }{4}\cos 2z_{1}\sin 2\tau
z_{2}+O(\varepsilon ^{3}).
\end{equation*}%
Now we have%
\begin{equation*}
d_{1}=\varepsilon d_{1}^{(1)}+\varepsilon
^{2}d_{1}^{(2)}+O(\varepsilon ^{3})
\end{equation*}%
with%
\begin{eqnarray*}
d_{1}^{(1)} &=&-\tau \cos x_{1}\sin \tau x_{2}, \\
d_{1}^{(2)} &=&\frac{\tau }{4}\left\{ \tau ^{2}+(1+\tau ^{2})\cos
2x_{1}\right\} \sin 2\tau x_{2},
\end{eqnarray*}%
and finally (\ref{q}) gives%
\begin{equation*}
\mathfrak{q}=\mu _{c}(1-\varepsilon \cos x_{1}\cos \tau
x_{2}+\varepsilon ^{2}\mathfrak{q}^{(2)})+O(\varepsilon ^{3}),
\end{equation*}%
with%
\begin{eqnarray*}
\mathfrak{q}^{(2)} &=&\frac{\mathfrak{a}^{(2)}}{\mu
_{c}}+Y^{(2)}+\partial _{x_{2}}d_{1}^{(2)}-(2+\tau ^{2}+\tau
^{4})\cos ^{2}x_{1}\cos ^{2}\tau x_{2}
\\
&=&\frac{\mu _{1}}{\mu _{c}}+\frac{1}{8}+\frac{1}{\mu _{c}^{2}}+\mathfrak{q}%
_{20}^{(2)}\cos 2x_{1}+\mathfrak{q}_{02}^{(2)}\cos 2\tau x_{2}+ \\
&&+\mathfrak{q}_{22}^{(2)}\cos 2x_{1}\cos 2\tau x_{2}, \\
\mathfrak{q}_{20}^{(2)} &=&\frac{3}{4(2-\mu
_{c})}+\frac{3}{8}-\frac{5}{4\mu _{c}}+\frac{1}{2\mu
_{c}^{2}}+\frac{1}{\mu _{c}^{3}}-\frac{1}{\mu _{c}^{4}},
\\
\mathfrak{q}_{02}^{(2)} &=&\frac{7}{8}+\frac{1}{4\mu _{c}^{2}},\text{ \ \ }%
\mathfrak{q}_{22}^{(2)}=\frac{1}{8}+\frac{1}{4\mu _{c}^{2}}.
\end{eqnarray*}%
Now, we obtain%
\begin{equation*}
(\mathfrak{q}\circ \mathcal{U}_{1}^{-1})(Z)=\mu _{0}(1-\varepsilon
\cos
z_{1}\cos \tau z_{2}+\varepsilon ^{2}\widetilde{\mathfrak{q}}%
^{(2)})+O(\varepsilon ^{3})
\end{equation*}%
with%
\begin{eqnarray*}
\widetilde{\mathfrak{q}}^{(2)} &=&\frac{\mu _{1}}{\mu _{c}}-\frac{1}{8}+%
\frac{5}{4\mu _{c}^{2}}+(\mathfrak{q}_{20}^{(2)}+\frac{\tau
^{2}}{4})\cos
2z_{1}+\frac{9}{8}\cos 2\tau z_{2}+ \\
&&+\frac{3}{8}\cos 2z_{1}\cos 2\tau z_{2},
\end{eqnarray*}%
and%
\begin{eqnarray*}
\{(\mathfrak{q}\circ \mathcal{U}_{1}^{-1})(Z)\}^{1/2} &=&\mu _{c}^{1/2}\{1-%
\frac{\varepsilon }{2}\cos z_{1}\cos \tau z_{2}+\varepsilon
^{2}(\frac{\mu
_{1}}{2\mu _{c}}-\frac{3}{32}+\frac{5}{8\mu _{c}^{2}})+ \\
&&+\varepsilon ^{2}(\frac{\mathfrak{q}_{20}^{(2)}}{2}+\frac{\tau ^{2}}{8}-%
\frac{1}{32})\cos 2z_{1}+\varepsilon ^{2}\frac{17}{32}\cos 2\tau z_{2}+ \\
&&+\varepsilon ^{2}\frac{5}{32}\cos 2z_{1}\cos 2\tau
z_{2}\}+O(\varepsilon ^{3}),
\end{eqnarray*}%
hence%
\begin{eqnarray*}
\frac{1}{2\pi }\int_{-\pi }^{\pi }\{(\mathfrak{q}\circ \mathcal{U}%
_{1}^{-1})(Z)\}^{1/2}dz_{1} &=&\mu _{c}^{1/2}\{1+\varepsilon
^{2}(\frac{\mu
_{1}}{2\mu _{c}}-\frac{3}{32}+\frac{5}{8\mu _{c}^{2}})+ \\
&&+\varepsilon ^{2}\frac{17}{32}\cos 2\tau z_{2}\}+O(\varepsilon
^{3}).
\end{eqnarray*}%
Then with (\ref{nu}), we obtain%
\begin{equation*}
\nu =\mu _{c}^{-1}\{1-\varepsilon ^{2}(\frac{\mu _{1}}{\mu _{c}}-\frac{3}{16}%
+\frac{5}{4\mu _{c}^{2}})+O(\varepsilon ^{3})\},
\end{equation*}%
as indicated in Lemma \ref{Lemcoef-nu}. Furthermore, in using (\ref{e_2}) , (%
\ref{d_2}), we notice that%
\begin{eqnarray*}
e_{2} &=&-\frac{17\varepsilon ^{2}}{32\tau }\sin 2\tau
z_{2}+O(\varepsilon
^{3}), \\
d_{2} &=&-\frac{\varepsilon }{2}\sin z_{1}\cos \tau
z_{2}+\varepsilon ^{2}(d_{21}\sin 2z_{1}+d_{22}\sin 2z_{1}\cos
2\tau z_{2})+O(\varepsilon
^{3}), \\
d_{21} &=&\frac{1}{4}\left\{ \frac{3}{4(2-\mu _{c})}+\frac{1}{16}-\frac{5}{%
4\mu _{c}}+\frac{3}{4\mu _{c}^{2}}+\frac{1}{\mu
_{c}^{3}}-\frac{1}{\mu _{c}^{4}}\right\} ,\text{ \ \
}d_{22}=\frac{5}{64}.
\end{eqnarray*}%
This, with (\ref{U_2}), leads to the principal part of the
bi-periodic
functions occurring in the diffeomorphism of the torus of Theorem \ref%
{thmChangeCoord}:
\begin{eqnarray*}
d(x,y) &=&-\frac{\varepsilon }{2}\sin x_{1}\cos \tau
x_{2}+\varepsilon ^{2}\sin 2x_{1}\{d_{21}-\frac{\tau
^{2}}{8}+(d_{22}+\frac{\tau ^{2}}{8})\cos
2\tau x_{2}\}+O(\varepsilon ^{3}) \\
e(x,y) &=&-\varepsilon \tau \cos x_{1}\sin \tau x_{2}+\frac{\varepsilon ^{2}%
}{4\tau }\sin 2\tau x_{2}\left\{ \frac{1}{\mu _{c}^{4}}-\frac{2}{\mu _{c}^{2}%
}-\frac{9}{8}+\frac{\tau ^{2}}{\mu _{c}^{2}}\cos 2x_{1}\right\}
+O(\varepsilon ^{3}).
\end{eqnarray*}%
Inverting $\mathcal{U}_{2}\circ \mathcal{U}_{1}$, and changing the
coordinates in%
\begin{eqnarray*}
\eta (X) &=&-\frac{\varepsilon }{\mu _{c}}\cos x_{1}\cos \tau
x_{2}+\varepsilon ^{2}\left\{ \frac{\mu _{c}^{2}+2\mu _{c}-2}{4\mu
_{c}^{2}(2-\mu _{c})}\cos 2x_{1}+\frac{\tau ^{2}}{4\mu _{c}}\cos
2\tau
x_{2}\right. \\
&&+\left. \frac{1}{4\mu _{c}^{3}}\cos 2x_{1}\cos 2\tau
x_{2}\right\} +O(\varepsilon ^{3})
\end{eqnarray*}%
leads to the results of Lemma \ref{Lemmcoordchange}.

\section{Distribution of numbers $\{\protect\omega
_{0}n^{2}\}$}\label{numbers}

Recall that for each $x\in \mathbb{R}$,
\begin{equation*}
\lbrack x]=\max \{N:N\in \mathbb{N},\quad N\leq x\},\quad
\{x\}=x-[x]\in \lbrack 0,1).
\end{equation*}%
In this section we consider the distribution of the numbers
\begin{equation}
\theta _{n}=\{\omega _{0}n^{2}-C\},\quad \sigma _{n}=\{\omega
_{0}n\},\quad n\geq 1  \label{71}
\end{equation}%
Due the famous Weil Theorem \cite{Weil}, for each polynomial
$f(x)=\alpha
_{k}x^{k}+..+\alpha _{0}$ with irrational $\alpha _{k}$, the numbers $%
\{f(n)\}$ are uniformly distributed in $[0,1]$. This means that for each $%
[\alpha ,\beta ]\subset \lbrack 0,1]$,
\begin{equation*}
\text{card~}\{n:0\leq n\leq N,\;\{f(n)\}\in \lbrack \alpha ,\beta
]\}=N(\beta -\alpha )+o(N).
\end{equation*}%
This result is the asymptotic relation in which the remainder
strongly depends on the choice of the interval. We consider the
simplest case
\begin{equation*}
f(x)=\omega _{0}x^{2}-C\text{~~with~~}\{f(n)\}=\theta _{n},\quad
\lbrack \alpha ,\beta ]=[0,\rho ]
\end{equation*}%
and deduce the rough, uniform in $\rho $ estimate which is
sufficient for our needs. The main result is the following .

\begin{proposition}
\label{distribution} Suppose that
\begin{equation*}
\frac{1}{|e^{2\pi i\omega _{0}l}-1|}\leq c_{1}l^{2}\text{~~for all
positive integers ~~}l.
\end{equation*}%
Then there is a constant $c$ depending only on $c_{1}$ such that
\begin{equation}
\frac{1}{N}\sum\limits_{1\leq n\leq N,\theta _{n}\in \lbrack
0,\rho ]}1\leq c\rho \text{~~for all~~}\rho \in
(0,1/4)\text{~~and~~}N>\rho ^{-78}. \label{72}
\end{equation}
\end{proposition}

The proof is based on the following lemma.

\begin{lemma}
\label{LC} Under the assumptions of Proposition
\ref{distribution}, there exists positive $c$ depending on $c_{1}$
only such that
\begin{equation}
\frac{1}{N}\sum\limits_{1\leq n\leq N,\sigma _{n}\in \lbrack
0,\varepsilon
]\cup \lbrack 1-\varepsilon ,1]}1\leq c\varepsilon \text{~~for all~~}%
\varepsilon \in (0,1/4)\text{~~and~~}N>\varepsilon ^{-9}.
\label{73}
\end{equation}
\end{lemma}

\begin{proof}
Fix an arbitrary positive $\varepsilon \in (0,1/4)$ and introduce
the function depending on parameter $\varepsilon $ and given by
the equalities
\begin{gather*}
\varphi _{\varepsilon }(x)=1\text{~~for~~}x\in \lbrack
0,\varepsilon ]\cup
\lbrack 1-\varepsilon ,1], \\
\varphi _{\varepsilon }(x)=0\text{~~for~~}x\in \lbrack
2\varepsilon
,1-2\varepsilon ], \\
\varphi _{\varepsilon }(x)=\frac{2\varepsilon -x}{\varepsilon }\text{~~for~~}%
x\in \lbrack \varepsilon ,2\varepsilon ], \\
\varphi _{\varepsilon }(x)=\frac{x-1+2\varepsilon }{\varepsilon }\text{%
~~for~~}x\in \lbrack 1-2\varepsilon ,1-\varepsilon ].
\end{gather*}%
We will assume that $\varphi _{\varepsilon }$ is extended 1-periodically onto $%
\mathbb{R}$. Obviously the extended function is absolutely
continuous and
\begin{equation}
\int\limits_{0}^{1}\varphi _{\varepsilon }(x)dx=3\varepsilon
,\quad \int\limits_{0}^{1}|\varphi _{\varepsilon }^{\prime
}(x)|^{2}dx=2\varepsilon ^{-1}.  \label{74}
\end{equation}%
It has the representation
\begin{equation*}
\varphi _{\varepsilon }(x)=\sum\limits_{l=-\infty }^{\infty
}\varphi _{\varepsilon ,l}e^{2\pi ilx},\quad \varphi _{\varepsilon
,l}=\int\limits_{0}^{1}e^{-2\pi ilx}\varphi _{\varepsilon }(x)dx.
\end{equation*}%
It is clear that
\begin{equation}
|\varphi _{\varepsilon ,l}|\leq \int\limits_{0}^{1}\varphi
_{\varepsilon }(x)dx=3\varepsilon .  \label{75}
\end{equation}%
Obviously
\begin{equation}
\frac{1}{N}\sum\limits_{1\leq n\leq N,\sigma _{n}\in \lbrack
0,\varepsilon ]\cup \lbrack 1-\varepsilon ,1]}1\leq
\frac{1}{N}\sum\limits_{1\leq n\leq N}\varphi _{\varepsilon
}(\sigma _{n}).  \label{76}
\end{equation}%
Represent $\varphi $ in the form
\begin{equation}
\varphi _{\varepsilon }(x)=\varphi _{\varepsilon
,0}+\sum\limits_{1\leq |l|\leq k}\varphi _{\varepsilon ,l}e^{2\pi
ilx}+Q_{k}(x),\quad Q_{k}(x)=\sum\limits_{k+1\leq |l|}\varphi
_{\varepsilon ,l}e^{2\pi ilx} \label{77}
\end{equation}%
We have, by equality \eqref{74},
\begin{gather}
|Q_{k}(x)|\leq 2\sum\limits_{k+1\leq l}|\varphi _{\varepsilon ,l}|\leq 2%
\sqrt{\sum\limits_{k+1\leq
l}l^{-2}}\sqrt{\sum\limits_{l}l^{2}|\varphi
_{\varepsilon ,l}|^{2}}=  \notag \\
\frac{1}{\pi }\sqrt{\sum\limits_{k+1\leq l}l^{-2}}\sqrt{\int\limits_{0}^{1}|%
\varphi _{\varepsilon }^{\prime }(x)|^{2}dx}\leq c\frac{1}{\sqrt{%
k\varepsilon }}.  \label{78}
\end{gather}%
Here $c$ is some absolute constant. Combining
\eqref{76}-\eqref{78} we obtain
\begin{equation}
\frac{1}{N}\sum\limits_{1\leq n\leq N,\sigma _{n}\in \lbrack
0,\varepsilon ]\cup \lbrack 1-\varepsilon ,1]}1\leq 3\varepsilon
+6\varepsilon \sum\limits_{1\leq l\leq k}\frac{1}{N}\left\vert
\sum\limits_{1\leq n\leq N}e^{2\pi il\sigma _{n}}\right\vert
+\frac{c}{\sqrt{k\varepsilon }}. \label{79}
\end{equation}%
Since the numbers
\begin{equation*}
e^{2\pi il\sigma _{n}}=e^{2\pi i\omega _{0}ln}
\end{equation*}%
form a geometric progression, we have
\begin{equation*}
\left\vert \sum\limits_{1\leq n\leq N}e^{2\pi il\sigma
_{n}}\right\vert \leq \frac{2}{|e^{2\pi i\omega _{0}l}-1|}\leq
2c_{1}l^{2}.
\end{equation*}%
Substituting this inequality into \eqref{79} we obtain
\begin{equation}
\frac{1}{N}\sum\limits_{1\leq n\leq N,\sigma _{n}\in \lbrack
0,\varepsilon
]\cup \lbrack 1-\varepsilon ,1]}1\leq 3\varepsilon +c\frac{\varepsilon }{N}%
\sum\limits_{1\leq l\leq k}l^{2}+\frac{c}{\sqrt{k\varepsilon
}}\leq 3\varepsilon +\frac{c\varepsilon
k^{3}}{N}+\frac{c}{\sqrt{k\varepsilon }}. \label{80}
\end{equation}%
Now set
\begin{equation*}
k=[\varepsilon ^{-3/7}N^{2/7}]
\end{equation*}%
and note that for $N\geq \varepsilon ^{-9}$,
\begin{equation*}
\frac{1}{\sqrt{\varepsilon k}}\leq \frac{\varepsilon
}{(1-\varepsilon ^{3})^{1/2}}\leq c\varepsilon ,\quad
\frac{\varepsilon k^{3}}{N}\leq \varepsilon ,
\end{equation*}%
which along with \eqref{80} yields
\begin{equation*}
\frac{1}{N}\sum\limits_{1\leq n\leq N,\sigma _{n}\in \lbrack
0,\varepsilon ]\cup \lbrack 1-\varepsilon ,1]}1\leq c\varepsilon
\text{~~for~~}N\geq \varepsilon ^{-9},
\end{equation*}%
and the lemma \ref{LC} follows.
\end{proof}

\textbf{Proof of Proposition} \ref{distribution}

The proof in main part imitates the proof of the Weil Theorem. Fix
$\rho \in (0,1/4)$ and consider the function $\varphi _{\rho }$ as
defined above. Obviously
\begin{equation}
\frac{1}{N}\sum\limits_{1\leq n\leq N,\theta _{n}\in \lbrack
0,\rho ]}1\leq \frac{1}{N}\sum\limits_{1\leq n\leq N}\varphi
_{\rho }(\theta _{n}). \label{716}
\end{equation}%
Combining \eqref{77}-\eqref{78} we obtain
\begin{equation}
\frac{1}{N}\sum\limits_{1\leq n\leq N,\theta _{n}\in \lbrack
0,\rho ]}1\leq 3\rho +6\rho \sum\limits_{1\leq l\leq
k}\frac{1}{N}\left\vert
\sum\limits_{1\leq n\leq N}e^{2\pi il\theta _{n}}\right\vert +\frac{c}{\sqrt{%
k\rho }}.  \label{719}
\end{equation}%
Next set
\begin{equation}
W_{l,N}=\left\vert \sum\limits_{1\leq n\leq N}e^{2\pi il\theta
_{n}}\right\vert .  \label{720}
\end{equation}%
We have from \eqref{719}
\begin{equation}
\frac{1}{N}\sum\limits_{1\leq n\leq N,\theta _{n}\in \lbrack
0,\rho ]}1\leq
3\rho +6\rho \sum\limits_{1\leq l\leq k}\frac{1}{N}W_{l,N}+\frac{c}{\sqrt{%
k\rho }}.  \label{7191}
\end{equation}%
Noting that
\begin{equation*}
e^{2\pi il\theta _{n}}=e^{2\pi il(\omega _{0}n^{2}-C)}
\end{equation*}%
we obtain
\begin{equation}
W_{l,N}=\left\vert \sum\limits_{1\leq n\leq N}e^{2\pi i\omega
_{0}ln^{2}}\right\vert .  \label{7200}
\end{equation}%
Next we have
\begin{equation*}
W_{l,N}^{2}=\sum\limits_{1\leq m,n\leq N}e^{2\pi
igw_{0}l(n^{2}-m^{2})}.
\end{equation*}%
Setting $r=m+n$, $q=n-m$ we arrive at
\begin{equation*}
W_{l,N}^{2}=\sum\limits_{2\leq r\leq N}\sum\limits_{|q|\leq
r}e^{2\pi i\omega _{0}lrq}+\sum\limits_{N<r\leq
2N}\sum\limits_{|q|\leq 2N-r}e^{2\pi i\omega _{0}lrq}
\end{equation*}%
Introduce the quantities
\begin{gather*}
\chi _{rl}=\left\vert \sum\limits_{|q|\leq r}e^{2\pi i\omega
_{0}lrq}\right\vert \text{~~for~~}2\leq r\leq N, \\
\chi _{rl}=\left\vert \sum\limits_{|q|\leq 2N-r}e^{2\pi i\omega
_{0}lrq}\right\vert \text{~~for~~}N<r\leq 2N
\end{gather*}%
Thus we get
\begin{equation}
W_{l,N}^{2}=\sum\limits_{2\leq r\leq 2N}\chi _{rl}  \label{721}
\end{equation}%
Obviously
\begin{equation}
\chi _{rl}\leq 2N.  \label{722}
\end{equation}%
On the other hand, since $\chi _{rl}$ is a geometric progression
in $q$,
\begin{equation}
\chi _{rl}\leq \frac{1}{\sin (\pi \omega _{0}rl)}.  \label{723}
\end{equation}%
Choose an arbitrary $\varepsilon \in (0,1/4)$ and denote by
$J_{l}$ the set of all $r$ such that
\begin{equation*}
2\leq r\leq 2N,\quad \sigma _{rl}\equiv \{\omega _{0}rl\}\in
\lbrack 0,\varepsilon ]\cup \lbrack 1-\varepsilon ,1].
\end{equation*}%
It is easy to see that for $r\in \lbrack 2,N]\setminus J_{l}$,
\begin{equation}
\chi _{rl}\leq \frac{1}{\sin (\pi \varepsilon )}\leq
\frac{c}{\varepsilon }. \label{724}
\end{equation}%
From this and \eqref{721}, \eqref{722} we conclude that
\begin{equation}
W_{l,N}^{2}\leq 2N\sum\limits_{r\in J_{l}}1+c\frac{2N}{\varepsilon
} \label{725}
\end{equation}%
On the other hand for fixed $l$, $r\in J_{l}$, and $p=rl$ we have
$\sigma _{p}=\{\omega _{0}p\}\in \lbrack 0,\varepsilon ]\cup
\lbrack 1-\varepsilon ,1]$. Hence, since $l\leq k$,
\begin{equation*}
\text{~card~}J_{l}\leq \text{card~}\{p:1\leq p\leq 2kN,\sigma
_{p}\in \lbrack 0,\varepsilon ]\cup \lbrack 1-\varepsilon ,1]\}.
\end{equation*}%
By lemma \ref{LC}, we have
\begin{equation*}
\frac{1}{2kN}\text{card~}\{p:1\leq p\leq 2kN,\sigma _{p}\in
\lbrack
0,\varepsilon ]\cup \lbrack 1-\varepsilon ,1]\}\leq c\varepsilon \text{%
~~for~~}kN\geq \varepsilon ^{-9},
\end{equation*}%
which gives
\begin{equation*}
\sum\limits_{r\in J_{l}}1\leq c2kN\varepsilon \text{~~for~~}kN\geq
\varepsilon ^{-9}.
\end{equation*}%
Substituting this inequality in \eqref{725} we obtain
\begin{equation*}
W_{l,N}^{2}\leq ckN^{2}\varepsilon +c\frac{N}{\varepsilon }\text{~~for~~}%
kN>\varepsilon ^{-9},
\end{equation*}%
or
\begin{equation*}
\frac{1}{N}W_{l,N}\leq c\sqrt{k\varepsilon }+\frac{c}{\sqrt{N\varepsilon }}%
\text{~~for~~}kN\geq \varepsilon ^{-9}.
\end{equation*}%
Substituting this result in \eqref{7191} we finally obtain
\begin{equation*}
\frac{1}{N}\sum\limits_{1\leq n\leq N,\theta _{n}\in \lbrack
0,\rho ]}1\leq
3\rho +c\rho \sqrt{k^{3}\varepsilon }+\frac{c\rho k}{\sqrt{N\varepsilon }}+%
\frac{c}{\sqrt{k\rho }}\text{~~for~~}kN>\varepsilon ^{-9}
\end{equation*}%
It follows from this
\begin{equation}
\frac{1}{N}\sum\limits_{1\leq n\leq N,\theta _{n}\in \lbrack
0,\rho ]}1\leq
3\rho +c\rho \sqrt{k^{3}\varepsilon }+c\rho k^{3/2}\varepsilon ^{4}+\frac{c}{%
\sqrt{k\rho }}\text{~~for~~}kN>\varepsilon ^{-9}.  \label{727}
\end{equation}%
Now choose
\begin{equation*}
k=[\frac{1}{\rho ^{3}}],\quad \varepsilon =k^{-3},\quad N\geq
k^{26}.
\end{equation*}%
Obviously
\begin{equation*}
k>\frac{1}{2\rho ^{3}}\geq 32,\quad \varepsilon <1/4,\quad Nk\geq
\varepsilon ^{-9},\quad k^{3/2}\varepsilon ^{4}\leq 1,\quad
k^{3}\varepsilon =1.
\end{equation*}%
From this and \eqref{727} we conclude that
\begin{equation*}
\frac{1}{N}\sum\limits_{1\leq n\leq N,\theta _{n}\in \lbrack
0,\rho ]}1\leq c\rho \text{~~for~~}N>\rho ^{-78},
\end{equation*}%
which completes the proof of the Proposition.

\section{Pseudodifferential operators}
\label{pseudodifferential}

In this section we collect basic facts from the theory of
pseudodifferential operators. We refer to the pioneering  paper
\cite{KohnNirenberg} and monographs \cite{MTaylor},
\cite{Petersen} for general theory. Note only that different maps
from functions $A(Y,k)$ to operators $\mathfrak
A=A(Y,-i\partial_Y)$
 give rise to different theories of pseudodifferential calculus.
 In these notes we assume that $Y$ is a coordinate on the
 $2D$-torus $\mathbb T^2$ and  that the dual variable $k$ belongs to the
 lattice $\mathbb Z^2$.
The first result constitutes the continuity properties of general
pseudodifferential operators.

\begin{proposition}
\label{boundedness} Let $|\mathfrak{A}|_{0,l}^{r}<\infty $ and
$0\leq s\leq l-3$, $r+s\geq 0$. Then there is a constant $c$
depending on $s$ only so
that for all $u\in H^{s+r}(%
\mathbb{R}
^{2}/\Gamma ),$%
\begin{equation}
\Vert \mathfrak{A}u\Vert _{s}\leq
c\Big(|\mathfrak{A}|_{0,l}^{r}\Vert u\Vert
_{r}+|\mathfrak{A}|_{0,3}^{r}\Vert u\Vert _{r+s}\Big). \label{041}
\end{equation}
\end{proposition}

The proof is based on the following estimate of the convolution of
non-negative sequences. Let us consider a non-negative sequences $\mathbf{a}%
^{j}=(a^{j}(n))_{n\in \mathbb{Z}^{2}}$, $1\leq j\leq m$, and $\mathbf{v}%
=(v(n))_{n\in \mathbb{Z}^{2}}$. Set
\begin{equation*}
|\mathbf{a}^{j}|_{s}=\sup\limits_{n\in \mathbb{Z}^{2}}(1+|n|)^{s}a^{j}(n),%
\quad |||\mathbf{v}|||_{s}^{2}=\sum\limits_{n\in \mathbb{Z}%
^{2}}(1+|n|)^{2s}v(n)^{2}.
\end{equation*}

\begin{lemma}
\label{convolution} Under the above assumptions, the convolution $\mathbf{w}=%
\mathbf{a}^{1}\ast ...\ast \mathbf{a}^{m}\ast \mathbf{v}$ has the
bound
\begin{equation}
|||\mathbf{w}|||_{s}\leq c(s)\sum\limits_{j}\Big(\prod_{p\neq j}|\mathbf{a}%
^{p}|_{3}\Big)|\mathbf{a}^{j}|_{s+3}\Vert |\mathbf{v}\Vert |_{0}+\Big(%
\prod_{j}|\mathbf{a}^{j}|_{3}\Big)|||\mathbf{v}|||_{s}.
\label{042}
\end{equation}
\end{lemma}

\begin{proof}
We begin with proving \eqref{042} for $m=1$. Recalling the formula
\begin{equation*}
{w}(n)=\sum_{k_{1}+...+k_{m+1}=n}a^{1}(k_{1})..a^{m}(k_{m})v(k_{m+1})
\end{equation*}%
and noting that for $k_{1}+k_{2}=n$,
\begin{equation*}
(1+|n|)^{s}\leq c(s)\big((1+|k_{1}|)^{s}+(1+|k_{2}|)^{s}\big)
\end{equation*}%
we obtain for $m=1$,
\begin{eqnarray*}
(1+|n|)^{s}w(n) &\leq &c(s)|\mathbf{a}^{1}|_{s+3}%
\sum_{k_{1}+k_{2}=n}(1+|k_{1}|)^{-3}v(k_{2})+ \\
&&+c(s)|\mathbf{a}^{1}|_{3}%
\sum_{k_{1}+k_{2}=n}(1+|k_{2}|)^{s}v(k_{2})(1+|k_{1}|)^{-3}.
\end{eqnarray*}%
From this and the classic inequality
\begin{equation*}
\sum_{n}\Big(\sum_{k_{1}+k_{2}=n}|a(k_{1})v(k_{2})|\Big)^{2}\leq \Big(%
\sum_{k}|a(k)|\Big)^{2}\sum_{k}|v(k)|^{2}.
\end{equation*}%
we obtain \eqref{042} in the case $m=1$. The general case
obviously follows from the mathematical induction principle and
the distributive property of the convolution.
\end{proof}

Let us turn to the proof of Proposition \ref{boundedness}. We have
\begin{equation*}
\widehat{\mathfrak{A}u}(n)=\frac{1}{2\pi }\sum_{k}\widehat{A}(n-k,k)\widehat{%
u}(k)\text{~~where~~}\widehat{A}(p,k)=\frac{1}{2\pi }\int\limits_{\mathbb{T}%
^{2}}A(Y,k)e^{-iYp}\,dY,
\end{equation*}%
which yields
\begin{equation}
|\widehat{\mathfrak{A}u}(n)|\leq c\sum_{k}|\widehat{A}(n-k,k)||\widehat{u}%
(k)|\leq \big[\mathbf{a}\ast \mathbf{v}\big](n),  \label{043}
\end{equation}%
with
\begin{equation*}
a(n)=\sup\limits_{k\in
\mathbb{Z}^{2}}(1+|k|)^{-r}|\widehat{A}(n,k)|\,\quad
v(k)=(1+|k|)^{r}|\widehat{u}(k)|.
\end{equation*}%
It is easy to see that
\begin{equation}
|\mathbf{a}|_{s}\leq c|\mathfrak{A}|_{0,s}^{r},\quad |||\mathbf{v}%
|||_{s}\leq c\Vert u\Vert _{{s+r}}.  \label{044}
\end{equation}%
Applying Lemma \ref{convolution} to \eqref{043}, using inequalities %
\eqref{044} and noting that
\begin{equation*}
\Big(\Vert \mathfrak{A}u\Vert _{s}\Big)^{2}\leq c\sum_{n}(1+|n|)^{2s}|%
\widehat{\mathfrak{A}u}|^{2}
\end{equation*}%
we obtain \eqref{041}, and the proposition follows.

The next proposition gives the representation for the composition
and commutators of pseudodifferential operators

\begin{proposition}
\label{p48} Let $\mathfrak{A}$ and $\mathfrak{B}$ be
pseudodifferential operators so that for some $r,\rho \in
\mathbb{R}^{1}$ and non-negative integers $m,l$,
\begin{equation*}
|\mathfrak{A}|_{m,l}^{r}+|\mathfrak{B}|_{m,l}^{\rho }<\infty .
\end{equation*}%
Let also
\begin{equation*}
l>|r|+5+s,\quad l>s+3,\quad m\geq 2.
\end{equation*}%
Then the composition $\mathfrak{A}\mathfrak{B}$ and the commutator $%
\mathfrak{A}\mathfrak{B}-\mathfrak{B}\mathfrak{A}$ have the
representations
\begin{gather}
\mathfrak{A}\mathfrak{B}=\sum\limits_{p=0}^{d}(\mathfrak{A}\mathfrak{B})_{p}+%
\mathfrak{D}_{d+1}^{(AB)},\quad d=0,1,  \label{50a} \\
\mathfrak{A}\mathfrak{B}-\mathfrak{B}\mathfrak{A}=\sum\limits_{p=1}^{d}[%
\mathfrak{A},\mathfrak{B}]_{p}+\mathfrak{D}_{d+1}^{[A,B]},\quad
d=0,1, \label{500}
\end{gather}%
in which $(\mathfrak{A}\mathfrak{B})_{p}$ and $[\mathfrak{A},\mathfrak{B}%
]_{p}$ are the pseudodifferential operators with symbols
\begin{gather}
(AB)_{0}(Y,k)=A(Y,k)B(Y,k),\quad (AB)_{1}(Y,k)=\frac{1}{i}\partial
_{k}A(Y,k)\partial _{Y}B(Y,k),  \label{49} \\
\lbrack A,B]_{0}(Y,k)=0,\quad \lbrack
A,B]_{1}(Y,k)=\frac{1}{i}\big(\partial _{k}A(Y,k)\partial
_{Y}B(Y,k)-\partial _{Y}A(Y,k)\partial _{k}B(Y,k)\big)
\label{490} \\
\text{~~for~~}k\neq 0,\text{~~
and~~}(AB)_{p}(Y,0)=[A,B]_{p}(Y,0)=0.  \notag
\end{gather}%
The reminders have the estimates
\begin{gather}
\Vert \mathfrak{D}_{d+1}^{(AB)}u\Vert _{s}\leq c\Big(|\mathfrak{A}%
|_{d+1,s}^{r}|\mathfrak{B}|_{d+1,|r|+d+4}^{\rho }+|\mathfrak{A}|_{d+1,3}^{r}|%
\mathfrak{B}|_{d+1,|r|+d+4+s}^{\rho }\Big)\Vert u\Vert _{r+\rho
-d-1}+
\notag \\
|\mathfrak{A}|_{d+1,3}^{r}|\mathfrak{B}|_{d+1,|r|+d+4}^{\rho
}\Vert u\Vert
_{s+r+\rho -d-1},  \label{51a} \\
\Vert \mathfrak{D}_{d+1}^{[A,B]}u\Vert _{s}\leq c\Big(|\mathfrak{A}%
|_{d+1,s+|\rho |+d+4}^{r}|\mathfrak{B}|_{d+1,|r|+d+4}^{\rho }+  \notag \\
|\mathfrak{A}|_{d+1,|\rho |+d+4}^{r}|1-\mathfrak{B}|_{d+1,|r|+d+4+s}^{\rho }%
\Big)\Vert u\Vert _{r+\rho -d-1}+  \notag \\
|\mathfrak{A}|_{d+1,|\rho
|+d+4}^{r}|1-\mathfrak{B}|_{d+1,|r|+d+4}^{\rho }\Vert u\Vert
_{s+r+\rho -d-1},  \label{510}
\end{gather}%
in which the constant $c$ depends on $s,r,\rho $ only.
\end{proposition}

\noindent The proof is based on the following lemma

\begin{lemma}
\label{l40} Let
\begin{equation}
\mathcal{R}_{d+1}(\eta ,\zeta ,k)=\widehat{A}(\eta ,\zeta
+k)-\sum\limits_{p=0}^{d}\sum\limits_{|\alpha |=p}\frac{1}{\alpha
!(i)^{|\alpha |}}\partial _{k}^{\alpha }\widehat{{A}}(\eta
,k)(i\zeta )^{\alpha }\text{~~for~~}k\neq 0,  \label{41a}
\end{equation}%
and $\mathcal{R}_{d+1}(\eta ,\zeta ,0)=\widehat{A}(\eta ,\zeta )$,
where
\begin{equation*}
\widehat{A}(\eta ,k)=\frac{1}{2\pi
}\int\limits_{\mathbb{T}^{2}}e^{-i\eta Y}A(Y,k)dY.
\end{equation*}%
Then for all $\eta ,\zeta ,k\in \mathbb{Z}^{2}$ and $0\leq s\leq
l$,
\begin{equation}
|\mathcal{R}_{d+1}(\eta ,\zeta ,k)|\leq
c(d,r)|A|_{d+1,l}^{r}(1+|\eta |)^{-l}(1+|\zeta
|)^{|r|+d+1}(1+|k|)^{r-d-1}.  \label{42a}
\end{equation}
\end{lemma}

\begin{proof}
It suffices to prove \eqref{42} for $k\neq 0$ only. If $|\zeta
|\leq |k|/2$, then the Taylor formula
\begin{equation*}
\mathcal{R}_{d+1}(\eta ,\zeta ,k)=\sum\limits_{|\alpha |=d+1}\frac{d+1}{%
\alpha !}\Big\{\int\limits_{0}^{1}[\partial _{k}^{\alpha
}\widehat{A}](\eta ,k+t\zeta )(1-t)^{d}\,dt\Big\}\zeta ^{\alpha }
\end{equation*}%
implies the estimate
\begin{gather}
|\mathcal{R}_{d+1}(\eta ,\zeta ,k)|\leq
c|\mathfrak{A}|_{d+1,s}^{r}(1+|\eta |)^{-s}(1+|\zeta
|)^{d+1}\int\limits_{0}^{1}(1+|k+t\zeta |)^{r-d-1}dt\leq
\notag \\
c|\mathfrak{A}|_{d+1,s}^{r}(1+|\eta |)^{-s}(1+|\zeta
|)^{d+1}(1+|k|)^{r-d-1}, \label{44a}
\end{gather}%
which obviously yields \eqref{42}. If $|\zeta |\geq |k|/2,$ we
have
\begin{gather*}
|\mathcal{R}_{d+1}(\zeta ,k)|\leq c|\mathfrak{A}|_{d+1,s}^{r}(1+|\eta |)^{-s}%
\Big[|k+\zeta |^{r}+\sum\limits_{0}^{d}(1+|\zeta |)^{p}|k|^{r-p}\Big] \\
\leq c|\mathfrak{A}|_{d+1,s}^{r}(1+|\eta |)^{-s}\Big[|\zeta
|^{r}+\sum\limits_{0}^{d}(1+|\zeta |)^{p}|k|^{r-p}\Big].
\end{gather*}%
Noting that for $0\leq p\leq d+1$ and $|\zeta |\geq |k|/2$,
\begin{equation*}
(1+|\zeta |)^{p}|k|^{r-p}\leq c(1+|\zeta |)^{d+1}(1+|k|)^{r-d-1}
\end{equation*}%
we obtain
\begin{equation*}
|\mathcal{R}_{d+1}(\zeta ,k)|\leq
c(s)|\mathfrak{A}|_{d+1,s}^{r}(1+|\eta |)^{-s}(1+|\zeta
|)^{|r|+d+1}(1+|k|)^{r-d-1},
\end{equation*}%
and the lemma follows.
\end{proof}

Let us turn to the proof of the proposition. Since $[\mathfrak{A},\mathfrak{B%
}]=-[\mathfrak{A},(1-\mathfrak{B})]$, it suffices to prove
\eqref{50} only. To this end note that, by the definition of
pseudodifferential operator,
\begin{equation*}
\widehat{\mathfrak{A}\mathfrak{B}u}(n)=\sum\limits_{p,k}\widehat{A}(n-p,p)%
\widehat{B}(p-k,k)\widehat{u}(k)=\sum\limits_{\zeta
,k}\widehat{A}(n-k-\zeta ,k+\zeta )\widehat{B}(\zeta
,k)\widehat{u}(k).
\end{equation*}%
Applying Lemma \ref{l40} to the Fourier transform
$\widehat{A}(\eta ,k)$ of the symbol $A$ we arrive at the identity
\begin{gather*}
\sum\limits_{{\small
\begin{array}{c}
k,\zeta \in \mathbb{Z}^{2} \\
k\neq 0%
\end{array}%
}}\sum\limits_{p=0}^{d}\sum\limits_{|\alpha |=p}\frac{1}{\alpha
!(i)^{|\alpha |}}[\partial _{k}^{\alpha }\widehat{A}](n-k-\zeta ,k)\Big[%
(i\zeta )^{\alpha }\widehat{B}(\zeta ,k)\Big]\widehat{u}(k)+ \\
+\sum\limits_{k,\zeta \in
\mathbb{Z}^{2}}\mathcal{R}_{d+1}(n-k-\zeta ,\zeta
,k)\widehat{B}(\zeta ,k)\widehat{u}(k)=\widehat{\mathfrak{A}\mathfrak{B}u}%
(n).
\end{gather*}%
Noting that $(i\zeta )^{\alpha }\widehat{B}(\zeta
,k)=\widehat{\partial _{Y}^{\alpha }B}(\zeta ,k)$, we obtain
\begin{equation*}
\sum\limits_{{\small
\begin{array}{c}
k,\zeta \in \mathbb{Z}^{2} \\
k\neq 0%
\end{array}%
}}\sum\limits_{|\alpha |=p}\frac{1}{\alpha !(i)^{|\alpha
|}}[\partial
_{k}^{\alpha }\widehat{A}](n-k-\zeta )\Big[(i\zeta )^{\alpha }\widehat{B}%
(\zeta
,k)\Big]\widehat{u}(k)=\widehat{(\mathfrak{A}\mathfrak{B})_{p}u}(n),
\end{equation*}%
which leads to representation \eqref{50} with the remainder
\begin{equation*}
\widehat{\mathfrak{D}_{m+1}^{(AB)}u}(n)=\sum\limits_{k,\zeta \in \mathbb{Z}%
^{2}}\mathcal{R}_{d+1}(n-k-\zeta ,\zeta ,k)\widehat{B}(\zeta ,k)\widehat{u}%
(k).
\end{equation*}%
In particular, we have the inequality
\begin{equation}
|\widehat{\mathfrak{D}_{m+1}^{(AB)}u}|\leq \mathbf{a}\ast
\mathbf{b}\ast \mathbf{v},  \label{50b}
\end{equation}%
in which the elements of the sequences $\mathbf{a}$, $\mathbf{b}$, $\mathbf{v%
}$ are given by
\begin{gather*}
a(n)=\sup\limits_{\zeta ,k}(1+|\zeta |)^{-|r|-d-1}(1+|k|)^{-r+d+1}|\mathcal{R%
}_{d+1}(n,\zeta ,k)|, \\
b(\zeta )=(1+|\zeta |)^{|r|+d+1}\sup\limits_{k}(1+|k|)^{-\rho }|\widehat{B}%
(\zeta ,k)|,v(k)=(1+|k|)^{r+\rho -d-1}|\widehat{u}(k)|.
\end{gather*}%
Applying Lemma \ref{convolution} to the right side of \eqref{50b}
and noting that by Lemma \ref{l40}
\begin{equation*}
|\mathbf{a}|_{s}\leq |\mathfrak{A}|_{d+1,s}^{r},\quad |\mathbf{b}|_{s}\leq |%
\mathfrak{B}|_{0,s+|r|+d+1}^{\rho },\quad |||\mathbf{v}|||_{s}\leq
c\Vert u\Vert _{\rho +r-d-1+s}
\end{equation*}%
we obtain \eqref{51} and the proposition follows.

It is useful to reformulate the above results in terms of infinite
matrices. To this end we introduce the Hilbert space
$\mathbf{F}_{s,t}$ which consists
of all operators $\mathcal{Y}:H^{s}(\mathbb{R}^{2}/\Gamma )\mapsto H^{t}(%
\mathbb{R}^{2}/\Gamma )$ having the representation
\begin{equation*}
\widehat{\mathcal{Y}u}(k)=\sum\limits_{p\in \mathbb{Z}^{2}}\mathcal{Y}_{kp}%
\widehat{u}(p)
\end{equation*}%
such that
\begin{equation}
\Vert \mathcal{Y}\Vert _{\mathbf{F}_{s,t}}^{2}:=\sup \big\{%
\sum\limits_{k}(1+|k|)^{2t}\Big(\sum\limits_{p}|\mathcal{Y}_{kp}||\widehat{u}%
(p)|\Big)^{2}:\quad \sum\limits_{k}(1+|k|)^{2s}|\widehat{u}(k)|^{2}\big\}%
<\infty .  \label{G11}
\end{equation}

\begin{corollary}
\label{matrix}

(i) Under the assumptions of Proposition \ref{boundedness}, operator $%
\mathfrak{A}$ has a matrix representation with $\mathfrak{\ A}_{kp}=\widehat{%
A}(k-p,p)$ and
\begin{multline}
\sum\limits_{k}(1+|k|)^{2s}\Big(\sum\limits_{p}|\mathfrak{A}_{kp}||\widehat{u%
}(p)|\Big)^{2}\leq  \notag \\
c(|\mathfrak{A}|_{0,l}^{r})^{2}\sum\limits_{k}(1+|k|)^{2r}|\widehat{u}%
(k)|^{2}+c(|\mathfrak{A}|_{0,3}^{r})^{2}\sum\limits_{k}(1+|k|)^{2r+2s}|%
\widehat{u}(k)|^{2}.  \label{G12}
\end{multline}%
In particular, $\Vert \mathfrak{A}\Vert _{\mathbf{F}_{r+s,s}}\leq c|%
\mathfrak{A}|_{0,l}^{r}$.

(ii) Under the assumptions of Proposition \ref{p48}, the operator $\mathfrak{%
D}_{d+1}^{(AB)}$ has a matrix representation so that
\begin{gather}
\sum\limits_{k}(1+|k|)^{2(s)}\Big(\sum\limits_{p}|\mathfrak{D}%
_{d+1,kp}^{(AB)}||\widehat{u}(p)|\Big)^{2}\leq c\Big(|\mathfrak{A}%
|_{d+1,s}^{r}|\mathfrak{B}|_{d+1,|r|+d+4}^{\rho }+  \notag \\
|\mathfrak{A}|_{d+1,3}^{r}|\mathfrak{B}|_{d+1,|r|+d+4+s}^{\rho }\Big)%
^{2}\sum\limits_{k}(1+|k|)^{2r+2\rho -2d-2}|\widehat{u}(k)|^{2}+  \notag \\
c(|\mathfrak{A}|_{d+1,3}^{r}|\mathfrak{B}|_{d+1,|r|+d+4}^{\rho
})^{2}\sum\limits_{k}(1+|k|)^{2s+2r+2\rho
-2d-2}|\widehat{u}(k)|^{2} \label{G13}
\end{gather}%
In particular,
\begin{equation*}
\Vert \mathfrak{D}_{d+1}^{(AB)}\Vert _{\mathbf{F}_{r+s+\rho -d-1,s}}\leq c%
\Big(|\mathfrak{A}|_{d+1,s}^{r}|\mathfrak{B}|_{d+1,|r|+d+4}^{\rho
}+\newline
|\mathfrak{A}|_{d+1,3}^{r}|\mathfrak{B}|_{d+1,|r|+d+4+s}^{\rho
}\Big)
\end{equation*}
\end{corollary}

\begin{proof}
Assertion (\textit{i}) is integral part of the proof of Proposition \ref%
{boundedness}. In order to prove (\textit{ii}) note that
\begin{equation*}
\mathfrak{D}_{d+1,kp}^{(AB)}=\sum\limits_{\zeta \in \mathbb{Z}^{2}}\mathcal{R%
}_{d+1}(p-k-\zeta ,\zeta ,k)\widehat{B}(\zeta ,k),
\end{equation*}%
where $\mathcal{R}_{d+1}$ is defined by formula \eqref{41a}. The
needed result follows from \eqref{50b} and Lemma
\ref{convolution}.
\end{proof}

\textbf{Proof of Proposition} \ref{101l}

We give the proof of representations \eqref{10200} and
\eqref{1021} only, and begin with proving \eqref{10200}. Note that
\begin{equation*}
A(Y,\xi )\equiv A_{r}(Y,\xi _{1},\xi _{2}^{2})+i\xi
_{2}A_{i}(Y,\xi _{1},\xi _{2}^{2}),
\end{equation*}%
where
\begin{equation*}
A_{r}(\cdot ,\cdot ,\rho )=\frac{1}{2}\Big(A(\cdot ,\cdot ,\sqrt{\rho }%
)+A(\cdot ,\cdot ,-\sqrt{\rho })\Big),A_{i}(\cdot ,\cdot ,\rho )=\frac{1}{2i%
\sqrt{\rho }}\Big(A(\cdot ,\cdot ,\sqrt{\rho })-A(\cdot ,\cdot ,-\sqrt{\rho }%
)\Big).
\end{equation*}%
Assuming $k_{1}\neq 0$ and noting that in this case $\xi
_{2}^{2}=1-\xi _{1}^{2}$, we arrive at the identity
\begin{equation*}
A(Y,\xi )\equiv \mathcal{A}_{r}(Y,\xi _{1})+i\xi
_{2}\mathcal{A}_{i}(Y,\xi _{1})\equiv \mathcal{A}(Y,\xi )\text{~~
where ~~}\mathcal{A}_{\beta }(Y,\xi _{1})=A_{\beta }(Y,\xi
_{1},1-\xi _{1}^{2}).
\end{equation*}%
From the Taylor formula we conclude that for $k_{1}\neq 0$,
\begin{gather}
\mathcal{A}(Y,\xi
)=\sum\limits_{j=0}^{j=2}\frac{1}{(i)^{j}j!}\Big[\partial
_{\xi _{1}}^{j}\mathcal{A}\Big](Y,0,\xi _{2})(i\xi _{1})^{j}+(\xi _{1})^{3}%
\mathcal{R}_{3}(Y,\xi )  \label{taylor} \\
\mathcal{R}_{3}(Y,\xi )=\frac{1}{6}\int\limits_{0}^{1}[\partial
_{\xi _{1}}^{3}\mathcal{A}](Y,s\xi _{1},\xi _{2})(1-s)^{3}\,ds.
\notag
\end{gather}%
Recalling symmetry property \eqref{104} we obtain
\begin{gather*}
\mathcal{A}_{r}+i\xi _{2}\mathcal{A}_{i}\Big|_{{\small \xi _{1}=0}}=\Big[%
\text{Re~}A+i\xi _{2}\text{Im~}A\Big](Y,0,1)=A_{0}(Y,\xi _{2}), \\
\partial _{\xi _{1}}(\mathcal{A}_{r}+i\xi _{2}\mathcal{A}_{i})\Big|_{{\small %
\xi _{1}=0}}=i\Big[\partial _{\xi _{1}}(\text{Im~}A-i\xi _{2}\text{Re~}A)%
\Big](Y,0,1)=-{\nu }A_{1}(Y,\xi _{2}) \\
\partial _{\xi _{1}}^{2}(\mathcal{A}_{r}+i\xi _{2}\mathcal{A}_{i})\Big|_{%
{\small \xi _{1}=0}}=\Big[(\partial _{\xi _{1}}^{2}-\partial _{\xi _{2}})(%
\text{Re~}A+i\xi _{2}\text{Im~}A)+i\xi _{2}\text{Im~}A\Big](Y,0,1)= \\
-2{\nu }^{2}A_{2}(Y,\xi _{2}),
\end{gather*}%
which being substituted into \eqref{taylor} leads to
\begin{equation}
A(Y,\xi )=A_{0}(Y,\xi _{2})-{\nu }A_{1}(Y,\xi _{2})i\xi _{1}+\nu
^{2}A_{2}(Y,\xi _{2})(i\xi _{1})^{2}+\mathcal{R}_{3}(Y,\xi )\xi
_{1}^{3}. \label{1012b}
\end{equation}%
Noting that
\begin{equation}
\begin{split}
ik_{1}i\xi _{1}& =-\frac{1}{\nu }+\frac{1}{\nu |\mathbb{T}^{-1}(k)|}L(k), \\
ik_{1}(i\xi _{1})^{2}& =\Big(\frac{1}{\nu }\Big)^{2}\frac{1}{ik_{1}}-\Big(%
\frac{1}{\nu
}\Big)^{2}\frac{1}{ik_{1}|\mathbb{T}^{-1}(k)|^{2}}\Big(\nu
k_{1}^{2}+|\mathbb{T}^{-1}(k)|\Big)L(k), \\
ik_{1}(\xi _{1})^{3}& =\frac{i}{\nu
^{2}|\mathbb{T}^{-1}(k)|}-\frac{i}{\nu
^{2}|\mathbb{T}^{-1}(k)|^{3}}(\nu
k_{1}^{2}+|\mathbb{T}^{-1}(k)|)L(k),
\end{split}
\label{1020k}
\end{equation}%
we get for $k_{1}\neq 0$,
\begin{equation}
ik_{1}A(Y,\xi )=\sum\limits_{j=0}^{2}A_{j}(Y,\xi
_{2})(ik_{1})^{1-j}+P_{A}(Y,\xi )L(k)+Q_{A}(Y,\xi ). \label{1020c}
\end{equation}%
Here the remainders are given by
\begin{gather}
P_{A}(Y,k)=\frac{1}{|\mathbb{T}^{-1}(k)|}A_{1}(Y,\xi _{2})+\frac{i}{k_{1}|%
\mathbb{T}^{-1}(k)|^{2}}\Big(A_{2}(Y,\xi _{2})-  \label{1020f} \\
\frac{k_{1}}{\nu ^{2}|\mathbb{T}^{-1}(k)|}\mathcal{R}_{3}(Y,\xi )\Big)\Big(%
\nu k_{1}^{2}+|\mathbb{T}^{-1}(k)|\Big),\quad Q_{A}(Y,k)=\frac{i}{\nu ^{2}|%
\mathbb{T}^{-1}(k)|}\mathcal{R}_{3}(Y,\xi ).  \notag
\end{gather}%
Next set $P_{A}(Y,k)=Q_{A}(Y,k)=0$ for $k_{1}=0$ and denote by $\mathfrak{P}%
_{A}$, $\mathfrak{Q}_{A}$ the pseudodifferential operators with the symbols $%
P_{A}$, $Q_{A}$. With this notation, decomposition \eqref{10200}
easy follows from \eqref{1020c}. It remains to note that formulae
\eqref{taylor} and \eqref{1020f} imply the estimate
\begin{equation*}
|\mathfrak{P}_{A}|_{0,l}^{-1}+|\mathfrak{Q}_{A}|_{0,l}^{-1}\leq |\mathfrak{A}%
|_{3,l},
\end{equation*}%
which along with Proposition \ref{boundedness} yields inequalities %
\eqref{1020e0} for $\mathfrak{P}_{A}$ and $\mathfrak{Q}_{A}$.

Let us turn to the proof of \eqref{1021}. We begin with the
observation that
for $k_{1}\neq 0$, the product of the symbols of the elementary operators $%
\mathfrak{A}_{j}$ and $\mathfrak{W}$ is equal to
\begin{equation}
A_{j}W\equiv S_{j}^{\prime }+\xi _{1}^{2}\text{Im~}\tilde{A_{j}}\text{Im~}%
\tilde{W},  \label{1021a}
\end{equation}%
where $S_{j}^{\prime }$ are the symbols of the elementary
operators associated with the functions $\tilde{A}_{j}\tilde{W}$.
Multiplying both sides of \eqref{1020c} by $W(Y,\xi _{2})$ and
using the identities
\begin{gather*}
ik_{1}\xi _{1}^{2}=-\Big(\frac{1}{\nu }\Big)^{2}\frac{1}{ik_{1}}+\Big(\frac{1%
}{\nu }\Big)^{2}\frac{1}{|\mathbb{T}^{-1}(k)|ik_{1}}\big(\nu k_{1}^{2}+|%
\mathbb{T}^{-1}(k)|\big)L(k), \\
\xi _{1}^{2}=\frac{1}{\nu }\frac{1}{|\mathbb{T}^{-1}(k)|}-\frac{1}{\nu }%
\frac{1}{|\mathbb{T}^{-1}(k)|^{2}}L(k)
\end{gather*}%
we arrive at
\begin{equation*}
ik_{1}S(Y,k)=\sum\limits_{j=0}^{2}(ik_{1})^{1-j}S_{j}^{\prime }(Y,\xi _{2})-%
\Big(\frac{1}{\nu }\Big)^{2}\text{Im~}\tilde{A_{0}}\text{Im~}\tilde{W}\frac{1%
}{ik_{1}}+L(k)U_{S}(Y,k)+V_{S}(Y,k),
\end{equation*}%
where the symbols $U_{S}(Y,k)$ and $V_{S}(Y,k)$ are defined by
\begin{gather*}
U_{S}=P_{A}W+\Big(\frac{1}{\nu }\Big)^{2}\frac{1}{|\mathbb{T}^{-1}(k)|ik_{1}}%
\big(\nu k_{1}^{2}+|\mathbb{T}^{-1}(k)|\big)\text{Im~}\tilde{A_{0}}\text{Im~}%
\tilde{W}- \\
\frac{1}{\nu |\mathbb{T}^{-1}(k)|^{2}}\big(\text{Im~}\tilde{A}%
_{1}+(ik_{1})^{-1}\text{Im~}\tilde{A}_{2}\Big)\text{Im~}\tilde{W}, \\
V_{S}=Q_{A}W+\frac{1}{\nu |\mathbb{T}^{-1}(k)|}\big(\text{Im~}\tilde{A}%
_{1}+(ik_{1})^{-1}\text{Im~}\tilde{A}_{2}\Big)\text{Im~}\tilde{W}.
\end{gather*}%
Recall that $\tilde{W}$ and $\tilde{A}$ are smooth function on tori $\mathbb{%
T}^{2}$ which do not depend on $k$. Setting $U_{S}(Y,k)=V_{S}(Y,k)=0$ for $%
k_{1}=0$, denoting by $\mathfrak{U}_{S}$ and $\mathfrak{V}_{S}$
the pseudodifferential operators with the symbols $U_{S}$, $V_{S}$
and arguing as before we obtain desired identity \eqref{1021}. It
remains to note that estimate \eqref{1020ee} follows from
Proposition \eqref{boundedness} which completes the proof.

\section{Dirichlet-Neuman operator\label{a03}}

In this section we deduce the basic decomposition for the
Dirichlet-Neumann operator and  prove Theorem \ref{thmDir-Neum}.
Let us denote the change of coordinates by
\begin{eqnarray}
x &=&x(y),\text{ \ }x=(X,x_{3}),\text{ \ }y=(Y,y_{3}),  \notag \\
X &=&X(Y),\text{ \ }x_{3}=y_{3}+\tilde{\eta}(Y),
\label{varchange1}
\end{eqnarray}%
where
\begin{equation*}
\tilde{\eta}(Y)=\eta (X(Y)),
\end{equation*}%
and $X(\cdot )$ is a  diffeomorphism of the torus of the form
(\ref{varchange2}) satisfying the Condition \ref{diffeomorphism}.
Let us notice that with the new coordinate $Y=(y_{1},y_{2})$ the
lattice $\Gamma ^{\prime }$ is generated by the two wave vectors
$(1,\pm 1).$ We still
denote the lattice of periods by $\Gamma .$ For a function $u(x)$ we define $%
\widetilde{u}(y)$ by $\widetilde{u}(y)=u(x(y)).$ The Jacobian
matrix of the
above diffeomorphism reads%
\begin{equation*}
\mathbb{B}_{1}(Y)=\left(
\begin{array}{ccc}
\partial _{y_{1}}X_{1} & \partial _{y_{2}}X_{1} & 0 \\
\partial _{y_{1}}X_{2} & \partial _{y_{2}}X_{2} & 0 \\
\partial _{y_{1}}\tilde{\eta} & \partial _{y_{2}}\tilde{\eta} & 1%
\end{array}%
\right)
\end{equation*}%
and the determinant satisfies $J=\det \mathbb{B}(Y)=\det
\mathbb{B}_{1}(Y).$ Now, we use the following identities for any
scalar function $u,$ and vector
function $V:$%
\begin{eqnarray*}
\nabla _{x}u(x(y)) &=&(\mathbb{B}_{1}^{\ast })^{-1}\nabla _{y}\widetilde{u}%
(y), \\
\nabla _{x}\cdot V(x(y)) &=&\frac{1}{J}\nabla _{y}\cdot (J\mathbb{B}_{1}^{-1}%
\tilde{V}(y)).
\end{eqnarray*}%
With these identities, the Dirichlet-Neumann operator
(\ref{Dir-Neu}) takes
the new following form%
\begin{eqnarray*}
\mathcal{A}\tilde{\varphi} &=&0,\text{ \ \ }y_{3}\in (-\infty ,0), \\
\tilde{\varphi}|_{y_{3}=0} &=&\tilde{\psi}(Y), \\
\nabla \tilde{\varphi} &\rightarrow &0\text{ as }y_{3}\rightarrow
-\infty ,
\end{eqnarray*}%
\begin{equation}
\mathcal{G}_{\eta }\psi =\frac{1}{J}(\mathbb{A}(Y)\nabla _{y}\tilde{\varphi}%
).\mathbf{e}_{3},  \label{Dir-Neu_1}
\end{equation}%
where%
\begin{eqnarray*}
\mathcal{A}\tilde{\varphi} &=&\nabla _{y}\cdot (\mathbb{A}(Y)\nabla _{y}%
\tilde{\varphi}), \\
\mathbb{A} &=&J(\mathbb{B}_{1}^{\ast }\mathbb{B}_{1})^{-1}\text{ \
(symmetric matrix)} \\
\det \mathbb{A} &=&J.
\end{eqnarray*}%
Notice that for computing the new expression of $\mathcal{G}_{\eta
}\psi ,$
we used the fact that $\Phi (x)=x_{3}-\eta (X)$ is such that $\widetilde{%
\Phi }(y)=y_{3},$ hence $\nabla _{x}\Phi =(\mathbb{B}_{1}^{\ast })^{-1}%
\mathbf{e}_{3}.$

We already defined the 2x2 matrix $\mathbb{G}(Y)$ of the first
fundamental
form of the free surface, and we have%
\begin{equation*}
\mathbb{B}_{1}^{\ast }\mathbb{B}_{1}(Y)=\left(
\begin{array}{ccc}
g_{11} & g_{12} & \partial _{y_{1}}\widetilde{\eta } \\
g_{12} & g_{22} & \partial _{y_{2}}\widetilde{\eta } \\
\partial _{y_{1}}\widetilde{\eta } & \partial _{y_{2}}\widetilde{\eta } & 1%
\end{array}%
\right) .
\end{equation*}%
We assume that elements of the matrix $\mathbb{A}(Y)$ are smooth $2\pi $%
-periodic functions, and for some $\rho $ and $l$, satisfying the
following inequalities for $9\leq \rho \leq l$,
\begin{equation}
\Vert \mathbb{A}_{0}-\mathbb{A}\Vert _{C^{\rho }}\leq c\varepsilon
,\quad
\Vert \mathbb{A}\Vert _{C^{l}}\leq E_{l},\text{~~ where~~}\mathbb{A}%
_{0}=\tau ^{-1}\text{\textrm{diag~}}\{1,\tau ^{2},1\},  \label{N3}
\end{equation}%
where, by construction, the estimates of $\mathbb{A}$ and $\mathbb{A}_{0}-%
\mathbb{A}$ would come from%
\begin{eqnarray}
||\widetilde{\eta }||_{C^{\rho
+1}}+||\widetilde{\mathcal{V}}(\cdot
)||_{C^{\rho +1}} &\leq &c\varepsilon ,  \label{estimA} \\
|\widetilde{|\eta }||_{C^{s+1}}+||\widetilde{\mathcal{V}}(\cdot
)||_{C^{s+1}} &\leq &c(l),\text{ \ \ }s\leq l,  \notag
\end{eqnarray}%
where $\widetilde{\mathcal{V}}$ is defined in (\ref{varchange2}).
By the
factorization theorem, there are first order pseudodifferential operators $%
\mathcal{G}^{\pm }$ so that
\begin{gather}
\mathcal{A}=a_{33}(\partial _{y_{3}}+\mathcal{G}^{+})(\partial _{y_{3}}+%
\mathcal{G}^{-}),  \label{NFAC} \\
e^{-ik\cdot Y}\mathcal{G}^{\pm }e^{ik\cdot Y}\rightarrow \pm \infty \text{%
~~as~~}|k|\rightarrow \infty .  \notag
\end{gather}%
\noindent It follows from \eqref{Dir-Neu_1} that
\begin{equation}
\mathcal{G}_{\eta }\psi
=\frac{1}{J}\big(\sum\limits_{j=1}^{2}a_{3j}\partial
_{y_{j}}-a_{33}\mathcal{G}^{-}\big)\widetilde{\psi }.  \label{N5a}
\end{equation}%
Hence the task now is to split $\mathcal{G}^{-}$ into a sum of
first and zero order pseudodifferential operators. The
corresponding result is given by the following

\begin{theorem}
\label{thmN1} Under the above assumptions there is $\varepsilon
_{0}$ depending on $\rho $ and $l$ only such that for all
$\varepsilon \in (0,\varepsilon _{0})$, the operator
$\mathcal{G}^{-}$ has the representation
\begin{equation}
\mathcal{G}^{-}=\mathcal{G}_{1}^{-}+\mathcal{G}_{0}^{-}+\mathcal{G}_{-1}^{-},
\label{N4}
\end{equation}%
in which the pseudodifferential operators
$\mathcal{G}_{0}^{-}$,$\mathcal{G}_{1}^{-}$, have the symbols
$G_{0}^{-}$ and $G_{1}^{-}$ defined by
\begin{gather}
G_{0}^{-}=\frac{a_{33}}{2D}\big(\,i\nabla_{k}G_{1}^{+}\nabla
_{Y}G_{1}^{-}-2bG_{1}^{-}+C_{1}\big),  \label{G_0formula}\\
G_{1}^{\pm}(Y,k)=\frac{1}{a_{33}}(i a_{31}k_{1}+i a_{32} k_{2}\pm
D),\label{symbol1}
\end{gather}
where\begin{equation}\begin{split}\label{Da} D(Y,k)=\Big\{
\sum\limits_{1\leq j,m\leq 2}a_{33}a_{jm}k_{j}k_{m}
-\bigl(a_{31}{k_{1}}+a_{32}k_{2}\bigr)%
^{2}\Big\} ^{1/2},\\
C_1(Y,k)=\frac{i}{a_{33}}\sum\limits_{j,m=1}^2(\partial_{y_j}a_{jm})k_m,
\quad 2b=\frac{1}{a_{33}}\sum\limits_{j=1}^2\partial_{y_j}a_{3j}.
\end{split}\end{equation}%
The zero-order pseudodifferential operator $\mathcal{G}_{0}^{-}$
satisfies the inequality
\begin{equation}
|{\mathcal{G}}_{0}^{-}|_{m,n-1}^{0}\leq c\Vert \mathbb{A}-\mathbb{A}%
_{0}\Vert _{C^{n}},\quad m\geq 0,\quad n\leq l,  \label{N8a}
\end{equation}%
and for any $u\in H^{r-1}(\mathbb{R}^{2}/\Gamma )$, $r<\rho -8$
and $1\leq s\leq l-8$, the rest term has the bound
\begin{equation}
\Vert \mathcal{G}_{-1}^{-}u\Vert _{r}\leq c\varepsilon \Vert u\Vert _{{r-1}%
},\quad \Vert \mathcal{G}_{-1}^{-}u\Vert _{s}\leq c(E_{l}\Vert u\Vert _{{r-1}%
}+\varepsilon \Vert u\Vert _{s-1}).  \label{symbol-1}
\end{equation}%
Moreover, operators $\mathcal{G}_{1}^{-},$ $\mathcal{G}_{0}^{-},$ $\mathcal{G%
}_{-1}^{-}$ verify the following symmetry properties%
\begin{equation}
\mathcal{G}_{j}^{-}u(\pm Y^{\ast })=\mathcal{G}_{j}^{-}u^{\ast }(\pm Y),%
\text{ \ \ }j=1,0,-1,\text{ \ }u^{\ast }(Y)=u(Y^{\ast }).
\label{symG_j}
\end{equation}
\end{theorem}

\begin{proof}
First we rewrite the operator $\mathcal{A}$ in the form
\begin{equation}
\mathcal{A}=a_{33}\frac{\partial ^{2}}{\partial y_{3}^{2}}+2a_{33}\frac{%
\partial }{\partial y_{3}}\mathcal{B}+a_{33}\mathcal{C},\quad \mathcal{C}=%
\mathcal{C}_{2}+\mathcal{C}_{1},\quad
\mathcal{B}=\mathcal{B}_{1}+b \label{N1}
\end{equation}%
with differential operators
\begin{equation}
\begin{split}
\mathcal{C}_{2}& =a_{33}^{-1}\sum\limits_{i,j=1}^{2}a_{ij}\partial
_{y_{i}}\partial _{y_{j}},\quad \mathcal{C}_{1}=a_{33}^{-1}\sum%
\limits_{i,j=1}^{2}(\partial _{y_{i}}a_{ij})\partial _{y_{j}}, \\
\mathcal{B}_{1}& =a_{33}^{-1}\sum\limits_{j=1}^{2}a_{3j}\partial
_{y_{j}},\quad 2b=a_{33}^{-1}\sum\limits_{j=1}^{2}\partial
_{y_{j}}a_{3j}.
\end{split}
\label{N2a}
\end{equation}%
Combining \eqref{NFAC} and \eqref{N1} we obtain
\begin{equation}
\mathcal{G}^{+}+\mathcal{G}^{-}=2\mathcal{B}\text{, \ \ \ \ }\mathcal{G}^{+}%
\mathcal{G}^{-}=\mathcal{C}.  \label{N5}
\end{equation}%
We find the solution of (\ref{N5}) in the form
\begin{equation}
\mathcal{G}^{\pm }=\mathcal{G}_{1}^{\pm }+\mathcal{G}_{0}^{\pm }\pm \mathcal{%
X},  \label{N6a}
\end{equation}%
with%
\begin{equation}
\mathcal{G}_{1}^{+}+\mathcal{G}_{1}^{-}=2\mathcal{B}_{1},\text{ \ }\mathcal{G%
}_{0}^{+}+\mathcal{G}_{0}^{-}=2b,  \label{identif}
\end{equation}%
where the symbol of pseudodifferential operator
$\mathcal{G}_{1}^{-}$ is
given by formula \eqref{symbol1}, $\mathcal{G}_{0}^{-}$ satisfying (\ref{N8a}%
) and $\mathcal{X}$ being unknown. It follows from these formulae
that
\begin{equation*}
\mathcal{G}_{1}^{\pm }=\pm (-\Delta
)^{1/2}+\tilde{\mathcal{G}}_{1}^{\pm }
\end{equation*}%
and for any integers $m\geq 0$ and $n\leq l$
\begin{equation}
|\tilde{\mathcal{G}}_{1}^{-}|_{m,n}^{1}+|{\mathcal{G}}_{0}^{-}|_{m,n-1}^{0}+|%
{\mathcal{B}}|_{m,n-1}^{1}\leq c\Vert
\mathbb{A}-\mathbb{A}_{0}\Vert _{C^{n}}.  \label{N8}
\end{equation}%
Representation \eqref{50a} from Proposition \ref{p48} yields the
decompositions
\begin{equation}
\begin{split}
\mathcal{G}_{1}^{+}\mathcal{G}_{1}^{-}& =\mathcal{G}_{1}^{(2)}+\mathcal{G}%
_{1}^{(1)}+\mathcal{R}_{1}, \\
\mathcal{G}_{0}^{+}\mathcal{G}_{1}^{-}& =\mathcal{G}_{01}^{(1)}+\mathcal{R}%
_{01},\text{ \ }\mathcal{G}_{1}^{+}\mathcal{G}_{0}^{-}=\mathcal{G}%
_{10}^{(1)}+\mathcal{R}_{10}.
\end{split}
\label{N7a}
\end{equation}%
Here the second order pseudodifferential operator
$\mathcal{G}_{1}^{(2)}$ have the symbol $G_{1}^{+}G_{1}^{-}$, and
satisfies the identity
\begin{equation}
\mathcal{G}_{1}^{(2)}=\mathcal{C}_{2},  \label{N9}
\end{equation}%
the first order pseudodifferential operators $\mathcal{G}_{1}^{(1)}$, $%
\mathcal{G}_{10}^{(1)}$, and $\mathcal{G}_{01}^{(1)}$ have the
symbols
\begin{equation}
\begin{split}
G_{1}^{(1)}& =-i\partial _{k}G_{1}^{+}\partial _{Y}G_{1}^{-}, \\
G_{10}^{(1)}& =G_{1}^{+}G_{0}^{-},\quad
G_{01}^{(1)}=G_{0}^{+}G_{1}^{-},
\end{split}
\label{N10}
\end{equation}%
and satisfy the identity
\begin{equation}
\mathcal{G}_{1}^{(1)}+\mathcal{G}_{10}^{(1)}+\mathcal{G}_{01}^{(1)}-\mathcal{%
C}_{1}=0.  \label{N11}
\end{equation}%
Now from (\ref{N5}) and (\ref{N9}) we have%
\begin{equation*}
G_{1}^{+}+G_{1}^{-}=2B_{1},\text{ \ }G_{1}^{+}G_{1}^{-}=C_{2}
\end{equation*}%
which leads to%
\begin{equation*}
G_{1}^{\pm }=B_{1}\pm (B_{1}^{2}-C_{2})^{1/2}
\end{equation*}%
and since $a_{33}>0$ for $\varepsilon $ small enough, and noticing from (\ref%
{Da}) that%
\begin{equation*}
(B_{1}^{2}-C_{2})^{1/2}=\frac{1}{a_{33}}D,
\end{equation*}%
the formula (\ref{symbol1}) follows. For obtaining $G_{0}^{-},$ we
use (\ref{N11}), (\ref{N10}) which leads to \eqref{G_0formula}

It is then clear that $G_{0}^{-}$ is a pseudodifferential operator
of zero
order which satisfies (\ref{N8a}). Now, the symmetry properties (\ref{symG_j}%
) for $\mathcal{G}_{1}^{-}$ and $\mathcal{G}_{0}^{-}$ follow from (\ref%
{symbol1}), (\ref{G_0formula}), the evenness of $a,$ $a_{33},$ $a_{11},$ $%
a_{22}$ in $y_{1}$ and $y_{2},$ the oddness of $b$ and $a_{12}$ in
$y_{1}$ and $y_{2},$ and the evenness in $y_{1},$ oddness in
$y_{2}$ of $a_{13}$, and the oddness in $y_{1},$ evenness in
$y_{2}$ of $a_{23}.$

Inequality \eqref{51a} from Proposition \ref{p48} along with
\eqref{N8} implies the estimates
\begin{eqnarray}
\Vert \mathcal{R}_{01}u\Vert _{s}+\Vert \mathcal{R}_{10}u\Vert
_{s}+\Vert \mathcal{R}_{1}u\Vert _{s} &\leq &c\Vert
\mathbb{A}-\mathbb{A}_{0}\Vert _{C^{6+s}}\Vert
\mathbb{A}-\mathbb{A}_{0}\Vert _{C^{6}}\Vert u\Vert _{0}+
\notag \\
&&+c\Vert \mathbb{A}-\mathbb{A}_{0}\Vert _{C^{6}}^{2}\Vert u\Vert
_{s}. \label{N12}
\end{eqnarray}%
Substituting identities \eqref{N7a},\eqref{N9} and \eqref{N11} into %
\eqref{N5} gives the equation for the operator $\mathcal{X}$
\begin{equation}
\mathcal{X}^{2}+(-\Delta )^{1/2}\mathcal{X}+\mathcal{X}(-\Delta )^{1/2}+%
\mathcal{U}\mathcal{X}+\mathcal{X}\mathcal{V}+\mathcal{W}=0,
\label{N14}
\end{equation}%
where%
\begin{equation*}
\begin{split}
\mathcal{U}& =(\widetilde{\mathcal{G}}_{1}^{+}+\mathcal{G}_{0}^{+}),\mathcal{%
V}=-(\widetilde{\mathcal{G}}_{1}^{-}+\mathcal{G}_{0}^{-}), \\
\mathcal{W}& =-(\mathcal{R}_{1}+\mathcal{R}_{10}+\mathcal{R}_{01}+\mathcal{G}%
_{0}^{+}\mathcal{G}_{0}^{-}).
\end{split}%
\end{equation*}%
Our task is to prove the existence of a "small" solution $\mathcal{X}$ to %
\eqref{N14}.

Introduce the Banach spaces of bounded operators
\begin{equation*}
\mathbf{X}_{s}=\mathbf{F}_{s-1,s}\cap \mathbf{F}_{s,s+1},\text{ \ \ \ }%
\mathbf{Y}_{s}=\mathbf{F}_{s,s-1},\quad
\mathbf{Z}_{s}=\mathbf{F}_{s,s}
\end{equation*}%
supplemented with the norms
\begin{equation*}
\Vert \mathcal{X}\Vert _{\mathbf{X}_{s}}=\Vert \mathcal{X}\Vert _{\mathbf{F}%
_{s,s+1}}+\Vert \mathcal{X}\Vert _{\mathbf{F}_{s-1,s}},\quad \Vert \mathcal{Y%
}\Vert _{\mathbf{Y}_{s}}=\Vert \mathcal{Y}\Vert
_{\mathbf{F}_{s,s-1}},\quad
\Vert \mathcal{Z}\Vert _{\mathbf{Z}_{s}}=\Vert \mathcal{Z}\Vert _{\mathbf{F}%
_{s,s}}.
\end{equation*}%
It follows from \eqref{N8}, Corollary \ref{matrix}, and
\eqref{N12}that for all $r\leq \rho -8$,
\begin{equation}
\Vert \mathcal{U}\Vert _{\mathbf{Y}_{r+1}}+\Vert \mathcal{V}\Vert _{\mathbf{Y%
}_{r}}+\Vert \mathcal{W}\Vert _{\mathbf{Z}_{s}}\leq c\varepsilon ,
\label{N15}
\end{equation}%
where the constant $c$ depends on $l$ and $\rho $ only. The rest
of the proof is based on the following lemma which is proved at
the end of this section.

\begin{lemma}
\label{N17l} Under the above assumptions, there exist $\varepsilon
_{0}>0$ and $c>0$ depending on $\rho $ and $l$ only such that for
$\varepsilon \in (0,\varepsilon _{0})$ equation \eqref{N14} has a
solution satisfying the inequalities
\begin{gather}
\Vert \mathcal{X}\Vert _{\mathbf{X}_{r}}\leq c\varepsilon \text{~~ when~~}%
0\leq r\leq \rho -8,  \notag \\
\Vert \mathcal{X}u\Vert _{s}\leq c(\varepsilon \Vert u\Vert
_{s-1}+E_{l}\Vert u\Vert _{0})\text{~~when~~}1\leq s\leq l-8.
\label{N18}
\end{gather}%
Moreover the operator $\mathcal{X}$ satisfies the following symmetry property%
\begin{equation*}
\mathcal{X}u(\pm Y^{\ast })=\mathcal{X}u^{\ast }(\pm Y).
\end{equation*}
\end{lemma}
\begin{proof} We start with the consideration of the simple linear
operator equation
\begin{equation}
(-\Delta )^{1/2}\mathcal{X}+\mathcal{X}(-\Delta )^{1/2}=\mathcal{Z}\text{%
~~with~~}\mathcal{Z}\in \mathbf{Z}_{s}.  \label{30a}
\end{equation}%
It is easy to see that for any $\mathcal{Z}\in \mathbf{Z}_{s}$ satisfying %
\eqref{30a} and having the matrix form with elements
$\mathcal{Z}_{kp}$ , operator $\mathcal{X}$ has the matrix
representation with elements
\begin{equation*}
\mathcal{X}_{kp}=\frac{1}{N(k)+N(p)}\mathcal{Z}_{kp},\text{~~where~~}%
N(k)=(k_{1}^{2}+\tau ^{2}k_{2}^{2})^{1/2}
\end{equation*}%
which obviously yields the estimate
\begin{equation}
\Vert \mathcal{X}\Vert _{\mathbf{X}_{s}}\leq c(\tau )\Vert
\mathcal{Z}\Vert _{\mathbf{Z}_{s}}.  \label{N31}
\end{equation}%
Hence the mapping $\mathcal{Z}\mapsto \mathcal{X}$ defines a
bounded linear operator $\Xi \in
\mathcal{L}(\mathbf{Z}_{s},\mathbf{X}_{s})$. Let us consider the
sequence of operators $\mathcal{X}_{n}$ defined by the equalities
\begin{equation*}
\mathcal{X}_{0}=0,\quad \mathcal{X}_{n+1}=-\Xi \Big(\mathcal{X}_{n}^{2}+%
\mathcal{U}\mathcal{X}_{n}+\mathcal{X}_{n}\mathcal{V}+\mathcal{W}\Big).
\end{equation*}%
Note that $\mathcal{U}$ is a pseudodifferential operator which symbol $%
U(Y,k) $ satisfying the inequalities
\begin{equation*}
\Vert U(\cdot ,k)\Vert _{C^{r}}\leq c\varepsilon ,\quad |\widehat{U}%
(p,q)|\leq c\varepsilon (1+|p|)^{-\rho +1}|q|.
\end{equation*}%
It is easy to see that for any $\mathcal{X}\in \mathbf{X}_{r}$,
the operator $\mathcal{U}\mathcal{X}$ has a matrix representation
with the matrix elements
\begin{equation*}
(\mathcal{U}\mathcal{X})_{kp}=\sum\limits_{q\in \mathbb{Z}^{2}}\widehat{U}%
(k-q,q)\mathcal{X}_{qp}.
\end{equation*}%
We have
\begin{equation*}
(1+|k|)^{r}||(\mathcal{U}\mathcal{X})_{kp}|\leq c\varepsilon
\sum\limits_{q}(1+|k-q|)^{r-\rho
+1}|\mathcal{X}_{qp}|(1+|q|)^{r+1},
\end{equation*}%
which gives
\begin{gather*}
\sum\limits_{k}(1+|k|)^{2r}\Big[\sum\limits_{p}|(\mathcal{U}\mathcal{X}%
)_{kp}||\widehat{u}(p)|\Big]^{2}\leq \\
c\varepsilon ^{2}\sum\limits_{k}\Big[\sum\limits_{p}\sum%
\limits_{q}(1+|q-k|)^{r-\rho +1}(1+|q|)^{r+1}|\mathcal{X}_{qp}||\widehat{u}%
(p)|\Big]^{2}\leq \\
c\varepsilon ^{2}\sum\limits_{k}\Big[\sum\limits_{q}(1+|q-k|)^{r-\rho +1}%
\Big((1+|q|)^{r+1}\sum\limits_{p}|\mathcal{X}_{qp}\Vert \widehat{u}(p)|\Big)%
\Big]^{2}\leq \\
c\varepsilon ^{2}\sum\limits_{q}\Big[(1+|q|)^{r+1}\sum\limits_{p}|\mathcal{X}%
_{qp}||\widehat{u}(p)|\Big]^{2}\leq c\varepsilon ^{2}\Vert
\mathcal{X}\Vert
_{\mathbf{F}_{r,r+1}}^{2}\sum\limits_{p}(1+|p|)^{2r}|\widehat{u}(p)|^{2}.
\end{gather*}%
Thus we get
\begin{equation*}
\Vert \mathcal{U}\mathcal{X}\Vert _{\mathbf{Z}_{r}}\leq c\Vert \mathcal{U}%
\Vert _{\mathbf{Y}_{r+1}}\Vert \mathcal{X}\Vert _{\mathbf{X}_{r}}.
\end{equation*}%
Repeating these arguments we obtain
\begin{equation*}
\Vert \mathcal{X}\mathcal{V}\Vert _{\mathbf{Z}_{r}}\leq c\Vert \mathcal{V}%
\Vert _{\mathbf{Y}_{r}}\Vert \mathcal{X}\Vert _{\mathbf{X}_{r}},
\end{equation*}%
and
\begin{equation*}
\Vert \mathcal{X}^{2}\Vert _{\mathbf{Z}_{r}}\leq c\Vert \mathcal{X}\Vert _{%
\mathbf{Z}_{r}}^{2}\leq c\Vert \mathcal{X}\Vert
_{\mathbf{X}_{r}}^{2},
\end{equation*}%
inequality \eqref{N31} yields the estimates
\begin{equation}
\Vert \mathcal{X}_{n+1}\Vert _{\mathbf{X}_{r}}\leq c\varepsilon
(\Vert
\mathcal{X}_{n}\Vert _{\mathbf{X}_{r}}+1)+\Vert \mathcal{X}_{n}\Vert _{%
\mathbf{X}_{r}}^{2},  \label{N33}
\end{equation}%
which holds true for all $r\in (1,\rho -8)$. On the other hand,
since
\begin{equation*}
\mathcal{X}_{n+1}^{2}-\mathcal{X}_{n}^{2}=(\mathcal{X}_{n+1}-\mathcal{X}_{n})%
\mathcal{X}_{n+1}+\mathcal{X}_{n}(\mathcal{X}_{n+1}-\mathcal{X}_{n})
\end{equation*}%
we have
\begin{equation}
\Vert \mathcal{X}_{n+1}-\mathcal{X}_{n}\Vert _{\mathbf{X}_{r}}\leq
c\Vert \mathcal{X}_{n}-\mathcal{X}_{n-1}\Vert
_{\mathbf{X}_{r}}(\varepsilon +\Vert
\mathcal{X}_{n}\Vert _{\mathbf{X}_{r}}+\Vert \mathcal{X}_{n-1}\Vert _{%
\mathbf{X}_{r}}).  \label{N34}
\end{equation}%
Here the constant $c$ depends on $\varrho $ and $\tau $ only. It
follows from \eqref{N33} that for all $\varepsilon \in
(0,\varepsilon _{0}(\rho ,\tau )$,\ the values $\Vert
\mathcal{X}_{n}\Vert _{\mathbf{X}_{r}}$ are less than
$c\varepsilon $. From this and \eqref{N34} we conclude that for all small $%
\varepsilon $ the sequence $\mathcal{X}_{n}$ converges in
$\mathbf{X}_{r}$. Repeating these arguments and using Corollary
\ref{matrix} gives the tame estimate \eqref{N18},
and the lemma follows once we observe that the symmetry property of $%
\mathcal{X}$ follows from the uniqueness of $\mathcal{X}$ and from
the
equivariance with respect to the required symmetry of the equation (\ref{N14}%
).\end{proof}

In order to complete the proof of Theorem  \ref{thmN1}, it remains to note that operator $\mathcal{G}_{-1}^{-}:=-\mathcal{X}$ with $%
\mathcal{X}$ given by Lemma \ref{N17l} satisfies \eqref{symbol-1}.
\end{proof}

\textbf{ Proof of Theorem \ref{thmDir-Neum}}
 It follows from formulae (\ref{symbol1}),  (\ref{N5a}), and (\ref%
{N4}) that we can write%
\begin{equation*}
\mathcal{G}_{\eta }\psi =\frac{1}{J}\{\mathcal{D}-a_{33}(\mathcal{G}%
_{0}^{-}+\mathcal{G}_{-1}^{-})\}\widetilde{\psi },
\end{equation*}%
where the first order pseudodifferential operator $\mathcal{D}$
has the
symbol $D(Y,k).$ Then we define%
\begin{equation*}
\mathcal{G}_{0}=-\frac{a_{33}}{J}\mathcal{G}_{0}^{-},\text{ \ }\mathcal{G}%
_{-1}=-\frac{a_{33}}{J}\mathcal{G}_{-1}^{-},
\end{equation*}%
and the symmetry properties follow from the evenness of $J,$
$a_{33}$ and from Theorem \ref{thmN1} and Lemma \ref{N17l}. The
zero order pseudodifferential operator $\mathcal{G}_{0}$ satisfies
(\ref{estimG0}) and
from Proposition \ref{boundedness} we have%
\begin{eqnarray}
||\mathcal{G}_{0}u||_{r} &\leq &c\varepsilon ||u||_{r},\text{ \ \
}0\leq
r\leq \rho -4  \notag \\
||\mathcal{G}_{0}u||_{s} &\leq &c(\varepsilon ||u||_{s}+E_{l}||u||_{0}),%
\text{ \ \ }0\leq s\leq l-4,  \label{estimG_0}
\end{eqnarray}%
while the operator $\mathcal{G}_{-1}$ satisfies
\begin{eqnarray}
||\mathcal{G}_{-1}u||_{r} &\leq &c\varepsilon \Vert u\Vert _{r-1},\text{ \ }%
1\leq r\leq \rho -8,  \notag \\
||\mathcal{G}_{-1}u||_{s} &\leq &c(\varepsilon \Vert u\Vert
_{s-1}+E_{l}\Vert u\Vert _{0}),\text{ \ }1\leq s\leq l-8.
\label{estimG1}
\end{eqnarray}%
We then deduce the estimates (\ref{estimG-1}) in using
(\ref{estimA}).

Now our task is to calculate the symbols  $G_j$. It is convenient
to introduce the scalar $\mathbf I(Y)$ and  linear form $\Pi(Y,k)$
defined by
$$
\mathbf{I}=J/\sqrt{\text{det~}\mathbb G},\quad \Pi=\mathbb
G^{-1}\nabla_Y\tilde \eta \cdot k.
$$
Recall the identity
\begin{eqnarray*}
&&a_{33}\sum\limits_{1\leq j,m\leq 2}a_{jm}k _{j}k _{m}-\Bigl(%
\sum\limits_{1\leq m\leq 2}a_{3m}k _{m}\Bigr)^{2} \\
&=&\det \mathbb{A}\Bigl((\mathbb{A}^{-1})_{22}k _{1}^{2}-2(\mathbb{A}%
^{-1})_{12}k _{1}\xi _{2}+(\mathbb{A}^{-1})_{11}k _{2}^{2}\Bigr).
\end{eqnarray*}%
Noting that $\det \mathbb{A}=J$,
$\mathbb{A}^{-1}=J^{-1}\mathbb{B}_{1}^{\ast }\mathbb{B}_{1}$, we
conclude from this that
\begin{equation*}
a_{33}\sum\limits_{1\leq j,m\leq 2}a_{jm}k_{j}k_{m}-\Bigl(%
\sum\limits_{1\leq m\leq 2}a_{3m}k_{m}\Bigr)^{2}=(\mathbb{B}_{1}^{\ast }%
\mathbb{B}_{1})_{22}k_{1}^{2}-2(\mathbb{B}_{1}^{\ast }\mathbb{B}%
_{1})_{12}k_{1}k_{2}+(\mathbb{B}_{1}^{\ast
}\mathbb{B}_{1})_{11}k_{2}^{2}.
\end{equation*}%
Hence, by the definition of the metric tensor $\mathbb{G}$,
\begin{equation*}
a_{33}\sum\limits_{1\leq j,m\leq 2}a_{jm}k _{j}\xi _{m}-\Bigl(%
\sum\limits_{1\leq m\leq 2}a_{3m}k_{m}\Bigr)^{2}=g_{22}k
_{1}^{2}-2g_{12}k _{1}k_{2}+g_{11}k_{2}^{2},
\end{equation*}%
which yields
\begin{equation}\label{dformula1}
a_{33}\sum\limits_{1\leq j,m\leq 2}a_{jm}k_{j}k_{m}-\Bigl(%
\sum\limits_{1\leq m\leq 2}a_{3m}k_{m}\Bigr)^{2}=(\det
\mathbb{G})\mathbb{G}^{-1}k\cdot k.
\end{equation}%
Substituting this relation into $D$ (see (\ref{Da})) finally gives
\begin{equation}
D(Y,k)=\sqrt{\det \mathbb{G}}\,\mathbf{G}_1(Y,k). \label{symbol6}
\end{equation}%
Noting that $\mathcal{G}_{1}=\frac{1}{J}\mathcal{D}$ we obtain the
needed formula \eqref{defG_1}. The calculation of $G_0$ is more
delicate task. Since $a_{33}=\text{det~}\mathbb G/J$, formula
\eqref{G_0formula} yields
\begin{equation*}
(2\mathbf I \mathbf G_1) G^-_0=i\nabla
_{k}G_{1}^{+}\nabla_{Y}G_{1}^{-}-
\frac{G^-_1}{a_{33}}\sum\limits_{j=1}^2\partial_{y_j}a_{3j}+
\frac{i}{a_{33}}\sum\limits_{j,m=1}^2\partial_{y_j}a_{jm}k_m.
\end{equation*}
It follows from the definition of the form   $\Pi$ that
\begin{equation}\label{divpi}
\sum\limits_{m=1}^2 a_{3m}k_m=-a_{33}\Pi,\quad \sum\limits_{j=1}^2
\partial_{y_j}a_{3j}=-\text{~div}_Y(a_{33}\nabla_k\Pi).
\end{equation}
From this, \eqref{symbol1}, and \eqref{symbol6} we conclude that
$$
G^-_1=-\mathbf I\mathbf G_1-i\Pi, \quad G^+_1=\mathbf I\mathbf
G_1-i\Pi,
$$
and hence
\begin{multline}\label{gminus1}
(2\mathbf I \mathbf G_1) G^-_0=i\nabla _{k} (\mathbf I\mathbf
G_1-i\Pi)\nabla_{Y}(-\mathbf I\mathbf G_1-i\Pi)
-\frac{1}{a_{33}}(\mathbf I\mathbf
G_1+i\Pi)\text{~div}_Y(a_{33}\nabla_k\Pi)+\\
\frac{i}{a_{33}}\sum\limits_{j,m=1}^2\partial_{y_j}a_{jm}k_m.
\end{multline}
Next differentiating both sides of \eqref{dformula1} with respect
to $k_j$ we arrive to
 $$
2a_{33}\sum\limits_{m=1}^2a_{jm}k_m-2a_{3j}\sum\limits_{m=1}^2a_{3m}k_m
=(\text{det~} \mathbb G)\partial_{k_j}\mathbf G_1^2,
$$
which along with \eqref{divpi} and the identity $\text{det~}
\mathbb G/a_{33}=J$ leads to
$$
\sum\limits_{m=1}^2a_{jm}k_m=a_{33}\Pi\partial_{k_j}\Pi+J\mathbf
G_1\partial_{k_j}\mathbf{G}_1.
$$
Substituting this expression into \eqref{gminus1} we finally
obtain
\begin{multline}\label{gminus2}
(2\mathbf I \mathbf G_1) G^-_0=i\nabla _{k} (\mathbf I\mathbf
G_1-i\Pi)\nabla_{Y}(-\mathbf I\mathbf G_1-i\Pi)
-\frac{1}{a_{33}}(\mathbf I\mathbf
G_1+i\Pi)\text{~div}_Y(a_{33}\nabla_k\Pi)+\\
\frac{i}{a_{33}}\text{div}_Y\Big(a_{33}\Pi\nabla_k\Pi+J\mathbf
G_1\nabla_k \mathbf G_1\Big).
\end{multline}
Let us calculate the real part of $G_0$. It is easy to see that
$$
(2\mathbf I \mathbf G_1) \text{~Re~}G^-_0=\mathbf I\nabla_k\mathbf
G_1\nabla_Y\Pi-\nabla_k\Pi\nabla_Y(\mathbf I\mathbf
G_1)-\frac{\mathbf I\mathbf
G_1}{a_{33}}\text{div}_Y(a_{33}\nabla_k\Pi),
$$
which along with  the equality $\mathbf I
a_{33}=\sqrt{\text{det~}\mathbb G}$ implies
$$
(2 \mathbf G_1) \text{~Re~}G^-_0= \nabla_k\mathbf
G_1\cdot\nabla_Y\Pi-\nabla_k\Pi\cdot\nabla_Y\mathbf G_1-\mathbf
G_1\text{~div}_Y(\nabla_k\Pi)-\mathbf
G_1\nabla_k\Pi\cdot\nabla_Y\ln \sqrt{\text{det~}\mathbb G}.
$$
Noting that
$$
\text{div}_Y(\nabla_k\Pi)+\nabla_k\Pi\cdot\nabla_Y\ln
\sqrt{\text{det~}\mathbb G}=\textbf{div}(\nabla_k\Pi),
$$
and recalling \begin{equation}\label{gtogminus}
G_0=-\frac{a_{33}}{J}G^-_0=-\frac{\text{det~}\mathbb G}{J^2}G^-_0,
\end{equation}
we obtain
$$
\text{Re~}G_0=\frac{\text{det~}\mathbb G}{2J^2\mathbf G_1} (
\nabla_k\Pi\cdot\nabla_Y\mathbf G_1-\nabla_k\mathbf
G_1\cdot\nabla_Y\Pi)+\frac{\text{det~}\mathbb
G}{2J^2}\textbf{div}(\nabla_k\Pi).
$$
From this and the identities
$$
\nabla_k\Pi=\mathbb G^{-1}\nabla_Y\tilde \eta, \quad \nabla_k
\mathbf G_1=\frac{1}{\mathbf G_1}\mathbb G^{-1}k,\quad \nabla_Y
\mathbf G_1=\frac{1}{2\mathbf G_1}\nabla_Y(\mathbb G^{-1}k\cdot k)
$$
we obtain the desired formula \eqref{cor11} for the real part of
$G_0$. Next \eqref{gminus2} yields
\begin{gather*}
(2\mathbf I \mathbf G_1) \text{~Im~}G^-_0=-\mathbf
I\nabla_k\mathbf G_1\cdot\nabla_Y(\mathbf I\mathbf G_1)
-\nabla_k\Pi\cdot\nabla_Y\Pi -\\
\frac{1}{a_{33}}\Pi\text{div}_Y(a_{33}\nabla_k\Pi)+\frac{1}{a_{33}}
\text{div}_Y\big[\Pi\nabla_k(a_{33}\Pi)+J\mathbf
G_1\nabla_k\mathbf G_1\big]=\\\frac{J}{a_{33}}\mathbf
G_1\text{div}_Y(\nabla_k\mathbf G_1)+\mathbf
G_1\Big(\frac{1}{a_{33}}\nabla_Y J-\mathbf I \nabla_Y \mathbf
I\Big)\cdot\nabla_k\mathbf G_1.
\end{gather*}
Since $ {J}/{a_{33}}=\mathbf I^2$ and
$$
\frac{1}{a_{33}}\nabla_Y J-\mathbf I \nabla_Y \mathbf
I=\frac{J^2}{(\text{det~} \mathbb
G)^{3/2}}\nabla_Y\sqrt{\text{det~}\mathbb G}=\mathbf I^2
\frac{1}{\sqrt{\text{det~} \mathbb
G}}\nabla_Y\sqrt{\text{det~}\mathbb G},
$$
we have
$$
2 \text{~Im~}G^-_0=\mathbf I\textbf {~div~}(\nabla_k \mathbf G_1).
$$
Recalling \eqref{gtogminus} we obtain \eqref{cor12} and the
theorem follows.

\textbf{Invariant form of} $\text{Re}\, G_0$. In the rest of the
section we prove the formula \eqref{geometry0}.  We start with
the calculation of the quadratic form
$$
Q(Y, \mathbb G \xi):=\tilde Q(Y,\xi)=\tilde Q_{11}\xi_1^2+2\tilde
Q_{12}\xi_1\xi_2+\tilde Q_{22}\xi_2^2.
$$
It follows from \eqref{cor13} and the identity $\mathbb
G\partial_{y_j}\mathbb G^{-1}\mathbb G=-\partial_{y_{j}}\mathbb G$
that
\begin{equation}\label{tildeq}
\tilde Q(Y,\xi)=\sum\limits_{j=1}^2(\xi_j\partial_{y_j}\mathbb
G)\mathbf q\cdot\xi-\frac{1}{2}\partial_{y_j}(\mathbb
G\xi\cdot\xi)q_j-\nabla_Y^2\tilde\eta \xi\cdot\xi,
\end{equation}
where the vector field
$
\mathbf q=\mathbb G^{-1}\nabla_Y\tilde \eta.
$
Thus we get
$$
\tilde
Q_{\alpha\beta}=\frac{1}{2}\sum\limits_{j=1}^2\Big(\partial_{y_\alpha}
g_{j\beta}+\partial_{y_\beta}
g_{j\alpha}-\partial_{y_j}
g_{\alpha\beta}\Big)q_j-\partial^2_{\alpha\beta}\tilde\eta,
$$
which along with the equality $
g_{\alpha\beta}=\partial_{y_\alpha}\mathbf r\cdot
\partial_{y_\beta}\mathbf r$  yields
\begin{equation}\label{qalbe}
\tilde
Q_{\alpha\beta}=\big(\sum\limits_{j=1}^2q_j\partial_{y_j}\mathbf
r-\mathbf e_3\Big)\cdot \partial^2_{y_{\alpha}y_{\beta}}\mathbf r.
\end{equation}
On the other hand, since  $\tilde\eta=\mathbf r\cdot\mathbf e_3$,
the expression for $\mathbf q$ reads
$$
\mathbf q=\frac{1}{\text{det~}\mathbb G}\left\{
\begin{array}{c}
( g_{22}\partial_{y_1}\mathbf r- g_{12}\partial_{y_2}\mathbf r)\cdot\mathbf e_3 \\
(-g_{12}\partial_{y_1}\mathbf r+ g
g_{11}\partial_{y_2}\mathbf r)\cdot\mathbf e_3
\end{array}%
\right\}.
$$
Now set
\begin{equation}\label{vectory}
    \mathbf a=\partial_{y_1}\mathbf{r},\quad \mathbf b=\partial_{y_2}\mathbf{r},
    \quad \mathbf c=\mathbf a\times\mathbf b.
\end{equation}
Noting that
\begin{gather*}
g_{11}=\mathbf a\cdot\mathbf a, \quad 
g_{12}=\mathbf a\cdot\mathbf b,\quad g_{22}=\mathbf
b\cdot\mathbf b,\\
\mathbf b\times \mathbf c=(\mathbf b\cdot\mathbf b)\mathbf a-
(\mathbf a\cdot\mathbf b)\mathbf b,\quad  \mathbf a\times \mathbf
c=(\mathbf a\cdot\mathbf b)\mathbf a- (\mathbf a\cdot\mathbf
a)\mathbf b,
\end{gather*}
we obtain
\begin{equation*}\label{Gq12}
(\text{det~}\mathbb G)q_1=(\mathbf b\times \mathbf c)\cdot \mathbf
e_3,\quad (\text{det~}\mathbb G)q_2=-(\mathbf a\times \mathbf
c)\cdot\mathbf e_3.
\end{equation*}
From this  and the identity
$$
[(\mathbf b\times\mathbf c)\cdot \mathbf e_3] \mathbf a- [(\mathbf
a\times\mathbf c)\cdot \mathbf e_3] \mathbf b-|\mathbf c|^2\mathbf
e_3=-(\mathbf c\cdot\mathbf e_3)\mathbf c,
$$
which holds true  for all $\mathbf a, \mathbf b\in \RR^3$ and
$\mathbf c=\mathbf a\times \mathbf b$, we conclude that
$$
(\text{det~}\mathbb G)(q_1\partial_{y_1}\mathbf
r+q_2\partial_{y_2}\mathbf r)-|\mathbf c|^2\mathbf e_3=-(\mathbf
c\cdot\mathbf e_3)\mathbf c.
$$
Noting that $|\mathbf c|^2=\text{det~}\mathbb G$ and $\mathbf
c\cdot\mathbf e_3=J$ we arrive at
$$
\sum\limits_{j=1}^2q_j\partial_{y_j}\mathbf r-\mathbf
e_3=-\frac{J}{\sqrt{\text{det~}\mathbb G}}\mathbf n,
$$
where $\mathbf n=\mathbf c/|\mathbf c|$ is the  unit normal
vector to $\Sigma$. Substituting this identity  into \eqref{qalbe}
gives $ \tilde Q_{\alpha\beta}=-(\mathbf n\cdot
\partial_{y_\alpha y_\beta}\mathbf r)J/\sqrt{\text{det~}\mathbb G}$
which leads to
$$
\tilde Q(Y,\xi)=-\frac{J}{\sqrt{\text{det~}\mathbb G}}(L\xi_1^2+2
M\xi_1\xi_2+N\xi_2^2).
$$
Since
$$
\mathbf G_1^2(Y,\mathbb G\xi)=\mathbb G\xi\cdot\xi:= E\xi_1^2+2
F\xi_1\xi_2+G\xi_2^2,
$$
we finally obtain
\begin{equation}\label{halfg0}
    \frac{\text{det~}\mathbb G}{2J^2}\left\{\frac{1}{\mathbf G_1^2(Y,\mathbb
    G\xi)}Q(Y,\mathbb G\xi)\right\}=-\frac{\sqrt{\text{det~}\mathbb G}}{J}
\frac{L\xi_1^2+2 M\xi_1\xi_2+N\xi_2^2}{2(E\xi_1^2+2
F\xi_1\xi_2+G\xi_2^2)}.
\end{equation}
Our next task is to express $ \text{\bf div}\,\mathbf q$ via
the geometric characteristics of $\Sigma$. First we do this in the
standard coordinates $Y=X$ with $\tilde\eta=\eta$. In this case
$$
\mathbb G^{-1}=\frac{1}{1+|\nabla\eta|^2} \left(
\begin{array}{cc}
1+\partial_{y_2}\eta^2,&-\partial_{y_1}\eta\ \partial_{y_2}\eta \\
-\partial_{y_1}\eta\ \partial_{y_2}\eta,&1+\partial_{y_1}\eta^2
\end{array}%
\right),\quad \text{det~}\mathbb G=1+|\nabla\eta|^2,\quad J=1,
$$
and $\mathbf q=(1+|\nabla\eta|^2)^{-1}\nabla \eta$, which  leads
to the formula \begin{equation}\label{meancurvature}
\frac{\text{det~} \mathbb G}{2J^2}\mathbf{div}\ \mathbf q=
\frac{\sqrt{\text{det~}\mathbb
G}}{J}\frac{1}{2}\text{div~}\Big(\frac{\nabla\eta}{\sqrt{1+|\nabla\eta|^2}}\Big)=
\frac{\sqrt{\text{det~}\mathbb G}}{J}\frac{\mathfrak k_1+\mathfrak
k_2}{2},
\end{equation}
where $\mathfrak k_i$ are the principal curvatures of $\Sigma$ at
the point $\mathbf r$. Next note that $\nabla\eta$ is a covariant
vector field on $\Sigma$, hence $\mathbb G^{-1}\nabla\eta$ is
a vector field on $\Sigma$. Since $\mathbf{div}$ is an
invariant operator on the space of vector fields on
$\Sigma$, the left side of \eqref{meancurvature} does not depend
on the choice of coordinates, hence
\begin{equation}\label{secondhalf}
\frac{\text{det~} \mathbb G}{2J^2}\mathbf{div}\ \mathbf
q=\frac{\sqrt{\text{det~}\mathbb G}}{J}\frac{\mathfrak
k_1+\mathfrak k_2}{2}=\frac{\sqrt{\text{det~}\mathbb
G}}{J}\frac{LG-2MF+NE}{2(EG-F^2)}.
\end{equation}
Combining \eqref{halfg0} and \eqref{secondhalf} gives the desired
identity \eqref{geometry0}. If we define by $\mathfrak n(\xi)$ the
normal curvature of $\Sigma$ in the direction $\xi$ at a point
$\mathbf r$, then \eqref{geometry0} becomes
$$
\frac{J}{\sqrt{\text{det~}\mathbb G}}\text{~Re~}G_0(Y, \mathbb
G\xi)= \frac{1}{2}(\mathfrak k_1+\mathfrak k_2-\mathfrak n(\xi)),
$$
which leads to
\begin{corollary} Assume that the manifold $\Sigma$
has a parametric representation $\mathbf r(Y)=(\mathbb
TY+\tilde{\mathcal V}(Y),\tilde \eta(Y))$,  $Y\in\RR^2$ so that
$\tilde{\mathcal V}$ and $\tilde \eta$, which are not defined yet,
satisfy all hypotheses of Theorem \ref{thmDir-Neum}. Assume also
that   the parametric form
$$
\mathfrak Gu:=\frac{J}{\sqrt{\text{\rm det~}\mathbb G}}\mathcal
G_\eta\check u\circ (\mathbb T+\mathcal V), \quad \check
u(X)=u(Y(X)),
$$
of the normal derivative operator is given for any bi- periodic
smooth function $u(Y)$. Then the manifold $\Sigma$ is defined by the
operator $\mathfrak G$ up to a translation and a rotation  of the
embedding space.
\end{corollary}
\begin{proof} Note that for all $k\in \ZZ^2$,
\begin{gather*}
\lim\limits_{n\to\infty}\frac{1}{n}e^{-i n k\cdot Y}\mathfrak
Ge^{in k\cdot Y}=\frac{J}{\sqrt{\text{\rm det~}\mathbb G}}G_1(Y,k)
=\sqrt{\mathbb G^{-1}k\cdot k},\\
\lim\limits_{n\to\infty}\text{~Re~}\big\{e^{-i n k\cdot
Y}\mathfrak Ge^{in k\cdot Y}-\frac{J}{\sqrt{\text{\rm det~}\mathbb
G}}G_1(Y,k)\Big\} =\frac{J}{\sqrt{\text{\rm det~}\mathbb
G}}\text{~Re~} G_{0}(Y,k).
\end{gather*}
Since $G_0$ is a homogeneous function of $k$,  it follows from
this that the right hand sides of these equalities are defined by
the operator  $\mathfrak G$ for all $k\in\RR^2$ and, in
particular, for  $k=\mathbb G\xi$ with an arbitrary $\xi\in\RR^2$.
Hence the first fundamental form $\mathbb G\xi\cdot \xi$ and the
difference $\mathfrak k_1+\mathfrak k_2-\mathfrak n(\xi)$ are
completely defined by  $\mathfrak G$ for all directions $\xi$ at
each point of $\Sigma$. Hence the principle curvatures of $\Sigma$
are also defined by the operator $\mathfrak G$. It remains to note
that, by the Bonnet Theorem, the first fundamental form and the
principal curvatures define $\Sigma$ up to a translation and
a rotation  of the embedding space.
\end{proof}

\section{Proof of Lemma \protect\ref{coef@}\label{app@}}

To be able to compute all terms in (\ref{@}), let us rewrite the system (\ref%
{basic1}), (\ref{basic2}) formally as a scalar equation for $\psi
$ and express the orders $\varepsilon $ and $\varepsilon ^{2}$ of
the differential
computed at $\psi _{\varepsilon }^{(N)},$ $N\geq 3.$ Then the operator $%
\mathfrak{L+H}$ is up to order $\varepsilon ^{2}$ closely linked
with the new form of this operator after applying the
diffeomorphism \ computed at Lemma \ref{Lemmcoordchange}.

Indeed we can formally solve (\ref{basic2}) with respect to $\eta
$ in
powers of $\psi $ as%
\begin{equation}
\eta =-\frac{1}{\mu }\partial _{x_{1}}\psi -\frac{1}{2\mu }(\nabla
\psi )^{2}+\frac{1}{2\mu ^{3}}(\partial _{x_{1}}^{2}\psi
)^{2}+O(||\psi ||^{3}), \label{formaleta}
\end{equation}%
and replace $\eta $ by this expression in (\ref{basic1}). We then
obtain a
new scalar \emph{formal equation} for $\psi ,$ under the form $\mathcal{E}%
(\psi ,\mu )=0,$ where%
\begin{equation*}
\mathcal{E}(\psi ,\mu )=\mathfrak{L}_{0}\psi +(\mu -\mu _{c})\mathfrak{L}%
_{1}\psi +\mathcal{E}_{2}(\psi ,\psi )+\mathcal{E}_{3}(\psi ,\psi
,\psi )+O(|\mu -\mu _{c}|^{2}||\psi ||+||\psi ||^{4}),
\end{equation*}%
where%
\begin{eqnarray*}
\mathfrak{L}_{0}\psi &=&\frac{1}{\mu _{c}}\partial _{x_{1}}^{2}\psi +%
\mathcal{G}^{(0)}\psi \\
\mathfrak{L}_{1}\psi &=&-\frac{1}{\mu _{c}^{2}}\partial
_{x_{1}}^{2}\psi ,
\end{eqnarray*}%
and $\mathcal{E}_{2}$ and $\mathcal{E}_{3}$ represent quadratic
and cubic terms in $\psi $, respectively. Let us write the formal
solution found at Theorem \ref{Lembifurc} for $\varepsilon
_{1}=\varepsilon _{2}=\varepsilon
/2,$ under the form ($N\geq 3)$%
\begin{eqnarray*}
\psi _{\varepsilon }^{(N)} &=&\varepsilon \psi _{1}+\varepsilon
^{2}\psi
_{2}+O(\varepsilon ^{3}), \\
\mu &=&\mu _{c}+\varepsilon ^{2}\mu _{1}+O(\varepsilon ^{3}),
\end{eqnarray*}%
where%
\begin{equation*}
\psi _{1}=\sin x_{1}\cos \tau x_{2},
\end{equation*}%
then we have the identities%
\begin{eqnarray}
\mathfrak{L}_{0}\psi _{1} &=&0,  \notag \\
\mathfrak{L}_{0}\psi _{2}+\mathcal{E}_{2}(\psi _{1},\psi _{1})
&=&0,
\label{ident} \\
\mathfrak{L}_{0}\psi _{3}+\mu _{1}\mathfrak{L}_{1}\psi _{1}+2\mathcal{E}%
_{2}(\psi _{1},\psi _{2})+\mathcal{E}_{3}(\psi _{1},\psi _{1},\psi
_{1}) &=&0.  \notag
\end{eqnarray}%
Now, we may observe that the operator (\ref{linequa}) we want to
invert acts on $\delta \phi =\delta \psi -\mathfrak{b}\delta \eta
.$ Since $\mathfrak{b}$ is $O(\varepsilon )$ and $\delta \eta $
may be expressed formally linearly in terms of $\delta \psi $ in
differentiating formally (\ref{formaleta}), we have formally
\begin{equation}
\partial _{\psi }\mathcal{E}(\psi ,\mu )(1+\mathcal{H}(\psi ,\mu ))=\kappa
\mathfrak{L}(\psi ,\mu )-\mathcal{R}(\mathcal{E},\mu )
\label{d_psi E}
\end{equation}%
where $\mathfrak{L}(\psi ,\mu )$ is the linear operator we want to invert, $%
\kappa $ is the function we introduced at Theorem
\ref{thmChangeCoord}, and
the operator $\mathcal{H}(\psi ,\mu )$ is such that formally%
\begin{equation*}
\delta \psi =(1+\mathcal{H}(\psi ,\eta ))\delta \phi ,
\end{equation*}%
and $\mathcal{R}(0,\mu )=0.$ Since we set%
\begin{eqnarray*}
\psi &=&\psi _{\varepsilon }^{(N)}+O(\varepsilon ^{N}),\text{ }N\geq 3 \\
\mu &=&\mu _{0}+\varepsilon ^{2}\mu _{1}(\tau _{0})+O(\varepsilon ^{3}),%
\text{ \ }\mu _{0}=\mu _{c}(\tau _{0}),
\end{eqnarray*}%
we have%
\begin{equation*}
\mathcal{E}(\psi _{\varepsilon }^{(2)},\mu _{0}+\varepsilon
^{2}\mu _{1})=O(\varepsilon ^{3}),\text{ \ \ \
}\mathcal{R}=O(\varepsilon ^{3}),
\end{equation*}%
hence
\begin{equation*}
\partial _{\psi }\mathcal{E}(\psi ,\mu )=\mathfrak{L}_{0}+\varepsilon
^{2}\mu _{1}\mathfrak{L}_{1}+2\mathcal{E}_{2}(\psi _{\varepsilon
}^{(2)},\cdot )+3\mathcal{E}_{3}(\psi _{\varepsilon }^{(2)},\psi
_{\varepsilon }^{(2)},\cdot )+O(\varepsilon ^{3}).
\end{equation*}

Making now the change of coordinates computed at Lemma \ref{Lemmcoordchange}%
, the new expressions of operators $\mathfrak{L}_{0}$, $\mathfrak{L}_{1},$ $%
\mathcal{E}_{2}(\psi _{\varepsilon }^{(2)},\cdot ),$
$\mathcal{E}_{3}(\psi
_{1},\psi _{1},\cdot )$ take the following form%
\begin{eqnarray*}
\text{new}\mathfrak{L}_{0} &=&\mathfrak{L}_{0}+\varepsilon \mathfrak{L}%
_{0}^{(1)}+\varepsilon ^{2}\mathfrak{L}_{0}^{(2)}+O(\varepsilon ^{3}), \\
\text{new}\mathfrak{L}_{1} &=&\mathfrak{L}_{1}^{(0)}+\varepsilon \mathfrak{L}%
_{1}^{(1)}+O(\varepsilon ^{2}), \\
\text{new}\mathcal{E}_{2}(\psi _{1},\cdot )
&=&\mathcal{E}_{2}^{(0)}(\psi ^{(0)},\cdot )+\varepsilon
\mathcal{E}_{2}^{(1)}(\psi ^{(0)},\cdot
)+O(\varepsilon ^{2}), \\
\text{new}\mathcal{E}_{3}(\psi _{1},\psi _{1},\cdot ) &=&\mathcal{E}%
_{3}^{(0)}(\psi ^{(0)},\psi ^{(0)},\cdot )+O(\varepsilon ),
\end{eqnarray*}%
and the functions $\psi _{1},\psi _{2}$ are transformed into%
\begin{eqnarray*}
\text{new }\psi _{1} &=&\psi ^{(0)}+\varepsilon \psi
_{1}^{(1)}+\varepsilon
^{2}\psi _{1}^{(2)}+O(\varepsilon ^{3}) \\
\text{new}\psi _{2} &=&\psi _{2}^{(0)}+\varepsilon \psi
_{2}^{(1)}+O(\varepsilon ^{2}).
\end{eqnarray*}%
Moreover, we have thanks to Lemma \ref{Lemmcoordchange}%
\begin{eqnarray*}
\kappa (Y) &=&1+\varepsilon \kappa _{1}(Y)+\varepsilon ^{2}\kappa
_{2}(Y)+O(\varepsilon ^{3}), \\
\text{new}\mathcal{H}(\psi ,\mu ) &=&\varepsilon
\mathcal{H}_{1}+\varepsilon ^{2}\mathcal{H}_{2}+O(\varepsilon
^{3}).
\end{eqnarray*}%
Hence, identities (\ref{ident}) lead to%
\begin{eqnarray}
\mathfrak{L}_{0}\psi ^{(0)} &=&0,  \notag \\
\mathfrak{L}_{0}^{(1)}\psi ^{(0)}+\mathfrak{L}_{0}\psi _{1}^{(1)}
&=&0,
\label{ident1} \\
\mathfrak{L}_{0}^{(2)}\psi ^{(0)}+\mathfrak{L}_{0}^{(1)}\psi _{1}^{(1)}+%
\mathfrak{L}_{0}\psi _{1}^{(2)} &=&0,  \notag
\end{eqnarray}%
\begin{eqnarray}
\mathfrak{L}_{0}\psi _{2}^{(0)}+\mathcal{E}_{2}^{(0)}(\psi
^{(0)},\psi
^{(0)}) &=&0,  \label{ident2} \\
\mathfrak{L}_{0}^{(1)}\psi _{2}^{(0)}+\mathfrak{L}_{0}\psi _{2}^{(1)}+2%
\mathcal{E}_{2}^{(0)}(\psi ^{(0)},\psi _{1}^{(1)})+\mathcal{E}%
_{2}^{(1)}(\psi ^{(0)},\psi ^{(0)}) &=&0,  \notag
\end{eqnarray}%
\begin{equation}
\mathfrak{L}_{0}\psi _{3}^{(0)}+\mu _{1}\mathfrak{L}_{1}^{(0)}\psi ^{(0)}+2%
\mathcal{E}_{2}^{(0)}(\psi ^{(0)},\psi _{2}^{(0)})+\mathcal{E}%
_{3}^{(0)}(\psi ^{(0)},\psi ^{(0)},\psi ^{(0)})=0.  \label{ident3}
\end{equation}%
We deduce from these formulae and from (\ref{d_psi E}), that the
linear
operator obtained after the change of coordinates satisfies%
\begin{equation*}
\mathfrak{L+H=L}_{0}+\varepsilon \mathfrak{H}^{(1)}+\varepsilon ^{2}%
\mathfrak{H}^{(2)}+O(\varepsilon ^{3})
\end{equation*}%
with%
\begin{equation}
\kappa _{1}\mathfrak{L}_{0}+\mathfrak{H}^{(1)}=\mathfrak{L}_{0}^{(1)}+2%
\mathcal{E}_{2}^{(0)}(\psi ^{(0)},\cdot
)+\mathfrak{L}_{0}\mathcal{H}_{1} \label{h(1)}
\end{equation}%
\begin{eqnarray*}
\kappa _{2}\mathfrak{L}_{0}+\kappa
_{1}\mathfrak{H}^{(1)}+\mathfrak{H}^{(2)}
&=&\mathfrak{L}_{0}^{(2)}+\mu _{1}\mathfrak{L}_{1}^{(0)}+2\mathcal{E}%
_{2}^{(1)}(\psi ^{(0)},\cdot )+2\mathcal{E}_{2}^{(0)}(\psi
_{2}^{(0)},\cdot
)+2\mathcal{E}_{2}^{(0)}(\psi _{1}^{(1)},\cdot )+ \\
&&+3\mathcal{E}_{3}^{(0)}(\psi ^{(0)},\psi ^{(0)},\cdot )+\left\{ \mathfrak{L%
}_{0}^{(1)}+2\mathcal{E}_{2}^{(0)}(\psi ^{(0)},\cdot )\right\} \mathcal{H}%
_{1}+\mathfrak{L}_{0}\mathcal{H}_{2}.
\end{eqnarray*}%
We now compute the terms under the integral in (\ref{@}). First we
observe
(thanks to (\ref{ident1}), (\ref{ident2}))%
\begin{eqnarray*}
\mathfrak{L}_{0}^{-1}\mathfrak{H}^{(1)}\psi ^{(0)} &=&\mathfrak{L}_{0}^{-1}\{%
\mathfrak{L}_{0}^{(1)}\psi ^{(0)}+2\mathcal{E}_{2}^{(0)}(\psi
^{(0)},\psi
^{(0)})\}+\mathcal{H}_{1}\psi ^{(0)} \\
&=&\mathfrak{L}_{0}^{-1}\{-\mathfrak{L}_{0}\psi _{1}^{(1)}-2\mathfrak{L}%
_{0}\psi _{2}^{(0)}\}+\mathcal{H}_{1}\psi ^{(0)} \\
&=&-(\psi _{1}^{(1)}+2\psi _{2}^{(0)})+\mathcal{H}_{1}\psi ^{(0)},
\end{eqnarray*}%
hence%
\begin{eqnarray*}
-\mathfrak{H}^{(1)}\mathfrak{L}_{0}^{-1}\mathfrak{H}^{(1)}\psi ^{(0)} &=&%
\mathfrak{L}_{0}^{(1)}(\psi _{1}^{(1)}+2\psi _{2}^{(0)})+2\mathcal{E}%
_{2}^{(0)}(\psi ^{(0)},\psi _{1}^{(1)})+ \\
&&+4\mathcal{E}_{2}^{(0)}(\psi ^{(0)},\psi _{2}^{(0)})+\mathfrak{L}_{0}%
\mathcal{H}_{1}(\psi _{1}^{(1)}+2\psi _{2}^{(0)})+ \\
&&-\kappa _{1}\mathfrak{L}_{0}(\psi _{1}^{(1)}+2\psi _{2}^{(0)})-\mathfrak{H}%
^{(1)}\mathcal{H}_{1}\psi ^{(0)}.
\end{eqnarray*}%
Moreover%
\begin{eqnarray*}
\mathfrak{H}^{(2)}\psi ^{(0)} &=&\mathfrak{L}_{0}^{(2)}\psi ^{(0)}+\mu _{1}%
\mathfrak{L}_{1}^{(0)}\psi ^{(0)}+2\mathcal{E}_{2}^{(1)}(\psi
^{(0)},\psi
^{(0)})+2\mathcal{E}_{2}^{(0)}(\psi _{2}^{(0)},\psi ^{(0)})+ \\
&&+2\mathcal{E}_{2}^{(0)}(\psi _{1}^{(1)},\psi ^{(0)})+3\mathcal{E}%
_{3}^{(0)}(\psi ^{(0)},\psi ^{(0)},\psi ^{(0)})-\kappa _{1}\mathfrak{H}%
^{(1)}\psi ^{(0)}+ \\
&&+\left\{ \mathfrak{L}_{0}^{(1)}+2\mathcal{E}_{2}^{(0)}(\psi
^{(0)},\cdot )\right\} \mathcal{H}_{1}\psi
^{(0)}+\mathfrak{L}_{0}\mathcal{H}_{2}\psi ^{(0)}
\end{eqnarray*}%
and in using again (\ref{ident1}), (\ref{ident2}), (\ref{ident3}), and (\ref%
{h(1)}) we obtain%
\begin{eqnarray*}
\mathfrak{H}^{(2)}\psi ^{(0)}-\mathfrak{H}^{(1)}\mathfrak{L}_{0}^{-1}%
\mathfrak{H}^{(1)}\psi ^{(0)} &=&-2\mu _{1}\mathfrak{L}_{1}^{(0)}\psi ^{(0)}-%
\mathfrak{L}_{0}(3\psi _{3}^{(0)}+2\psi _{2}^{(1)}+\psi _{1}^{(2)})+ \\
&&+\mathfrak{L}_{0}\{\mathcal{H}_{1}(\psi _{1}^{(1)}+2\psi _{2}^{(0)})+%
\mathcal{H}_{2}\psi ^{(0)}-\mathcal{H}_{1}^{2}\psi ^{(0)}\}.
\end{eqnarray*}%
Hence (\ref{@}) leads to%
\begin{equation*}
@=-2\mu _{1}\int_{\mathbb{T}^{2}}(\mathfrak{L}_{1}^{(0)}\psi
^{(0)})\psi ^{(0)}dY.
\end{equation*}%
Since
\begin{equation*}
\mathfrak{L}_{1}^{(0)}=-\frac{1}{\mu _{0}^{2}}\partial
_{y_{1}}^{2}
\end{equation*}%
we finally obtain the result of Lemma \ref{coef@}.

\section{Fluid particles dynamics}\label{particles}
The kinematic and dynamic boundary conditions (\ref{equation2},
\ref{equation3})
  give  two equations
for two unknowns $\eta$ and $\psi$.
 Assume for the moment that we  know $\eta$ i.e. the free surface
$\Sigma$.  The question is \emph{ can we restore $\psi$ without
solving PDE equations}? The answer is yes: it suffices to solve a
problem of moving a heavy single mass point along the free
surface, or equivalently to find the corresponding geodesic flow
on the surface with an appropriate metric.


Let us begin with the consideration of the motion of a single mass
point along the surface $\Sigma=\{x_3=\eta(X)\}$. Assuming that
gravity $\mu$ acts in the $-e_3=(0,0,-1)$ direction we can write
the governing equations in the form
$$
\ddot{x}+\mu e_3 =\lambda \mathbf{n},\quad x_3=\eta(X),
$$
where $\lambda$ is the Lagrange multiplier and $\mathbf{n}$ is a
normal vector to $\Sigma$. Choosing components of $X$ as
generalized coordinates we rewrite equivalently these equation in
the Lagrange form with the  Lagrangian
$$
\mathbf L(X,\dot X)=\mathbf T(X,\dot X)-\mathbf U(X)=\frac{1}{2}
\mathbb G(X)\dot X\cdot \dot X-\mu \eta(X),
$$
where $\mathbb GdX\cdot dX$ is   the first fundamental form of the
free surface $\Sigma$ with
$$
\mathbb G(X)= I+\nabla_X\eta(X)\otimes \nabla_X\eta(X).
$$
More precisely, we have
$$
\frac{d}{dt}\partial_{\dot X}\mathbf{L}(X,\dot X)-
\partial_{ X}\mathbf{L}(X,\dot X)=0.
$$
If we define the  moments  and  Hamiltonian  by
$$
y=\mathbb G \dot X, \quad \mathbf H(X,y)=\frac{1}{2}\mathbb
G^{-1}y\cdot y+\mu \eta(X),
$$
then the governing equations can be rewritten in the Hamilton form
\begin{equation}\label{Ham} \dot X=\partial_y \mathbf H(X,y),\quad \dot y=-\partial_X
\mathbf H(X,y). \end{equation} Next, suppose that a
  $C^1$-generating function $S(X)$ satisfies the Hamilton-Jacobi
equation
\begin{equation}\label{jacobi}
\mathbf H(X, \nabla_X S(X))=h=\text{~~const~~},
\end{equation}
 and the periodicity conditions
\begin{equation}\label{generatorper}
  S(X+2\pi
e_1)-2\pi = S(X+\frac{2\pi}{\tau} e_2)-\frac{2\pi}{\tau}= S(X).
\end{equation}
Suppose also that $X(t)$ is a solution of the equations
\begin{equation}\label{flow}
\dot X=\mathbb G^{-1}(X)\nabla_X S(X),
\end{equation}
then, it is known that $(X(t), y(t))$, with $y(t)=\nabla_X
S(X(t))$, is a solution of the Hamiltonian system \eqref{Ham}, and
the surface ${y=\nabla_XS(X)}$ is an invariant manifold of
$\eqref{Ham}$, the flow being defined by (\ref{flow}).

Finally note that due the periodicity conditions, the mapping
$X\to (X,\nabla_XS(X))$ defines an embedding of the torus
$\RR^2/\Gamma$ into $\RR^2/\Gamma\times \RR^2$. Therefore,
$\{y=\nabla_X S(X)\}$ is an invariant torus of \eqref{Ham}
 lying on the energy surface $H=h$. In coordinates $X$ the
 Hamiltonian flow on the torus is given by \eqref{flow}

Let us turn to the diamond wave problem. Set
$$\varphi^*(x)=\mathbf u_0\cdot X+\varphi(x),\text{~~and~~}
\psi^*(X)=\varphi^*(x_1,x_2, \eta(X)),
$$
and recall  $|\mathbf u_0|=1$. In these notations   kinematic
condition \eqref{equation2}  and dynamic condition
\eqref{equation3} can be rewritten in the equivalent form
\begin{gather}\label{equation2a}
\nabla_X\varphi^*(x_1,x_2,x_3)=\mathbb
G(X)^{-1}\nabla_X\psi^*(X)\text{~~for~~} x_3=\eta(X),\\
\label{equation3a} \frac{1}{2}\mathbb
G^{-1}\nabla_X\psi^*\cdot\nabla_X\psi^*+\mu\eta=\frac{1}{2}.
\end{gather}
By construction the equation $\dot x=\nabla\varphi^*(x)$
determines the trajectories of liquid particles. From
\eqref{equation2a}, such a particle moving along the free surface
satisfies
 \begin{equation}\label{equation4a}
 \dot X=G^{-1}(X)\nabla_X\psi^*(X).
\end{equation}
On the other hand, equation \eqref{equation3a} reads
$$
\mathbf H(X, \nabla_X \psi^*)=1/2.
$$
Hence $S(X)=\psi^*(X)$ is a generating function for the
Hamiltonian system \eqref{Ham}.

From this we conclude that trajectories of liquid particles
$X(t)$, which are defined by $\dot X= \mathbb G^{-1}\nabla_X
\psi^*(X)$, along with $y(t)=\nabla_X S(X(t))$ serve as solutions
of \eqref{Ham} and belongs to the invariant torus $\{y=\nabla_X
S(X)\}$. In other words,  they coincide with projections $(X,y)\to
X$ of solutions to \eqref{Ham} belonging to the invariant torus
$\{y=\nabla \psi^*(X)\}\subset \{\mathbf H=1/2\}$. Moreover, since
by \eqref{Z} we have $V=\mathbb G^{-1}\nabla_X \psi^*$, they also
coincide with the integral curves of the vector field $V$.

Finally  note that, by the  Maupertuis principle, the projections
$(X,y)\to X$ of solutions of \eqref{Ham} belonging to the energy
surface $\mathbf H=1/2$, coincide with the geodesics on the
manifold $\Sigma$ endowed with the Jacobi metric
\begin{equation}\label{metrictensor} ds^2=\big(1/2-\mu \eta(X)\big)\,\mathbb
G(X)\, dX  \cdot  dX\equiv 2(1/2-\mathbf U)\mathbf T(X,dX).
\end{equation}
Hence the integral curves of the vector field $V$ form the
geodesic flow associated with  the metric \eqref{metrictensor}.

\begin{corollary} Suppose that $\eta$ is an arbitrary bi-periodic smooth
function so that the hamiltonian system \eqref{Ham} has an
invariant torus $\{y=\nabla_X S(X)\}$ with a smooth generating
function $S(X)$ satisfying the periodicity conditions
\eqref{generatorper}.  Then the solution $\varphi^*$ of the Cauchy
problem
$$
\varphi^*(x)=S(X),\quad \partial_n \varphi^*(x)=0\text {~~for
~~}x_3=\eta(X),
$$
for the Laplace equation,  satisfies  kinematic and dynamic
conditions \eqref{equation2},\eqref{equation3}. (The local
existence of such a solution follows from the Cauchy-Kowalewski
theorem, and the existence and boundedness in the lower half plane
is true only for the "good" choice of generating function.)
\end{corollary}

%


\begin{thebibliography}{99}
\bibitem{AntonineBarucq} X. Antoine, H. Barucq, A.Bendali.  Bayliss-Turkel-
like Radiation Conditions on Surfaces of Arbitrary Shape. J. Math.
Anal. Appl. \textbf{229} (1999), 184-211.
\bibitem{Arnold} V.I. Arnold. Proof of a theorem of A.N. Kolmogorov
on the invariance of quasi-periodic motions under small
perturbations of the Hamiltonian. Russ. Math. Surv.(1963),
\textbf{18}, 9-36.
\bibitem{B95} J. Bourgain. Construction of periodic solutions of
nonlinear wave equations in higher dimension. Geom. Funct. Anal.
\textbf{5} (1995), 629-639.
\bibitem{B98} J. Bourgain. Quasi-periodic solutions of Hamiltonian
perturbations of $2D$ linear Schr\"odinger equations. Ann. of
Math. \textbf{148} (1998), 363-439

\bibitem{B-D-M} T.Bridges, F.Dias, D.Menasce. Steady three-dimensional
water-wave patterns on a finite-depth fluid. J.Fluid Mech. (2001), \textbf{%
436}, 145-175.
\bibitem{Cassels} J.W.S. Cassels. An Introduction to Diophantine
approximations.  Cambridge University Press 1957.

\bibitem{C00} W. Craig. Probl\`emes de petits diviseurs dans les
\'equations aux d\'eriv\'ees partielles. Panoramas et Synth\`eses
\textbf{9}. Soci\'et\'e Math\'ematique de France, Paris  2000.

\bibitem{craig} W.Craig, D.Nicholls. Traveling gravity water waves in two
and three dimensions. EJMB/Fluids \textbf{21} (2002) 615-641.


\bibitem{Craig-Nicholls2000} W.Craig, D.Nicholls. Travelling two and
three-dimensional capillary gravity water waves. SIAM J. Math. Anal. \textbf{%
32} (2000) 323-359.

\bibitem{Craig-Sch-Sul} W.Craig, U.Schanz, C.Sulem. The modulational regime
of three-dimensional water waves and the Davey-Stewartson system.
Ann. Inst. Henri Poincar\'{e}, \textbf{14}, 5 (1997), 615-667.

\bibitem{CW93} W. Craig, E. Wayne. Newton's method and periodic
solutions of nonlinear wave equation. Commu. Pure Applied Math.
\textbf{ XLVI} (1993) 1409-1501.


\bibitem{Deimling} K. Deimling. Nonlinear Functional
Analysis, Springer-Verlag, Heidelberg, 1985
\bibitem{Dias-Io} F.Dias, G.Iooss. Water waves as a spatial dynamical
system. Handbook of Mathematical Fluid Dynamics, chap 10, p.443
-499. S.Friedlander, D.Serre Eds., Elsevier 2003.

\bibitem{dia-khar} F.Dias, C.Kharif. Nonlinear gravity and capillary-gravity
waves. Annu. Rev. Fluid Mech. (1999), \textbf{31}, 301-346.
\bibitem{FS83} J. Fr\"ohlich, T. Spencer. Absence of diffusion in
the Anderson tight binding model for large disorder or low energy.
Comm. Math. Phys. \textbf{88} (1983), 151-184.
\bibitem{Fuchs} R.Fuchs. On the theory of short-crested oscillatory waves.
U.S. Natl. Bur. Stand. Circ. \textbf{521} (1952), 187-200.

\bibitem{Gro} M.D.Groves. An existence theory for three-dimensional periodic
travelling gravity-capillary water waves with bounded transverse
profiles. Physica D \textbf{152-153} (2001), 395-415.

\bibitem{Grov-Harag} M.D.Groves, M.Haragus. A bifurcation theory for
three-dimensional oblique travelling gravity-capillary water
waves. J.Nonlinear Sci. \textbf{13} (2003), 397-447.

\bibitem{Gr-Mi} M.Groves, A.Mielke. A spatial dynamics approach to
three-dimensional gravity-capillary steady water waves. Proc. Roy.
Soc. Edin. A \textbf{131} (2001), 83-136.

\bibitem{Ham-Hend-Seg} J.Hammack, D.Henderson, H.Segur. Progressive waves
with persistent, two-dimensional surface patterns in deep water.
J.Fluid Mech. \textbf{532} (2005), 1-52.

\bibitem{Harag-Kirch} M.Haragus-Courcelle, K.Kirchg\"{a}ssner.
Three-dimensional steady capillary-gravity waves. Ergodic theory,
Analysis and efficient simulation of dynamical systems (Ed.
B.Fiedler), p.363-397. Berlin, Springer-Verlag 2001.
\bibitem{Hormander}L. H\"ormander. Pseudo-differential operators
and non-elliptic boundary problems. Ann. of Math. \textbf{83}
(1966), 129-209.
\bibitem{Io} G.Iooss.  Capillary and Capillary-Gravity periodic travelling waves
for two superposed fluid layers, one being of infinite depth . J.
Math. Fluid Mech. \textbf{1}, (1999), 24-61.

\bibitem{IPT} G.Iooss, P.Plotnikov, J.Toland. Standing waves on an
infinitely deep perfect fluid under gravity. Arch. Rat. Mech. Anal. \textbf{%
177} (2005), 3, 367-478.

\bibitem{Kirch} K.Kirchg\"{a}ssner. Wave solutions of reversible systems and
applications. J.Diff. Eqns. \textbf{45} (1982), 113-127.

\bibitem{KohnNirenberg}J.J. Kohn, L. Nirenberg. An algebra of
pseudodifferential operators. Comm. Pure. Apll. Math. \textbf{18}
(1965), (269-305).

\bibitem{Lannes} D.Lannes. Well-posedness of the water-waves equations.
J.Amer. Math. Soc. \textbf{18} (2005), 605-654.

\bibitem{LeviCivita} T.Levi-Civita. D\'{e}termination rigoureuse des ondes
permanentes d'ampleur finie. Math. Annalen \textbf{93} (1925),
264-314.
\bibitem{Moser}J. Moser. Minimal foliation on a torus.
Topics in calculus of variations (Montecatini Terme). Lecture
Notes in Math, \textbf{1365} (1989), 62-99.
\bibitem{Nekrasov} A.I.Nekrasov. On waves of permanent type. Izv.
Ivanovo-Voznesensk. Politekhn. Inst., \textbf{3} (1921), 52-65.
\bibitem{PIP2}P.I.  Plotnikov. Solvability of the problem of spatial gravitational waves on the
surface of an ideal fluid. Dokl. Akad. Nauk SSSR.
\textbf{251}(1980), 170-171.

\bibitem{PIP1} L.V. Ovsiannykov, N.I. Makarenko, V.I. Nalimov, V. Yu. Liapidevskii,
P.I. Plotnikov, I.V. Sturova, V.I. Bukreev, V.A. Vladimirov.
Nonlinear problems in the theory of surface and internal waves
(Russian). (1985), 165-199, Nauka, Novosibirsk. Plotnikov, P.I.;

\bibitem{Petersen} B.E. Petersen. An introduction to the Fourier
transform and pseudodifferential operators.  Pitman Advanced
Publishing programm, Boston, London, Melbourne 1983.


\bibitem{PIP3} P.I. Plotnikov, L.N. Yungerman. Periodic solutions of a weakly linear wave
equation with an irrotational ratio of the period to the interval
length. Diff.Equations \textbf{24} (1988), 1059-1065.

\bibitem{Plot-Tol} P.Plotnikov, J.Toland. Nash-Moser theory for standing
waves. Arch. Rat. Mech. Anal. \textbf{159 }(2001), 1--83.

\bibitem{Reed-Shin} J.Reeder, M.Shinbrot. Three-dimensional, nonlinear wave
interaction in water of constant depth. Nonlinear Anal., T.M.A.,
\textbf{5} (1981), 3, 303-323.

\bibitem{Rob-Schw} A.Roberts, L.Schwartz. The calculation of nonlinear
short-crested gravity waves. Phys. Fluids \textbf{26}, 9 (1983),
2388-2392.
\bibitem{Siegel} C. L. Siegel. Vorlesungen \" Uber Himmelsmechanik
(German, Russian). Inostarannaya Literatura, Moskow 1959.

\bibitem{Sret} L.Sretenskii. Spatial problem of determination of steady
waves of finite amplitude (russian). Dokl. Akad. Nauk SSSR (N.S.)
\textbf{89} (1953), 25-28.

\bibitem{Stokes} G.G.Stokes. \ On the theory of oscillatory waves. Trans.
Camb. Phil. Soc. \textbf{8} (1847), 441-455.

\bibitem{Struik} D. Struik. D\'{e}termination rigoureuse des ondes
irrotationnelles p\'{e}riodiques dans un canal \`{a} profondeur
finie. Math. Ann. 95 (1926), 595-634.

\bibitem{MTaylor} M.E. Taylor. Pseudodifferential operators.
 Princeton, New Jersey 1981.


\bibitem{Weil} H.Weil. Uber die gleichverteilung der zahlen rood eins. Math.
Ann., \textbf{77}, (1916) 313-352.

\bibitem{V.E.Zakharov} V.E.Zakharov. Stability of periodic waves of finite
amplitude on the surface of a deep fluid. Zh. Prikl. Mekh. Tekh.
Fiz. \textbf{9} (1968), 86-94, J.Appl. Mech. Tech. Phys.
\textbf{9} (1968) 190-194.
\end{thebibliography}
\end{document}